\documentclass[12pt,twoside]{article}
\usepackage[hmargin=0.8in,vmargin=0.9in]{geometry}
\usepackage{fancyhdr}
\usepackage{graphicx}
\usepackage{subfigure}
\usepackage{amssymb}
\usepackage{amsmath}
\usepackage{amsfonts}
\usepackage{theorem}
\usepackage{mathrsfs}
\usepackage{mathtools}
\usepackage{bm}
\usepackage{color}
\usepackage{setspace}
\usepackage{exscale}
\usepackage{relsize}
\usepackage{epstopdf}
\usepackage{float}
\usepackage{cite}
\usepackage{nicefrac}

\usepackage{booktabs,multirow} 
\usepackage{array} 
\usepackage{paralist} 
\usepackage{verbatim} 
\usepackage{subfigure} 

\theoremstyle{plain}			
\newtheorem{thm}{Theorem}[section]

{\theorembodyfont{\rmfamily}}
\newenvironment{acknowledgment}{{\flushleft \bf Acknowledgment:}}{}

\usepackage{tikz}
\usetikzlibrary{positioning}
\usepackage{xcolor}

\allowdisplaybreaks[1]

\numberwithin{equation}{section}
\numberwithin{figure}{section}
\numberwithin{table}{section}

\newcommand\eref[1]{(\ref{#1})}

\newcommand*\xbar[1]{%
  \hbox{%
    \vbox{%
      \hrule height 0.5pt 
      \kern0.4ex
      \hbox{%
        \kern-0.05em
        \ensuremath{#1}%
        \kern-0.00em
      }%
    }%
  }%
}

\newcommand{\mF}{\bm{F}}

\newcommand{\mH}{\bm{H}}
\newcommand{\mU}{\bm{U}}
\newcommand{\mE}{\bm{E}}
\newcommand{\mV}{\bm{V}}

\newcommand{\mW}{\bm{W}}

\newcommand{\dt}{\Delta t}
\newcommand{\dx}{\Delta x}

\newcommand{\hf}{{\frac{1}{2}}}

\newcommand{\jph}{{j+\frac{1}{2}}}
\newcommand{\jmh}{{j-\frac{1}{2}}}

\newcommand{\cV}{{\cal V}}
\newcommand{\cU}{{\cal U}}

\newtheorem{rmk}[thm]{Remark}

\title{Experimental Convergence Rate Study for Three Shock-Capturing Schemes and Development of Highly Accurate Combined Schemes}
\author{Shaoshuai Chu\thanks{Department of Mathematics, Southern University of Science and Technology, Shenzhen, 518055, China;
{\tt 11930702@mail.sustech.edu.cn}},~Olyana A. Kovyrkina\thanks{Lavrentyev Institute of Hydrodynamics of SB RAS and Novosibirsk State
University, Novosibirsk, 630090, Russia; {\tt olyana@ngs.ru}},~Alexander Kurganov,\thanks{Corresponding author; Department of Mathematics,
Shenzhen International Center for Mathematics, and Guangdong Provincial Key Laboratory of Computational Science and Material Design,
Southern University of Science and Technology, Shenzhen, 518055, China; {\tt alexander@sustech.edu.cn}}\\and Vladimir V.  Ostapenko
\thanks{Lavrentyev Institute of Hydrodynamics of SB RAS and Novosibirsk State University, Novosibirsk, 630090, Russia;
{\tt ostapenko\_vv@ngs.ru}}}

\begin{document}

\date{}
\maketitle

\begin{abstract}
We study experimental convergence rates of three shock-capturing schemes for hyperbolic systems of conservation laws: the second-order
central-upwind (CU) scheme, the third-order Rusanov-Burstein-Mirin (RBM), and the fifth-order alternative weighted essentially
non-oscillatory (A-WENO) scheme. We use three imbedded grids to define the experimental pointwise, integral, and $W^{-1,1}$ convergence
rates. We apply the studied schemes to the shallow water equations and conduct their comprehensive numerical convergence study. We verify
that while the studied schemes achieve their formal orders of accuracy on smooth solutions, after the shock formation, a part of the
computed solutions is affected by shock propagation and both the pointwise and integral convergence rates reduce there. Moreover, while the
$W^{-1,1}$ convergence rates for the CU and A-WENO schemes, which rely on nonlinear stabilization mechanisms, reduce to the first order, the
RBM scheme, which utilizes a linear stabilization, is clearly second-order accurate. Finally, relying on the conducted experimental
convergence rate study, we develop two new combined schemes based on the RBM and either the CU or A-WENO scheme. The obtained combined
schemes can achieve the same high-order of accuracy as the RBM scheme in the smooth areas while being non-oscillatory near the shocks.
\end{abstract}



\section{Introduction}
We focus on measuring the experimental convergence rates of three different high-order schemes for the one-dimensional (1-D) hyperbolic
systems of conservation laws
\begin{equation}
\mU_t+\mF(\mU)_x=\bm0,
\label{1.1}
\end{equation}
where $x$ is a spatial variable, $t$ is the time, $\mU\in\mathbb R^d$ is a vector of unknowns, and $\mF: \mathbb R^d\to\mathbb R^d$ are
nonlinear fluxes.

It is well-known since the pioneering work of Godunov \cite{God} that linear monotone schemes can be at most first-order accurate. In order
to overcome this Godunov's order barrier theorem, many nonlinear numerical methods had been introduced; see, e.g.,
\cite{Hes,LeV02,Shu_Acta,Tor} and references therein. Most of these methods utilize nonlinear mechanisms to detect nonsmooth parts of the
computed solution, where the numerical fluxes are corrected to prevent oscillations at shocks and other ``rough'' areas. However, when
the computed solutions are discontinuous, these schemes are typically only formally high-order as they produce ${\cal O}(1)$ errors in the
vicinities of discontinuities. Moreover, as it was shown in \cite{CC,ES,Ost} the local convergence rate of such schemes typically reduces
to the first or even lower order in the area of influence of the shock waves, which are smeared by the numerical viscosity. This occurs due
to numerical inaccuracy in the Rankine-Hugoniot conditions at the shock and hence to the error propagation behind the shock.

In order to quantify the influence shock waves exert on the accuracy of the computed smooth parts of the solution, one can measure not only
the local pointwise convergence rates, but also the integral ones, which, in fact, correspond to the convergence in the negative norms. In
particular, one can measure $W^{-1,1}_{\rm loc}$ convergence rates as it was done in \cite{KO10,KO14,KO21,LNOT} for a variety of high-order
schemes, for which the order reduction in the shock influence areas had been observed. At the same time, there are high-order schemes, for
instance, the third-order Rusanov-Burstein-Mirin (RBM) scheme \cite{Burstein70,Rusanov68} and the compact third-order weak approximation
scheme \cite{Ost00}, which are not based on nonlinear numerical flux correction and which can achieve higher experimental convergence rates
in the smooth parts of the solution even in presence of propagating shock waves; see \cite{KKO,KO10,KO14}.

In this paper, we study the experimental convergence rates of the semi-discrete second-order finite-volume (FV) central-upwind (CU) scheme
from \cite{KNP}, the third-order finite-difference (FD) RBM scheme and the semi-discrete fifth-order FD alternative weighted essentially
non-oscillatory (A-WENO) scheme from \cite{Kurganov20}, which is based on the CU numerical fluxes from \cite{Kurganov07} and \mbox{WENO-Z}
interpolation. The studied schemes are briefly described in Appendices \ref{appa}--\ref{appc}. Both the CU and A-WENO schemes use nonlinear
stabilization mechanisms, which, as we will demonstrate, lead to a substantial convergence rate reduction, both local and integral ones. On
contrary, the RBM scheme uses a linear stabilization mechanism, which leads to a better pointwise and integral convergence in the smooth
parts of the solution affected by the propagating shock waves.

We also study the experimental $W^{-1,1}$ convergence rates for the aforementioned three schemes. As it was shown in \cite{NT92,NTT}, the
convergence in the $W^{-1,1}$-norm together with the TV boundedness of the approximate solutions yield the $L^\infty_{\rm loc}$ convergence
away from shocks, and the $L^\infty_{\rm loc}$ convergence rates are the same as the $W^{-1,1}$ ones.

We conduct the experimental convergence rate study on two Cauchy problems for the Saint-Venant system of shallow water equations with smooth
periodic initial data. The solutions of the studied problems develop shocks, which propagate at a variable speed and the domain of the shock
influence grows in time. Thus, these examples serve as excellent benchmarks. We demonstrate that when the solutions are still smooth (that
is, before the breakdown), all of the studied three schemes achieve their formal order of accuracy. After the shock formation, the orders
inside the shock influence areas reduce to the first one for the CU and A-WENO schemes and to the second one for the RBM scheme. At larger
time, when the entire domain is influenced by the propagating shock waves, the CU and A-WENO schemes are globally first-order, while the RBM
is second-order accurate.

Relying on the conducted convergence study, we develop two new combined schemes based on the RBM as a basic scheme and either the CU or
A-WENO as an internal scheme. The technique for constructing combined schemes has been recently proposed in \cite{KO18,LNOT19,ZKO}. In
combined schemes, the numerical solution is first constructed in the entire computational domain according to the basic nonmonotone scheme
that maintains increased accuracy inside the shock influence areas. In the ``rough'' parts of this solution containing large gradients,
nonphysical oscillations develop and thus the nonmonotone solution has to be corrected there. This is done by solving the internal
initial-boundary value problems in the ``rough'' areas by a high-order non-oscillatory internal scheme. A somewhat similar scheme adaption
approach is used in hybrid schemes (see, e.g., \cite{ConK,DKL,Karni02}), in which different methods are used for the evolution of the
computed solution in the ``smooth'' and ``rough'' areas. We, however, should point out at a fundamental difference between the hybrid and
combined methods: in the latter ones, an internal non-oscillatory scheme does not affect the solution obtained by the basic nonmonotone
scheme, and thus combined schemes are capable of accurately localizing shock waves and capturing them in a non-oscillatory manner while
preserving the high accuracy away from the shocks. At the same time, hybrid schemes, like other schemes whose stabilization mechanism is
based on a certain nonlinear limiting procedure, suffer from the convergence rate reduction in the areas of shock influence.

A major novelty and advantage of the new combined schemes is in the way the ``rough'' parts of the computed solution are detected. Unlike
the combined schemes from \cite{KO18,LNOT19,ZKO}, where the smoothness was determined using the size of the solution gradients, here we use
the method, which is based on the weak local residual (WLR) introduced in \cite{Karni05,Karni02}. Using the fact that the magnitudes of the
WLR in the smooth areas and near the shocks are significantly different, we can accurately determine ``rough'' areas, and this leads to a
very good agreement between the numerical solutions obtained by the basic and internal schemes at the boundary of the internal computational
domain.  In order to further investigate the developed combined schemes, we perform numerical tests to experimentally check
their rates of pointwise and integral convergence. From the numerical results reported in \S\ref{sec7}, one can see that while the pointwise
convergence of the combined schemes is about the same as of the RBM scheme, the integral convergence rates are reduced and this reduction is
attributed to combining the schemes of a different nature.

The paper is organized as follows. In \S\ref{sec2a}, we describe the construction of three imbedded grids required to conduct the
experimental convergence rate study. In \S\ref{sec3f}, we introduce the way experimental pointwise and integral rates of convergence are
going to be measured. The studied schemes are applied to the Saint-Venant system, briefly described in \S\ref{sec4}. The obtained numerical
results are presented and analysed in \S\ref{sec5}. In \S\ref{sec7}, we introduce our new combined schemes and test them on the same
examples studied in \S\ref{sec5}. Finally, in \S\ref{sec8}, we give some concluding remarks.

\section{Construction of Three Imbedded Grids}\label{sec2a}
\begin{figure}[ht!]
\centerline{\includegraphics[trim=0.1cm 0.1cm 0.1cm 0.1cm, clip, width=14cm]{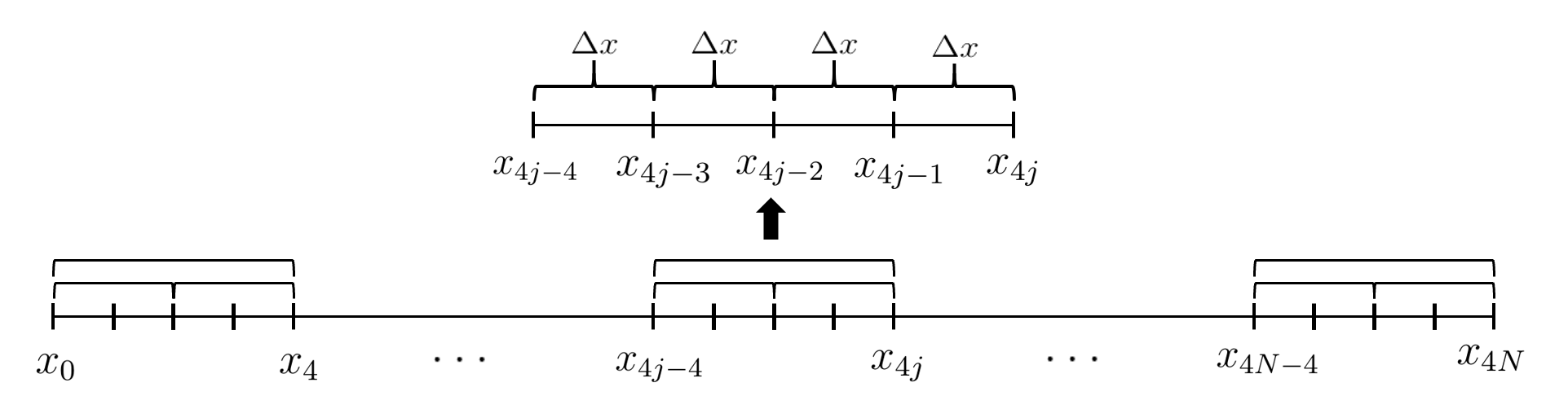}}
\caption{\sf Structure of the three imbedded meshes.\label{Fig1.1}}
\end{figure}

We cover the computational domain with three imbedded uniform grids consisting of cells of sizes $\dx$, $2\dx$, and $4\dx$, respectively;
see Figure \ref{Fig1.1}. We first introduce the finest of these three meshes. To this end, we split the computational domain $[a,b]$ into
$4N$ uniform cells denoted by $C^{4N}_\jph:=[x_j,x_{j+1}]$ for $j=0,\ldots,4N-1$. Here, $x_{j+1}-x_j=\dx=\nicefrac{(b-a)}{4N}$. The second
mesh consists of $2N$ uniform cells denoted by $C^{2N}_{2j+1}:=[x_{2j},x_{2j+2}]$ for $j=0,\ldots,2N-1$. The coarsest mesh consists of $N$
uniform cells denoted by $C^N_{4j+2}:=[x_{4j},x_{4j+4}\big]$ for $j=0,\ldots,N-1$. We note that the centers of the cells $C^{4N}_\jph$,
$C^{2N}_{2j+1}$, and $C^N_{4j+2}$ are $x_\jph=(x_j+x_{j+1})/2$, $x_{2j+1}$, and $x_{4j+2}$, respectively.

We will study both the FV and FD methods. For the FV CU scheme, we will evolve in time the cell averages of the solution on the three
aforementioned grids and we will denote these quantities by
\begin{equation}
\begin{aligned}
&\xbar\mV^{\,4N}_\jph(t)\approx\frac{1}{\dx}\!\int\limits_{C^{4N}_\jph}\mU(x,t)\,{\rm d}x,\quad
\xbar\mV^{\,2N}_{2j+1}(t)\approx\frac{1}{2\dx}\!\int\limits_{C^{2N}_{2j+1}}\mU(x,t)\,{\rm d}x,\\
&\xbar\mV^{\,N}_{4j+2}(t)\approx\frac{1}{4\dx}\!\int\limits_{C^N_{4j+2}}\mU(x,t)\,{\rm d}x.
\end{aligned}
\label{2.1}
\end{equation}
For the FD RBM and A-WENO schemes, the evolved in time quantities are the corresponding point values denoted by
\begin{equation}
\mV^{4N}_\jph(t)\approx\mU(x_\jph,t),\quad\mV^{2N}_{2j+1}(t)\approx\mU(x_{2j+1},t),\quad\mV^N_{4j+2}(t)\approx\mU(x_{4j+2},t).
\label{2.2}
\end{equation}
\begin{rmk}
We note that the time-dependence of most of the indexed quantities will be omitted in the rest of the paper except for Appendix \ref{appb},
where the fully discrete RBM scheme will be described.
\end{rmk}

\section{Experimental Rates of Convergence}\label{sec3f}
In this section, we introduce the ways experimental rates of convergence can be estimated.

\subsection{Pointwise Convergence}\label{sec2}
We begin with the pointwise convergence. We denote the local rates of convergence by $r_{4j}$ and compute them at the end points of the
cells of the coarsest mesh, that is, at $x=x_{4j}$ (see Figure \ref{Fig1.1}) using the Runge formula:
\begin{equation}
r_{4j}:=\log_\hf\Bigg|\frac{\cV^{2N}_{4j}-\cV^{4N}_{4j}}{\cV^N_{4j}-\cV^{2N}_{4j}}\Bigg|,\quad j=0,\ldots,N.
\label{3.1}
\end{equation}
Here, $\cV$ denotes a certain component of the vector $\mV$, and $\cV^{mN}_{4j}$, $m=1$, 2, 4 are the approximations of $\cU$ (which denotes
the same component of the vector $\mU$) at the point $x=x_{4j}$ computed by the three imbedded uniform grids introduced in \S\ref{sec2a}. We
compute these point values in the following manner:
\begin{equation*}
\cV^{mN}_{4j}=\hf\left[\big(\cV^{mN}_{4j}\big)^-+\big(\cV^{mN}_{4j}\big)^+\right],\quad m=1,2,4,
\end{equation*}
where $\big(\cV^{mN}_{4j}\big)^\pm$ are the corresponding right- and left-sided values at the cell interface $x=x_{4j}$. For the CU scheme,
we compute $\big(\cV^{mN}_{4j}\big)^\pm$ using a piecewise linear reconstruction \eref{A5}--\eref{A7} described in Appendix \ref{appa}. For
the third-order RBM and fifth-order A-WENO schemes, these values are computed using the fifth-order WENO-Z interpolant (see
\eref{C2}--\eref{C5} in Appendix \ref{appc}) applied to the local characteristic variables. Note that even though the WENO-Z interpolant is
not used in the RBM scheme, we apply it as a post-processing procedure to generate the required point values. As the WENO-Z interpolant has
a higher order of accuracy than the RBM scheme, this post-processing procedure does not affect the estimated convergence rates of the RBM
scheme.

\subsection{Integral Convergence}\label{sec3}
We now show how to compute the experimental integral rates of convergence, which at the point $x=x_{4j}$ is defined by analogy with
\eref{3.1} as
\begin{equation}
r^{\rm INT}_{4j}:=\log_\hf\Bigg|\frac{I^{2N}_{4j}-I^{4N}_{4j}}{I^N_{4j}-I^{2N}_{4j}}\Bigg|,\quad j=1,\ldots,N,
\label{4.1}
\end{equation}
where $I^{mN}_{4j}$, $m=1$, 2, 4 are computed recursively as follows. First, we set $I^{mN}_0=0$ and then for a certain component $\cU$ of
the vector $\mU$ and for $m=1$, 2 and 4, we compute
\begin{equation}
I^{mN}_{4j}=I^{mN}_{4j-4}+J_{mN}(x_{4j-4},x_{4j}),\quad j=1,\ldots,N,
\label{4.2}
\end{equation}
where
\begin{equation}
J_{mN}(x_{4j-4},x_{4j})\approx\int\limits_{x_{4j-4}}^{x_{4j}}\cU\,{\rm d}x.
\label{4.2a}
\end{equation}
In order to estimate the integral rate of convergence for the CU, RBM, and A-WENO schemes, the integral in \eref{4.2a} needs to be evaluated
for $m=1$, 2, and 4.

For the FV CU scheme, the numerical integrals $J_{mN}$ is computed in a straightforward manner:	
\begin{equation*}
\begin{aligned}
J_N(x_{4j-4},x_{4j})&=4\dx\,\xbar\cV^{\,N}_{4j-2}, \\
J_{2N}(x_{4j-4},x_{4j})&=J_{2N}(x_{4j-4},x_{4j-2})+ J_{2N}(x_{4j-2},x_{4j})=
2\dx\left(\,\xbar\cV^{\,2N}_{4j-3}+\,\xbar\cV^{\,2N}_{4j-1}\right),\\
J_{4N}(x_{4j-4},x_{4j})&= \sum_{i=1}^4J_{4N}(x_{4j-i},x_{4j-i+1})= \dx\sum_{i=1}^4\xbar\cV^{\,4N}_{4j-i+1/2},
\end{aligned}
\end{equation*}	
where the values $\,\xbar\cV^N_{\,4j+2}$, $\,\xbar\cV^{\,2N}_{2j+1}$, and $\,\xbar\cV^{\,4N}_\jph$ are the corresponding components of the
cell averages defined in \eref{2.1}.

For the RBM and A-WENO schemes, we calculate $J_{mN}$ by approximating the integral in \eref{4.2a} using the following six-order accurate
quadrature. For $m=1$, we first use five points stencil centered at $x=x_{4j-2}$ to construct the corresponding fourth-degree interpolating
polynomial. We then integrate this polynomial over the cell $C^N_{4j-2}$ to obtain
\begin{equation}
J_N(x_{4j-4},x_{4j})=4{\cal L}[\cV^N_{4j-2};4],
\label{4.3}
\end{equation}
where
\begin{equation}
{\cal L}[\cV^M_k; \ell]:=\frac{\dx}{5760}\Big(17\cV^M_{k-2\ell}+308\cV^M_{k-\ell}+5178\cV^M_k+308\cV^M_{k+\ell}-17\cV^M_{k+2\ell}\Big).
\label{4.3a}
\end{equation}
For $m=2$ and 4, we split the cell $C^N_{4j-2}=[x_{4j-4},x_{4j}]$ into two or four subintervals of sizes $2\dx$ and $\dx$, respectively, and
then use the rescaled versions of the quadrature \eref{4.3}--\eref{4.3a} on each of the subintervals. This results in
\begin{align}
J_{2N}(x_{4j-4},x_{4j})&=J_{2N}(x_{4j-4},x_{4j-2})+J_{2N}(x_{4j-2},x_{4j})=
2\left(\mathcal{L}[\cV^{2N}_{4j-3};2] +\mathcal{L}[\cV^{2N}_{4j-1};2]\right),\notag\\
J_{4N}(x_{4j-4},x_{4j})&=\sum_{i=1}^4J_{4N}(x_{4j-i},x_{4j-i+1})=\sum_{i=1}^4\mathcal{L}[\cV^{4N}_{4j-i+\hf};1].\label{4.3c}
\end{align}
Note that the values $\cV^N_{\,4j+2}$, $\cV^{\,2N}_{2j+1}$, and $\cV^{\,4N}_\jph$ used in \eref{4.3}--\eref{4.3c} are the corresponding
components of the point values defined in \eref{2.2}.

We would like to stress that the quantities $I^{mN}_{4j}$ are, in fact, approximated anti-derivatives of $\cU$, that is,
\begin{equation}\label{4.4a}
I^{mN}_{4j}\approx\int\limits_{x_0}^{x_{4j}}\cU\,{\rm d}x,
\end{equation}
which immediately follows from the recursive relation \eref{4.2}. Therefore, the integral convergence rates introduced in \eref{4.1} are the
pointwise convergence rates for the anti-derivatives. It is also instructive to evaluate the corresponding global $L^1$ convergence rate for
the anti-derivatives, which is, in fact, the $W^{-1,1}$ convergence rate and which is defined by
\begin{equation}
r^{\rm INT}:=\log_\hf\bigg(\frac{\|I^{2N}-I^{4N}\|_{L^1}}{\|I^N-I^{2N}\|_{L^1}}\bigg),
\label{4.4}
\end{equation}
where the discrete $L^1$-norm of any grid function $\psi$ is defined by
\begin{equation*}
\|\psi\|_{L^1}:=4\dx\left(|\psi_4|+\ldots+|\psi_{4N}|\right).
\end{equation*}

\section{Application to the Saint-Venant System}\label{sec4}
In this section, we apply the studied CU, RBM, and A-WENO schemes, briefly described in Appendices \ref{appa}, \ref{appb}, and \ref{appc},
respectively, to the Saint-Venant system of shallow water equations. In the case of a flat bottom topography, the Saint-Venant system reads
as \eref{1.1} with
\begin{equation}
\bm U=\begin{pmatrix}h\\q\end{pmatrix}\quad\mbox{and}\quad\bm F(\bm U)=\begin{pmatrix}q\\\dfrac{q^2}{h}+\dfrac{g}{2}h^2\end{pmatrix},
\label{5.1}
\end{equation}
where $h$ is the liquid depth, $u$ is the liquid velocity, $q=hu$ is the discharge, and $g$ is the acceleration due to gravity. In the
numerical examples reported in \S\ref{sec5} and \S\ref{sec7}, we have taken $g=10$.

The eigenvalues of the Jacobian $A=\nicefrac{\partial\mF}{\partial\mU}$ are
\begin{equation*}
\lambda_1(A)=u-c,\quad\lambda_2(A)=u+c,\quad c:=\sqrt{gh}.
\end{equation*}

As mentioned in \S\ref{sec2}, in order to compute the right- and left-sided values of $\mU$ at the cell interfaces, we will use either a
piecewise linear reconstruction (for the CU scheme) or the WENO-Z interpolation using the local characteristic variables, which are obtained
using the local characteristic decomposition (for the A-WENO scheme). In order to implement the local characteristic decomposition, we
proceed as follows. Consider, for instance, the finest mesh with the cells $C^{4N}_\jmh$ and introduce the Roe averages
\begin{equation}
\widehat h_j=\hf\left(h_\jmh+h_\jph\right),\quad
\widehat u_j=\frac{\sqrt{h_\jmh}\,u_\jmh+\sqrt{h_\jph}\,u_\jph}{\sqrt{h_\jmh}+\sqrt{h_\jph}},
\label{5.3}
\end{equation}
where $u_j=q_j/h_j$, $\widehat c_j=\sqrt{g\widehat h_j}$, the Roe matrix $\widehat A_j=A(\widehat\mV_j)$,
$\widehat\mV_j=(\widehat h_j,\widehat q_j)^\top$, and the matrices
\begin{equation*}
R_j=\begin{pmatrix}1&1\\\widehat u_j-\widehat c_j&\widehat u_j+\widehat c_j\end{pmatrix}\quad\mbox{and}\quad
R^{-1}_j=\frac{1}{2\widehat c_j}\begin{pmatrix}\widehat c_j+\widehat u_\jph&-1\\\widehat c_j-\widehat u_\jph&1\end{pmatrix}
\end{equation*}
satisfy $R^{-1}_j\widehat A_jR_j={\rm diag}\big(\lambda_1(\widehat A_j),\lambda_2(\widehat A_j)\big)$. We then introduce the local
characteristic variables in the neighborhood of $x=x_j$:
\begin{equation*}
\bm\Gamma_k=R^{-1}_j\mV_k,\quad k=j-\frac{5}{2},\ldots,j+\frac{5}{2},
\end{equation*}
and apply the WENO-Z interpolation to these six values of $\bm\Gamma$ to obtain the point values $\bm\Gamma^\pm_j$. At the end, we switch
back to the original variables and obtain
\begin{equation*}
\mV^\pm_j=R_j\bm\Gamma^\pm_j.
\end{equation*}
Equipped with the reconstructed values of $h^\pm_j$ and $q^\pm_j$, we estimate the one-sided local speeds of propagation $a^+_j$ and
$a^-_j$ used in \eref{A3} and \eref{C2a}, \eref{C2c} as follows:
\begin{equation}
a^+_j=\max\Big(u^+_j+c^+_j,\,u^-_j+c^-_j\,\Big),\quad a^-_j=\min\Big(u^+_j-c^+_j,\,u^-_j-c^-_j\,\Big),
\label{5.7}
\end{equation}
where $u^\pm_j=q^\pm_j/h^\pm_j$.

We note that for the coarser meshes, the formulae analogous to \eref{5.3}--\eref{5.7} can be obtained in a straightforward way.

\section{Numerical Examples}\label{sec5}
In this section, we compute the experimental convergence rates introduced in \S\ref{sec2} and \S\ref{sec3} using the studied CU, RBM, and
A-WENO schemes. While the RBM scheme is fully discrete, the CU and A-WENO schemes are semi-discrete ones and the corresponding
semi-discretizations result in the time-dependent ODE systems \eref{A2} and \eref{C1}, respectively. We numerically solve these ODE systems
using the three-stage third-order strong stability preserving (SSP) Runge-Kutta solver; see, e.g.,\cite{Gottlieb11,Gottlieb12}. The time
steps for all of the studied schemes are supposed to be selected adaptively based on the corresponding CFL-like stability restrictions. It
is well-known, however, that when the experimental rate of convergence is measured, it is better to use constant time steps, which will be
specified in every reported numerical example. We would also like to stress that in order to achieve the fifth order of accuracy when the
A-WENO scheme is used, we take very small time steps, being proportional to $(\dx)^{5/3}$ rather than to $\dx$ as done when the lower-order
CU and RBM schemes are used.

In all of the examples reported in this section, we use the periodic boundary conditions.

\subsubsection*{Example 1---Test with One Shock}
In the first example originally introduced in \cite{KO14}, we consider the following smooth 10-periodic initial conditions:
\begin{equation}
u(x,0)=2\sin\Big(\frac{\pi x}{5}+\frac{\pi}{4}\Big),\quad h(x,0)=\frac{\big(u(x,0)+10\big)^2}{4g},
\label{6.1}
\end{equation}
which correspond to the following initial values of the invariants $w_1=u-2c$ and $w_2=u+2c$:
$$
w_1(x,0)\equiv-10,\quad w_2(x,0)=2u(x,0)+10.
$$
One can show that the solution of the studied initial value problem \eref{1.1}, \eref{5.1}, \eref{6.1} develops one shock discontinuity per
each period at about $t\approx0.54$.

\begin{figure}[ht!]
\centerline{\includegraphics[trim=3.6cm 1.1cm 2.8cm 0.5cm, clip, width=5.3cm]{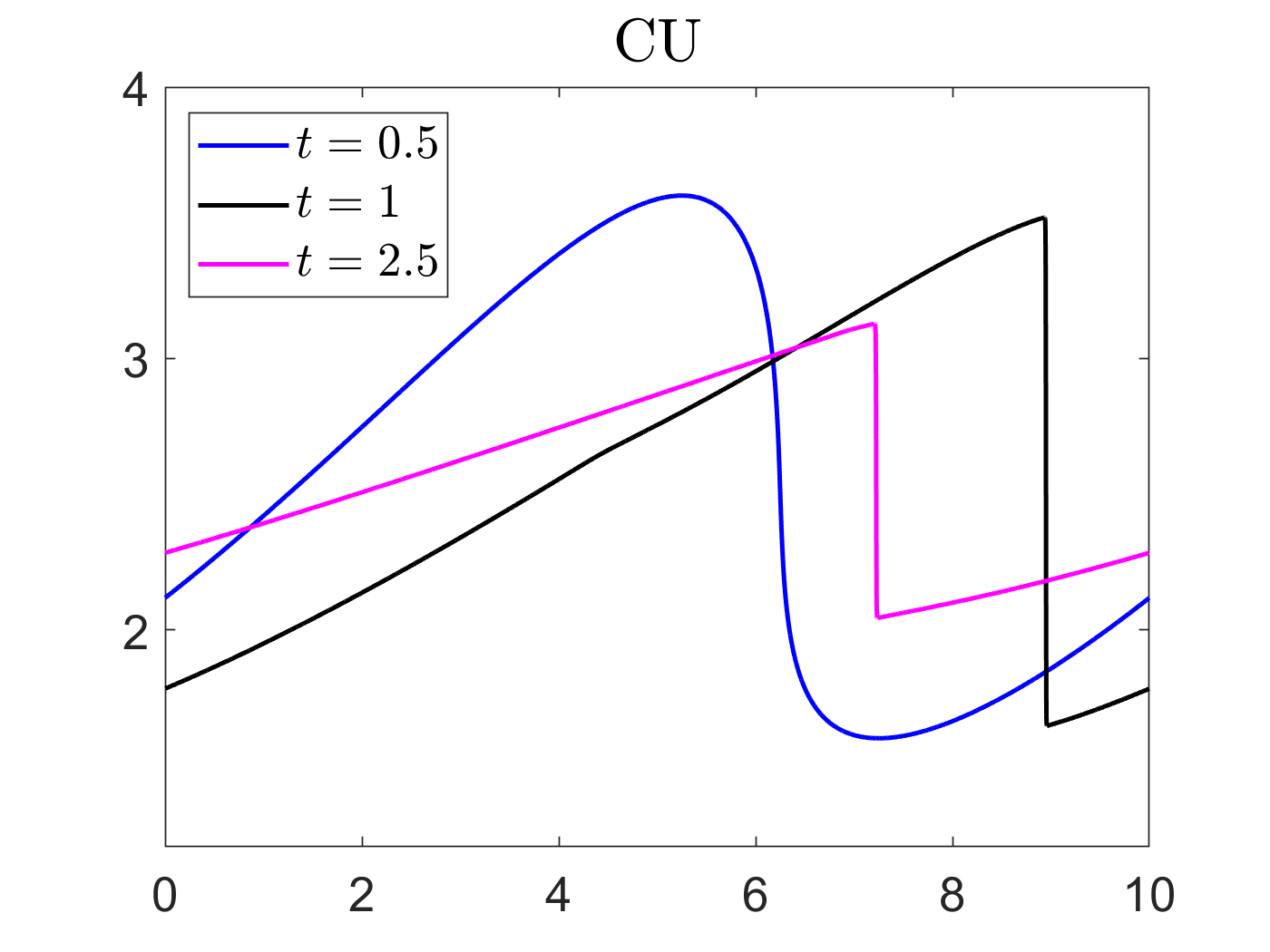}\hspace*{0.5cm}
\includegraphics[trim=3.6cm 1.1cm 2.8cm 0.5cm, clip, width=5.3cm]{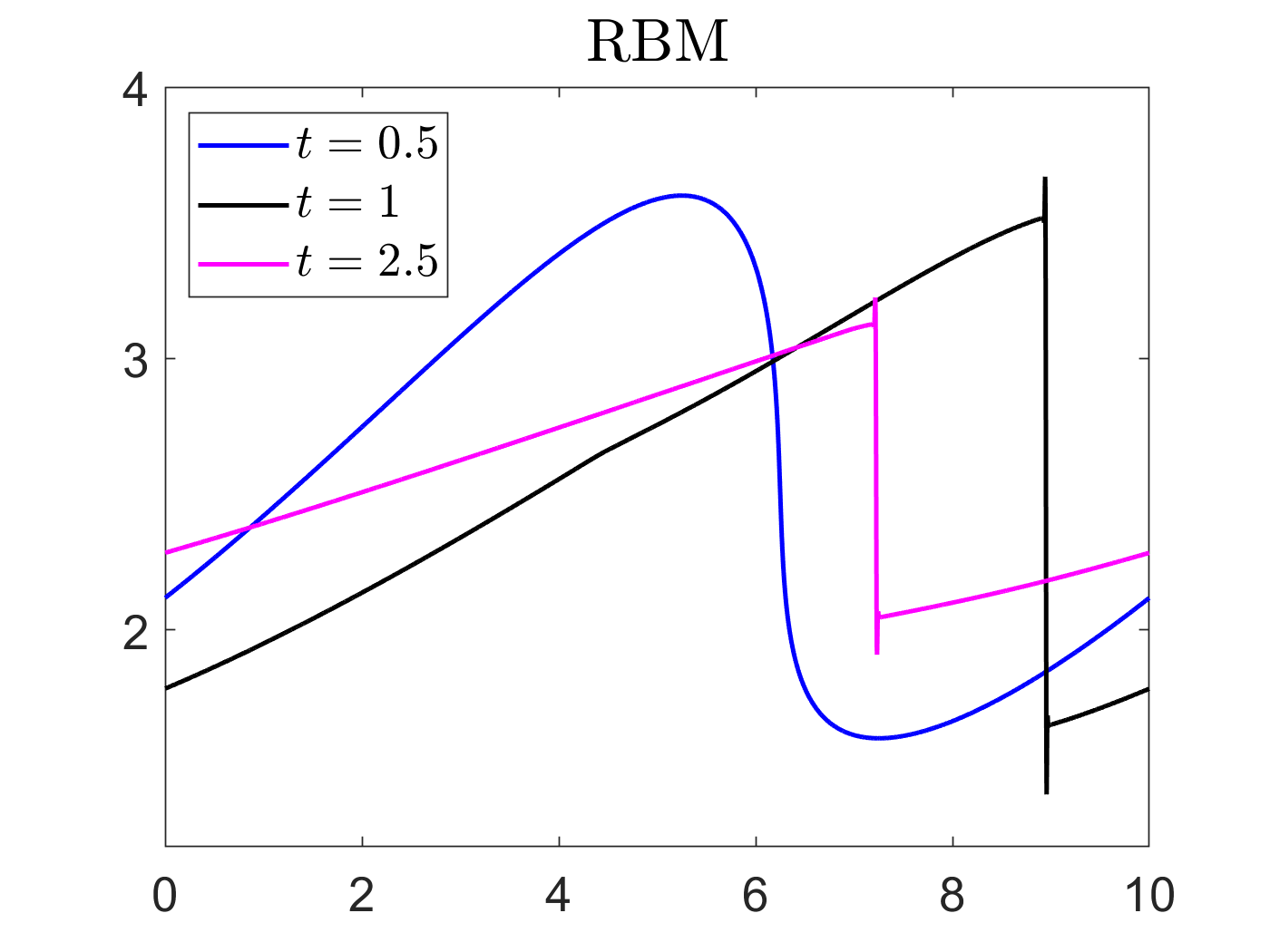}\hspace*{0.5cm}
\includegraphics[trim=3.6cm 1.1cm 2.8cm 0.5cm, clip, width=5.3cm]{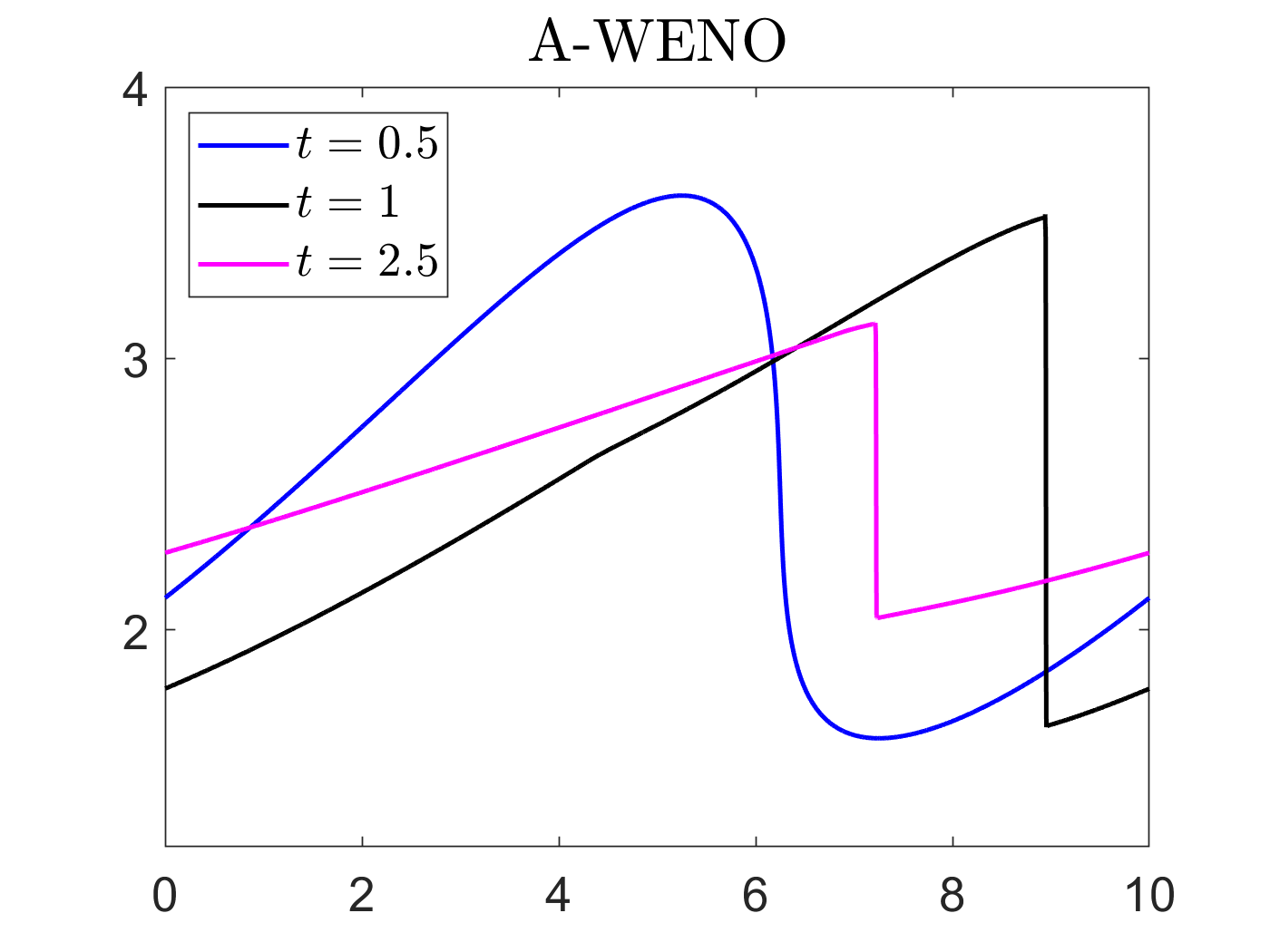}}
\caption{\sf Example 1: Water depth $h$ computed by the CU (left), RBM (middle), and A-WENO (right) schemes.\label{Fig4.1}}
\end{figure}

We compute the solutions by the studied CU, RBM, and A-WENO schemes at times $t=0.5$, 1, and 2.5 on the computational domain $[0,10]$ using
1000, 2000, 4000, and 8000 uniform cells. The water depths $h$ computed using 4000 uniform cells are presented in Figure \ref{Fig4.1}. As
one can see, all of the three studied schemes can capture the shock position correctly. At the same time, one can notice that there are
${\cal O}(1)$ oscillations in the immediate vicinity of the shock in the RBM solution.

\begin{figure}[ht!]
\centerline{\includegraphics[trim=0.9cm 0.4cm 1.1cm 0.2cm, clip, width=5.4cm]{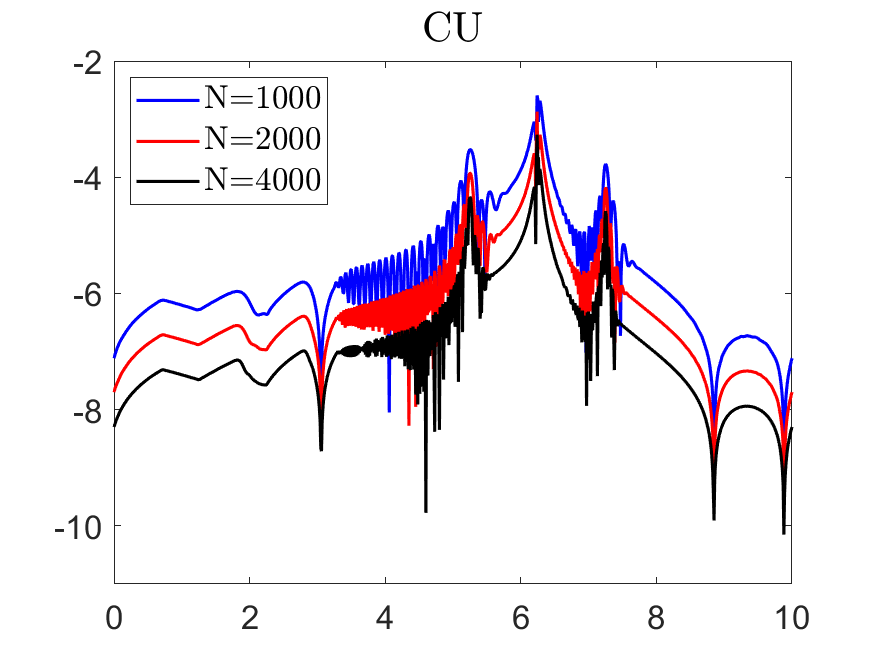}\hspace*{0.5cm}
\includegraphics[trim=0.9cm 0.4cm 1.1cm 0.2cm, clip, width=5.4cm]{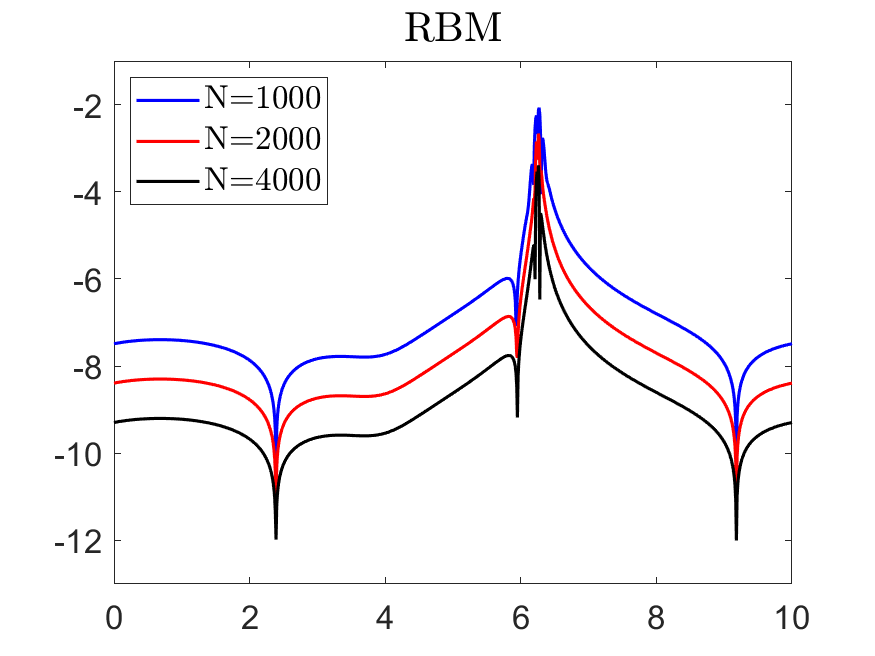}\hspace*{0.5cm}
\includegraphics[trim=0.9cm 0.4cm 1.1cm 0.2cm, clip, width=5.4cm]{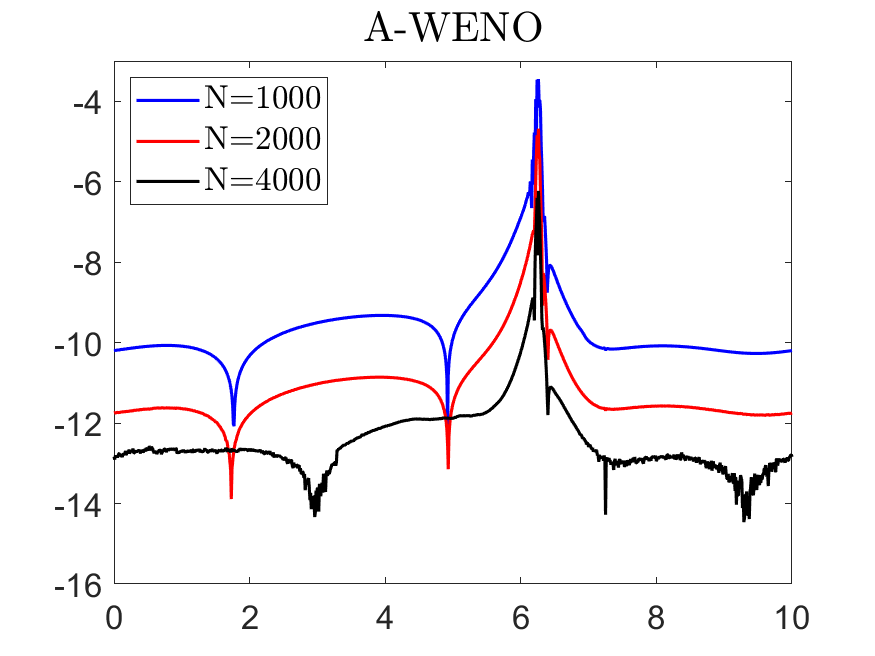}}
\caption{\sf Example 1: $\log_{10}|\frac{h^N-h^{\rm ref}}{h^{\rm ref}}|$, where $h^{\rm ref}$ is the reference solution  for the CU (left), RBM (middle) and A-WENO
(right) schemes at $t=0.5$.\label{Fig4.2}}
\end{figure}

\begin{figure}[ht!]
\centerline{\includegraphics[trim=0.9cm 0.4cm 1.1cm 0.2cm, clip, width=5.4cm]{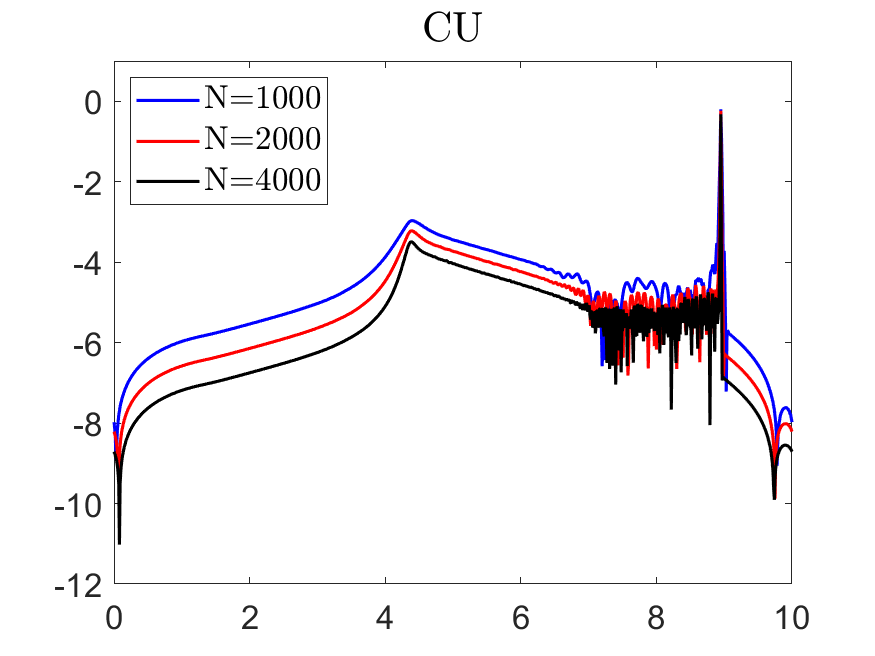}\hspace*{0.5cm}
\includegraphics[trim=0.9cm 0.4cm 1.1cm 0.2cm, clip, width=5.4cm]{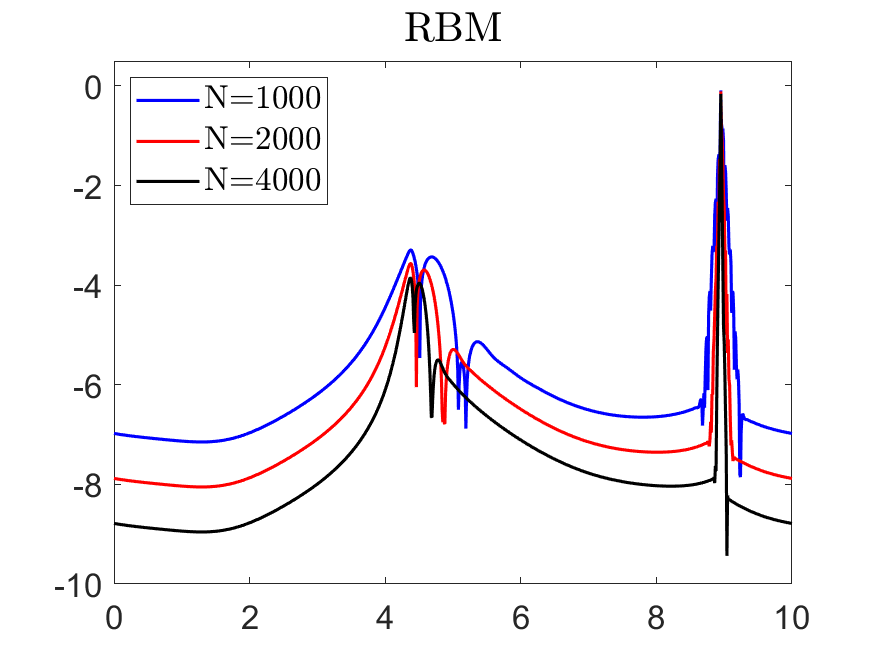}\hspace*{0.5cm}
\includegraphics[trim=0.9cm 0.4cm 1.1cm 0.2cm, clip, width=5.4cm]{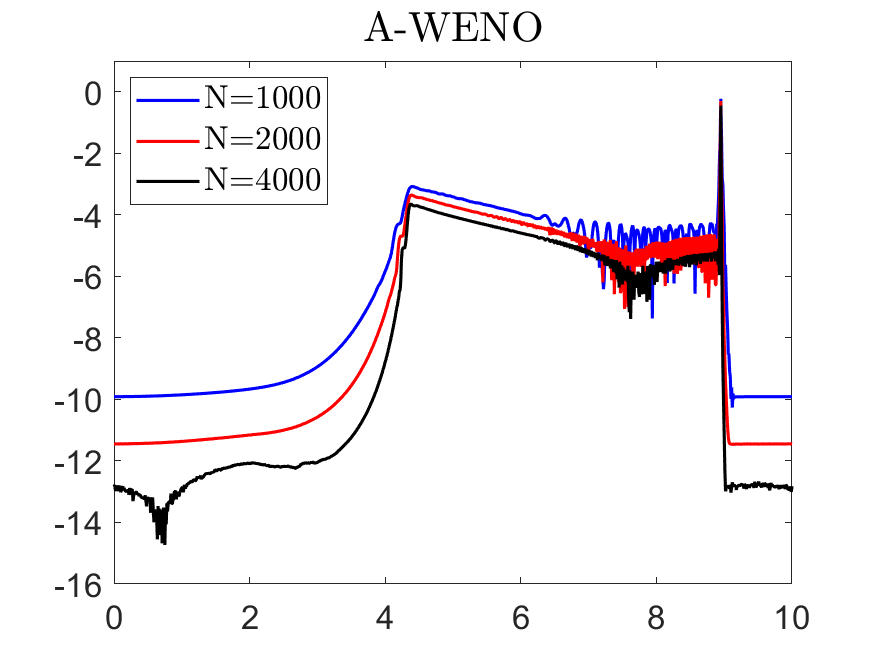}}
\caption{\sf Example 1: Same as in Figure \ref{Fig4.2}, but for $t=1$.\label{Fig4.3}}
\end{figure}

\begin{figure}[ht!]
\centerline{\includegraphics[trim=0.9cm 0.4cm 1.1cm 0.2cm, clip, width=5.4cm]{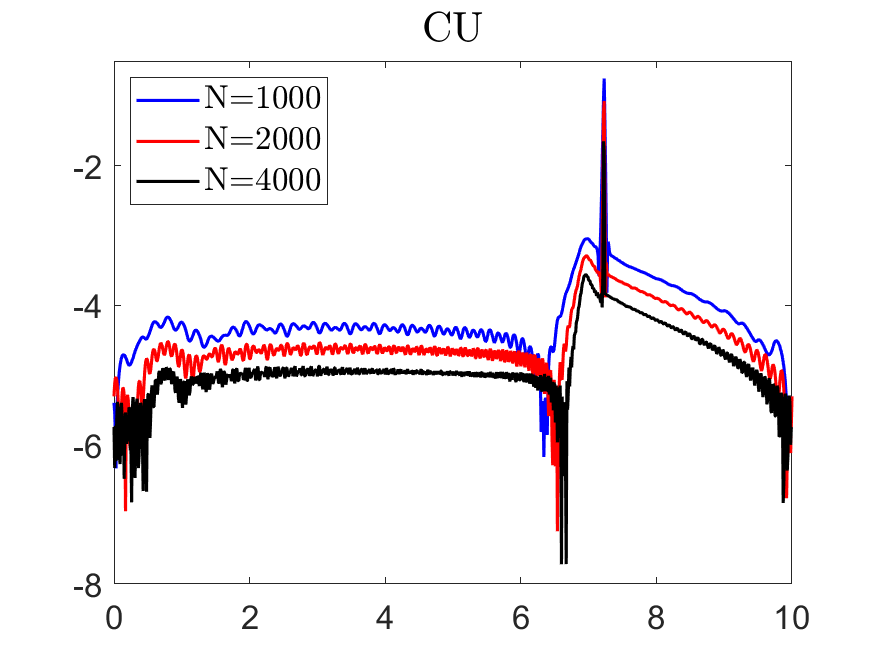}\hspace*{0.5cm}
\includegraphics[trim=0.9cm 0.4cm 1.1cm 0.2cm, clip, width=5.4cm]{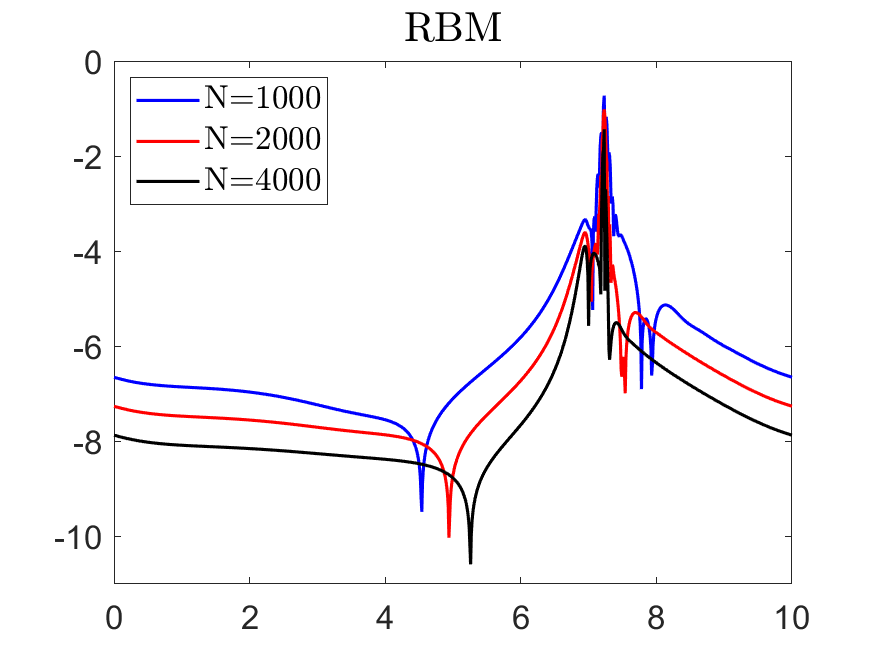}\hspace*{0.5cm}
\includegraphics[trim=0.9cm 0.4cm 1.1cm 0.2cm, clip, width=5.4cm]{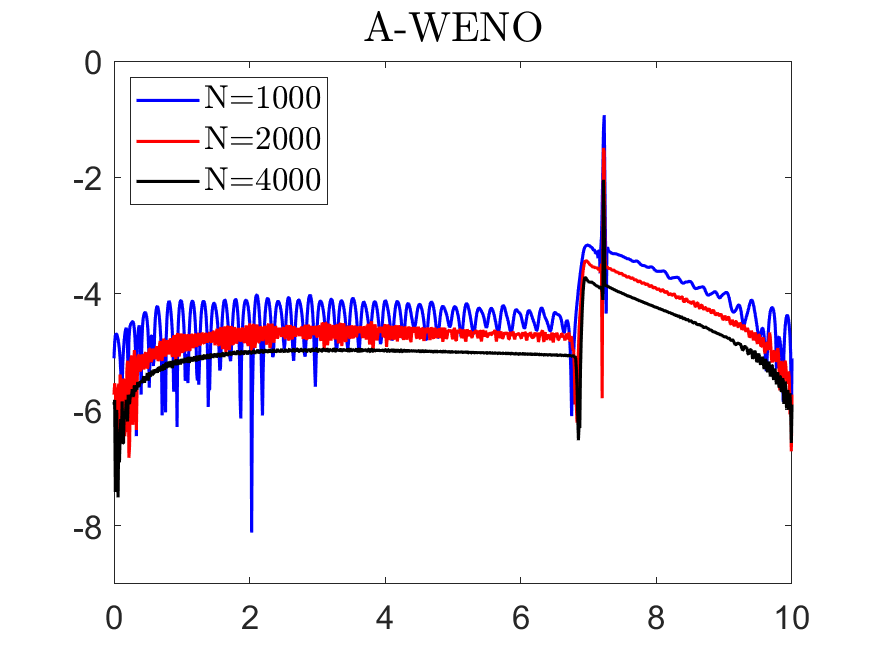}}
\caption{\sf Example 1: Same as in Figure \ref{Fig4.2}, but for $t=2.5$.\label{Fig4.4}}
\end{figure}

We now study the error in the above computations. First, we compute the corresponding reference solutions by the same schemes but using a
very fine mesh with 40000 uniform cells, and present the relative differences between the computed and reference water depths measured in
the logarithmic scale in Figures \ref{Fig4.2}--\ref{Fig4.4}. As one can see from Figure \ref{Fig4.1}, at $t=0.5$ the water depth is smooth
but it contains a high gradient area around $x=6.2$. The corresponding errors plotted in Figure \ref{Fig4.2} are small and decay when the
mesh is refined, that is, when $N$ increases from 1000 to 4000. The decay, however, is monotone only for the RBM and A-WENO schemes, while
the CU errors exhibit somewhat oscillatory behavior near the local extrema ($x\approx5.3$ and $x\approx7.1$), where the order of the
generalized minmod reconstruction utilized in the CU scheme reduces to the first one (this occurs due to a well-known clipping of extrema
phenomenon). One can also observe the fastest convergence of the fifth-order A-WENO scheme as expected, and also the saturation  of errors
when the calculations start getting into the range of machine errors.

The shock is formed at $t\approx0.54$ at about $x\approx6.7$ and then propagates to $x\approx8.9$ by the time $t=1$. From Figure
\ref{Fig4.3}, one can see that there are ${\cal O}(1)$ errors around the shock, and these errors propagate and affect the area behind the
shock for the CU and A-WENO schemes. One can also observe an oscillatory behavior of the errors in this area.

At the final time $t=2.5$, the formed shocks affects the entire computational domain. For the CU and especially A-WENO schemes, the error
behavior is very oscillatory, while the RBM errors still decay in quite monotone way and same low frequency oscillations are contained
within the area near the shock; see Figure \ref{Fig4.4}.

\begin{figure}[ht!]
\centerline{\includegraphics[trim=1.7cm 0.6cm 1.7cm 0.3cm, clip, width=5.5cm]{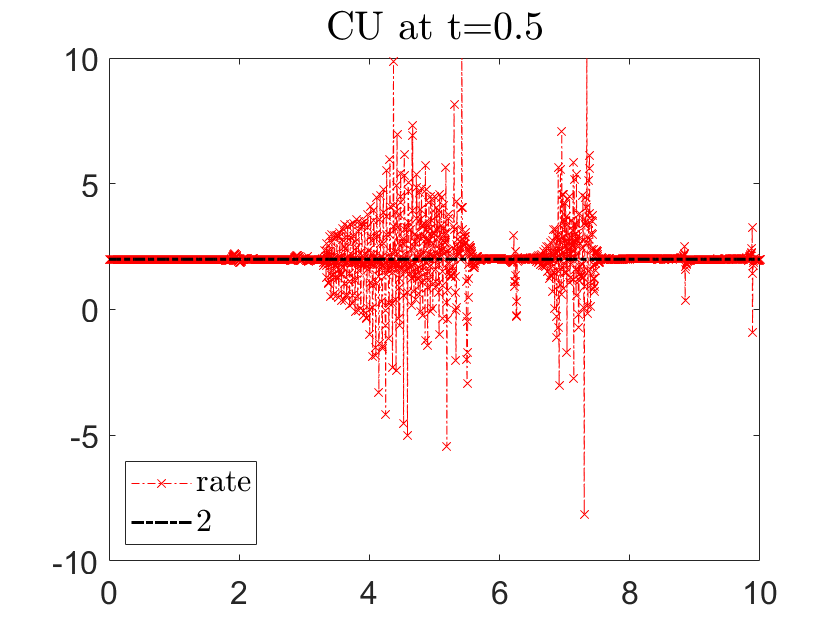}\hspace*{0.3cm}
\includegraphics[trim=1.7cm 0.6cm 1.7cm 0.3cm, clip, width=5.5cm]{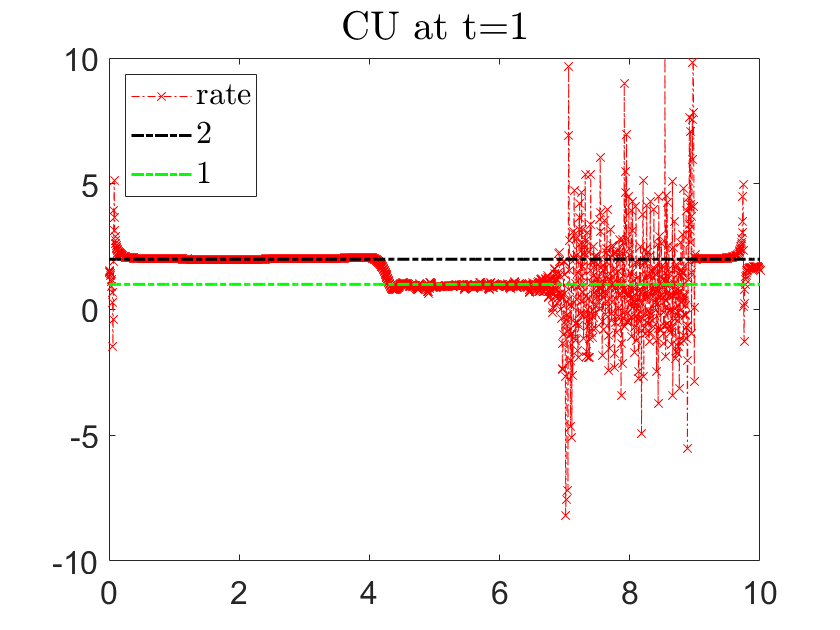}\hspace*{0.3cm}
\includegraphics[trim=1.7cm 0.6cm 1.7cm 0.3cm, clip, width=5.5cm]{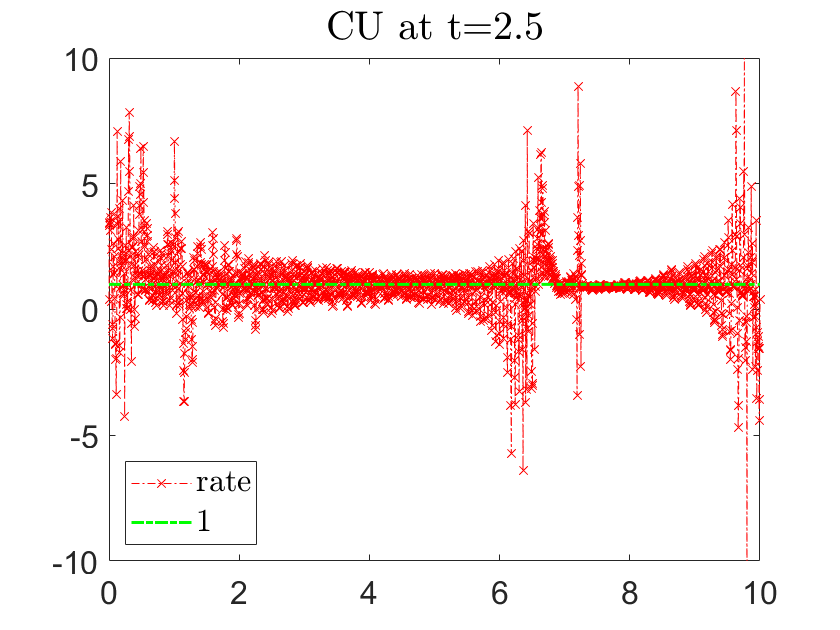}}
\vskip8pt
\centerline{\includegraphics[trim=1.7cm 0.6cm 1.7cm 0.3cm, clip, width=5.8cm]{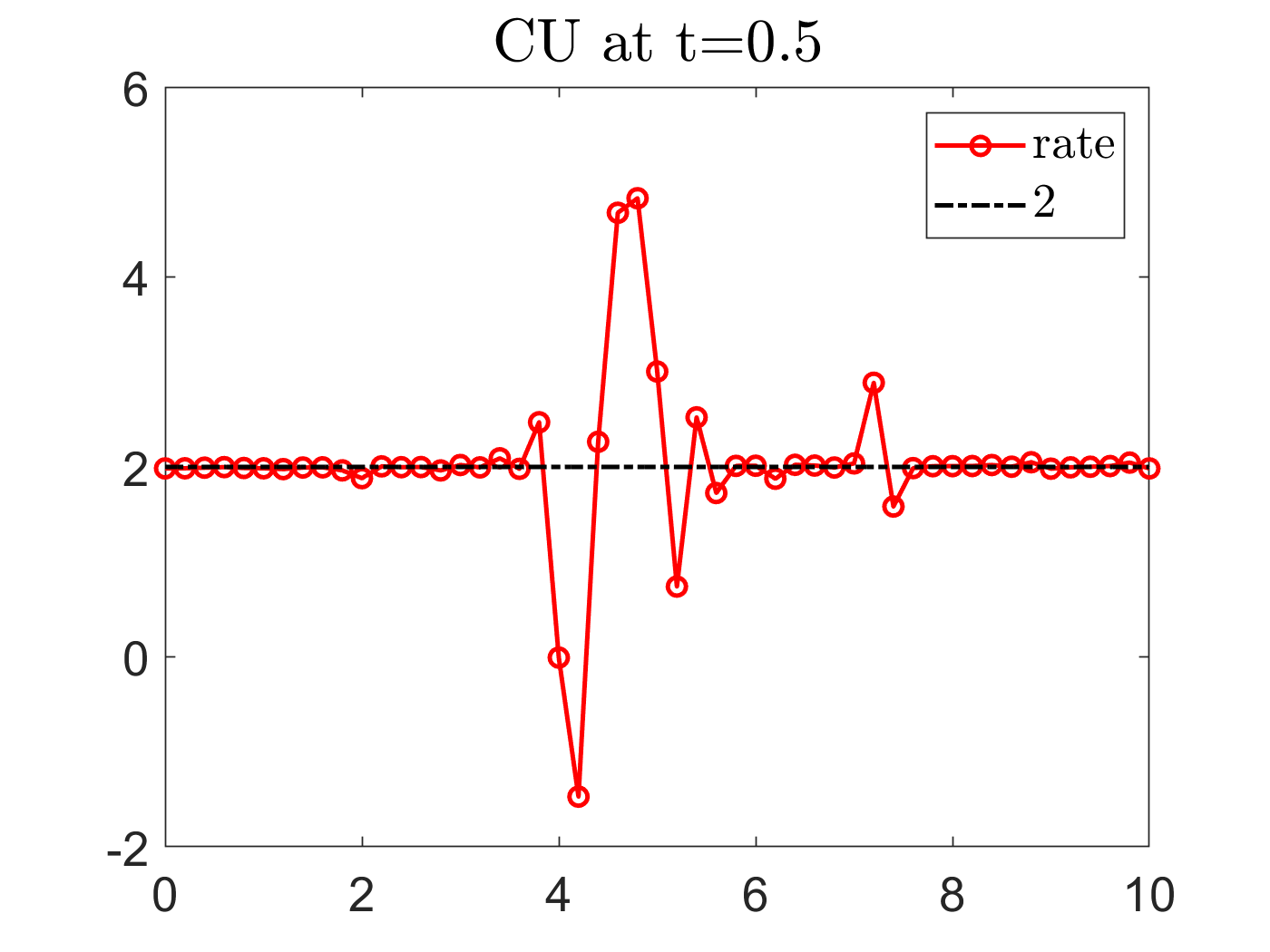}\hspace*{0.0cm}
\includegraphics[trim=1.7cm 0.6cm 1.7cm 0.3cm, clip, width=5.8cm]{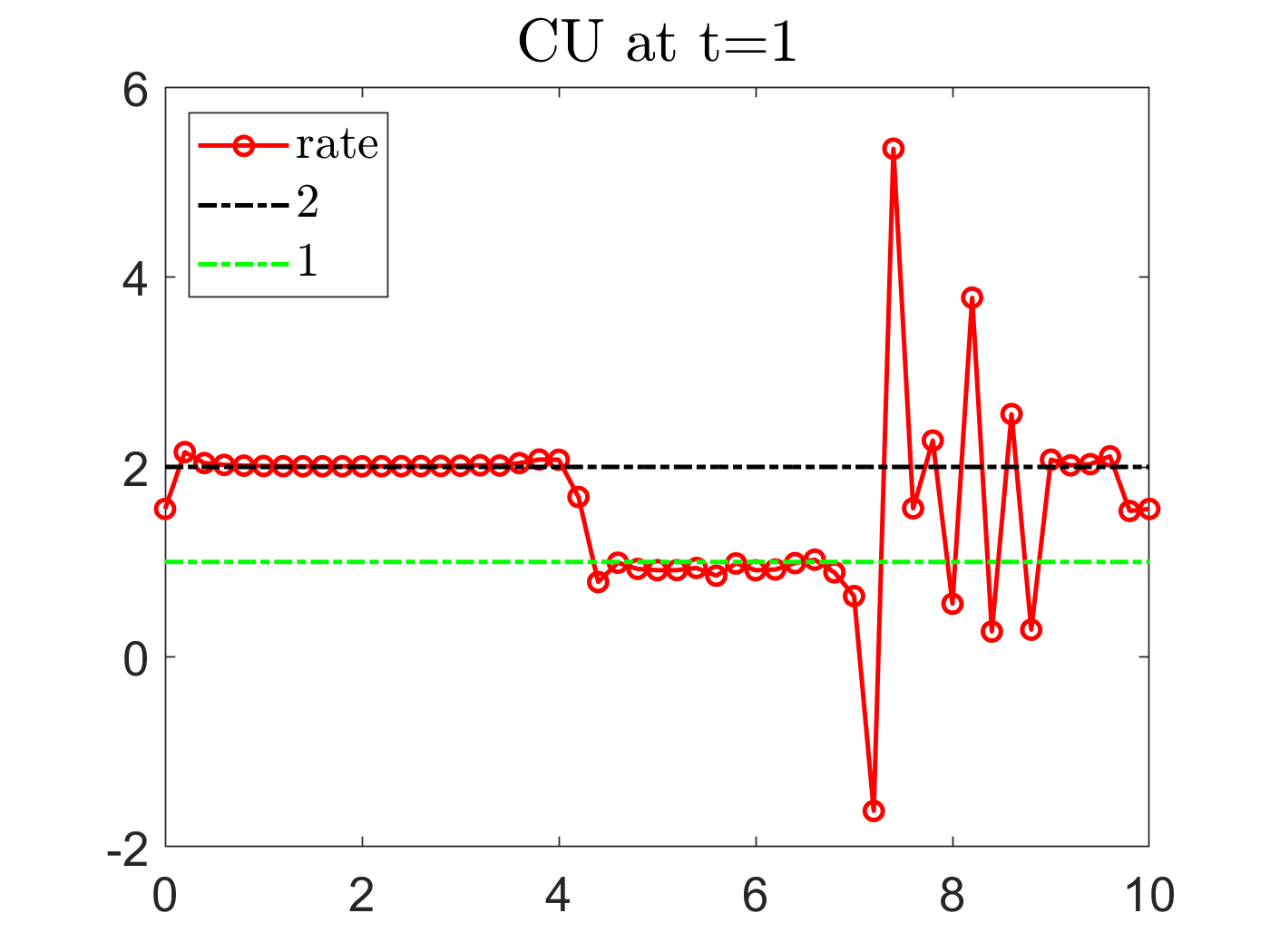}\hspace*{0.0cm}
\includegraphics[trim=1.7cm 0.6cm 1.7cm 0.3cm, clip, width=5.8cm]{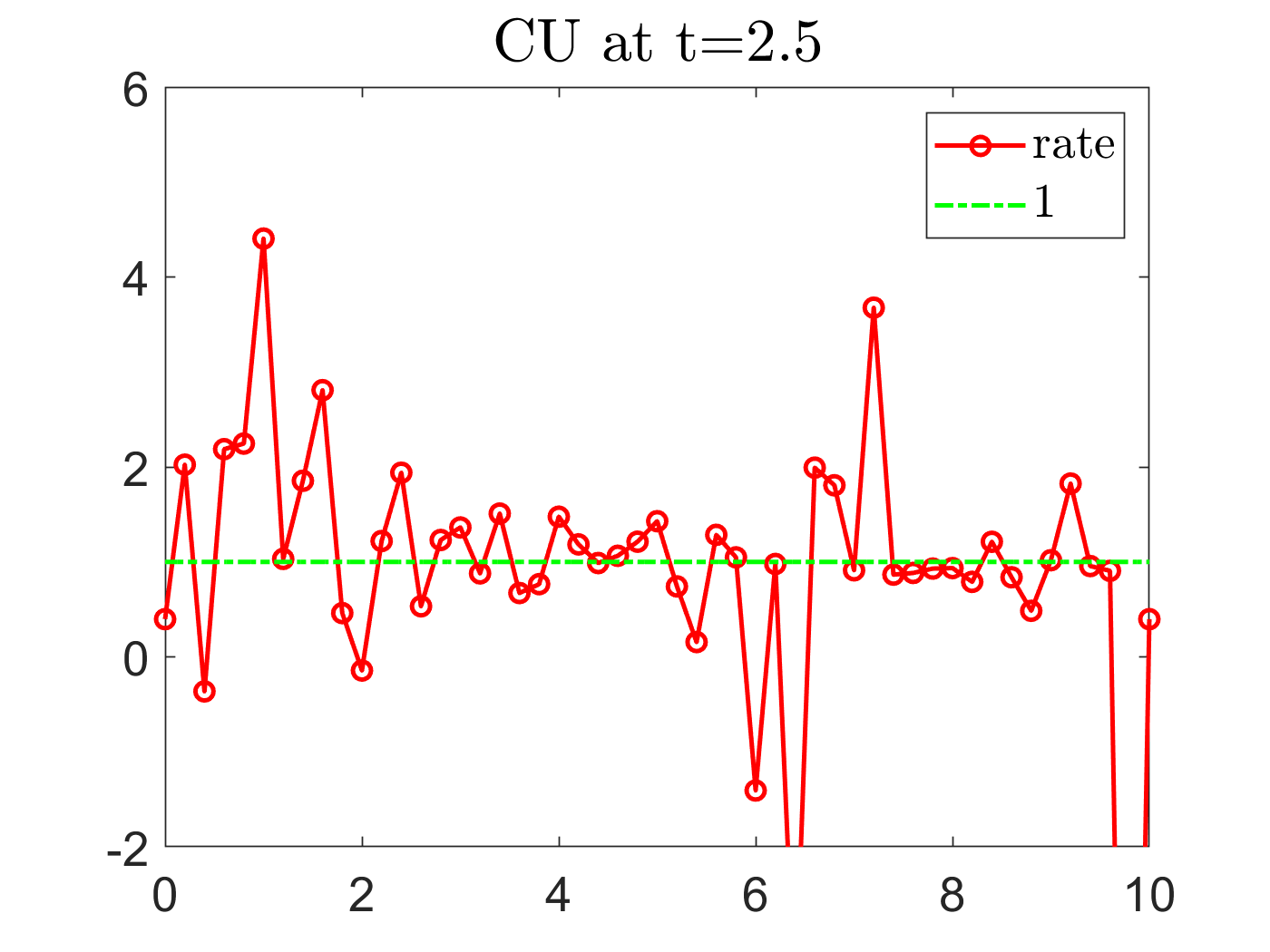}}
\caption{\sf Example 1: Experimental rates of pointwise convergence for the CU scheme: $r_{4j}$ for $j=0,\ldots,2000$ (top row) and
$r_{160k}$ for $k=0,\ldots,50$ (bottom row).\label{Fig4.4a}}
\end{figure}

\begin{figure}[ht!]
\centerline{\includegraphics[trim=1.7cm 0.6cm 1.7cm 0.3cm, clip, width=5.8cm]{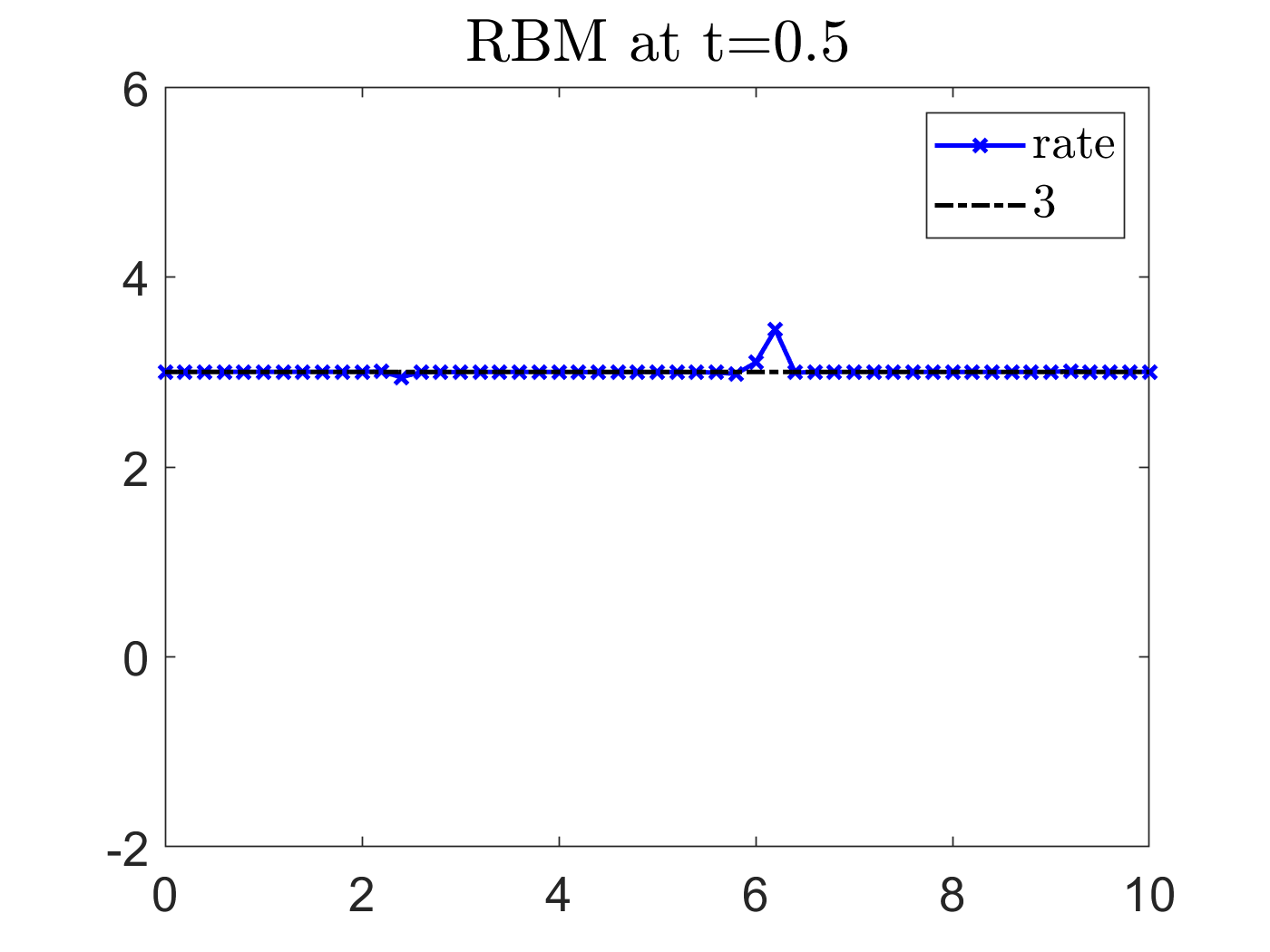}\hspace*{0.0cm}
\includegraphics[trim=1.7cm 0.6cm 1.7cm 0.3cm, clip, width=5.8cm]{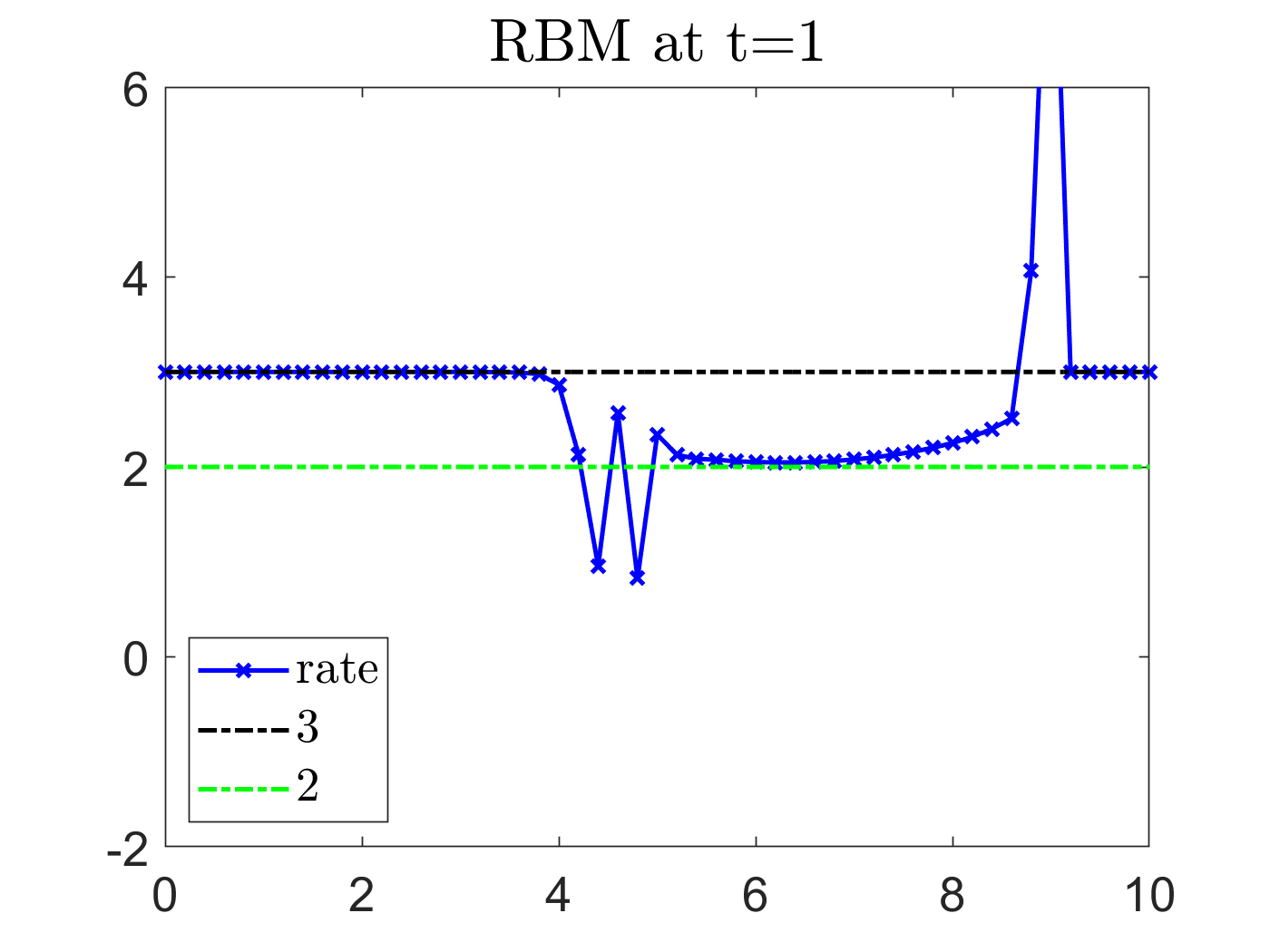}\hspace*{0.0cm}
\includegraphics[trim=1.7cm 0.6cm 1.7cm 0.3cm, clip, width=5.8cm]{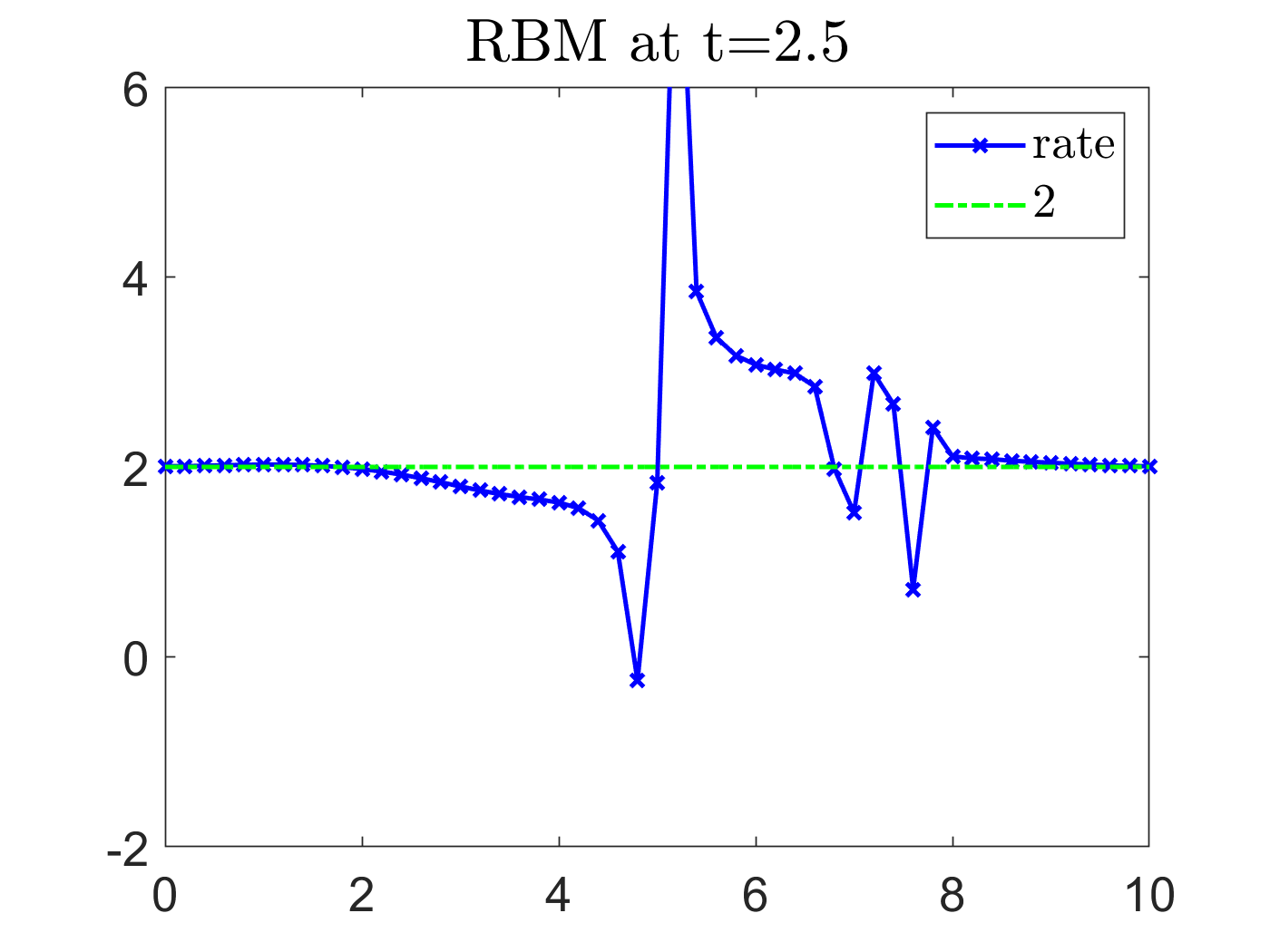}}
\vskip8pt
\centerline{\includegraphics[trim=1.7cm 0.6cm 1.7cm 0.3cm, clip, width=5.8cm]{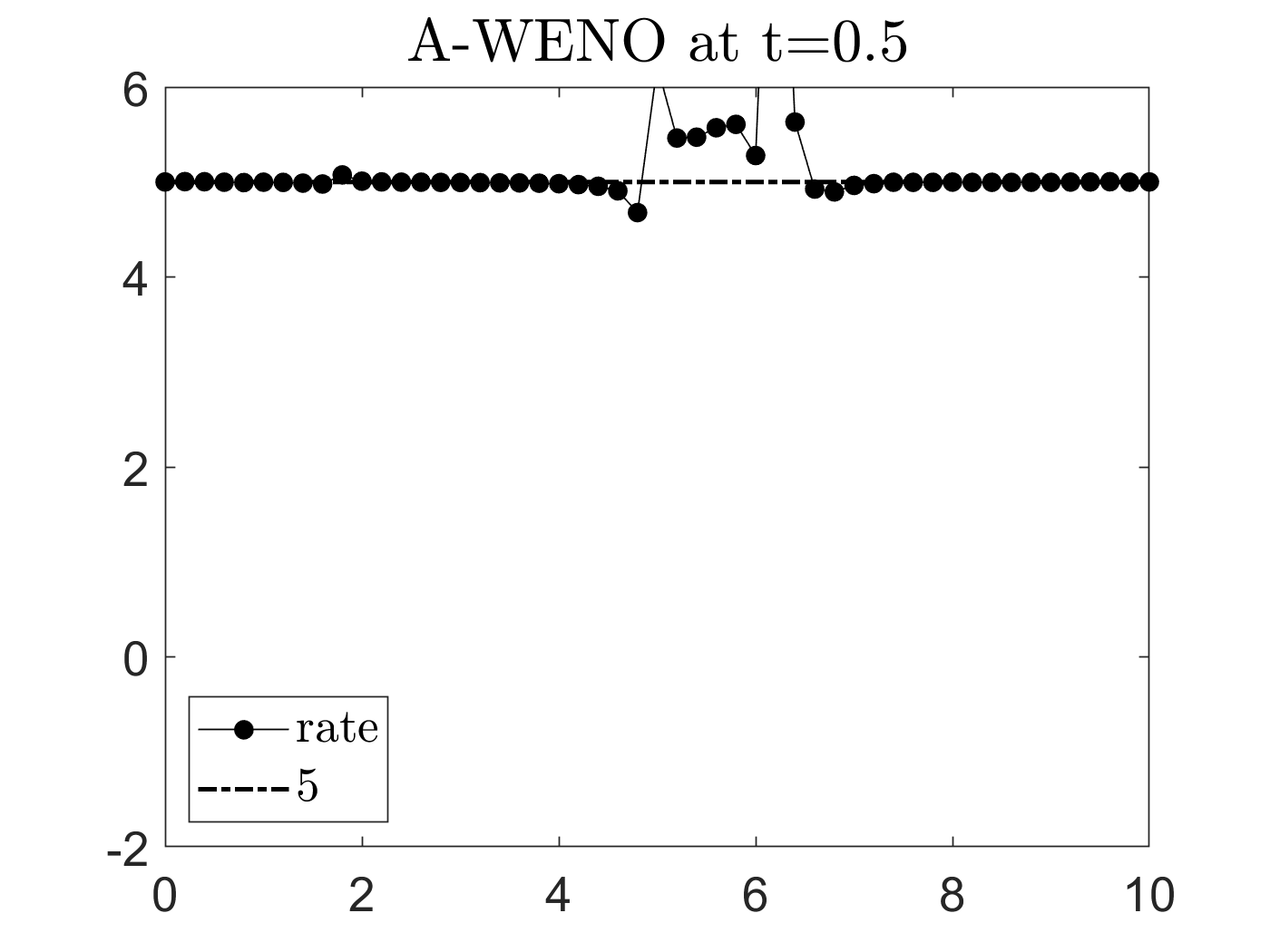}\hspace*{0.0cm}
\includegraphics[trim=1.7cm 0.6cm 1.7cm 0.3cm, clip, width=5.8cm]{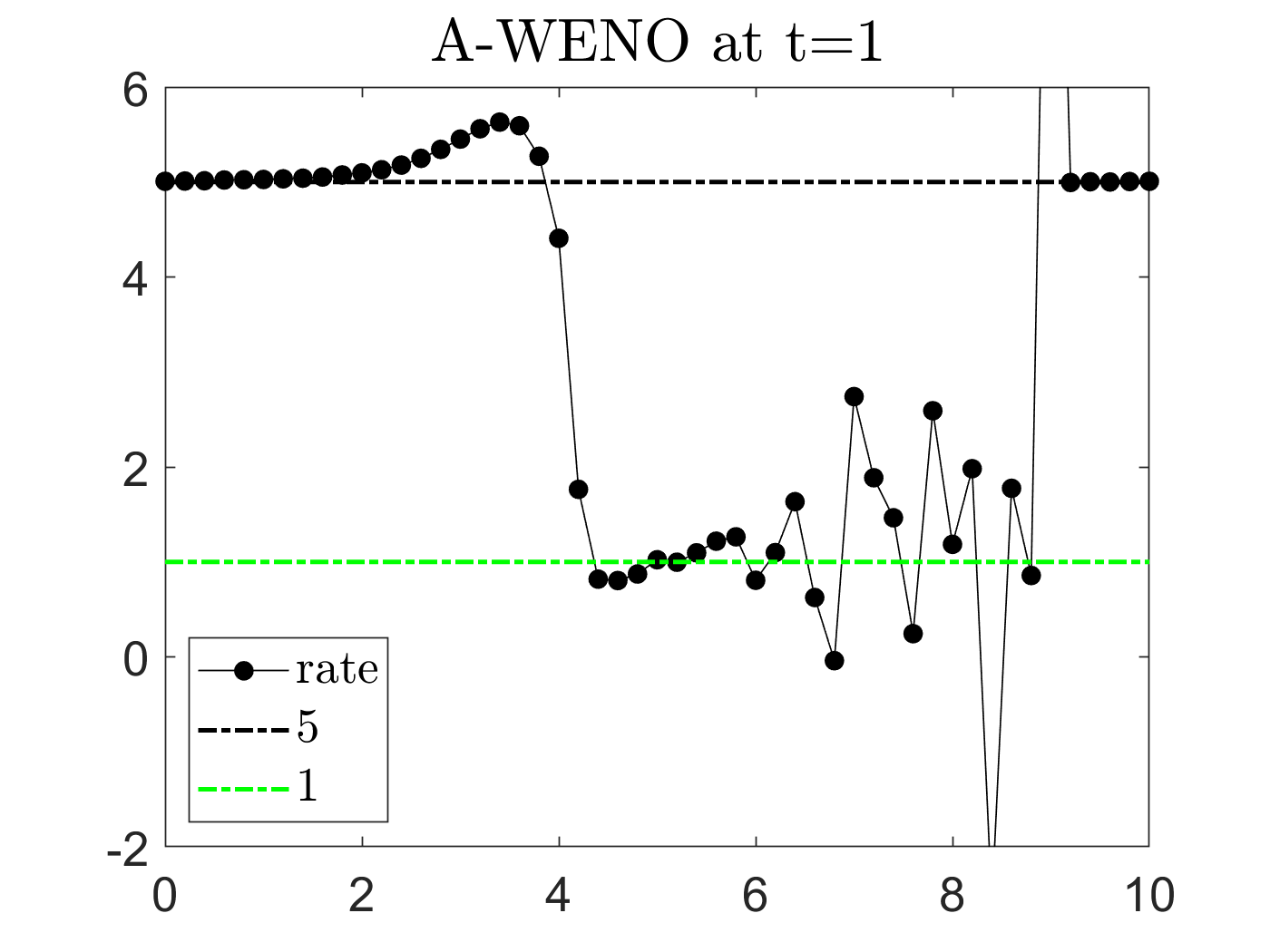}\hspace*{0.0cm}
\includegraphics[trim=1.7cm 0.6cm 1.7cm 0.3cm, clip, width=5.8cm]{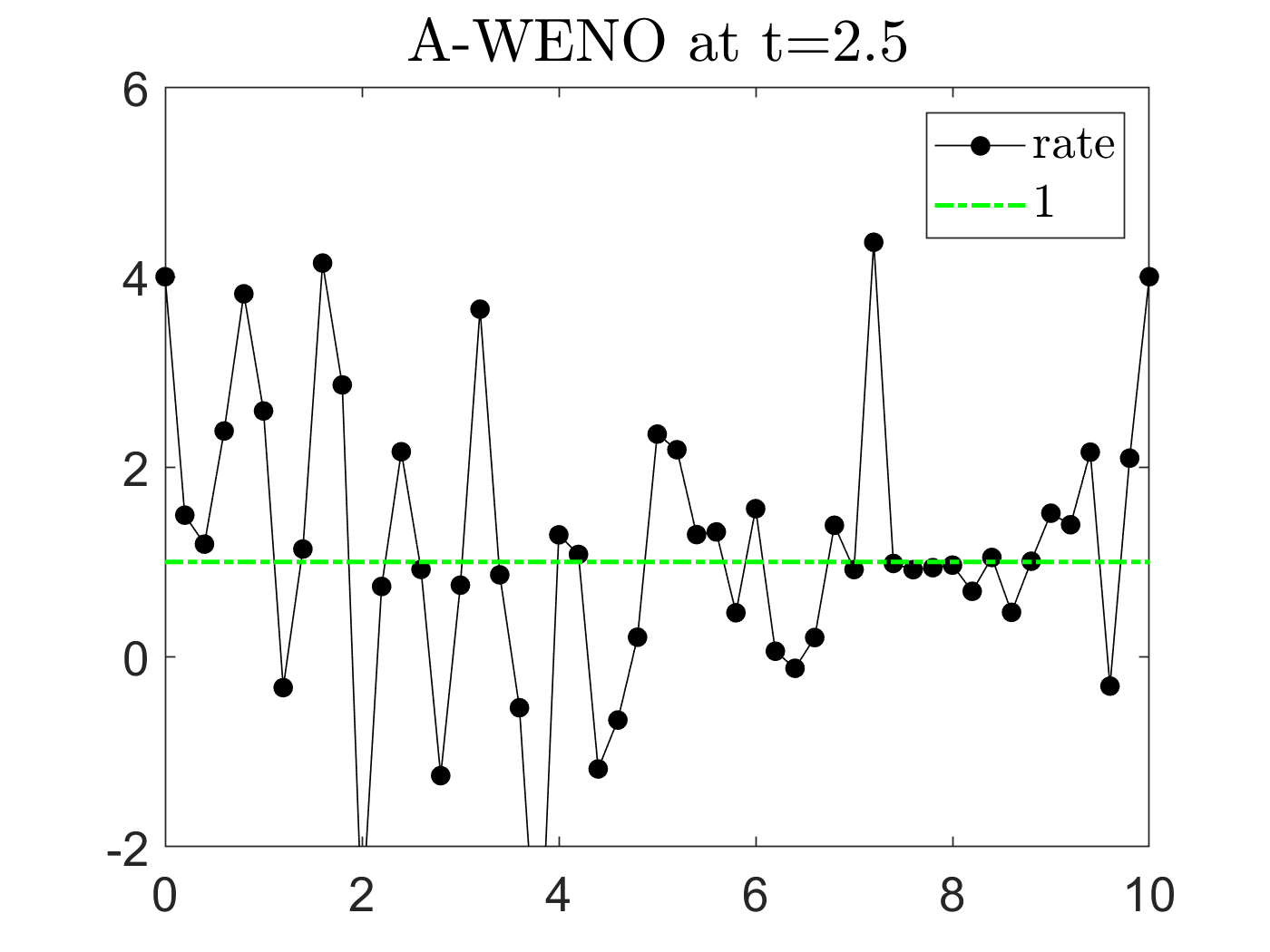}}
\caption{\sf Example 1: Experimental rates of pointwise convergence $r_{160k}$ for $k=0,\ldots,50$ for the RBM scheme (top row) and
$r_{80k}$ for $k=0,\ldots,50$ for the A-WENO scheme (bottom row).\label{Fig4.4aa}}
\end{figure}

In order to better quantify the pointwise convergence, we compute the experimental convergence rates \eref{3.1} for the water depth
computed by the studied schemes. To this end, we use three imbedded grids with $N=2000$ for the CU and RBM schemes and $N=1000$ for the
A-WENO scheme. Note that we use slightly coarser grids in the A-WENO computations as this is a fifth-order scheme and the roundoff errors
start affecting the rate calculations when the mesh is too fine. In Figure \ref{Fig4.4a}, we present the results obtained for the CU scheme.
In the top row, we plot $r_{4j}$ for $j=0,\ldots,2000$ and one can see that even when the solution is still smooth (at $t=0.5$) the rates
are enormously oscillatory and it is clear that the convergence is not monotone even away from the high gradient area. In order to better
see the convergence pattern, we plot in the bottom row of Figure \ref{Fig4.4a} the same quantities, but at every 40th grid point, that is,
we plot $r_{160k}$ for $k=0,\ldots,50$ there. As one can see, the pointwise convergence pattern for the CU scheme is more clear now: it
seems to be second order at $t=0.5$, first order in the zone behind the shock propagation at $t=1$ and first order throughout the entire
computational domain at $t=2.5$. We also plot $r_{160k}$ for the RBM scheme and $r_{80k}$ for the A-WENO scheme ($k=0,\ldots,50$) in Figure
\ref{Fig4.4aa}, where the following patterns can be observed: The corresponding rates for the RBM and A-WENO schemes are close to third and
fifth orders at $t=0.5$, second and first orders in the zone behind the shock at $t=1$, and second and first orders throughout the entire
computational domain at $t=2.5$.

\begin{figure}[ht!]
\centerline{\includegraphics[trim=1.7cm 0.6cm 1.7cm 0.3cm, clip, width=5.8cm]{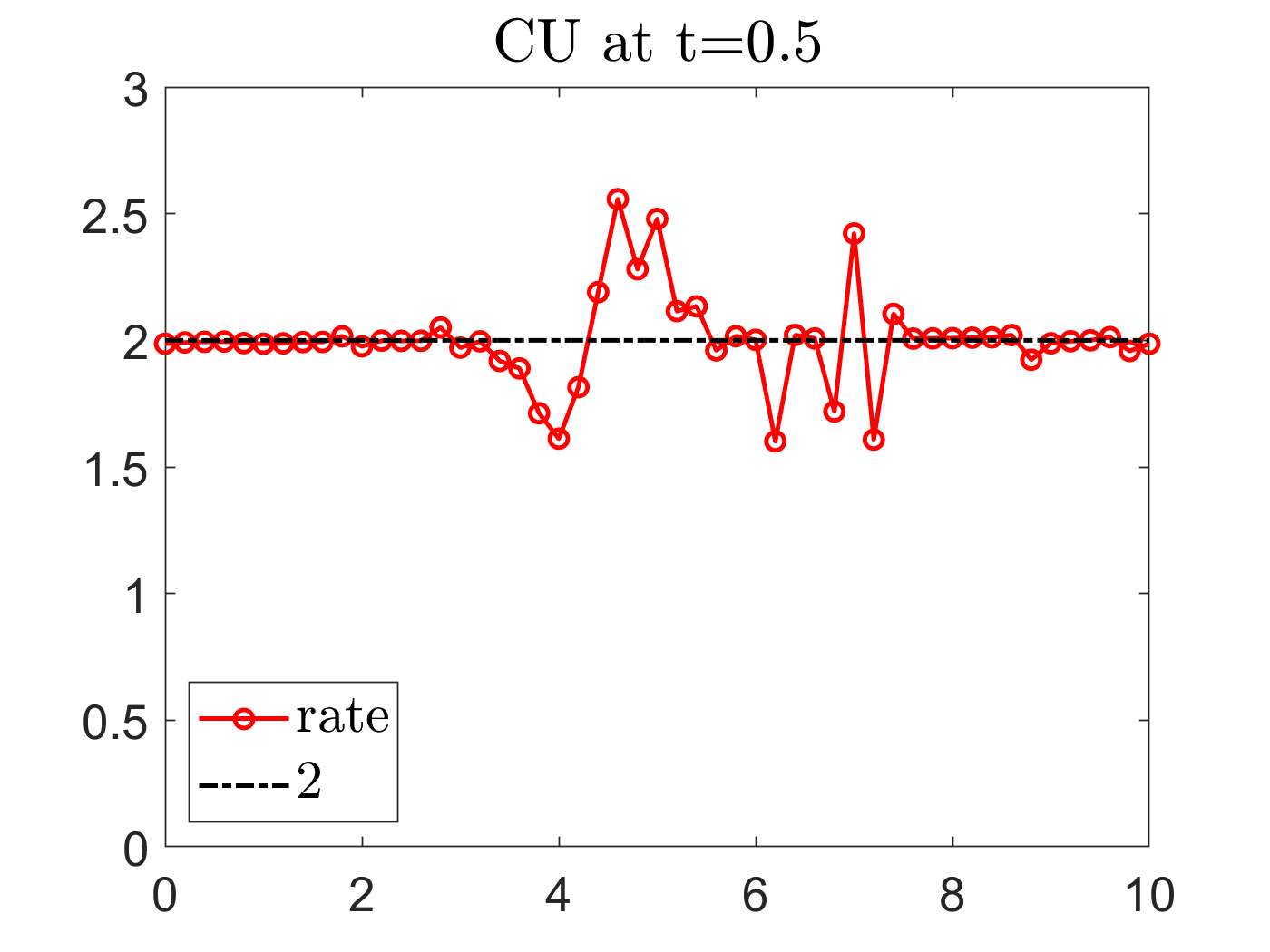}\hspace*{0.0cm}
\includegraphics[trim=1.7cm 0.6cm 1.7cm 0.3cm, clip, width=5.8cm]{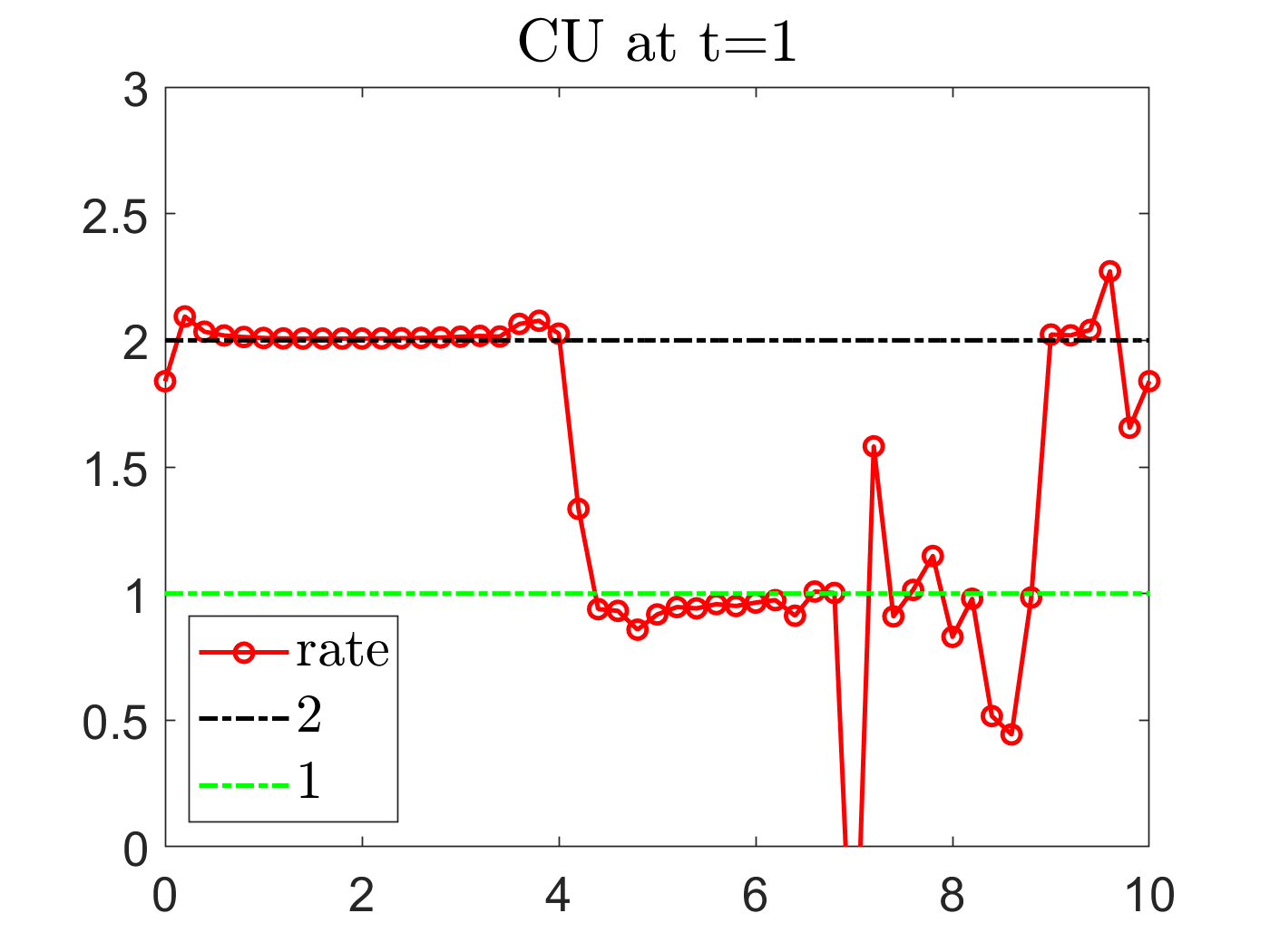}\hspace*{0.0cm}
\includegraphics[trim=1.7cm 0.6cm 1.7cm 0.3cm, clip, width=5.8cm]{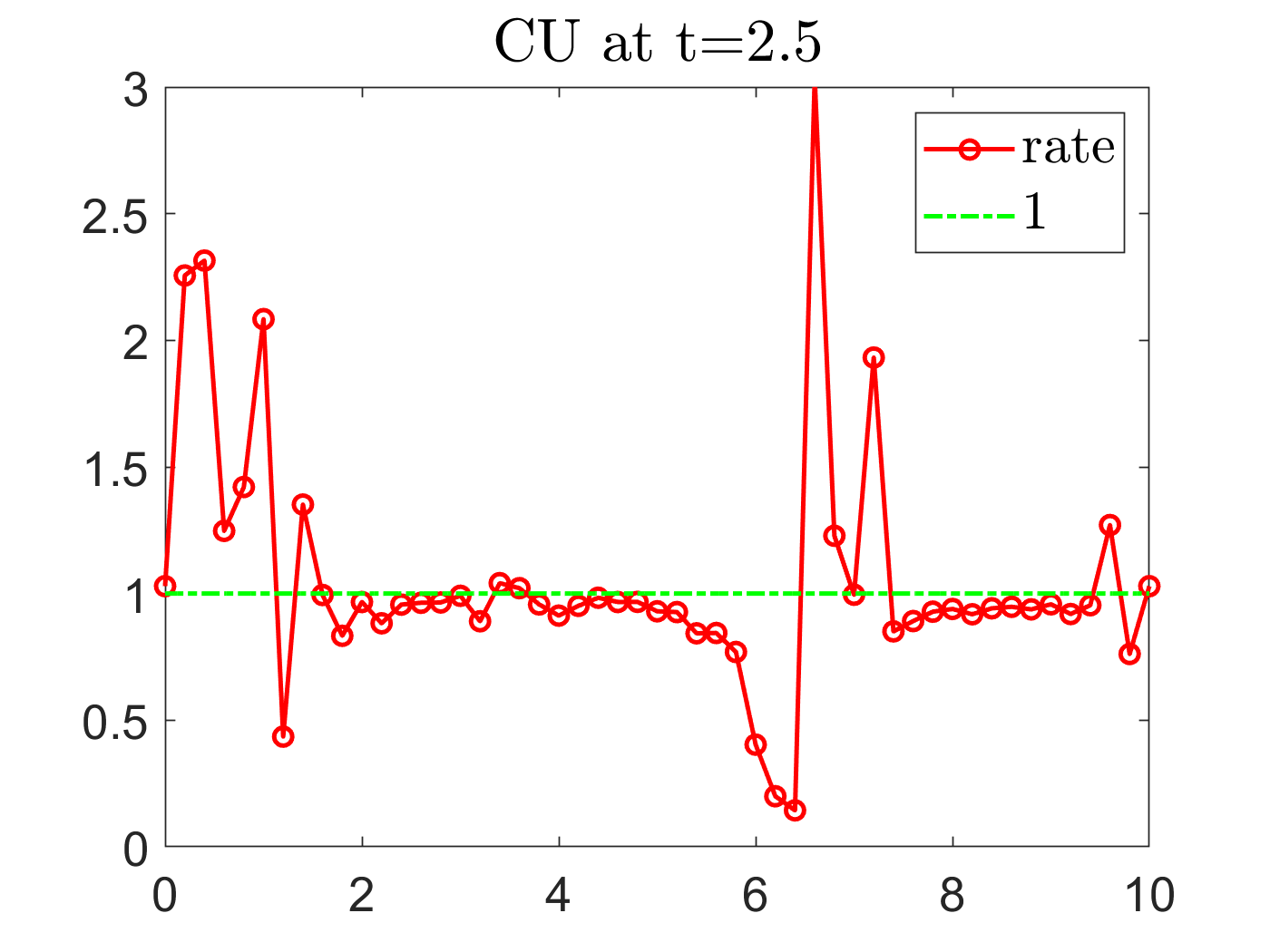}}
\vskip8pt
\centerline{\includegraphics[trim=1.7cm 0.6cm 1.7cm 0.3cm, clip, width=5.8cm]{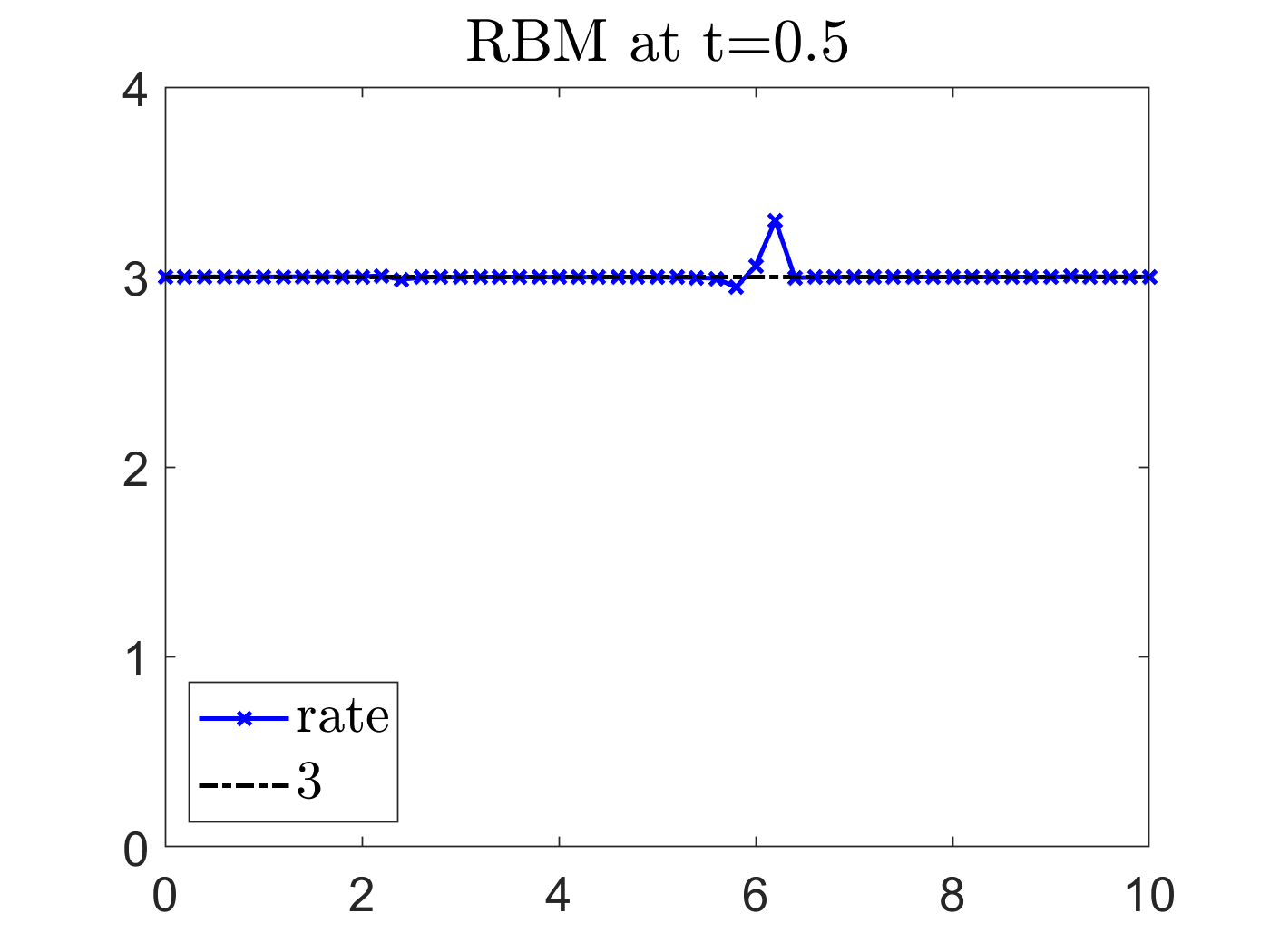}\hspace*{0.0cm}
\includegraphics[trim=1.7cm 0.6cm 1.7cm 0.3cm, clip, width=5.8cm]{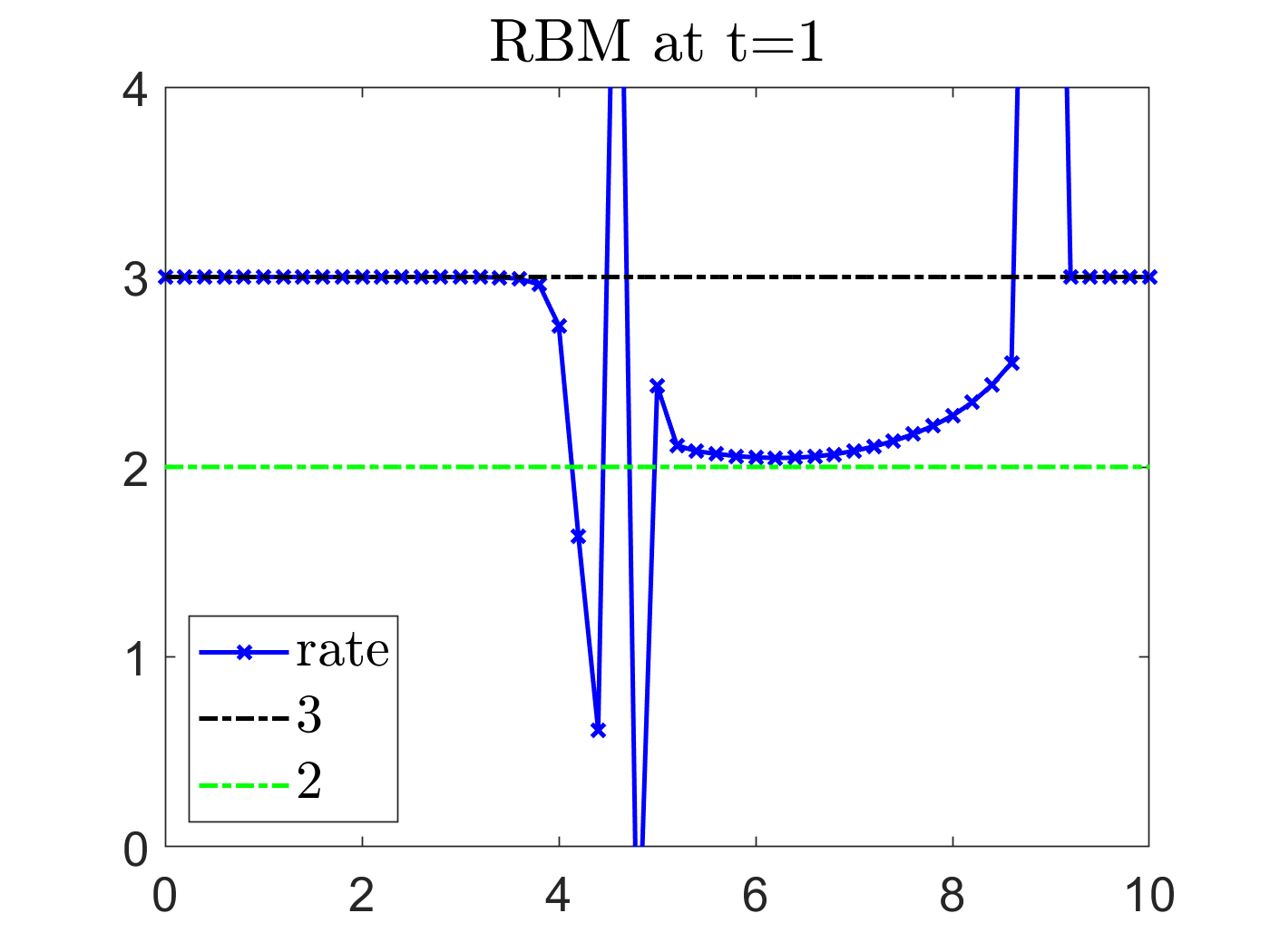}\hspace*{0.0cm}
\includegraphics[trim=1.7cm 0.6cm 1.7cm 0.3cm, clip, width=5.8cm]{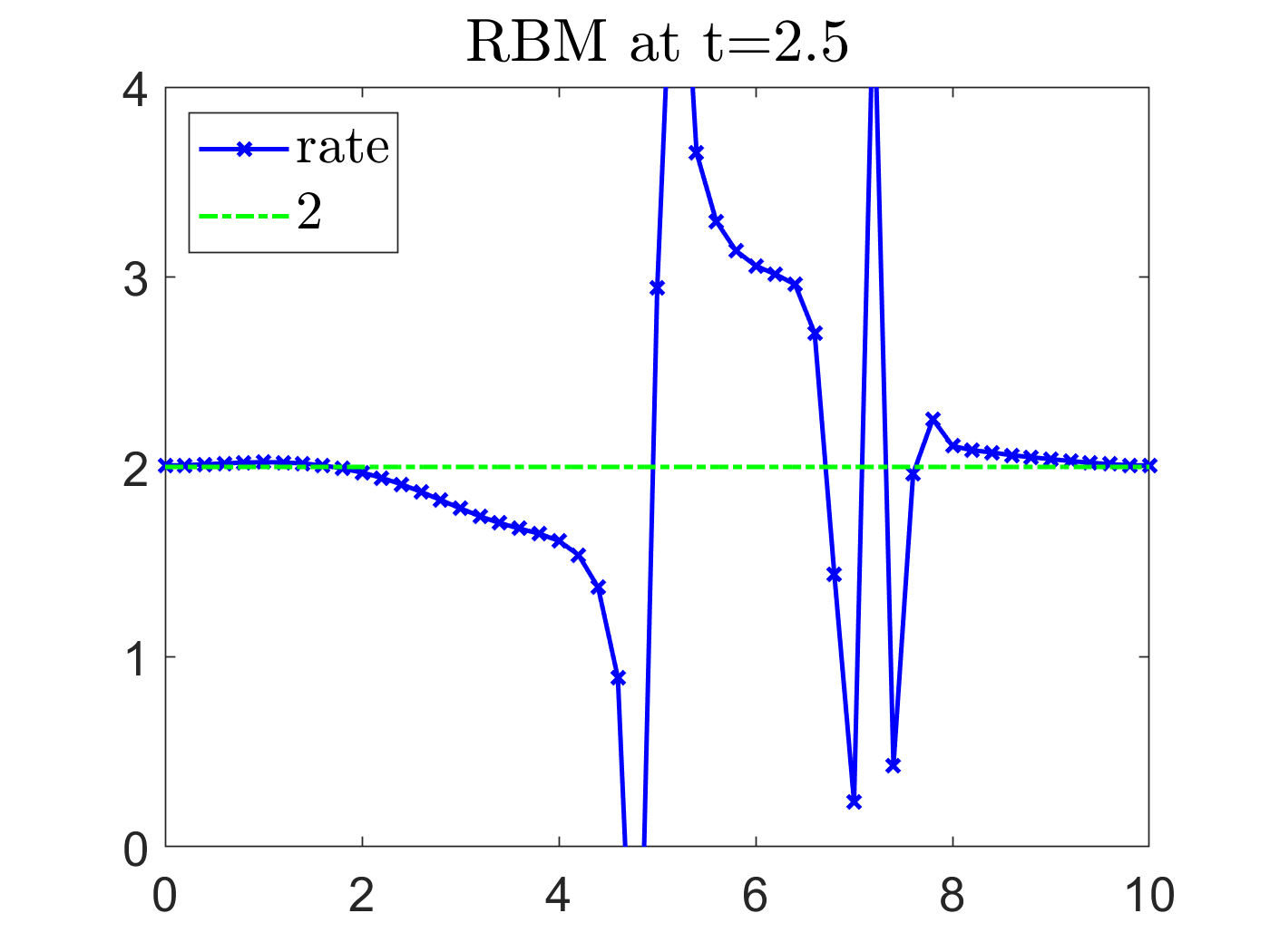}}
\vskip8pt
\centerline{\includegraphics[trim=1.7cm 0.6cm 1.7cm 0.3cm, clip, width=5.8cm]{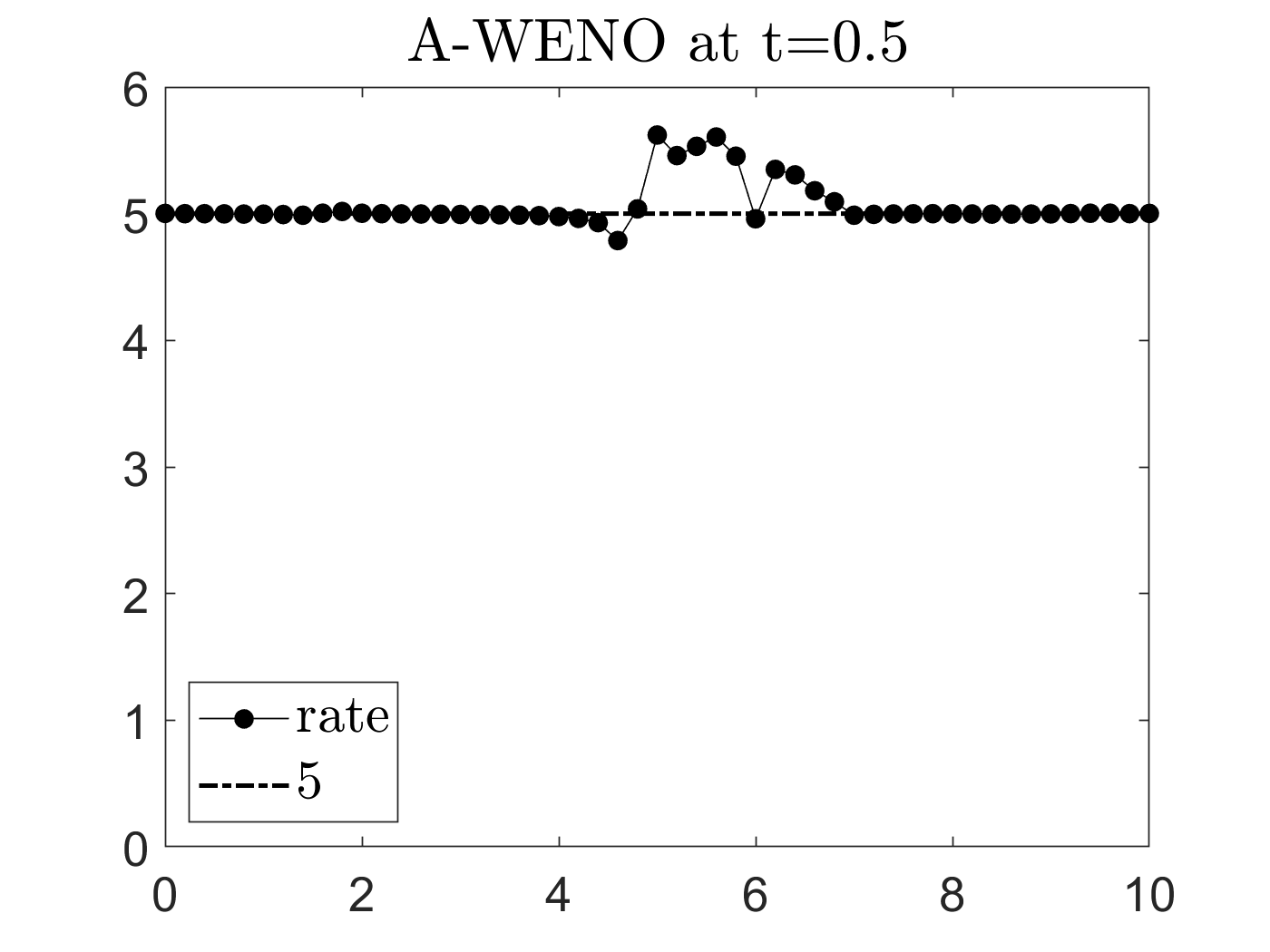}\hspace*{0.0cm}
\includegraphics[trim=1.7cm 0.6cm 1.7cm 0.3cm, clip, width=5.8cm]{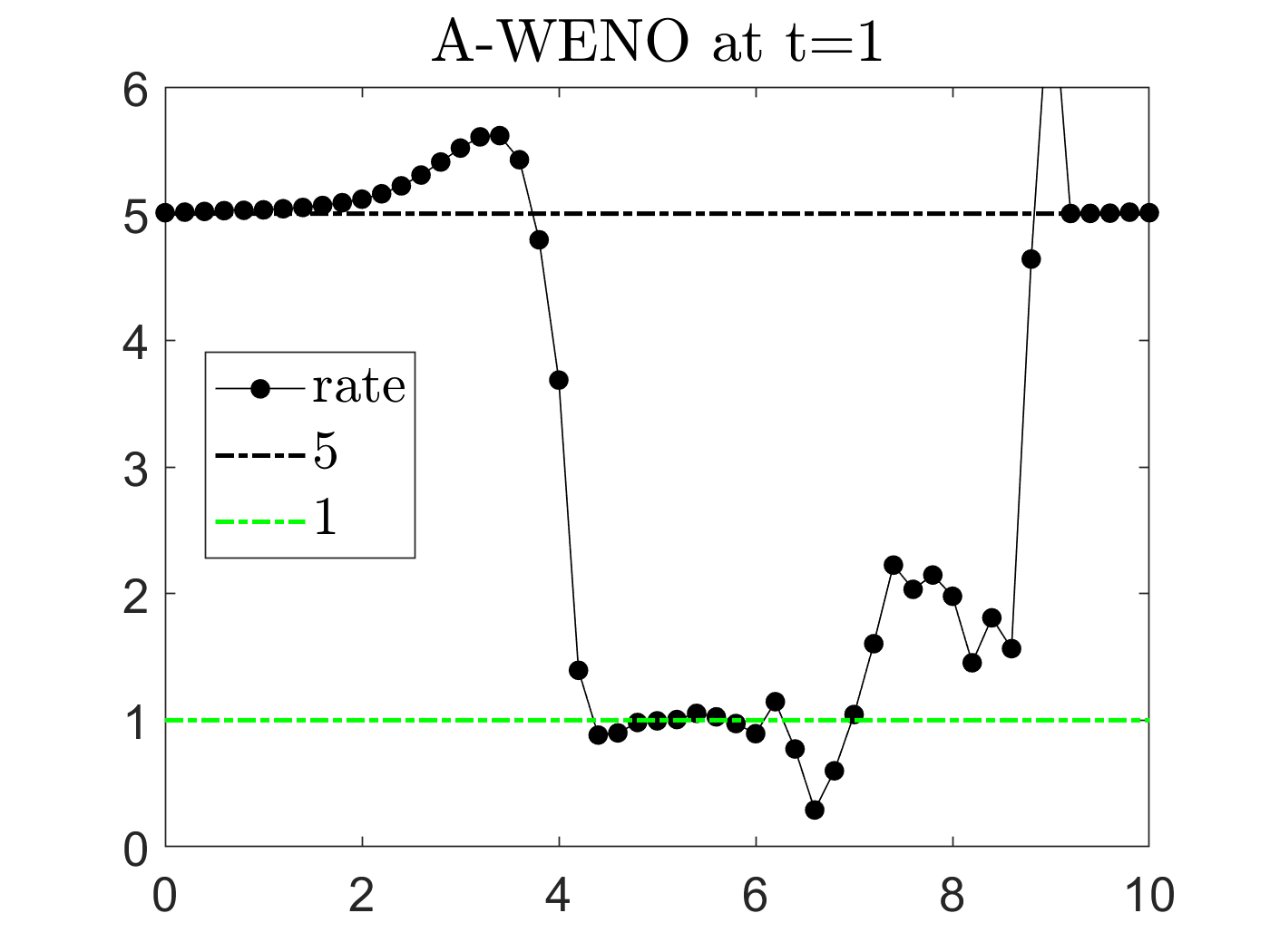}\hspace*{0.0cm}
\includegraphics[trim=1.7cm 0.6cm 1.7cm 0.3cm, clip, width=5.8cm]{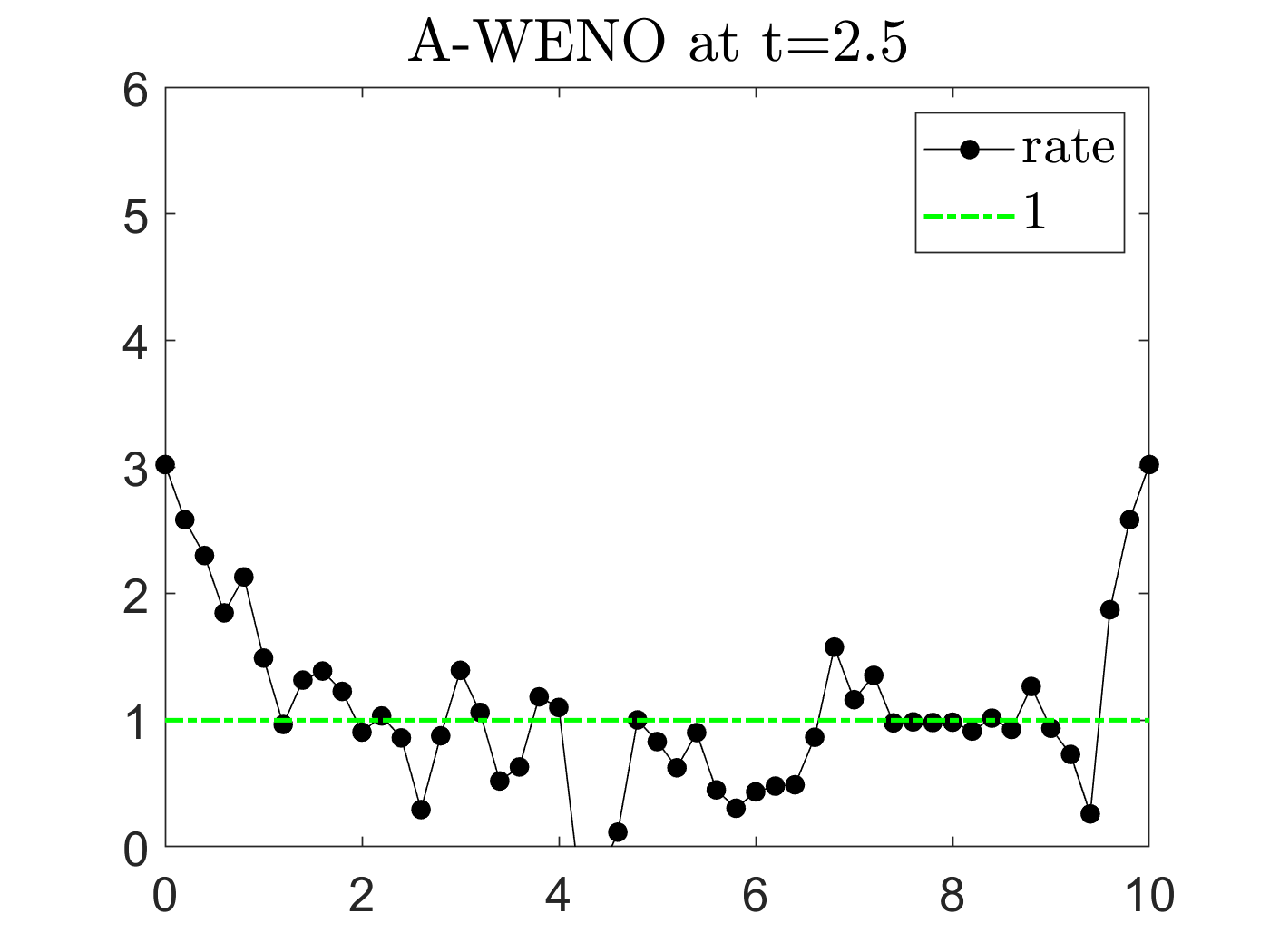}}
\caption{\sf Example 1: Average experimental rates of pointwise convergence for the CU (top row), RBM (middle row), and A-WENO (bottom row)
schemes.\label{Fig4.5}}
\end{figure}

An alternative way to quantify the nonmonotone pointwise convergence is to measure the average experimental pointwise convergence rates
defined by averaging of several, say, 25 convergence rates obtained using \eref{3.1} at the grid point $x=x_{4j}$ and its 24 neighbors on
the left and on the right. Namely, we measure
\begin{equation}
r^{\rm AVE}_{4j}=\frac{1}{25}\sum_{i=-12}^{12}r_{4(j+i)},\quad j=0,\ldots,N,
\label{equ5.6}
\end{equation}
which can be easily calculated for all $j$ thanks to the periodicity of the computed solutions. We compute these average rates for the water
depth for $N=2000$ for the CU and RBM schemes and $N=1000$ for the A-WENO scheme and then plot $r^{\rm AVE}_{160k}$ for the CU and RBM
schemes and $r^{\rm AVE}_{80k}$ for the A-WENO scheme for $k=0,\ldots,50$ in Figure \ref{Fig4.5}. As one can see, the behavior of the
average rates in the CU results shown in the top row of Figure \ref{Fig4.5} is less oscillatory than the behavior of the pointwise rates
reported in Figure \ref{Fig4.4a}. Thus, the average rates give one a somewhat better sense on the way the computed solution converges. From
Figure \ref{Fig4.5}, one can clearly see that at $t=0.5$ (when the numerical solution is still smooth), the average convergence rates for
the CU, RBM and A-WENO schemes are about second-, third- and fifth-order, respectively. At $t=1$ (after the formation of the shock), the
rates for the CU and A-WENO schemes reduce to the first-order in the area behind the shock, while the rates for the RBM scheme only reduce
to the second-order in the most part of the affected area. At the final time $t=2.5$, the average convergence rates for the CU and A-WENO
schemes reduce to the first-order in the entire domain, while they stay close to the second-order for the RBM scheme.

\begin{figure}[ht!]
\centerline{\includegraphics[trim=1.7cm 0.6cm 1.7cm 0.3cm, clip, width=5.8cm]{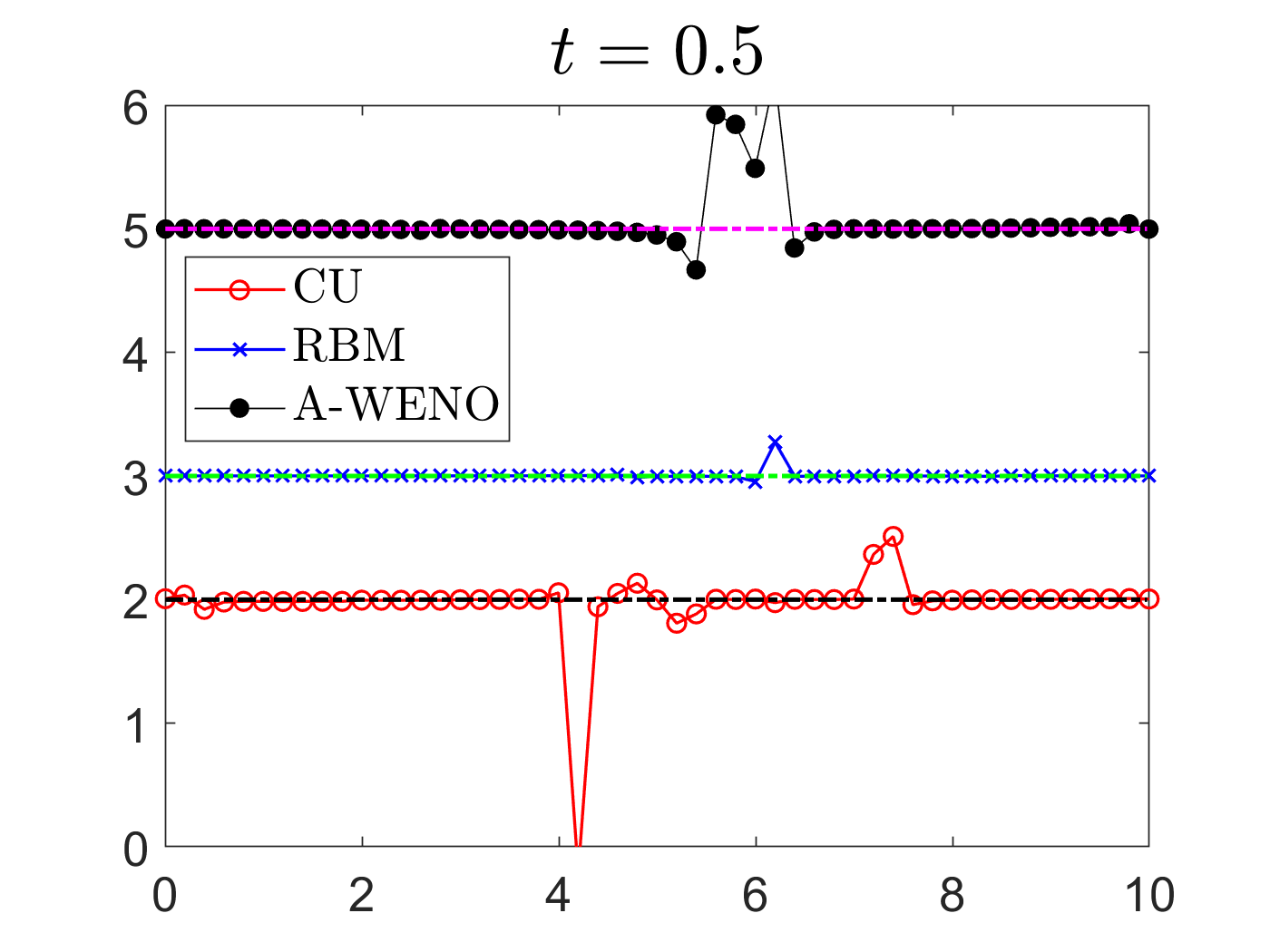}\hspace*{0.0cm}
\includegraphics[trim=1.7cm 0.6cm 1.7cm 0.3cm, clip, width=5.8cm]{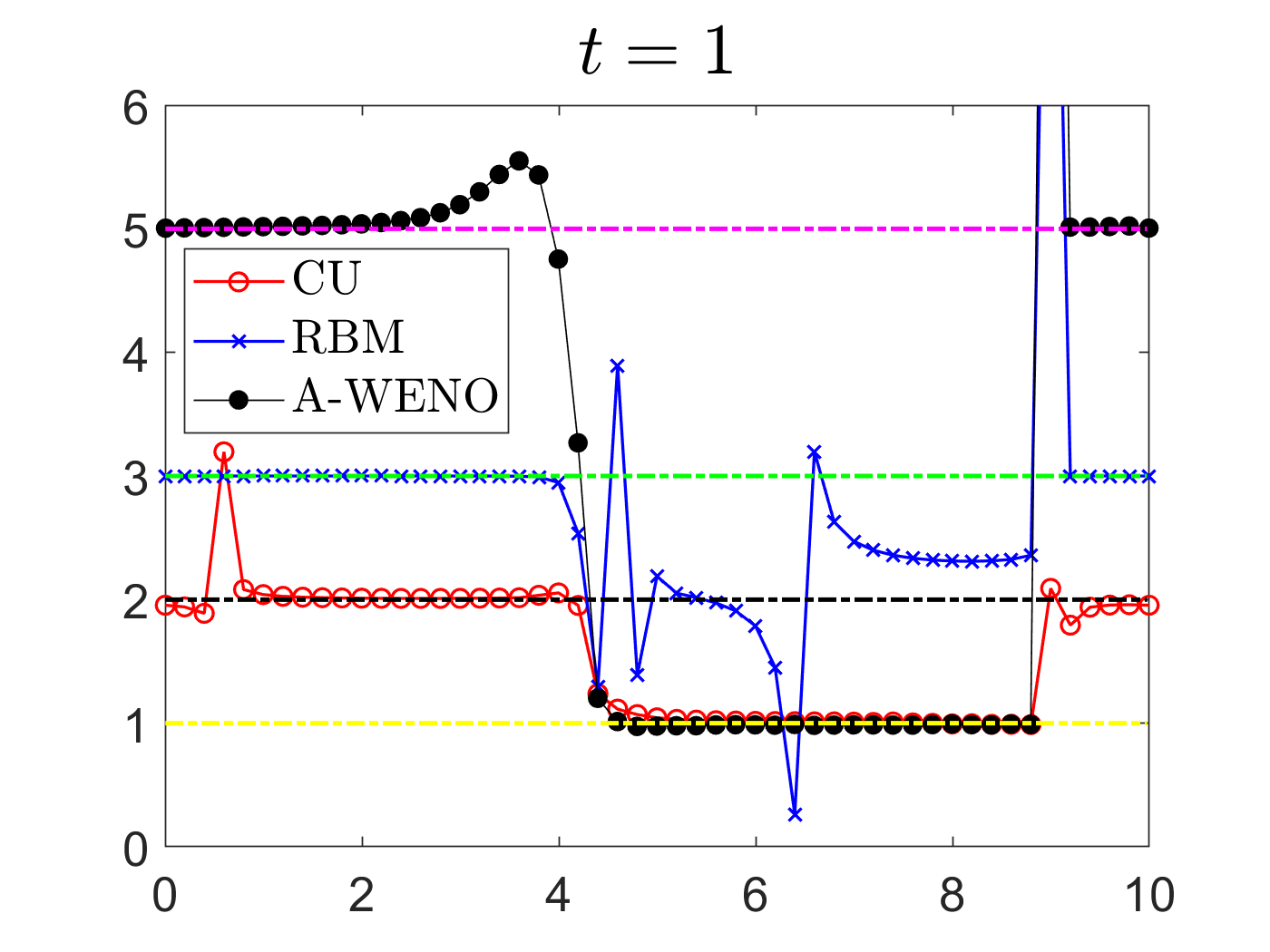}\hspace*{0.0cm}
\includegraphics[trim=1.7cm 0.6cm 1.7cm 0.3cm, clip, width=5.8cm]{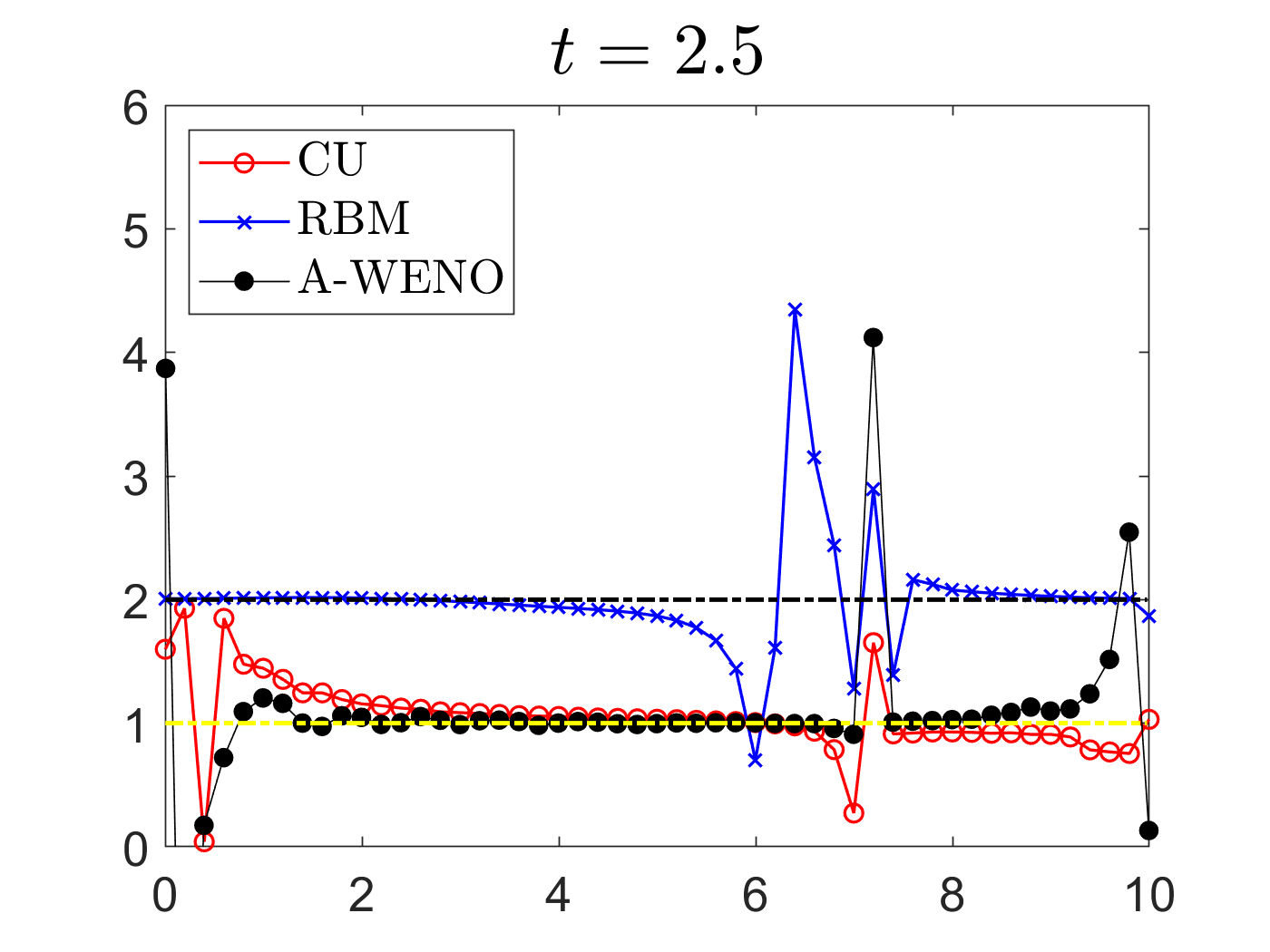}}
\caption{\sf Example 1: Experimental integral rates of convergence for the CU, RBM, and A-WENO schemes at different times.\label{Fig4.6}}
\end{figure}

Finally, we measure the experimental integral rates of convergence \eref{4.1}. We use three imbedded grids with $N=2000$ for the CU and RBM
schemes and with $N=1000$ for the A-WENO scheme. In Figure \ref{Fig4.6}, we present $r^{\rm INT}_{160k}$ for the CU and RBM schemes and
$r^{\rm INT}_{80k}$ for the A-WENO scheme for $k=0,\ldots,50$. Once again, one can see that when the solution is still smooth (at $t=0.5$),
the experimental integral rates of convergence correspond to the formal orders of accuracy for each of the studied schemes, while after the
shock develops and propagates (at $t=1$) the rates reduce in the area behind the shock. At a much later time $t=2.5$, the corresponding
integral rates of convergence reduce in the entire computational domain to the first order for both CU and the A-WENO schemes and to about
the second order for the RBM scheme. This means that the patterns we have observed in the pointwise convergence study manifest themselves in
the computation of the integrals rates as well.

It is also instructive to measure the $W^{-1,1}$ convergence rates given by \eref{4.4} as they would allow us to conclude about the
``actual'' convergence rates for each of the studied schemes. The obtained results are reported in Table \ref{Tab1}. We can clearly see that
the CU scheme is second-order at $t=0.5$ (when the solution is smooth), and only first-order at larger times. The same is true for the
A-WENO scheme, for which the formal fifth order of accuracy is achieved at $t=0.5$ and then it reduces to the first one at larger times. On
contrary, the rate of the convergence for the RBM scheme, which is third-order when the solution is still smooth, reduces to the second one
after the shock formation and it remains almost second-order even at larger times, when the shock propagation affects the entire
computational domain.
\begin{table}[ht!]
\centering
\begin{tabular}{|c|c|c|c|c|c|c|c|}\hline
\multicolumn{2}{|c|}{$~$}&\multicolumn{2}{c|}{CU Scheme}&\multicolumn{2}{c|}{RBM Scheme}&\multicolumn{2}{c|}{A-WENO Scheme}\\\hline
{$t$}&$N$&$||I^N-I^{2N}||_{L^1}$&$r^{\rm INT}$&$||I^N-I^{2N}||_{L^1}$&$r^{\rm INT}$&$||I^N-I^{2N}||_{L^1}$&$r^{\rm INT}$\\\hline
\multirow{3}*{0.5}&1000&1.95e-4&2.01&3.01e-5&2.84&2.45e-8 &4.78\\
\cline{2-8}&2000&4.85e-5&2.01&4.21e-6&2.96&8.92e-9 &4.99\\
\cline{2-8}&4000&1.20e-5&--- &5.41e-7&--- &2.80e-10&--- \\\hline\hline
\multirow{3}*{1}&1000&5.26e-3&1.04&3.00e-4&2.01&3.89e-3&0.99\\
\cline{2-8}&2000&2.56e-3&1.03&7.42e-5&2.01&1.95e-3&1.00\\
\cline{2-8}&4000&1.26e-3&--- &1.85e-5&--- &9.77e-4&--- \\\hline\hline
\multirow{3}*{2.5}&1000&1.74e-3&1.07&2.45e-4&1.90&1.57e-3&1.01\\
\cline{2-8}&2000&8.28e-4&1.01&6.55e-5&1.91&7.80e-4&1.00\\
\cline{2-8}&4000&4.11e-4&--- &1.74e-5&--- &3.89e-4&--- \\
\hline
\end{tabular}
\caption{\sf Example 1: $W^{-1,1}$ convergence rates at different times.\label{Tab1}}
\end{table}

\subsubsection*{Example 2---Test with Two Interacting Shocks}
In the second example originally introduced in \cite{OK}, we consider the following smooth 10-periodic initial conditions:
\begin{equation}
h(x,0)=2\cos\bigg(\frac{\pi x}{5}\bigg)+3,\qquad q(x,0)=0.
\label{6.3}
\end{equation}
One can show that the studied initial value problem \eref{1.1}, \eref{5.1}, \eref{6.3} develops two shock discontinuities per each period
at about $t\approx0.55$.

\begin{figure}[ht!]
\centerline{\includegraphics[trim=3.6cm 1.1cm 2.8cm 0.5cm, clip, width=5.3cm]{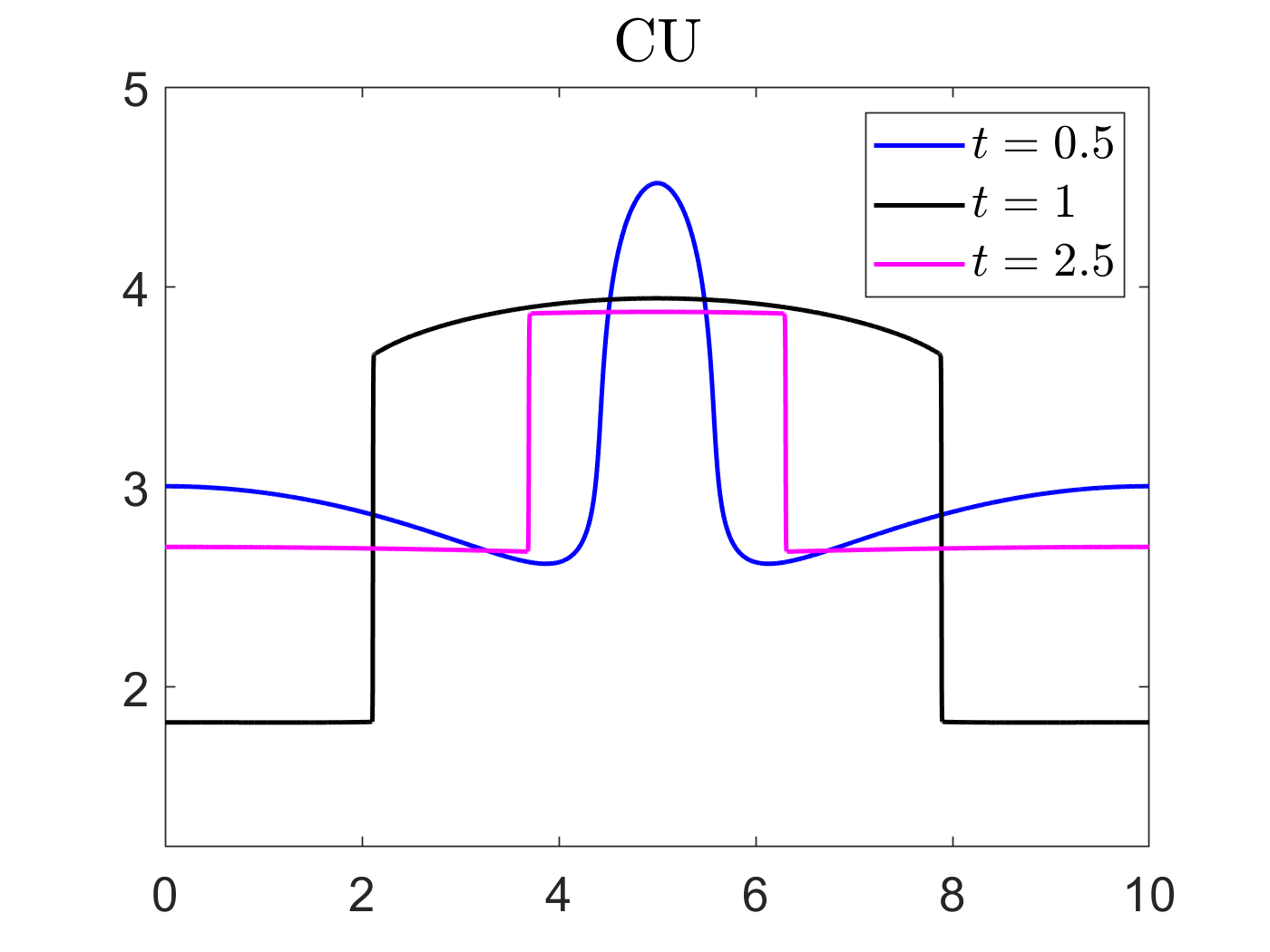}\hspace*{0.5cm}
\includegraphics[trim=3.6cm 1.1cm 2.8cm 0.5cm, clip, width=5.3cm]{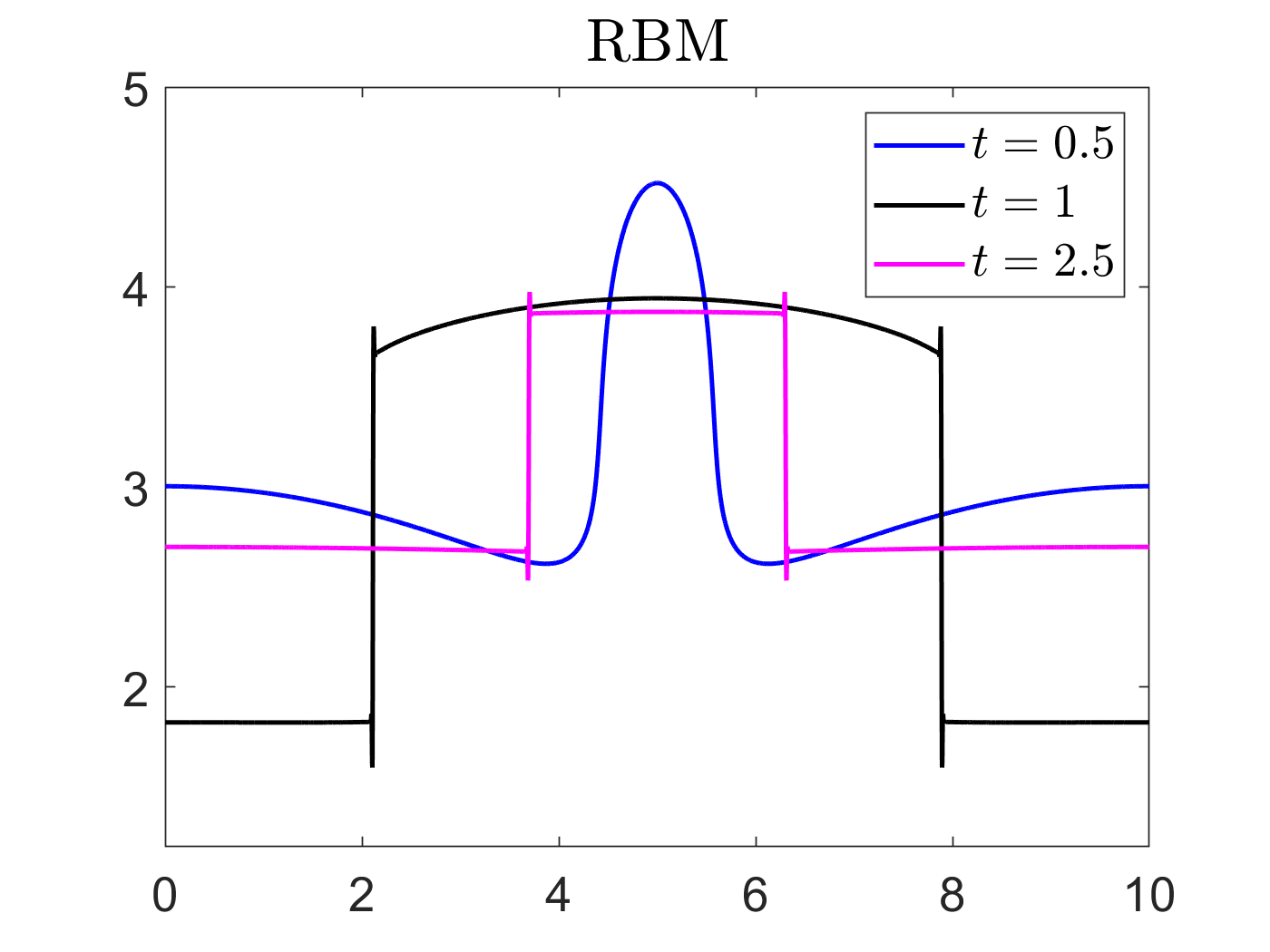}\hspace*{0.5cm}
\includegraphics[trim=3.6cm 1.1cm 2.8cm 0.5cm, clip, width=5.3cm]{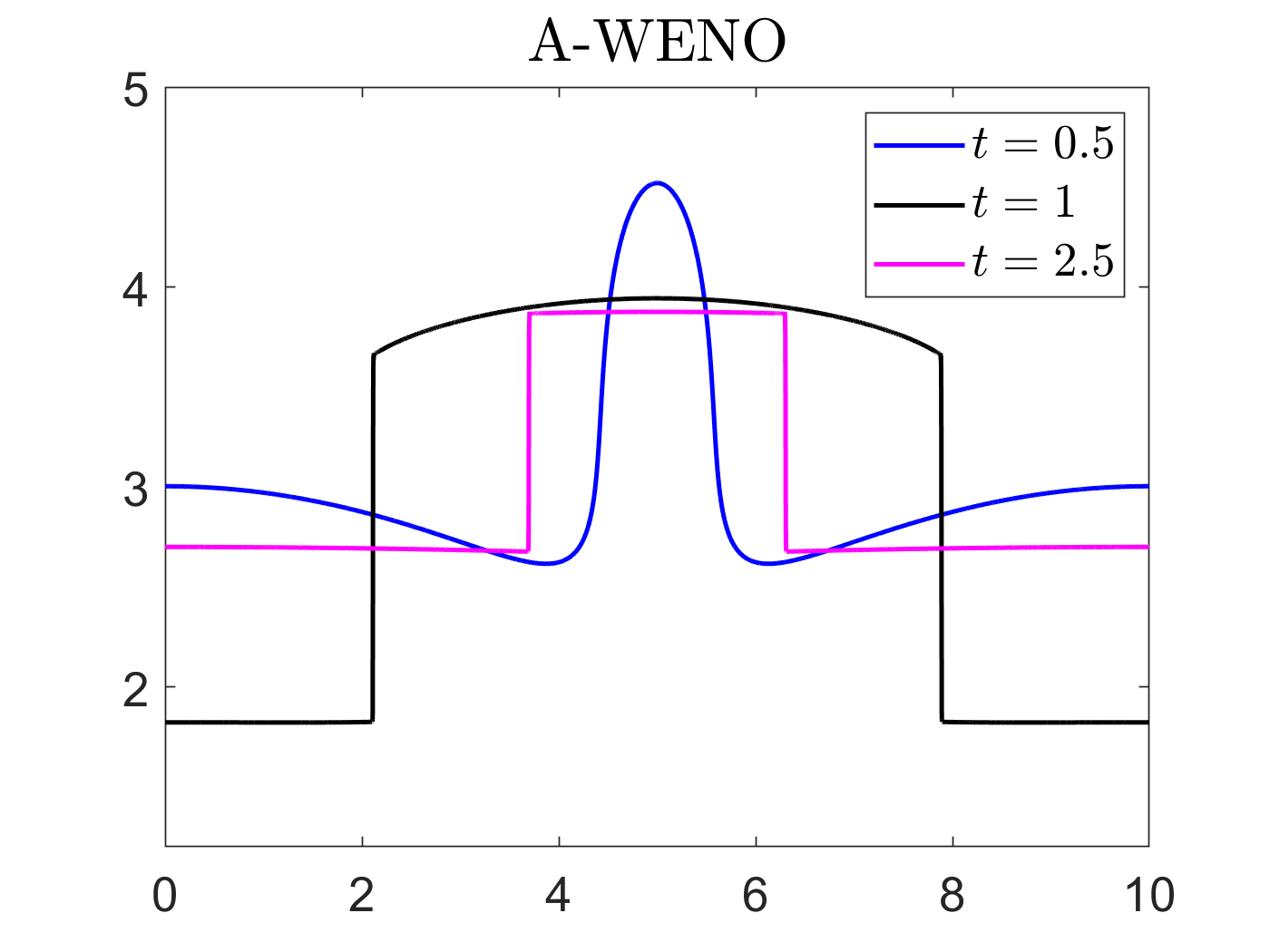}}
\caption{\sf Example 2: Water depth $h$ computed by the CU (left), RBM (middle), and A-WENO (right) schemes.\label{Fig4.7}}
\end{figure}

We compute the solutions by the studied CU, RBM and A-WENO schemes at times $t=0.5$, 1, and 2.5 on the computational domain $[0,10]$ using
1000, 2000, 4000, and 8000 uniform cells. The water depths $h$ computed using 4000 uniform cells are presented in Figure \ref{Fig4.7}. As
one can see, all of the three studied schemes can capture the shock position correctly. At the same time, one can notice that there are
${\cal O}(1)$ oscillations in the neighborhoods of the shocks in the RBM solution.

\begin{figure}[ht!]
\centerline{\includegraphics[trim=1.7cm 0.6cm 1.7cm 0.3cm, clip, width=5.8cm]{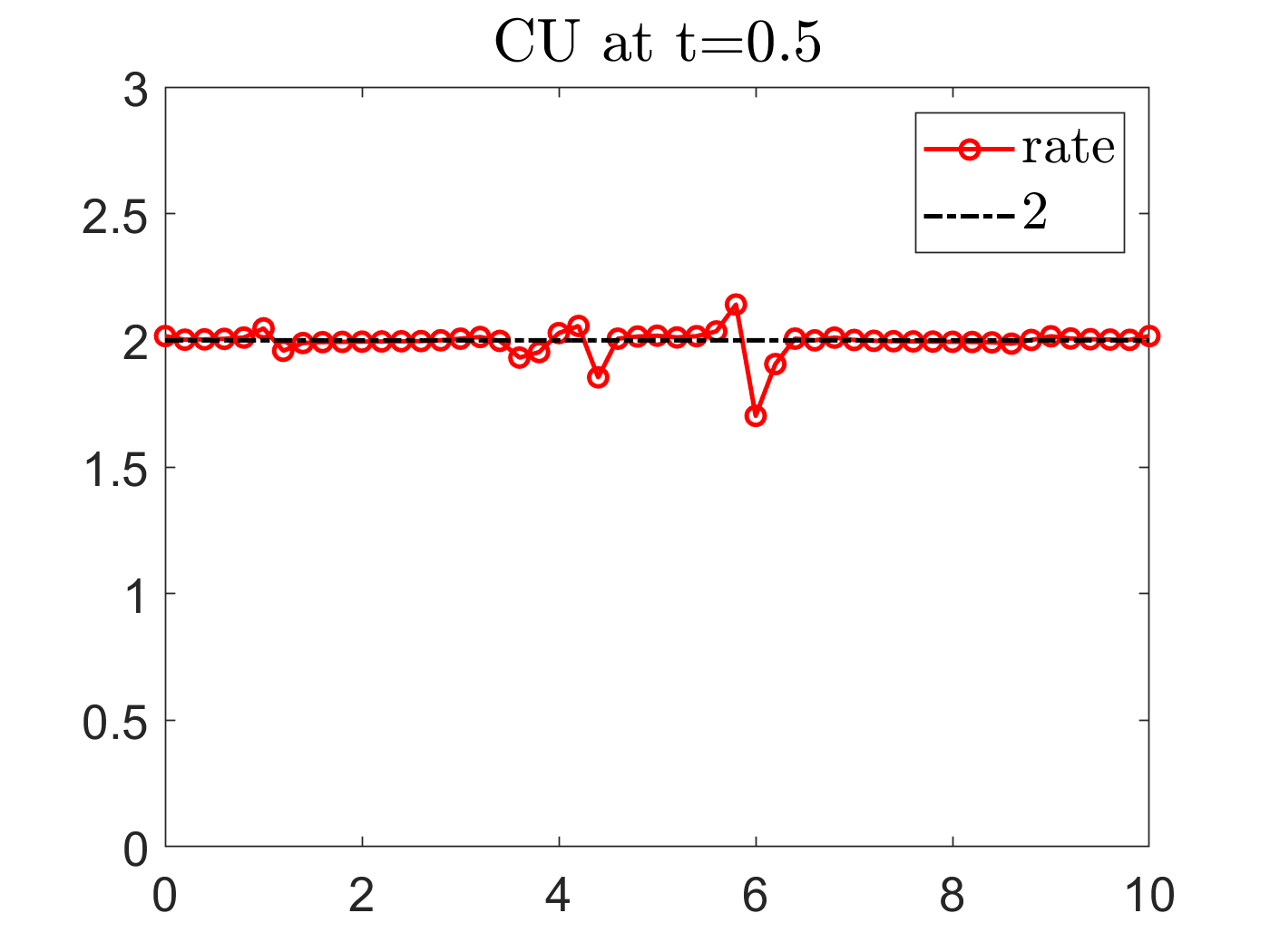}\hspace*{0.0cm}
\includegraphics[trim=1.7cm 0.6cm 1.7cm 0.3cm, clip, width=5.8cm]{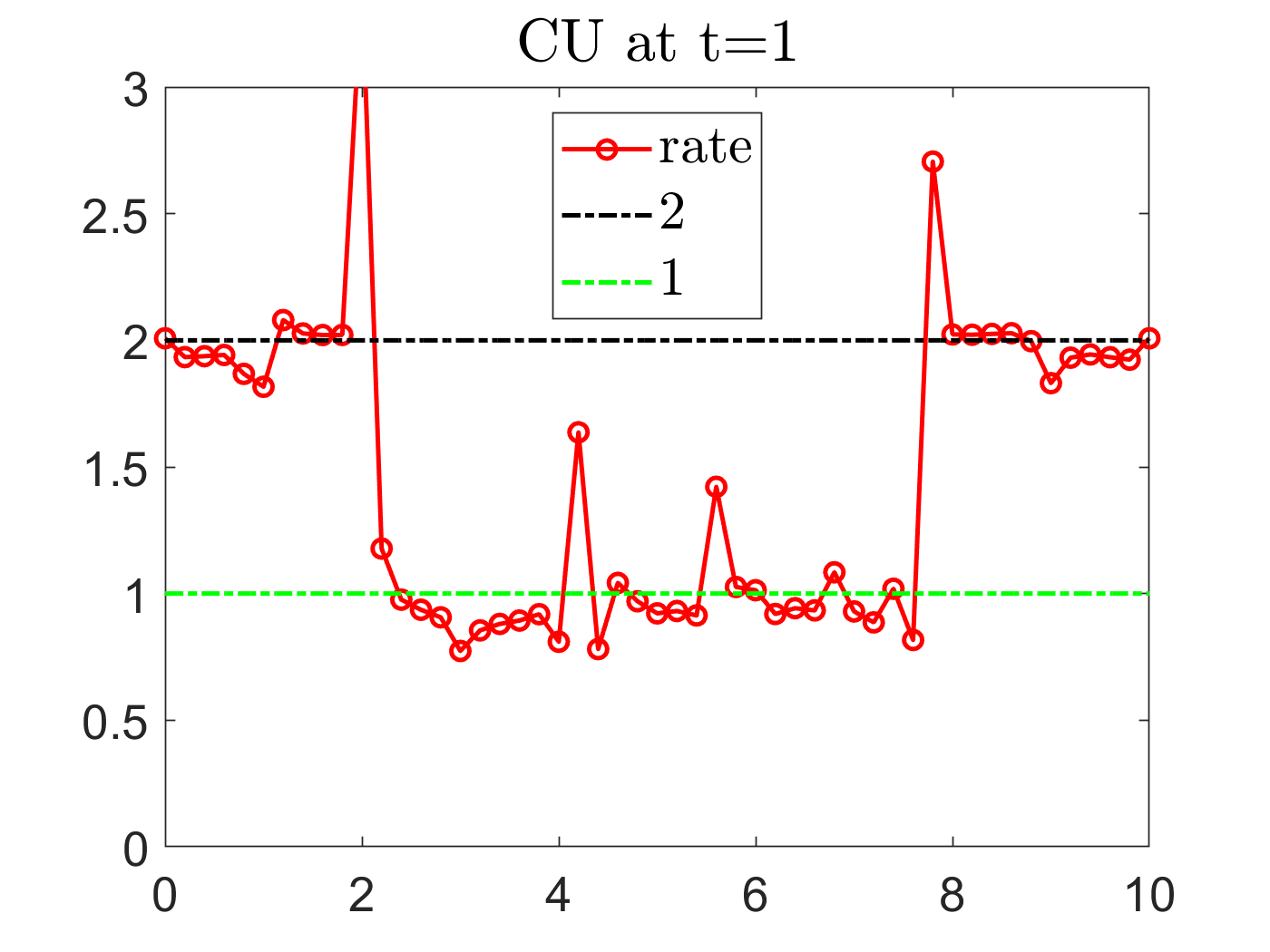}\hspace*{0.0cm}
\includegraphics[trim=1.7cm 0.6cm 1.7cm 0.3cm, clip, width=5.8cm]{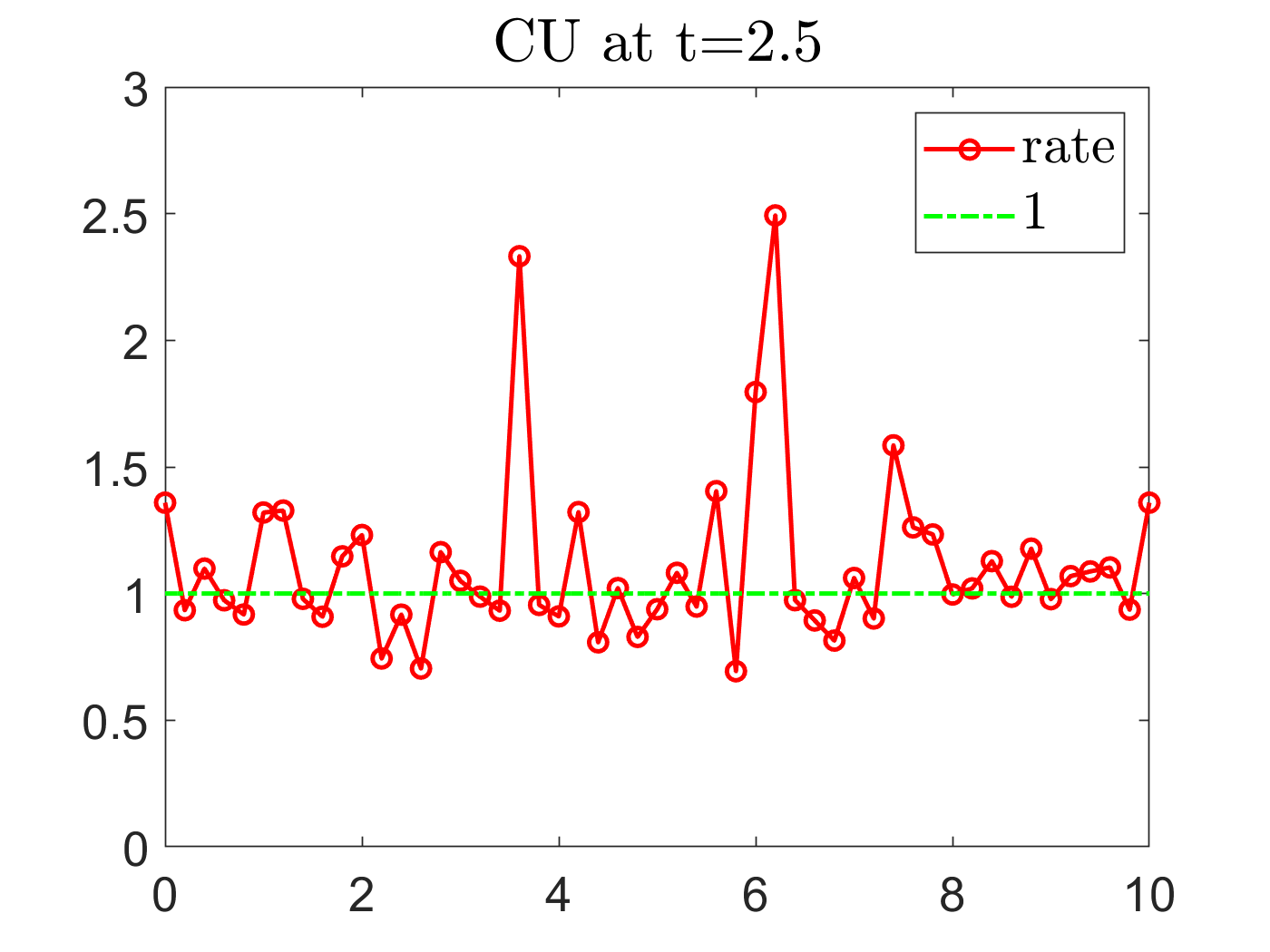}}
\vskip8pt
\centerline{\includegraphics[trim=1.7cm 0.6cm 1.7cm 0.3cm, clip, width=5.8cm]{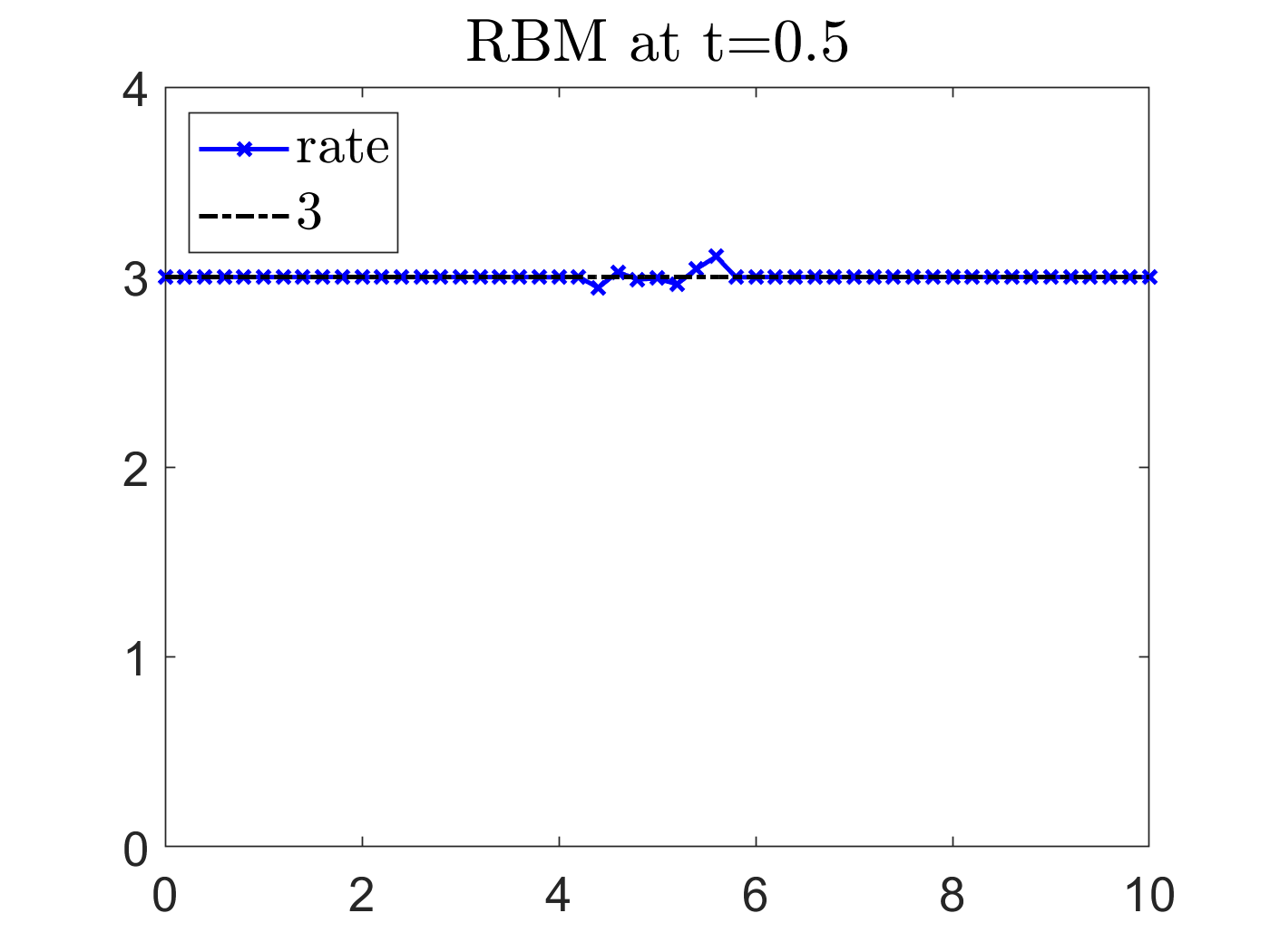}\hspace*{0.0cm}
\includegraphics[trim=1.7cm 0.6cm 1.7cm 0.3cm, clip, width=5.8cm]{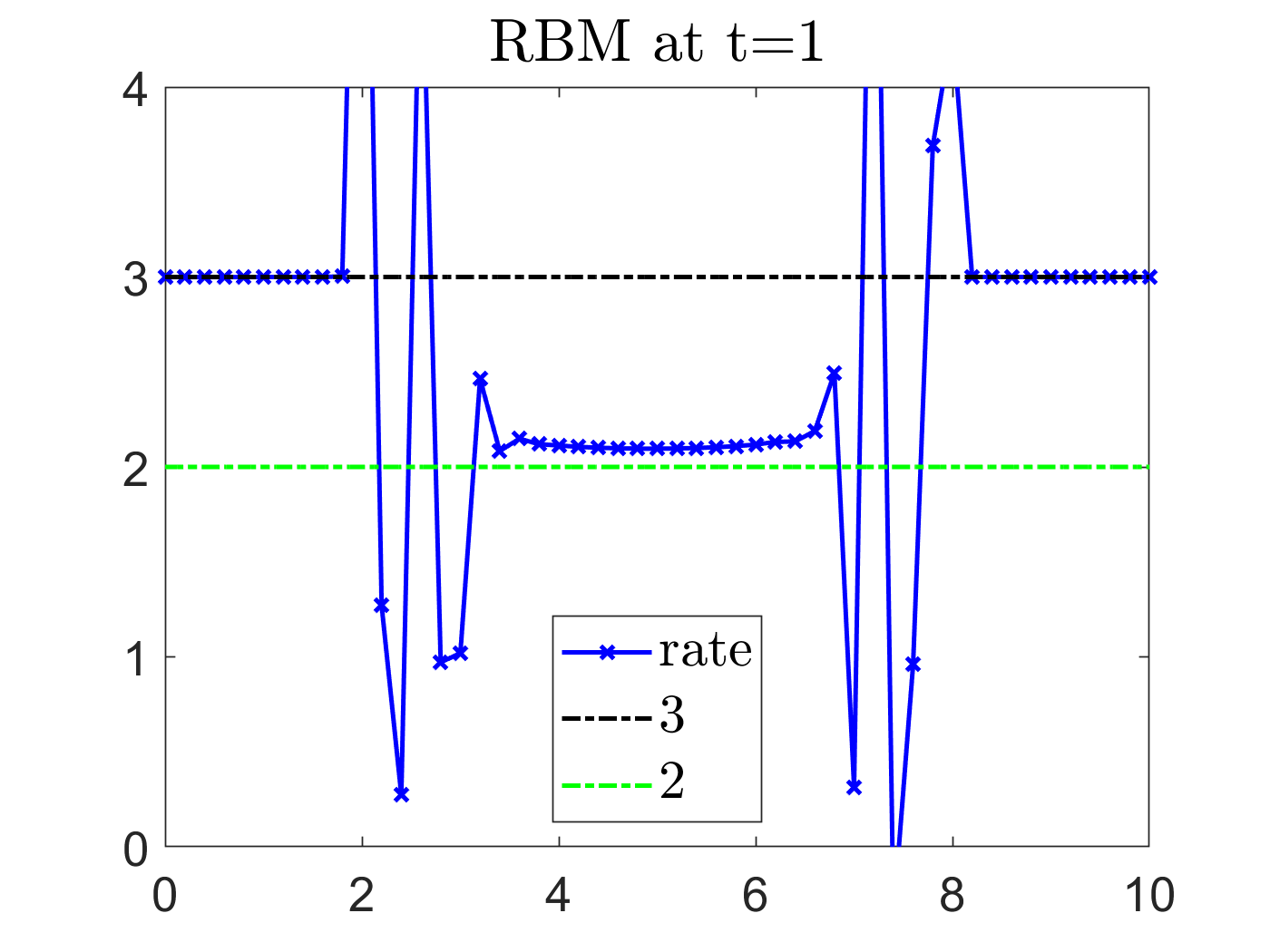}\hspace*{0.0cm}
\includegraphics[trim=1.7cm 0.6cm 1.7cm 0.3cm, clip, width=5.8cm]{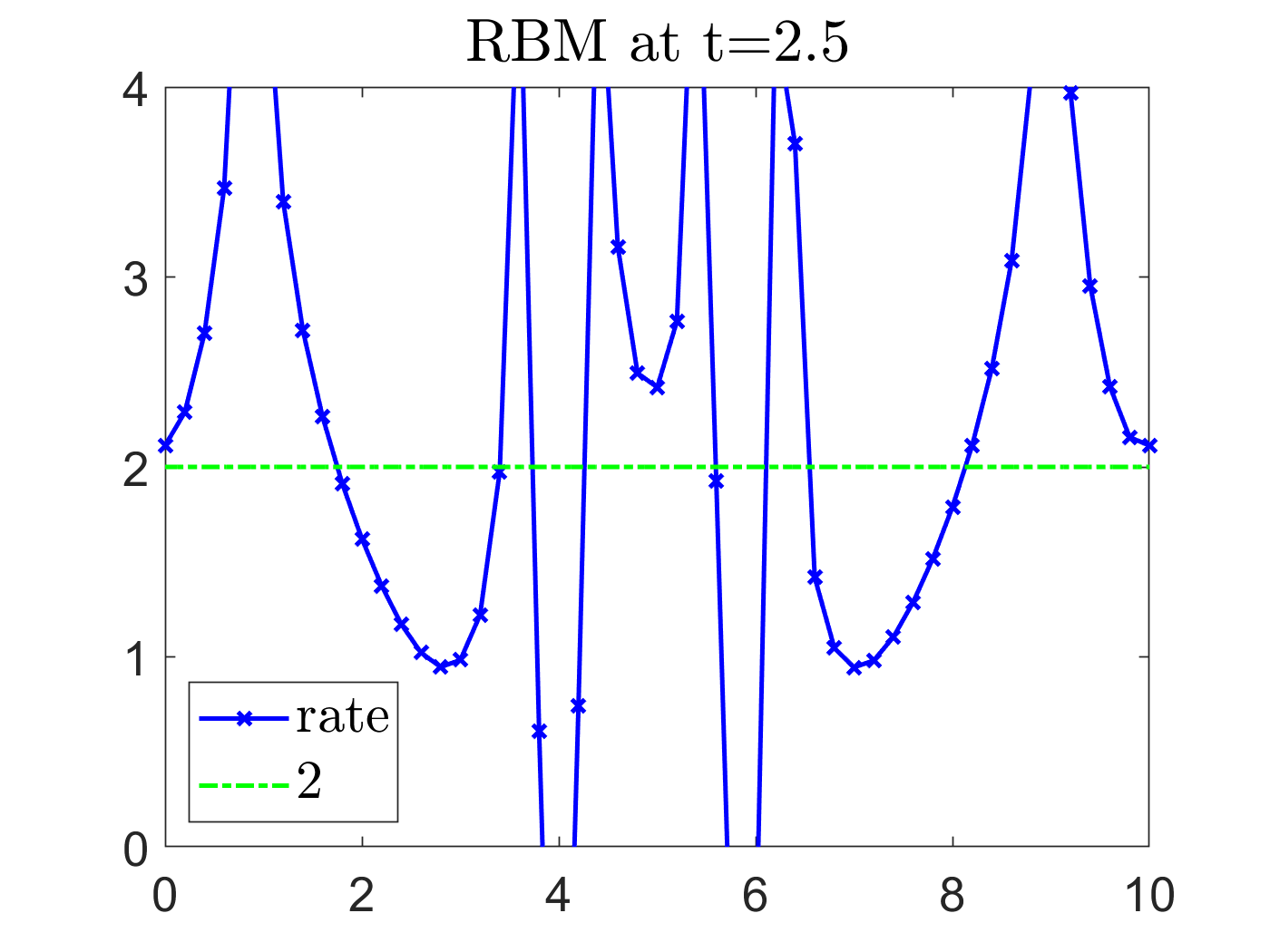}}
\vskip8pt
\centerline{\includegraphics[trim=1.7cm 0.6cm 1.7cm 0.3cm, clip, width=5.8cm]{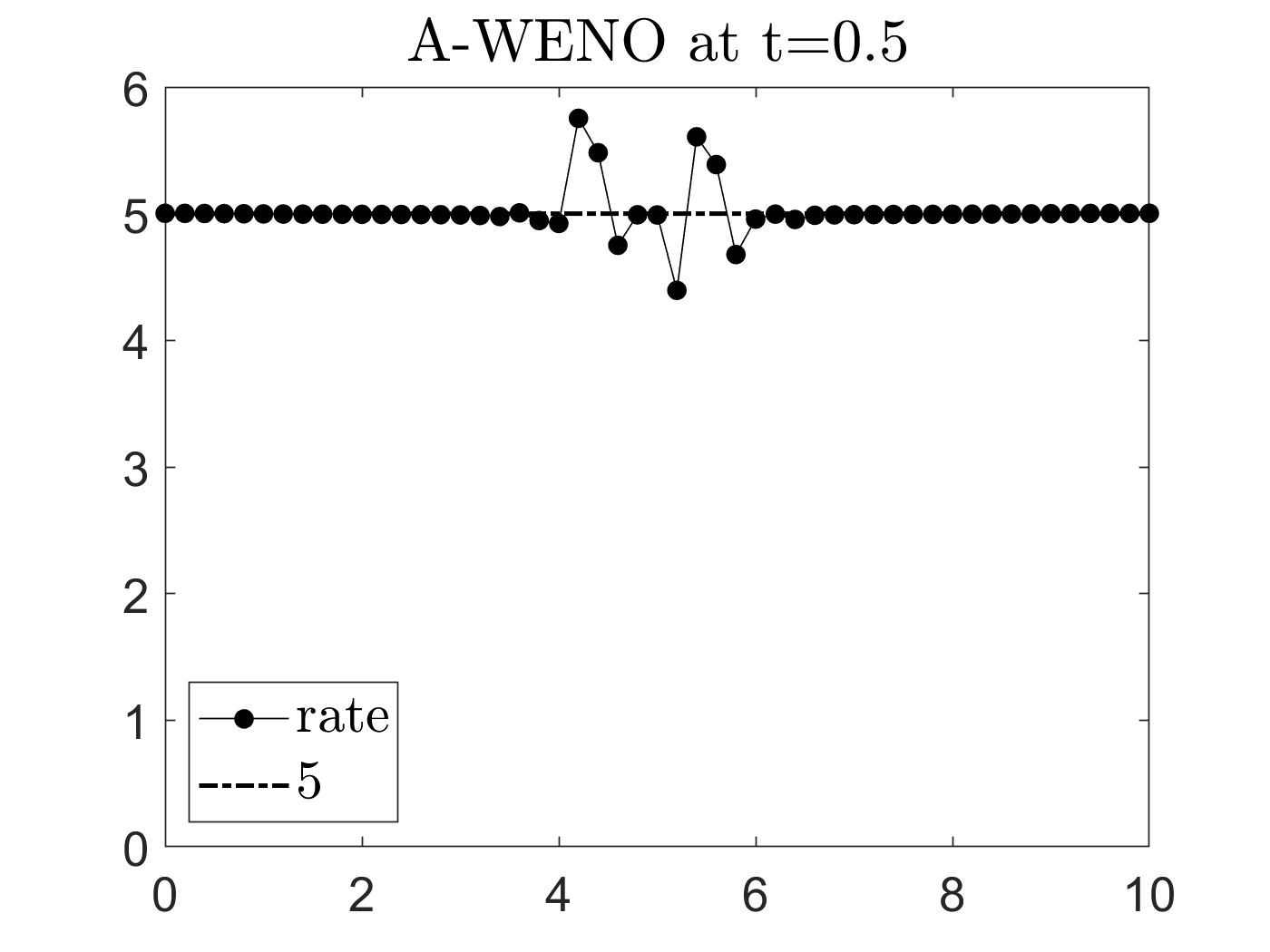}\hspace*{0.0cm}
\includegraphics[trim=1.7cm 0.6cm 1.7cm 0.3cm, clip, width=5.8cm]{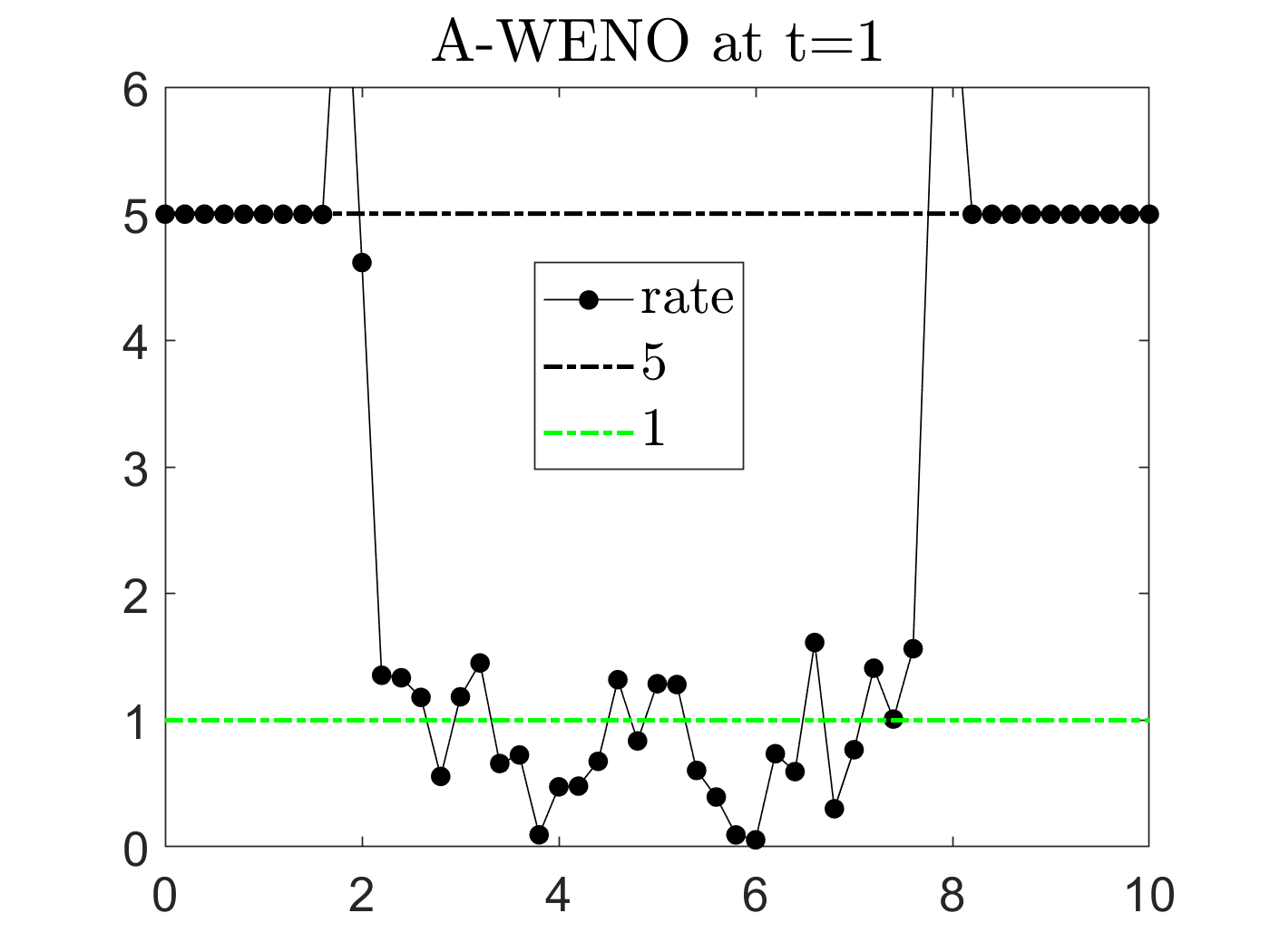}\hspace*{0.0cm}
\includegraphics[trim=1.7cm 0.6cm 1.7cm 0.3cm, clip, width=5.8cm]{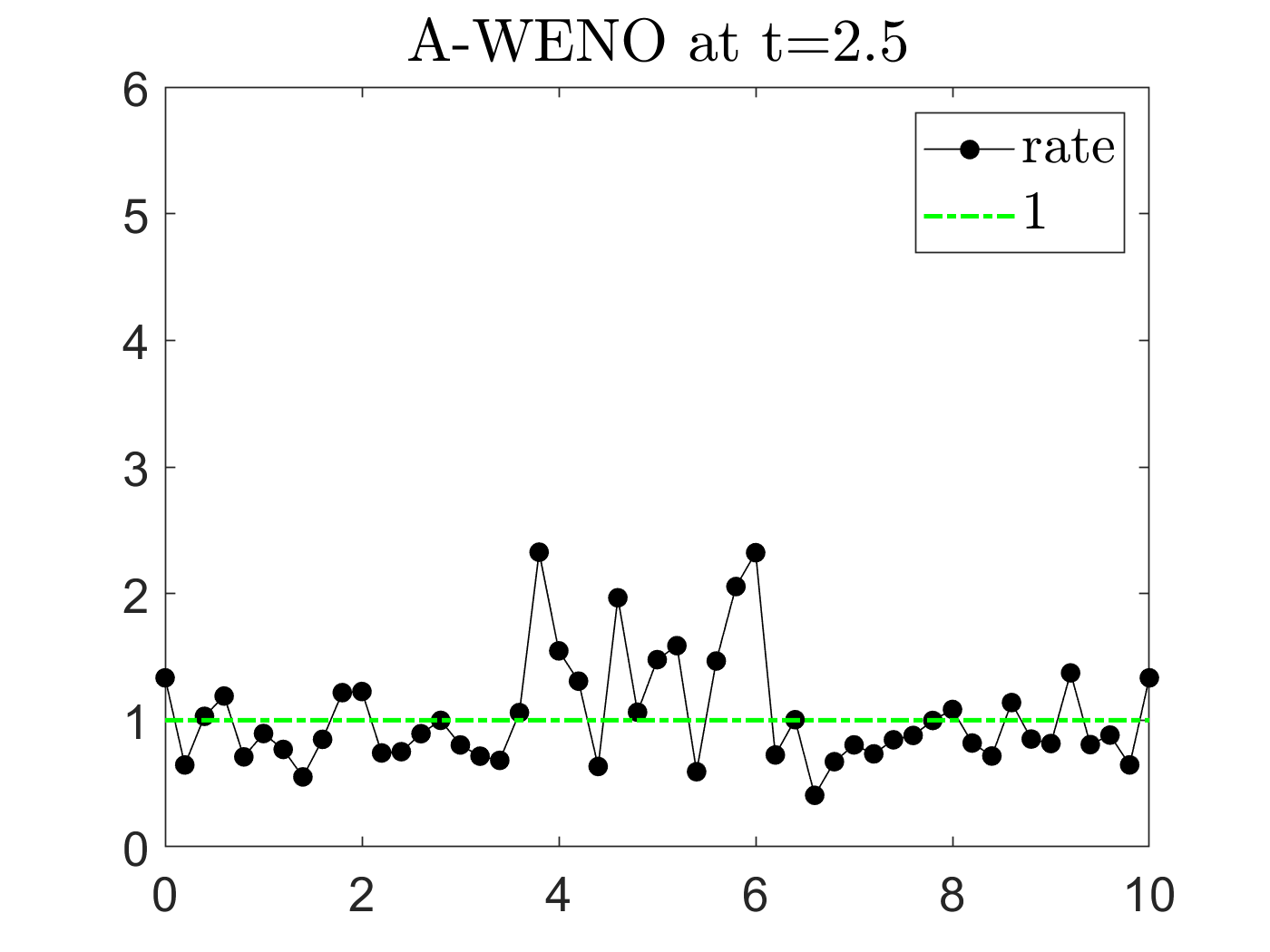}}
\caption{\sf Example 2: Average experimental rates of pointwise convergence for the CU (top row), RBM (middle row), and A-WENO (bottom row)
schemes.\label{Fig4.11}}
\end{figure}

As in Example 1, we check the pointwise convergence of the water depth $h$ by computing the average experimental pointwise convergence rates
defined in \eref{equ5.6} using three imbedded grids with $N=2000$ for the CU and RBM schemes and with $N=1000$ for the A-WENO scheme. We
plot $r^{\rm AVE}_{160k}$ for the CU and RBM schemes and $r^{\rm AVE}_{80k}$ for the A-WENO scheme for $k=0,\ldots,50$ in Figure
\ref{Fig4.11}. One can clearly see that at $t=0.5$ (when the numerical solution is still smooth), the average convergence rates for the CU,
RBM, and A-WENO schemes are about second-, thirda,- and fifth-order, respectively. At $t=1$ (after the formation of the shocks), the rates
for the CU and A-WENO schemes reduce to the first-order in the area between the shocks, while the rates for the RBM scheme only reduce to
the second-order in the most part of the affected area. At the final time $t=2.5$, the average convergence rates for the CU and A-WENO
schemes reduce to the first-order in the entire domain, while they stay around to the second-order for the RBM scheme.

\begin{figure}[ht!]
\centerline{\includegraphics[trim=1.7cm 0.6cm 1.7cm 0.3cm, clip, width=5.8cm]{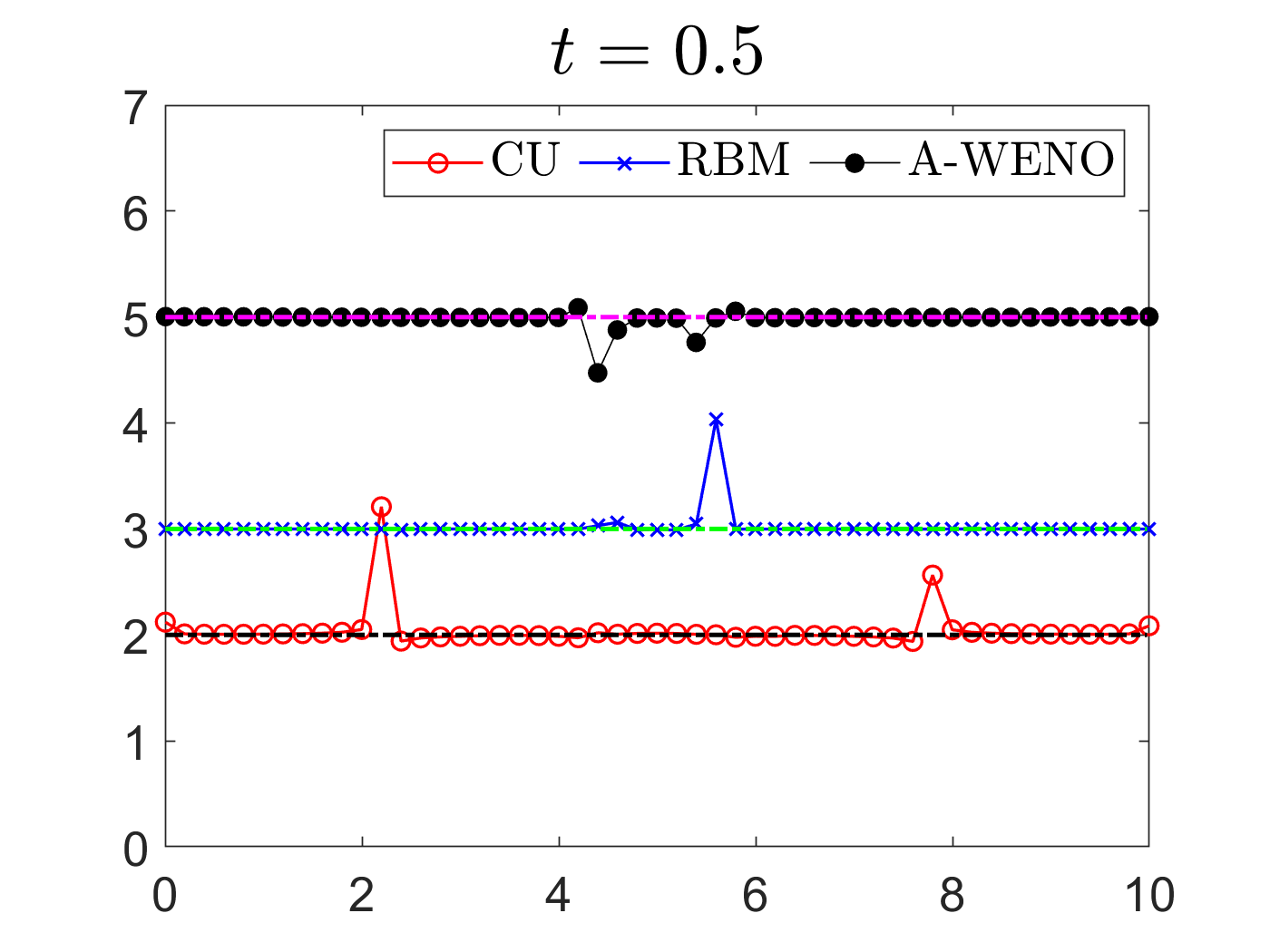}\hspace*{0.0cm}
            \includegraphics[trim=1.7cm 0.6cm 1.7cm 0.3cm, clip, width=5.8cm]{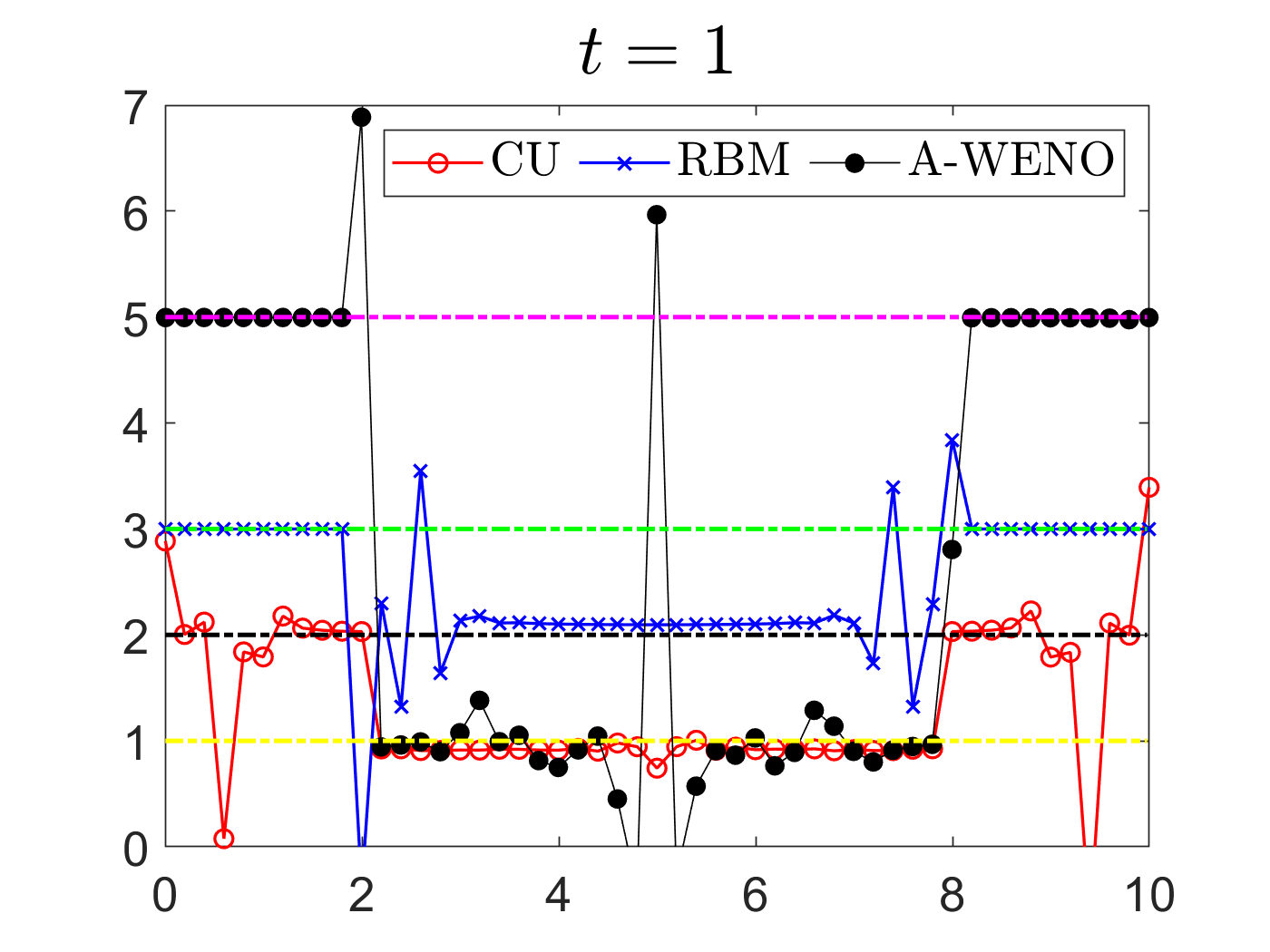}\hspace*{0.0cm}
            \includegraphics[trim=1.7cm 0.6cm 1.7cm 0.3cm, clip, width=5.8cm]{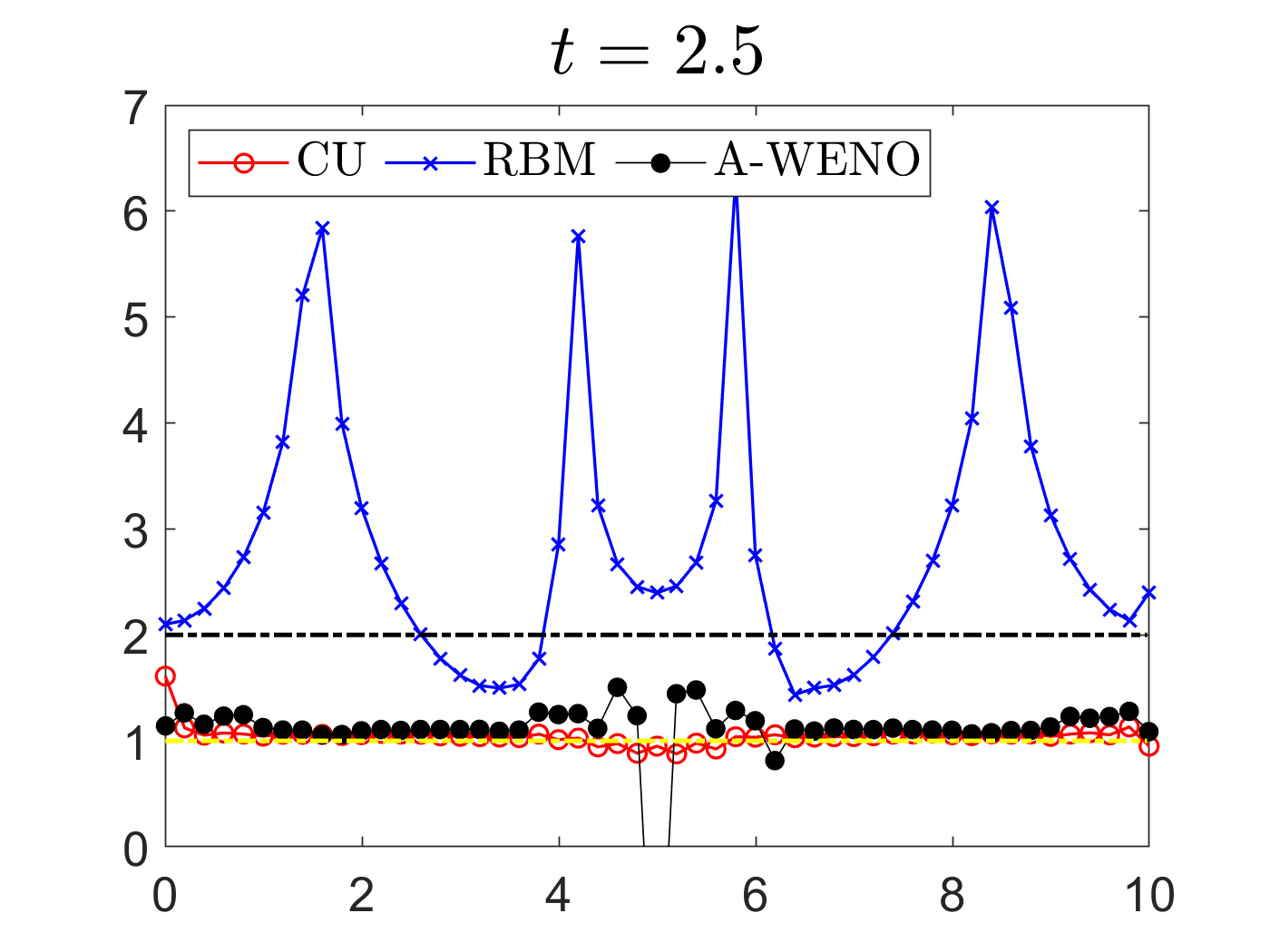}}
\caption{\sf Example 2: Experimental integral rates of convergence for the CU, RBM, and A-WENO schemes at different times.\label{Fig4.12}}
\end{figure}

We then compute the integral rates of convergence using \eref{4.1} and present $r^{\rm INT}_{160k}$ for the CU and RBM schemes and
$r^{\rm INT}_{80k}$ for the A-WENO scheme for $k=0,\ldots,50$ in Figure \ref{Fig4.12}. Once again, one can see that when the solution is
still smooth (at $t=0.5$), the experimental integral rates of convergence correspond to the formal orders of accuracy for each of the
studied schemes, while after the shocks develop and propagate (at $t=1$) the rates reduce in the area between the shocks to the first-order
for the CU and A-WENO schemes and to the second-order for the RBM scheme. At a much later time $t=2.5$, the corresponding reduction of the
integral rates of convergence occurs in the entire computational domain.

Finally, we compute the $W^{-1,1}$ convergence rates given by \eref{4.4} and report the results in Table \ref{Tab2}, where one can clearly
see that the obtained convergence rates are the same as in Example 1, that is, all of the studied three schemes achieve their formal order
of accuracy when the solution is smooth ($t=0.5$) and later on the rates drop to either first (CU and A-WENO) or second (RBM) order.
\begin{table}[ht!]
\centering
\begin{tabular}{|c|c|c|c|c|c|c|c|}\hline
\multicolumn{2}{|c|}{$~$}&\multicolumn{2}{c|}{CU Scheme}&\multicolumn{2}{c|}{RBM Scheme}&\multicolumn{2}{c|}{A-WENO Scheme}\\\hline
$t$&$N$&$||I^N-I^{2N}||_{L^1}$&$r^{\rm INT}$&$||I^N-I^{2N}||_{L^1}$&$r^{\rm INT}$&$||I^N-I^{2N}||_{L^1}$&$r^{\rm INT}$\\\hline
\multirow{3}*{0.5}&1000&1.22e-4&2.02&3.89e-5&2.95&7.81e-8&4.93\\
\cline{2-8}&2000&2.99e-5&2.01&5.02e-6&2.99&2.56e-9&5.00\\
\cline{2-8}&4000&7.41e-6&--- &6.31e-7&--- &8.00e-10&---\\\hline\hline
\multirow{3}*{1}&1000&1.53e-3&0.93&6.06e-4&2.00&1.19e-3&0.99\\
\cline{2-8}&2000&8.00e-4&0.95&1.52e-4&1.99&5.97e-4&1.01\\
\cline{2-8}&4000&4.15e-4&--- &3.82e-5& ---&2.96e-4&--- \\\hline\hline
\multirow{3}*{2.5}&1000&2.16e-3&1.11&4.57e-4&1.88&2.29e-3&1.13\\
\cline{2-8}&2000&9.98e-4&1.07&1.24e-4&2.01&1.05e-3&1.06\\
\cline{2-8}&4000&4.76e-4&--- &3.07e-5&--- &5.01e-4&--- \\
\hline
\end{tabular}
\caption{\sf Example 2: $W^{-1,1}$ convergence rates at different times.\label{Tab2}}
\end{table}

\subsubsection*{Example 3---Test with an Isolated Shock}
In the third example, we consider the following initial conditions:
\begin{equation*}
\big(h(x,0),q(x,0)\big)=\big(h_0(x),q_0(x)\big)=
\begin{cases}
(1,0),&x<5,\\
\Big(\dfrac{-5+3\sqrt{5}}{10},\dfrac{3\sqrt{5}-15}{10}\Big),&\mbox{otherwise},
\end{cases}
\end{equation*}
which correspond to an isolated shock moving to the right with the constant velocity 1. In this case, the exact solution is given by
$h(x,t)=h_0(x-t)$ and $q(x,t)=q_0(x-t)$.

\begin{figure}[ht!]
\centerline{\includegraphics[trim=0.8cm 0.3cm 1.0cm 0.1cm, clip, width=5.3cm]{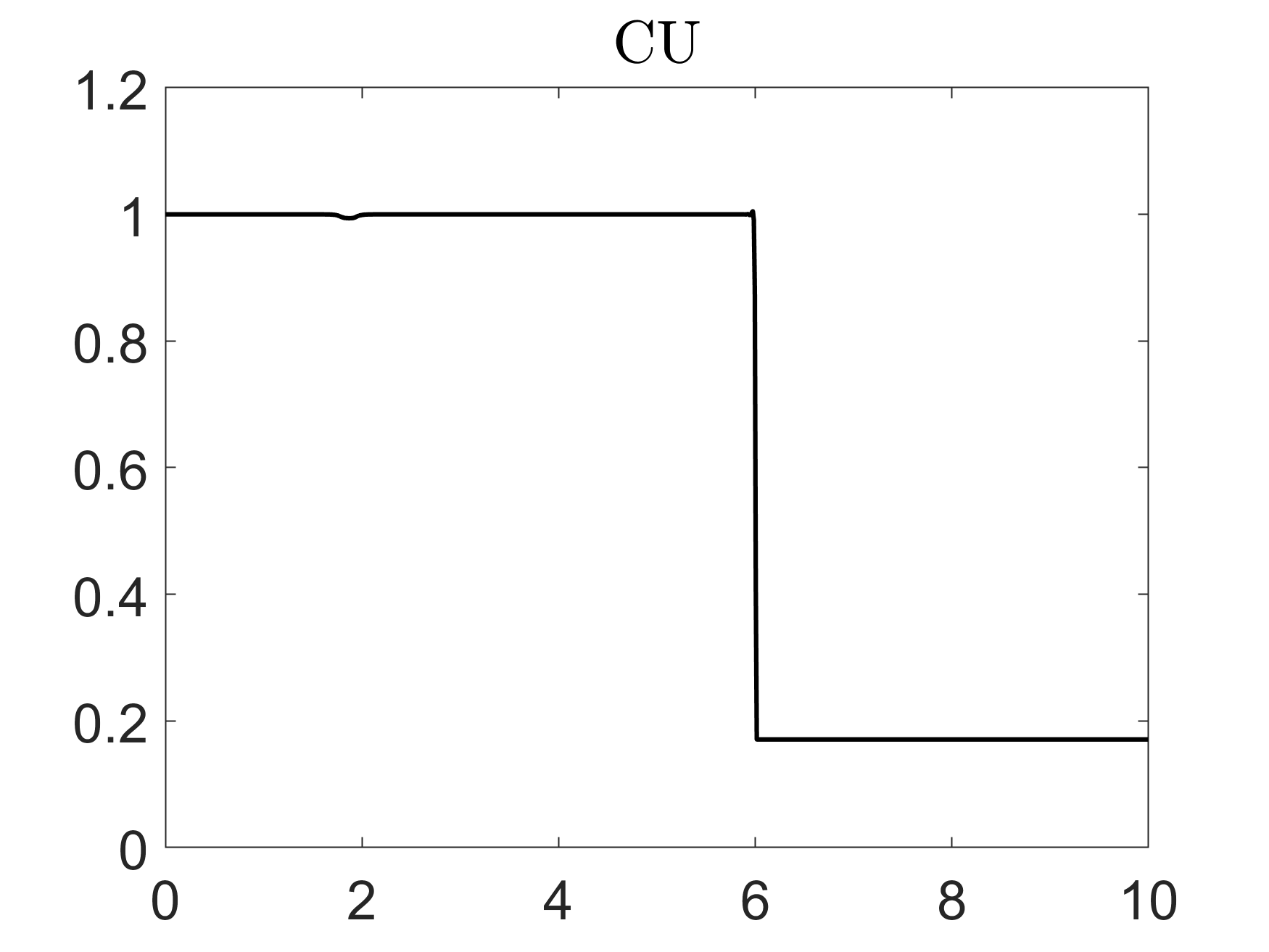}\hspace*{0.5cm}
            \includegraphics[trim=0.8cm 0.3cm 1.0cm 0.1cm, clip, width=5.3cm]{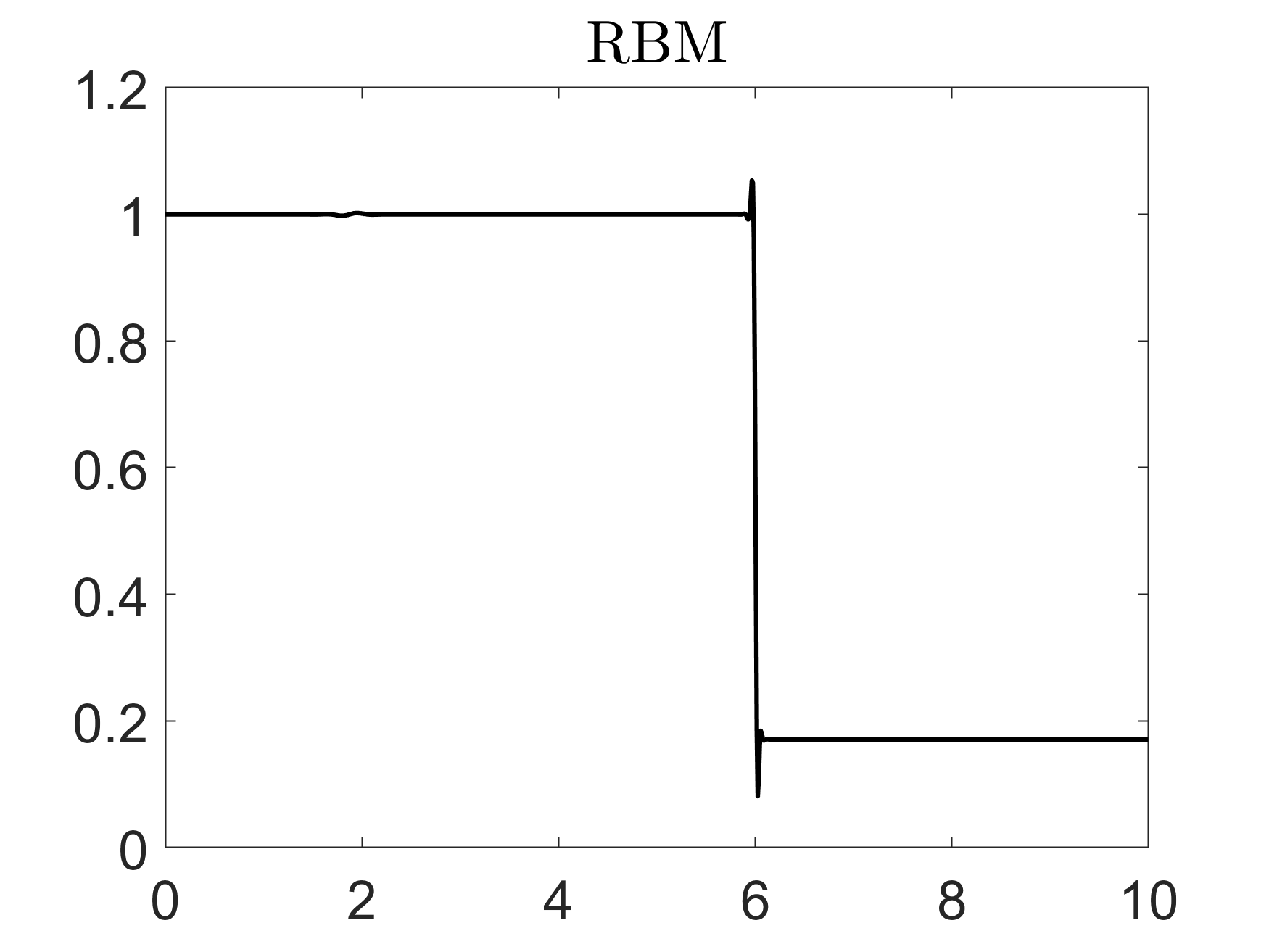}\hspace*{0.5cm}
            \includegraphics[trim=0.8cm 0.3cm 1.0cm 0.1cm, clip, width=5.3cm]{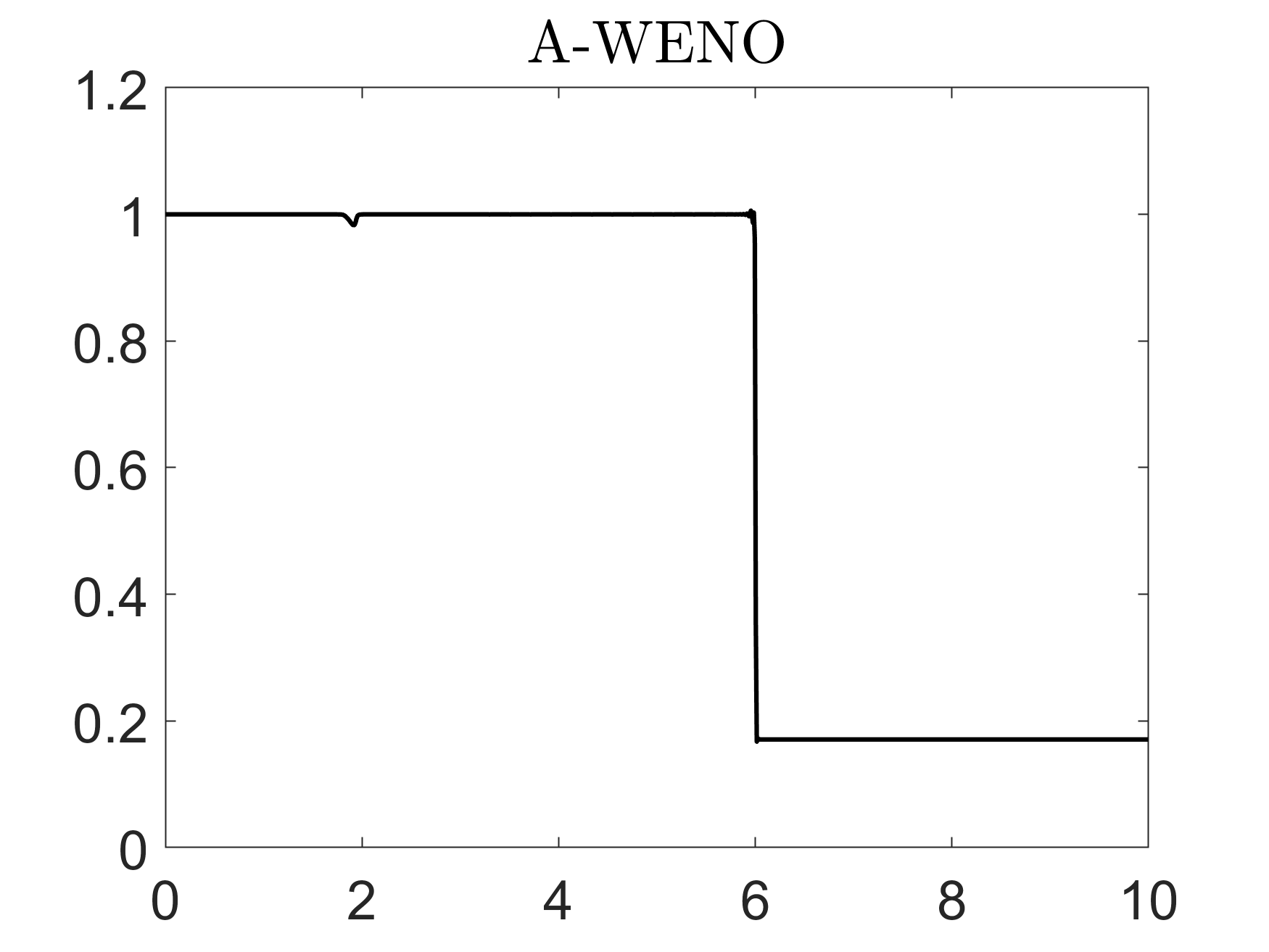}}
\caption{\sf Example 3: Water depth $h$ computed by the CU (left), RBM (middle), and A-WENO (right) schemes.\label{Fig5.12}}
\end{figure}

We compute the solutions by the studied CU, RBM, and A-WENO schemes until the final time $t=1$ on the computational domain $[0,10]$ using
1000, 2000, 4000, and 8000 uniform cells subject to the free boundary conditions. The water depths $h$ computed on 1000 uniform cells are
presented in Figure \ref{Fig5.12}. As one can see, all of the three studied schemes can capture the shock position correctly and there are
${\cal O}(1)$ oscillations in the neighborhoods of the shock in the RBM solution.

\begin{figure}[ht!]
\centerline{\includegraphics[trim=0.9cm 0.4cm 1.1cm 0.2cm, clip, width=5.4cm]{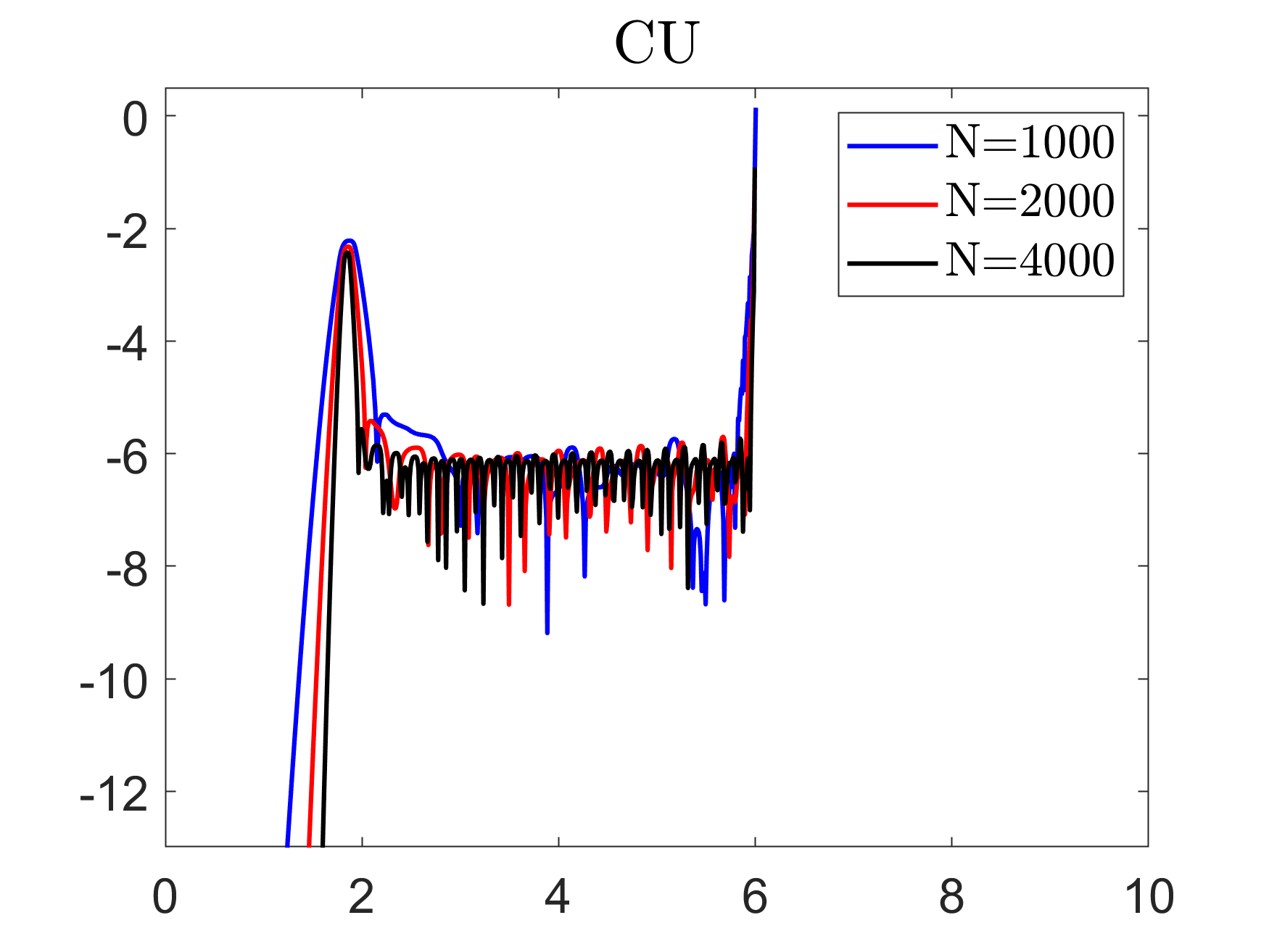}\hspace*{0.5cm}
\includegraphics[trim=0.9cm 0.4cm 1.1cm 0.2cm, clip, width=5.4cm]{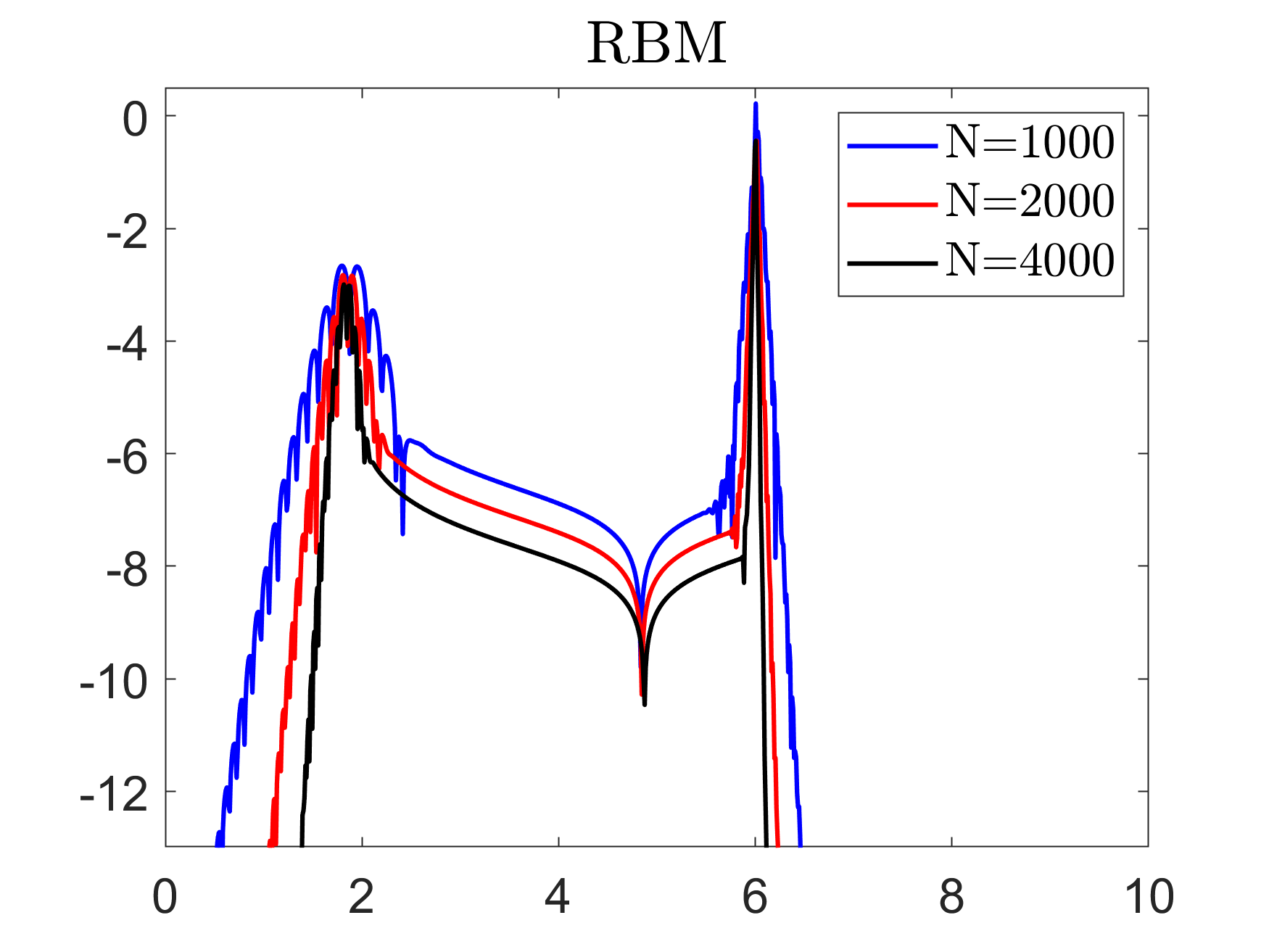}\hspace*{0.5cm}
\includegraphics[trim=0.9cm 0.4cm 1.1cm 0.2cm, clip, width=5.4cm]{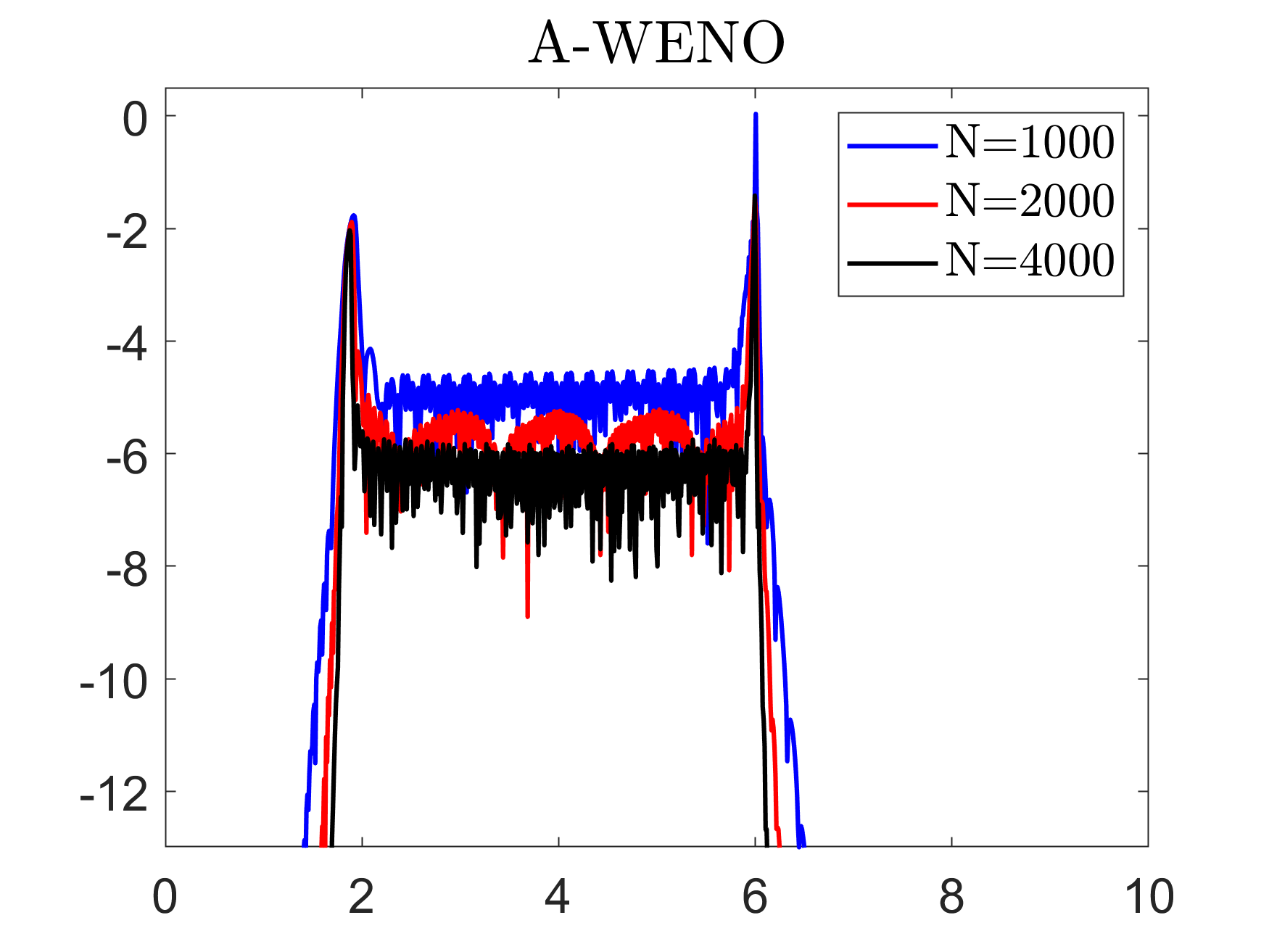}}
\caption{\sf Example 1: $\log_{10}|\frac{h^N-h^{\rm Exact}}{h^{\rm Exact}}|$, where $h^{\rm Exact}$ is the exact solution for the CU (left),
RBM (middle), and A-WENO (right) schemes.\label{FigNew}}
\end{figure}

We now study the error in the above computations. In Figure \ref{FigNew}, we present the relative differences between the computed and exact
water depths measured in the logarithmic scale. As one can see, the errors of the three studied schemes exhibit very different behavior. The
CU errors oscillate and do not seem to decay in the shock influence region as the mesh is refined, which means that only weak convergence
may be possible there. At the same time, the A-WENO errors do decay in that area while heavily oscillating, and the RBM errors decay in a
rather monotone way everywhere except for the vicinities of the shock (around $x=6$), ``postshock oscillations'' (around $x=2$, that is,
at the left edge of the shock influence region), and the initial breaking point $x=5$.

\begin{figure}[ht!]
\centerline{\includegraphics[trim=0.7cm 0.4cm 1.1cm 0.2cm, clip, width=5.3cm]{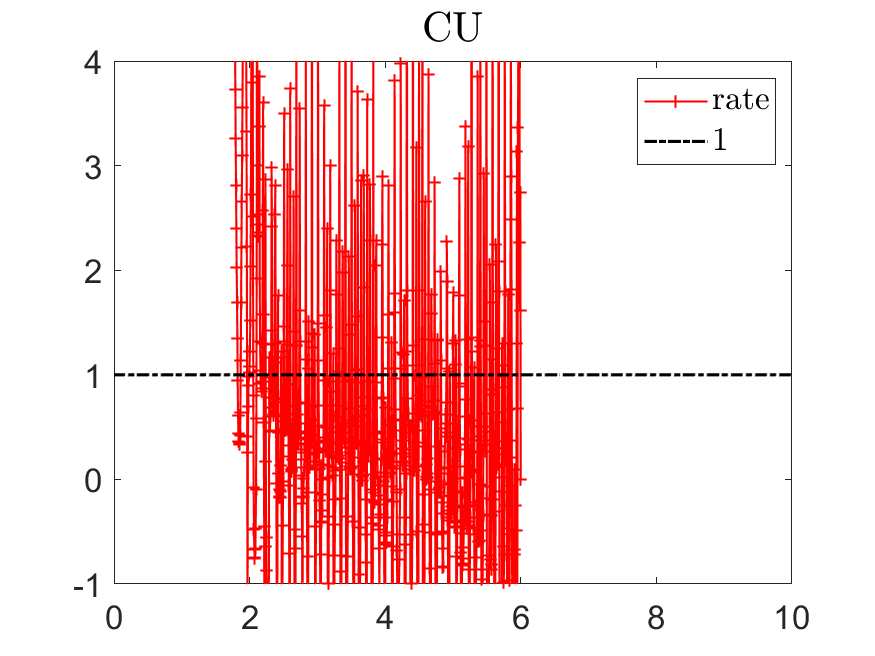}\hspace*{0.5cm}
            \includegraphics[trim=0.7cm 0.4cm 1.1cm 0.2cm, clip, width=5.3cm]{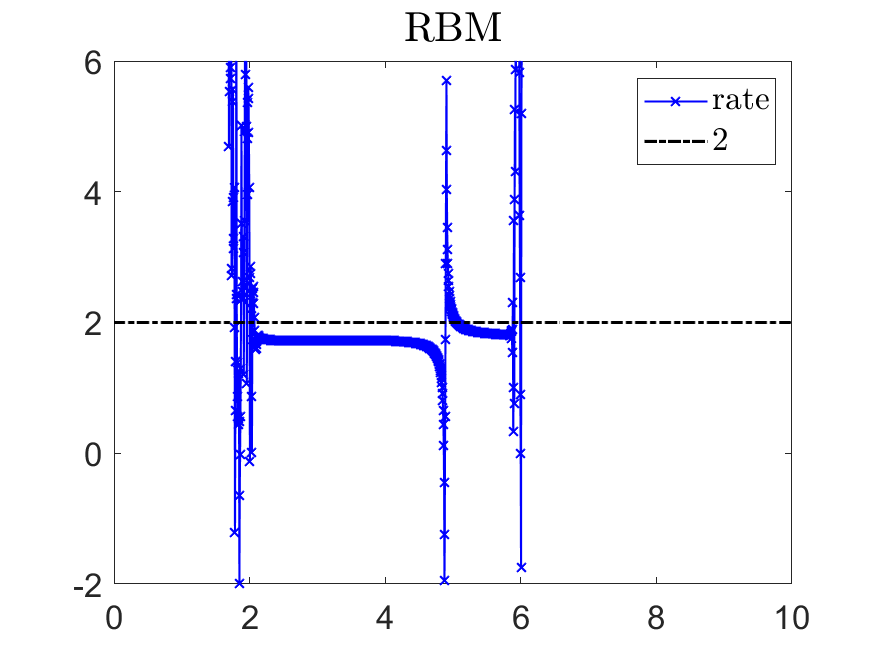}\hspace*{0.5cm}
            \includegraphics[trim=0.7cm 0.4cm 1.1cm 0.2cm, clip, width=5.3cm]{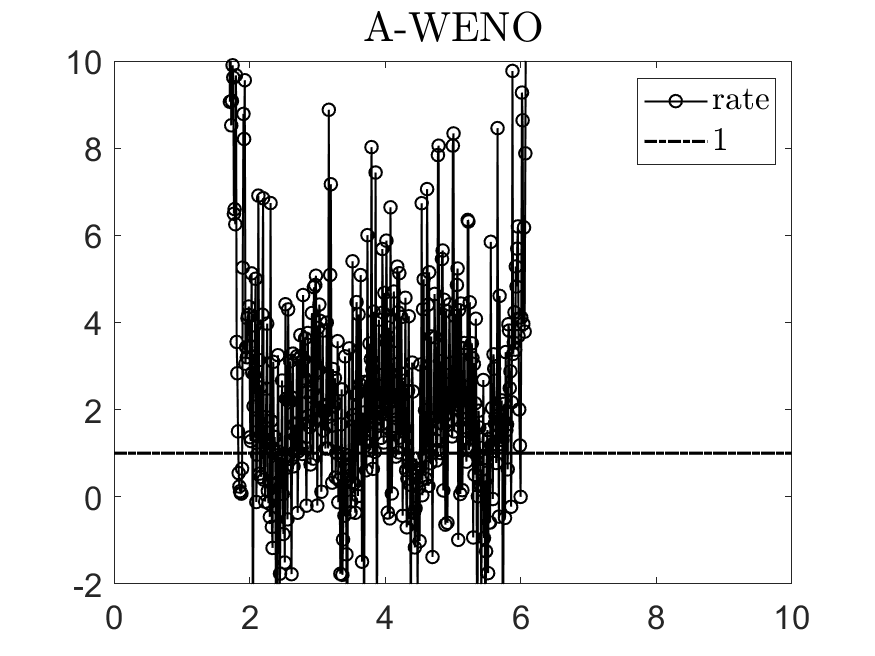}}
\caption{\sf Example 3: Experimental rates of pointwise convergence for the CU (left), RBM (middle), and A-WENO (right) schemes.
\label{Fig5.13a}}
\end{figure}

\begin{figure}[ht!]
\centerline{\includegraphics[trim=0.7cm 0.4cm 1.1cm 0.2cm, clip, width=5.3cm]{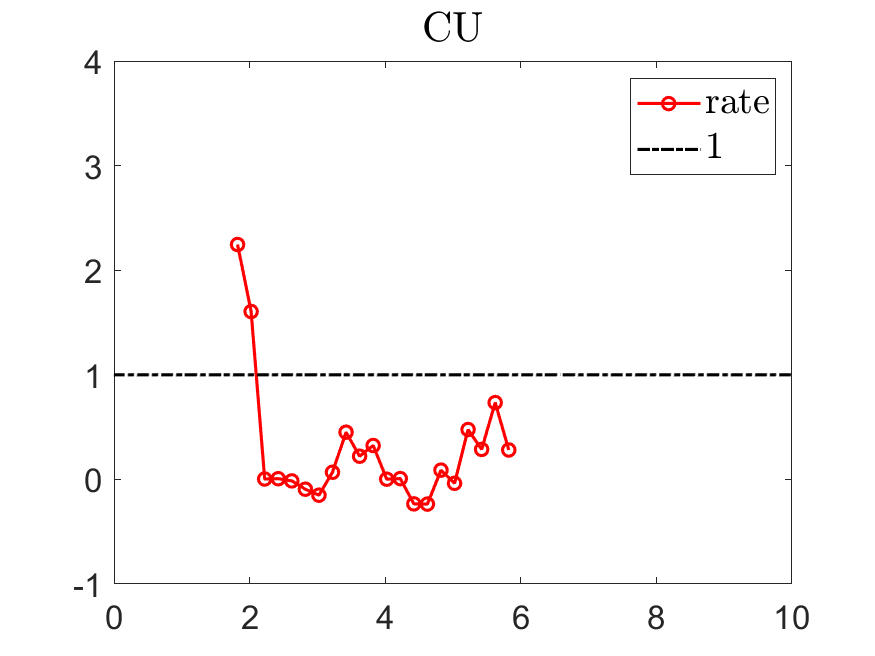}\hspace*{0.5cm}
            \includegraphics[trim=0.7cm 0.4cm 1.1cm 0.2cm, clip, width=5.3cm]{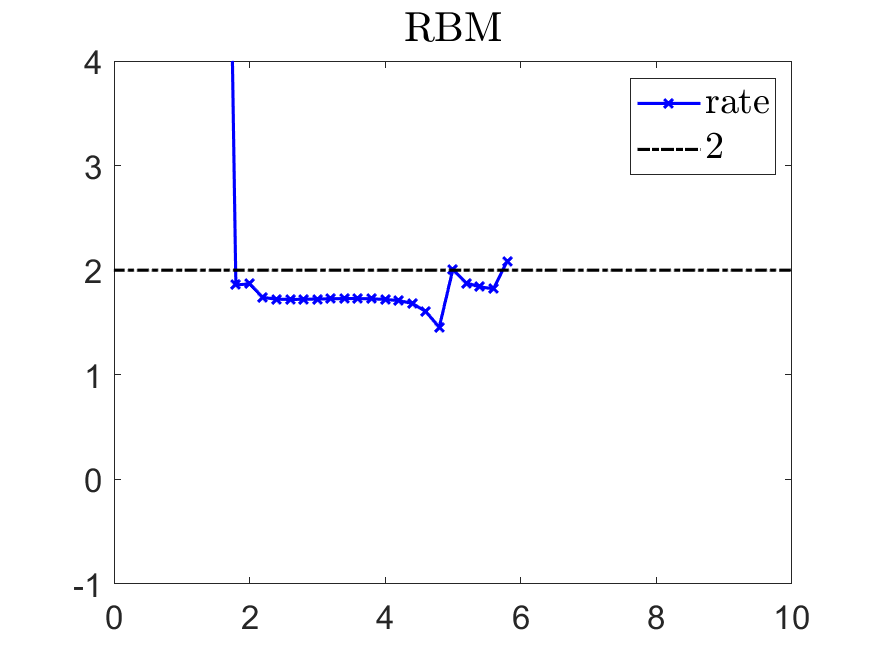}\hspace*{0.5cm}
            \includegraphics[trim=0.7cm 0.4cm 1.1cm 0.2cm, clip, width=5.3cm]{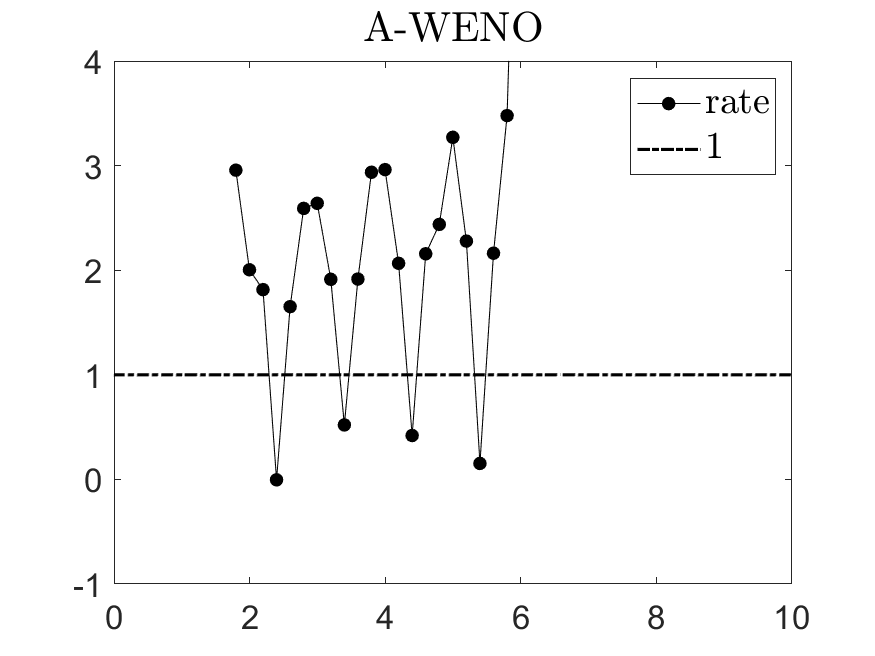}}
\caption{\sf Example 3: Average experimental rates of pointwise convergence for the CU (left), RBM (middle), and A-WENO (right) schemes.
\label{Fig5.13}}
\end{figure}

We then check the pointwise convergence of the water depth $h$ by computing the average experimental pointwise convergence rates defined in
\eref{equ5.6} with the local rates of convergence $r_{4j}$ computed by
\begin{equation}
r_{4j}:=\log_\hf\Bigg|\frac{h^{2N}_{4j}-h^{\rm Exact}_{4j}}{h^N_{4j}-h^{\rm Exact}_{4j}}\Bigg|,\quad j=0,\ldots,N,
\label{5.5a}
\end{equation}
with $N=2000$ for the CU and RBM schemes and with $N=1000$ for the A-WENO scheme. Here, $h_{4j}^{\rm Exact}$ is the corresponding values of
the exact solution. In Figure \ref{Fig5.13a}, we first plot the experimental convergence rates \eref{5.5a} for the water depth computed by
the studied schemes. As one can see, the rates are extremely oscillatory for the CU and A-WENO schemes. One can also see that the rates peak
both at the shock and ``postshock oscillations'' region, that is, around $x=6$ and $x=2$, respectively. We also stress that outside of the
shock influence region, the computed solutions remain constant and thus no convergence rates can be measured there. In order to better see
the convergence pattern, we plot $r^{\rm AVE}_{160k}$ for the CU and RBM schemes and $r^{\rm AVE}_{80k}$ for the A-WENO scheme for
$k=0,\ldots,50$ in Figure \ref{Fig5.13}. One can see that the CU scheme does not seem to convergence pointwise at all, the rate of
convergence for the RBM scheme is slightly below 2, while the rates for the A-WENO scheme seem to oscillate at about 2.

\begin{figure}[ht!]
\centerline{\includegraphics[trim=0.7cm 0.4cm 1.1cm 0.2cm, clip, width=5.3cm]{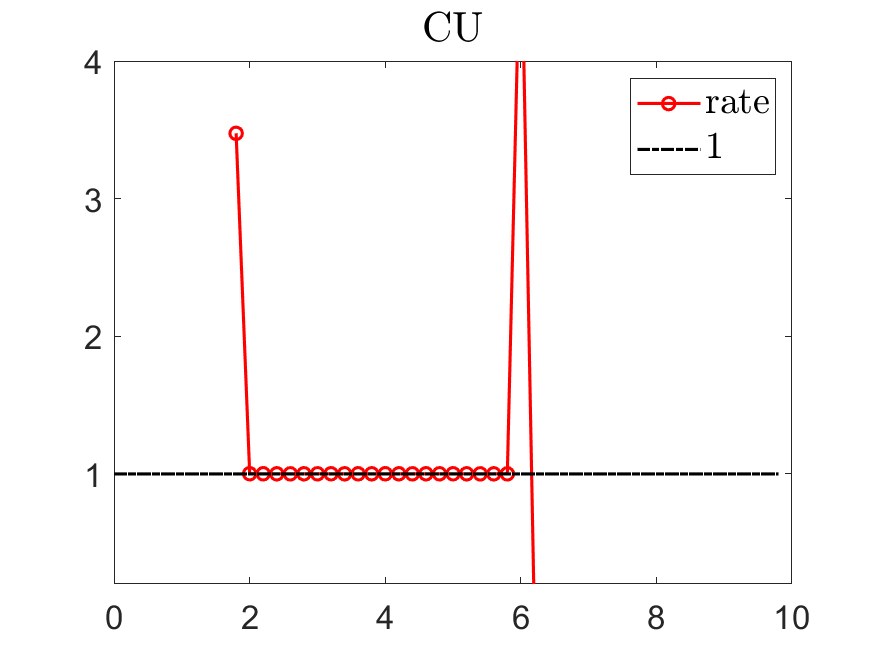}\hspace*{0.5cm}
            \includegraphics[trim=0.7cm 0.4cm 1.1cm 0.2cm, clip, width=5.3cm]{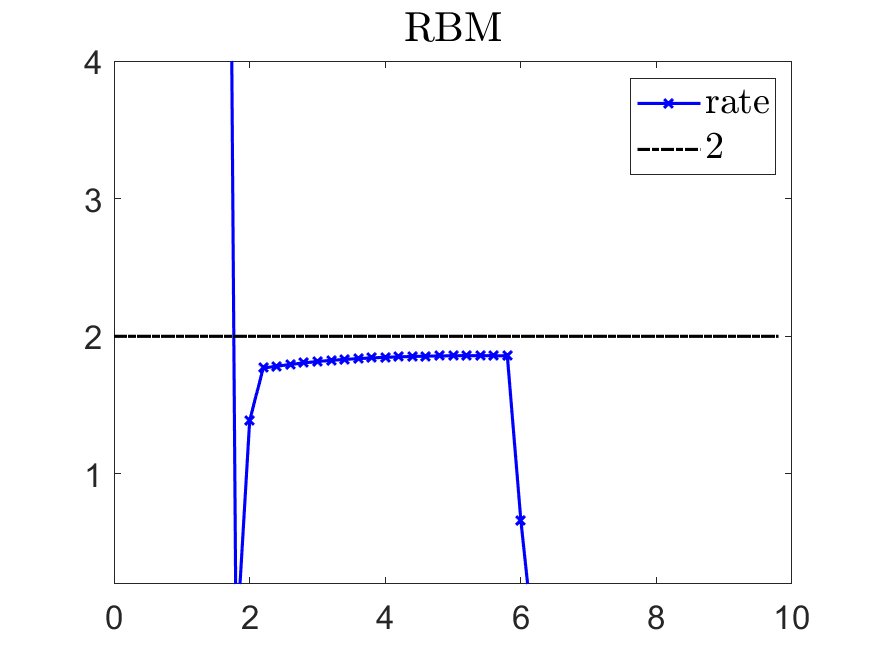}\hspace*{0.5cm}
            \includegraphics[trim=0.7cm 0.4cm 1.1cm 0.2cm, clip, width=5.3cm]{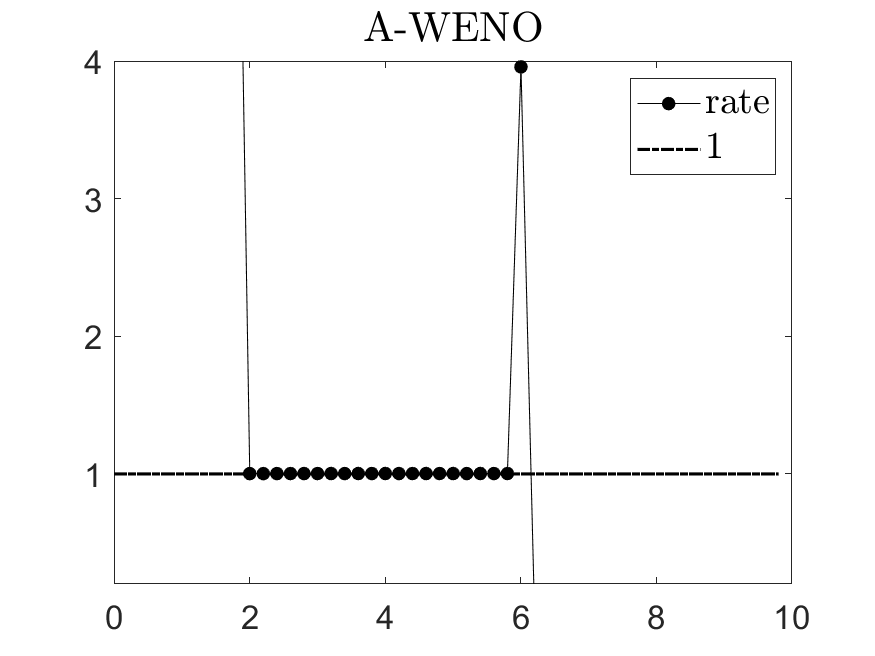}}
\caption{\sf Example 3: Experimental integral rates of convergence for the CU (left), RBM (middle), and A-WENO (right) schemes.
\label{Fig5.14}}
\end{figure}

We then compute the integral rates of convergence, which are now defined by
\begin{equation}
r^{\rm INT}_{4j}:=\log_\hf\Bigg|\frac{I^{2N}_{4j}-I^{\rm Exact}_{4j}}{I^N_{4j}-I^{\rm Exact}_{4j}}\Bigg|,\quad j=1,\ldots,N,
\label{5.5ab}
\end{equation}
where $I_{4j}^{\rm Exact}$ is the integral in \eref{4.4a} computed using the exact solution. In Figure \ref{Fig5.14}, we present
$r^{\rm INT}_{160k}$ for the CU and RBM schemes and $r^{\rm INT}_{80k}$ for the A-WENO scheme for $k=0,\ldots,50$. As one can see, unlike
the pointwise convergence rates, the integral ones  exhibit very uniform behavior (in the area where the computed solutions are no constant)
for all of the three studied schemes: it is first order for the CU or A-WENO schemes ans almost second order for the RBM scheme. Such a
uniform behavior of the integral rates is attributed to a very simple structure of the smooth parts of the solution in this example.

Finally, we compute the $W^{-1,1}$ convergence rates, which are now defined by
\begin{equation*}
r^{\rm INT}:=\log_\hf\bigg(\frac{\|I^{2N}-I^{\rm Exact}\|_{L^1}}{\|I^N-I^{\rm Exact}\|_{L^1}}\bigg),
\end{equation*}
where as in \eref{5.5ab}, $I^{\rm Exact}$ represents the quantities computed by evaluating the integrals in \eref{4.4a} using the exact
solution. The obtained results are reported in Table \ref{Tab3}, where one can clearly see that the studied three schemes drop to either
first (CU and A-WENO) or second (RBM) order. The $W^{-1,1}$ convergence rates for all of the studied schemes are the same as in Examples 1
and 2 at larger times, that is, after the shock formation: the CU and A-WENO schemes are first-order, while the RBM scheme is second-order
accurate.
\begin{table}[ht!]
\centering
\begin{tabular}{|c|c|c|c|c|c|c|c|}\hline
&\multicolumn{2}{c|}{CU Scheme}&\multicolumn{2}{c|}{RBM Scheme}&\multicolumn{2}{c|}{A-WENO Scheme}\\\hline
$N$&$||I^N-I^{Exact}||_{L^1}$&$r^{\rm INT}$&$||I^N-I^{Exact}||_{L^1}$&$r^{\rm INT}$&$||I^N-I^{Exact}||_{L^1}$&$r^{\rm INT}$\\\hline
2000&2.39e-3&---  &3.00e-5&---  &2.85e-3&---\\ \hline
4000&1.19e-3&1.00 &7.47e-6&2.01 &1.42e-3&1.01\\ \hline
8000&5.97e-4&1.00 &1.86e-6&2.01 &7.12e-4&1.00\\
\hline
\end{tabular}
\caption{\sf Example 3: $W^{-1,1}$ convergence rates for the studied three schemes.\label{Tab3}}
\end{table}

\section{Combined Schemes}\label{sec7}
In this section, we introduce two new combined schemes for the hyperbolic system of conservation laws \eref{1.1}. The proposed schemes are
based on the RBM as a basic scheme and either the CU or A-WENO as an internal scheme.

We develop the combined schemes on a particular spatial grid and use the notation introduced in \S\ref{sec2a} and Appendices
\ref{appa}--\ref{appc} for the finest mesh of size $\dx$. We note that the CU scheme is a FV scheme and thus the evolved quantities in the
CU scheme are the cell averages $\xbar\mV_\jph(t)$. However, as the CU scheme is second-order, the point values $\mV_\jph(t)$ are within the
accuracy of the scheme from the corresponding cell averages and thus one may omit the $\xbar{(\cdot)}$ notation and view the CU scheme from
Appendix \ref{appa} as a FD scheme.

The key idea in the construction of the combined schemes (see, e.g., \cite{KO18,LNOT19,ZKO}) is to have two copies of the computed
solutions: one of them is the RBM solution, which is computed throughout the entire computational domain and which is highly accurate in the
smooth areas, and the second one is the non-oscillatory solution computed with the help of either the CU or A-WENO scheme near the shock
regions only.

We therefore start with the evolution of the RBM solution. As the RBM data from the areas near shocks will be used in the evolution of the
CU/A-WENO solutions, we need to have the RBM and CU/A-WENO solutions at the same time levels. Recall that the CU/A-WENO solutions are
evolved using the three-stage third-order SSP Runge-Kutta method, which can be written in the operator form as follows. Let us denote the
computed values of the CU/A-WENO solutions by $\mW_\jph(t)$. Let $\mW$ stands for $\mW(t):=\{\mW_{\ell+\hf}(t)\}$, where $\ell$ are the
indices at which the RBM solution, denoted by $\mV_{\ell+\hf}(t)$, has been detected as being inside the shock neighborhood, and several of
the nearby points depending on the size of the CU/A-WENO stencil (at these points we will not have the values $\mW_{\ell+\hf}(t)$, but will
simply set $\mW_{\ell+\hf}(t)=\mV_{\ell+\hf}(t)$). Let ${\cal P}[\mW(t)]$ be the nonlinear operator representing the RHS of either \eref{A2}
or \eref{C1}. Then, the solution is evolved from the time level $t^n$ to the time level $t^{n+1}:=t^n+\dt^n$ according to
\begin{equation}
\begin{aligned}
&\mW^{\rm I}(t^{n+1})=\mW(t^n)+\dt^n{\cal P}[\mW(t^n)],\\
&\mW^{\rm II}(t^{n+\hf})=\frac{3}{4}\mW(t^n)+\frac{1}{4}\left(\mW^{\rm I}(t^{n+1})+\dt^n{\cal P}[\mW^{\rm I}(t^{n+1})]\right),\\
&\mW(t^{n+1})=\frac{1}{3}\mW(t^n)+\frac{2}{3}\left(\mW^{\rm II}(t^{n+\hf})+\dt^n{\cal P}[\mW^{\rm II}(t^{n+\hf})]\right),
\end{aligned}
\label{7.1}
\end{equation}
where $t^{n+\hf}:=t^n+\dt^n/2$, and $\mW^{\rm I}(t^{n+1})$ and $\mW^{\rm II}(t^{n+\hf})$ are lower-order intermediate solutions at the time
levels $t^{n+1}$ and $t^{n+\hf}$, respectively. The time step $\dt^n$ can be selected using the CFL number $1/2$, namely, by setting
\begin{equation}
\dt^n\le\frac{\dx}{2a^n},\quad
a^n=\max_j\left\{\max\left(\lambda_d\big(A(\mV_\jph(t^n))\big),-\lambda_1\big(A(\mV_\jph(t^n))\big)\right)\right\},
\label{7.2}
\end{equation}
where $a^n$ is the upper bound on the magnitude of the local speeds of propagation at the time level $t=t^n$.

We first detect nonsmooth points of the computed solution at the time level $t=t^n$. To this end, we use the WLR, introduced
in \cite{Karni05,Karni02} and given by
\begin{equation*}
\begin{aligned}
\mE_\jph^{n}=&\frac{1}{12}\Big\{\left[\mV_{j+\frac{3}{2}}^{n+1}-\mV_{j+\frac{3}{2}}^{n-1}+4\left(\mV_\jph^{n+1}-\mV_\jph^{n-1}\right)+
\mV_\jmh^{n+1}-\mV_\jmh^{n-1}\right]\dx\\
&+\left[\mF_{j+\frac{3}{2}}^{n+1}-\mF_\jmh^{n+1}+4\left(\mF_{j+\frac{3}{2}}^{n}-\mF_\jmh^{n}\right)+
\mF_{j+\frac{3}{2}}^{n-1}-\mF_\jmh^{n-1}\right]\dt^n\Big\},
\end{aligned}
\end{equation*}
where $\mV_\jph^i:=\mV_\jph(t^i)$ and $\mF_\jph^i:=\mF(\mV_\jph(t^i))$, $i=n-1$, $n$, $n+1$. As it was shown in \cite{Karni02,Karni05}, the
magnitude of $\mE_\jph^{n}$ is proportional to $\dx$ near the shock while being of order $(\dx)^5$ in the smooth parts of the solution. In
order to take advantage of the big discrepancy in these values of the WLR, we select an $i$-th component of $\mE$ ($i=1$ for the
Saint-Venant system \eref{5.1}), introduce the quantity
\begin{equation*}
\varepsilon^{n}_\jph=\max\left(\big|(E^{(i)})_\jmh^{n}\big|,\big|(E^{(i)})_\jph^{n}\big|,
\big|(E^{(i)})_{j+\frac{3}{2}}^{n}\big|\right),
\end{equation*}
and determine the local smoothness of $\mV^n_\jph$ based on the size of $ \varepsilon^n_\jph$. More precisely, we say that if
$\varepsilon^n_\jph>\mu(\dx)^3$, where $\mu$ is a tunable constant, then the point $x=x_\jph$ is in the ``rough'' part of
computed solution denoted by $\Omega^n:=\big\{x_\jph~\big|~\varepsilon^n_\jph>\mu(\dx)^3\big\}$.

In order to be able to use \eref{7.1}, we would need the RBM solution not only at the time levels $t=t^n$ and $t^{n+1}$, but also at the
intermediate time level $t=t^{n+\hf}$.  To this end, we compute the RBM solution $\mV(t^{n+\hf})$ at the time level $t^{n+\hf}$
using a CFL number 0.25, but only in the immediate vicinity of the detected ``rough" areas. Equipped with the RBM solutions $\mV(t^n)$,
$\mV(t^{n+\hf})$, and $\mV(t^{n+1})$, we evolve the CU/A-WENO solutions from $t^n$ to $t^{n+1}$ inside $\Omega^n$ and in the nearby points
$x_\jph$ such that either $x_\jmh\in\Omega^n$ or $x_{j+\frac{3}{2}}\in\Omega^n$. As we have already mentioned, there might be points
$x_{\ell+\hf}\in\Omega^n$, at which the values $\mW^n_{\ell+\hf}$ are unavailable at the time level $t^n$, $t^{n+\hf}$, and $t^{n+1}$. We
then set $\mW^n_{\ell+\hf}=\mV^n_{\ell+\hf}$ at these points. There might be also points $x_\jph\notin\Omega^n$, at which $\mW^n_\jph$ are
available. At these points the computed solution is smooth and hence we replace the CU/A-WENO solution values with more accurate RBM ones by
setting $\mW^n_\jph=\mV^n_\jph$.
\begin{rmk}
The value of $\mu$ can be selected at each example experimentally. As in \cite{Karni02}, we tune it on a coarser mesh and then use the same
$\mu$ at finer meshes. This strategy has been proved to be robust; see, e.g., \cite{ConK,DKL,Karni02}.
\end{rmk}
\begin{rmk}
At the final computational time, the solution computed by the combined scheme will consist of the ``pure'' RBM solution values everywhere
except for the shock areas where it will be replaced with the CU/A-WENO solutions. It should be observed, however, that this CU/A-WENO part
of the solution is calculated using the RBM data at each evolution step.
\end{rmk}

\subsection{Numerical Examples}

In this section, we apply the developed combined schemes based on the RBM--CU and RBM--A-WENO schemes to Examples 1--3 from \S\ref{sec5}.
These combined schemes will be referred to as the RBM--CU and RBM--A-WENO schemes below.

\subsubsection*{Example 4---Test with One Shock}
\begin{figure}[ht!]
\centerline{\includegraphics[trim=1.3cm 0.3cm 1.0cm 0.8cm, clip, width=5.3cm]{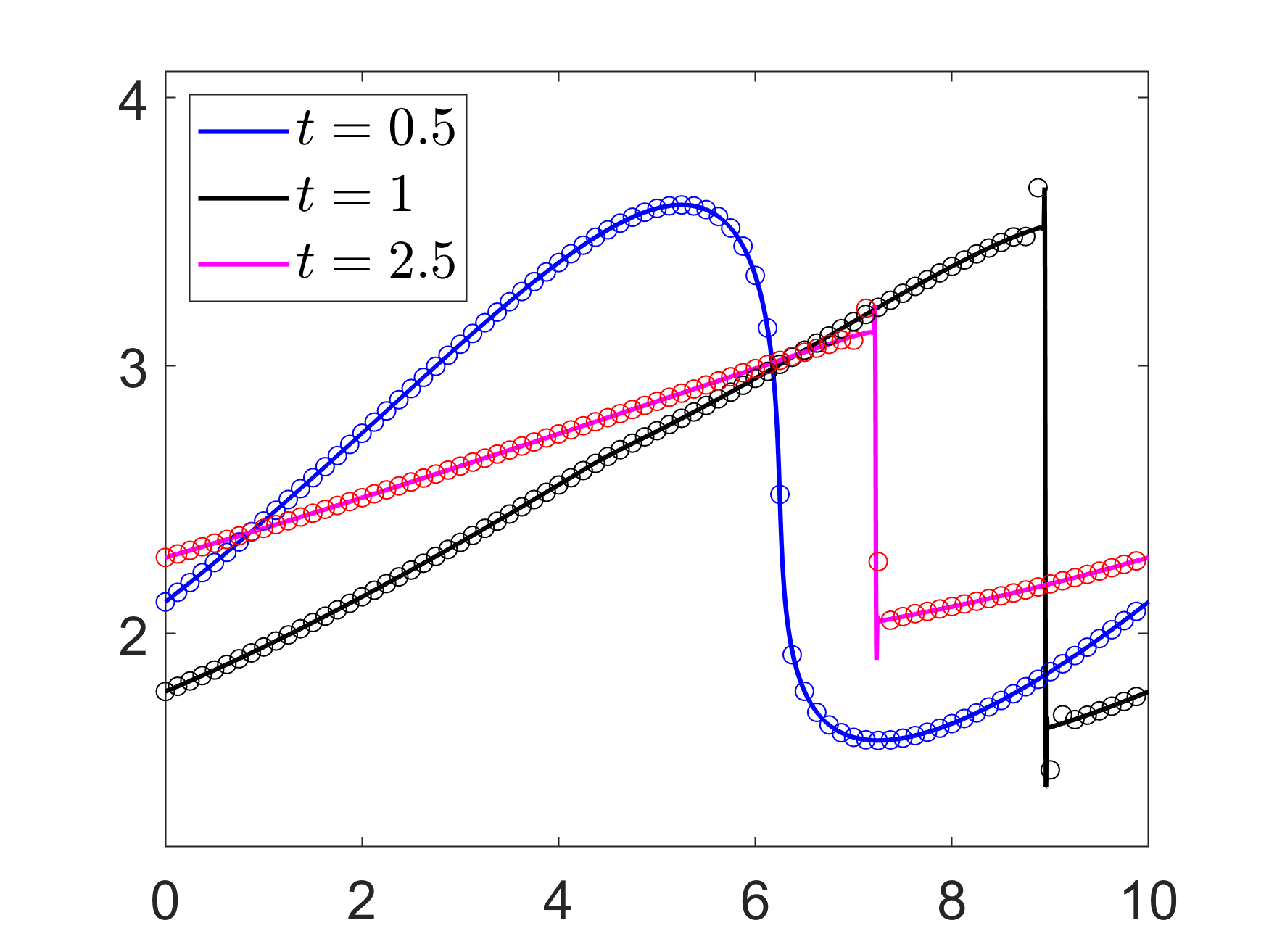}\hspace*{0.5cm}
            \includegraphics[trim=1.3cm 0.3cm 1.0cm 0.8cm, clip, width=5.3cm]{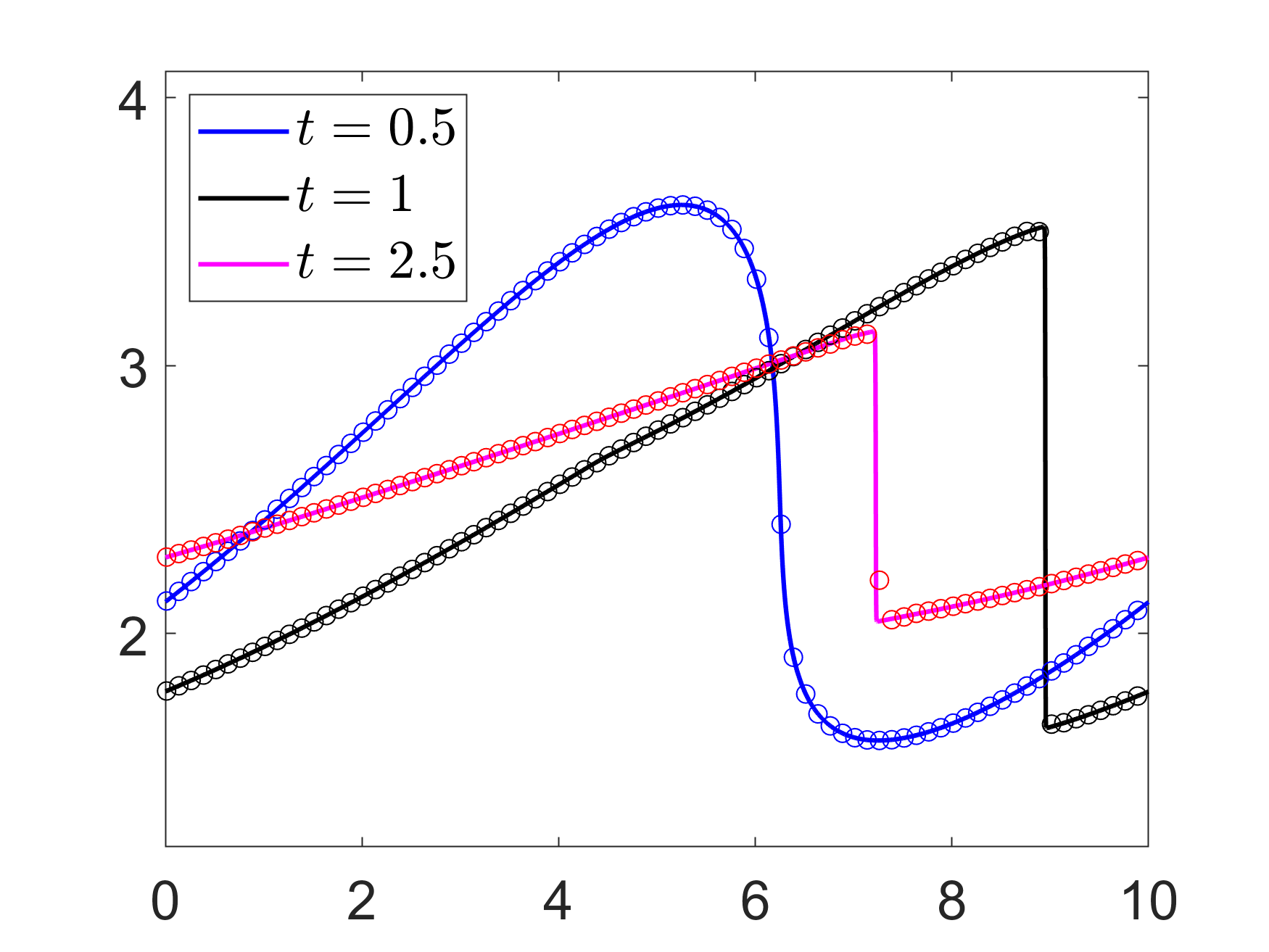}\hspace*{0.5cm}
            \includegraphics[trim=1.3cm 0.3cm 1.0cm 0.8cm, clip, width=5.3cm]{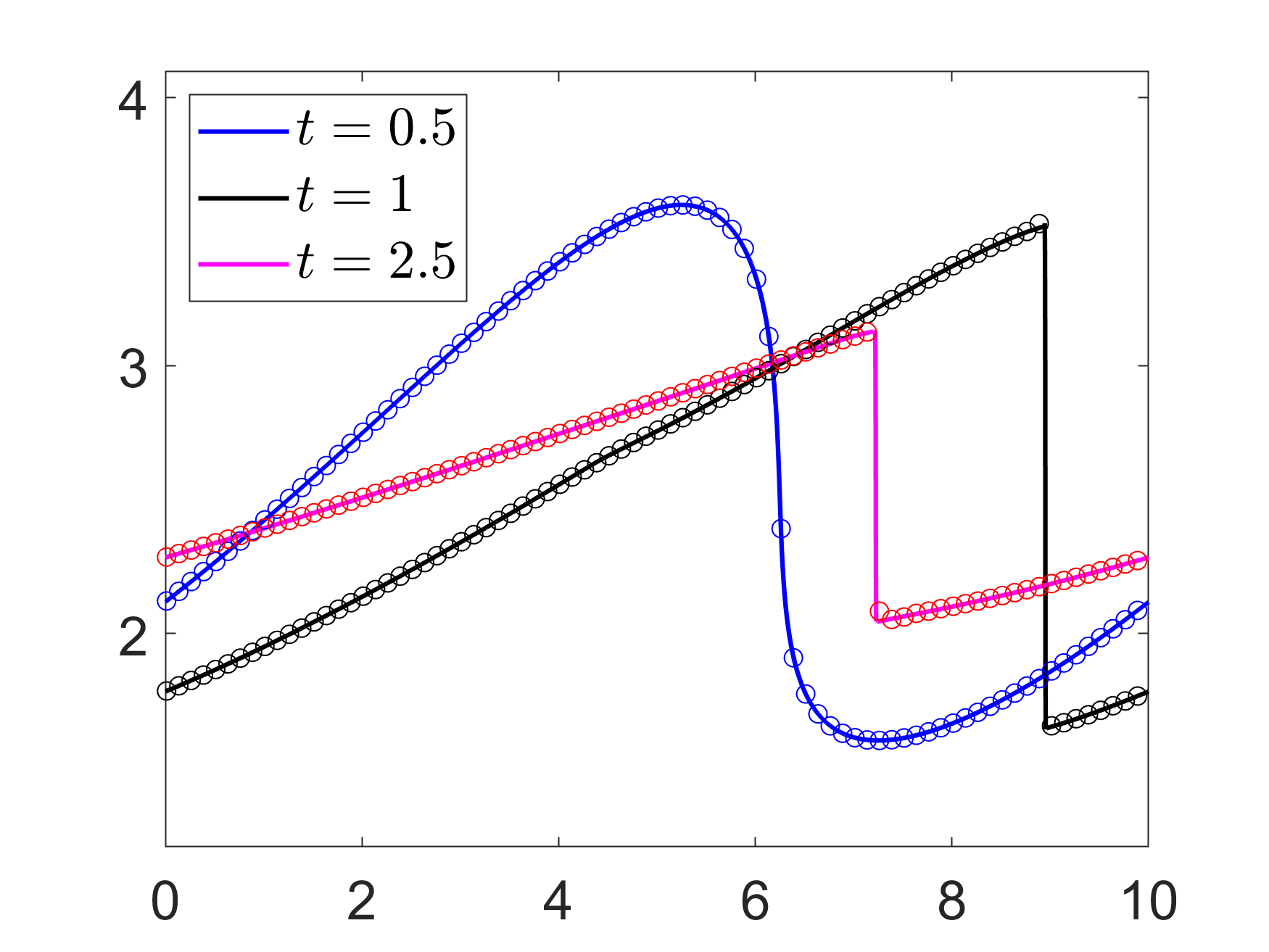}}
\caption{\sf Example 4: Water depth $h$ computed by the RBM (left), RBM--CU (middle), and RBM--A-WENO (right) schemes on the coarse (circles) and fine (solid lines) grids.\label{fig71}}
\end{figure}

In this example, we use the same settings as in Example 1 and compute the solutions by the studied RBM--CU and RBM--A-WENO schemes at times
$t=0.5$, 1, and 2.5 on the computational domain $[0,10]$ using 1000, 2000, 4000, and 8000 uniform cells. Here, we first tune the coefficient
$\texttt{C}$ on a coarse mesh with 400 uniform cells and then use it for the finer meshes. We take $\mu=0.2$ for both the RBM--CU and
RBM--A-WENO schemes. The results computed on the grids with 400 and 4000 uniform cells, are presented in Figures \ref{fig71} along with the
``pure'' RBM solutions. In order to improve the visibility, we plot the values $h_{5\jph}$, $j=0,\ldots,80$ for the numerical results
computed on the coarser mesh. As one can see, both combined schemes produce non-oscillatory solutions and the transitions between the RBM
and CU/A-WENO parts in the combined solutions are smooth.

\begin{figure}[ht!]
\centerline{\includegraphics[trim=1.3cm 0.3cm 1.0cm 0.1cm, clip, width=5.3cm]{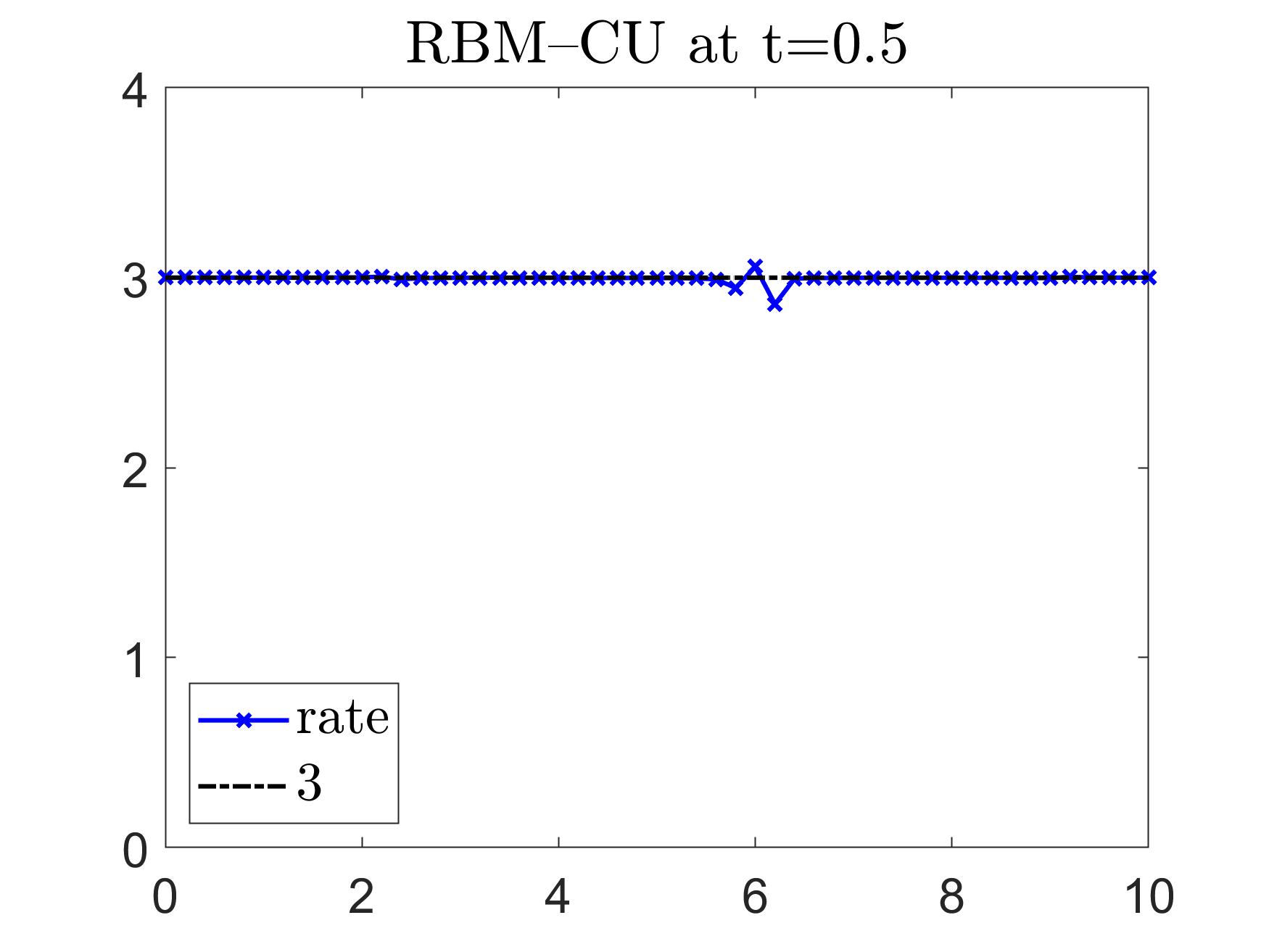}\hspace*{0.5cm}
            \includegraphics[trim=1.3cm 0.3cm 1.0cm 0.1cm, clip, width=5.3cm]{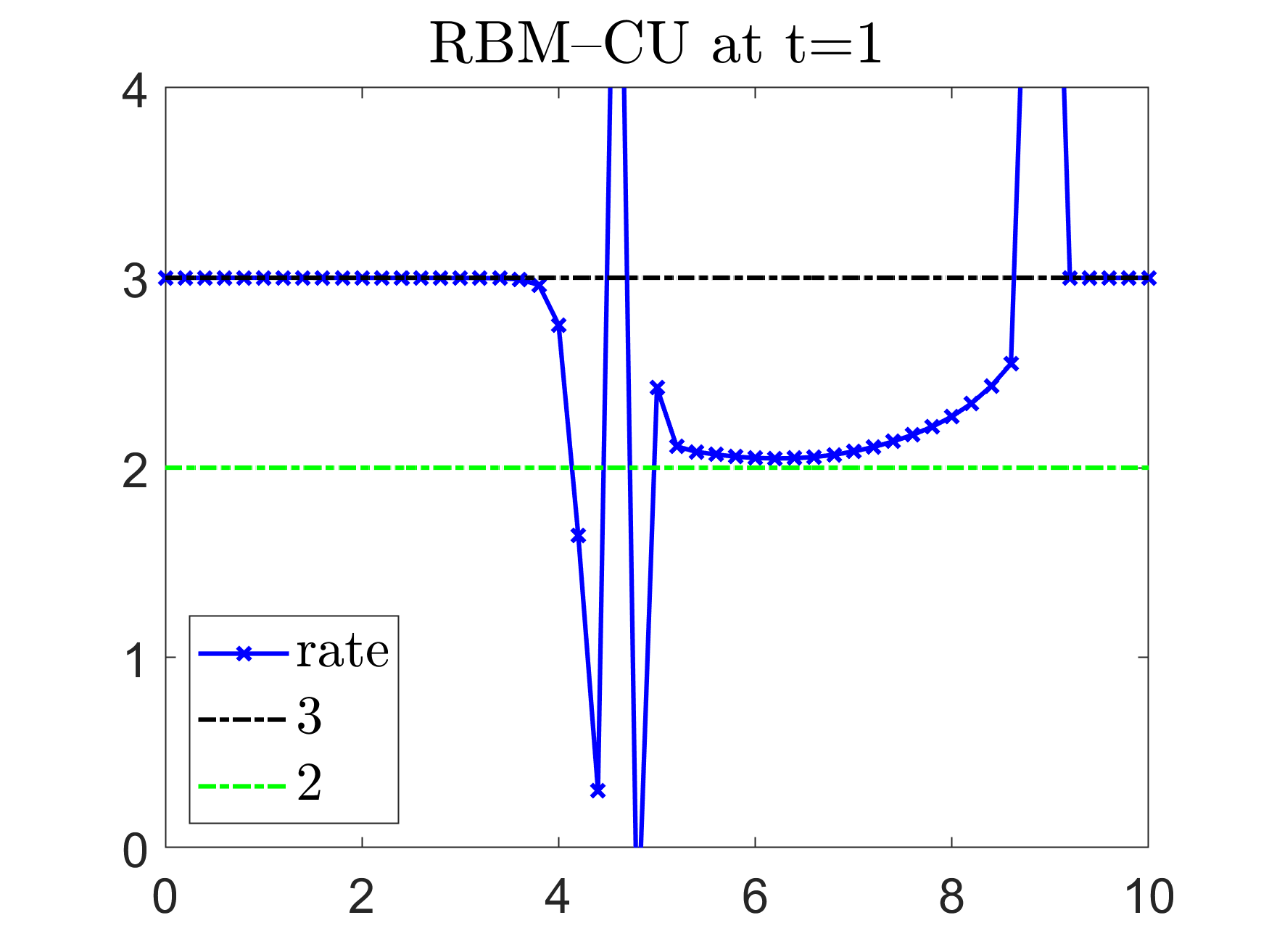}\hspace*{0.5cm}
            \includegraphics[trim=1.3cm 0.3cm 1.0cm 0.1cm, clip, width=5.3cm]{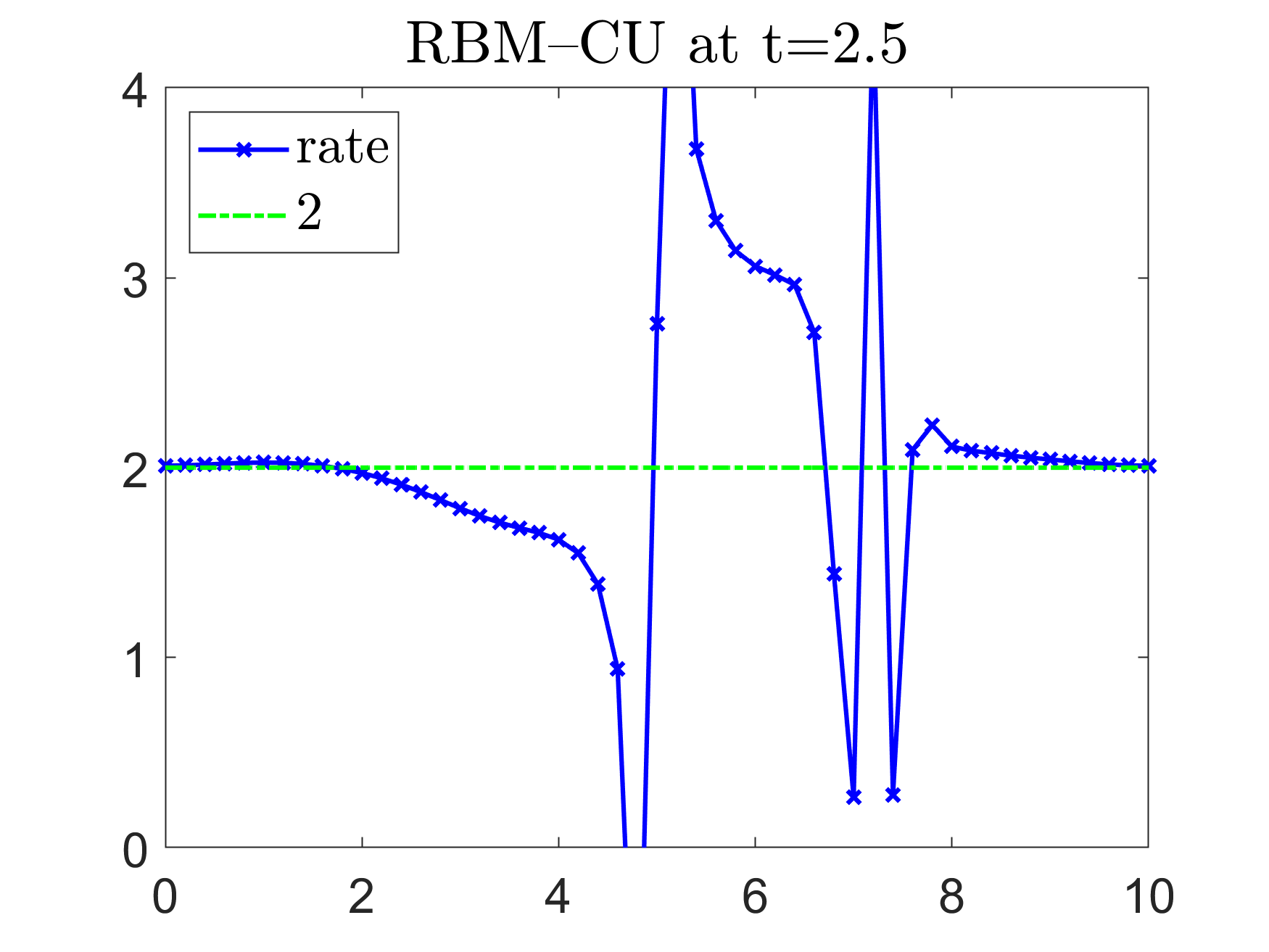}}
\vskip8pt
\centerline{\includegraphics[trim=1.3cm 0.3cm 1.0cm 0.1cm, clip, width=5.3cm]{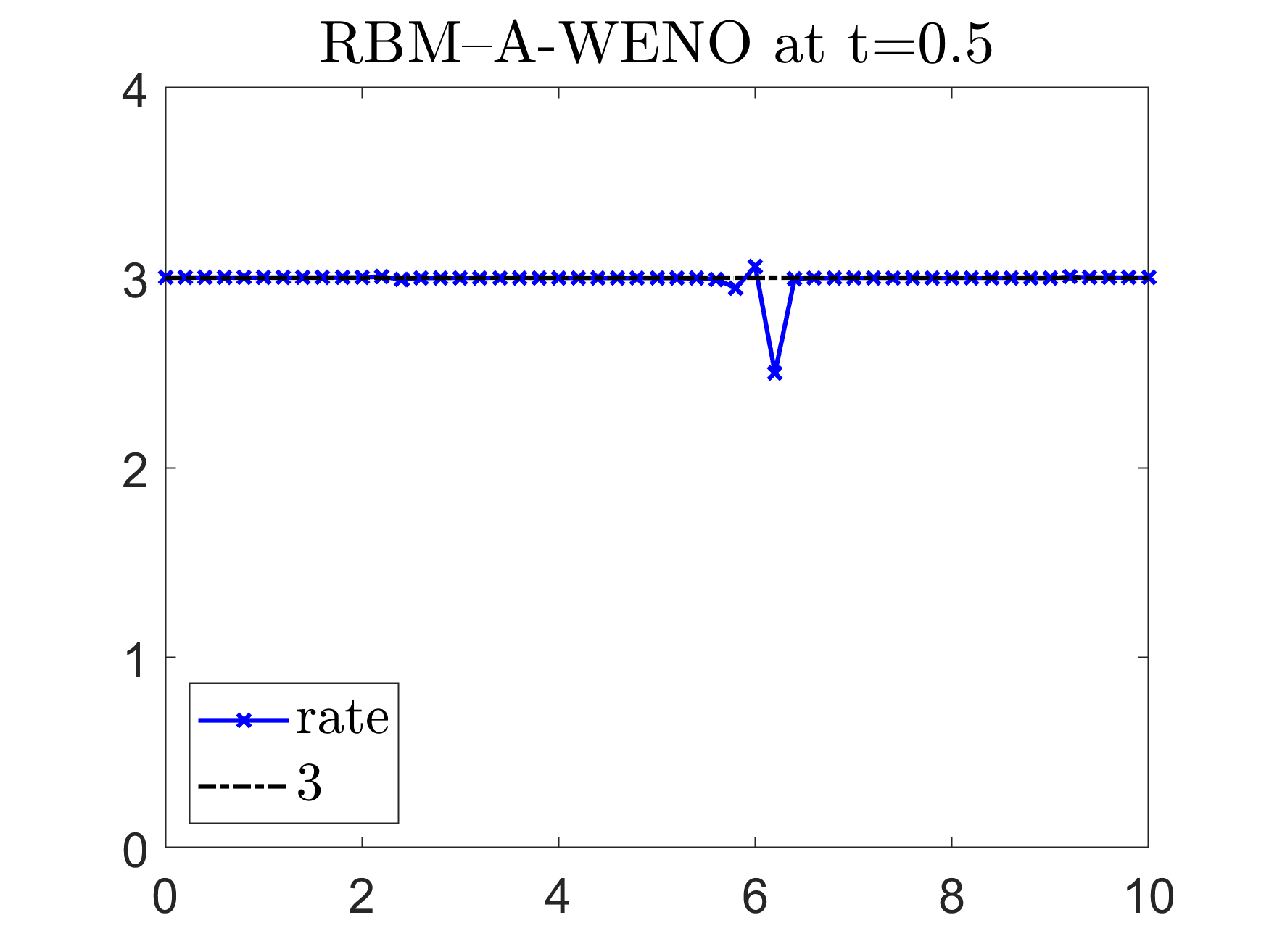}\hspace*{0.5cm}
            \includegraphics[trim=1.3cm 0.3cm 1.0cm 0.1cm, clip, width=5.3cm]{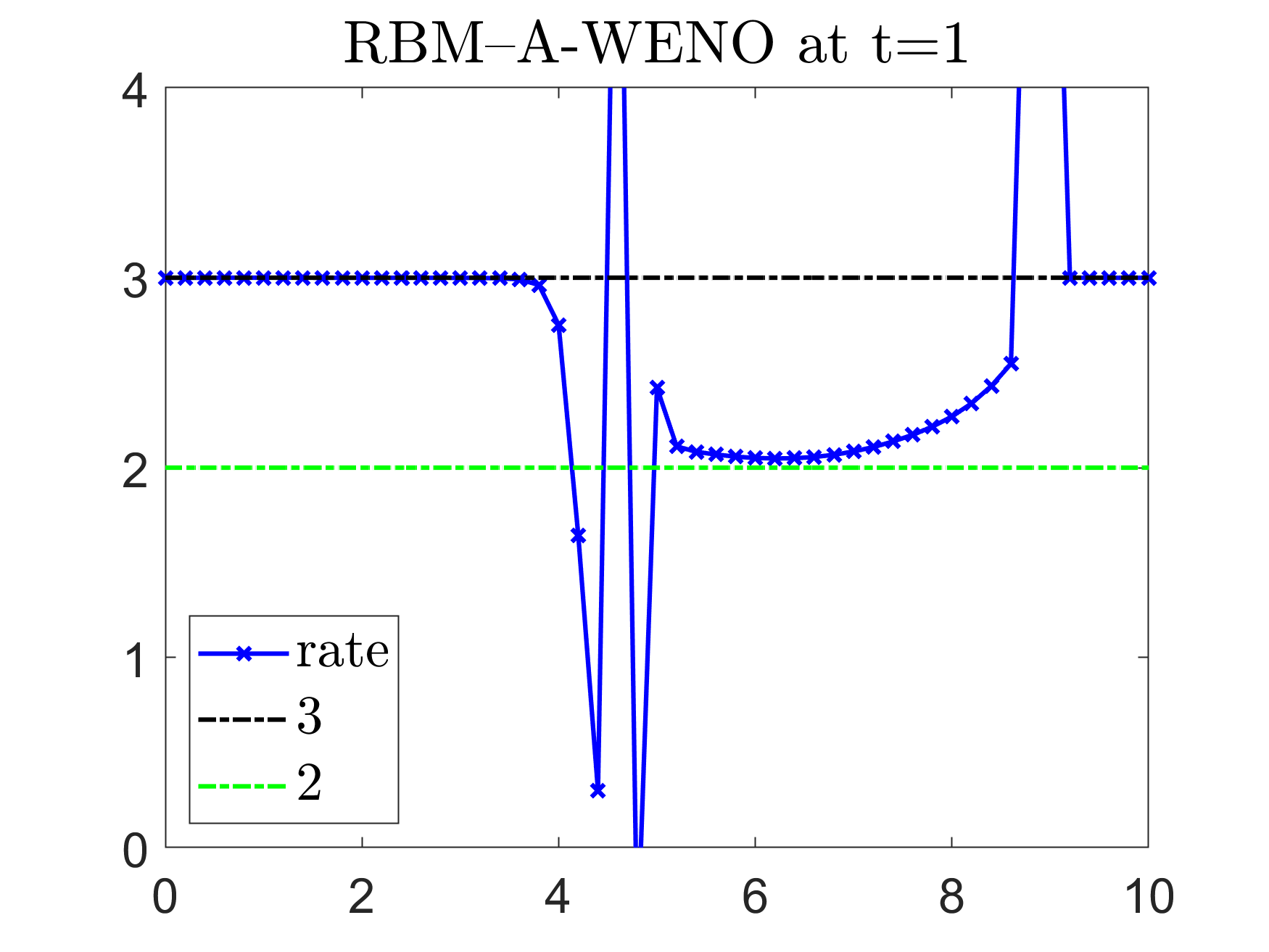}\hspace*{0.5cm}
            \includegraphics[trim=1.3cm 0.3cm 1.0cm 0.1cm, clip, width=5.3cm]{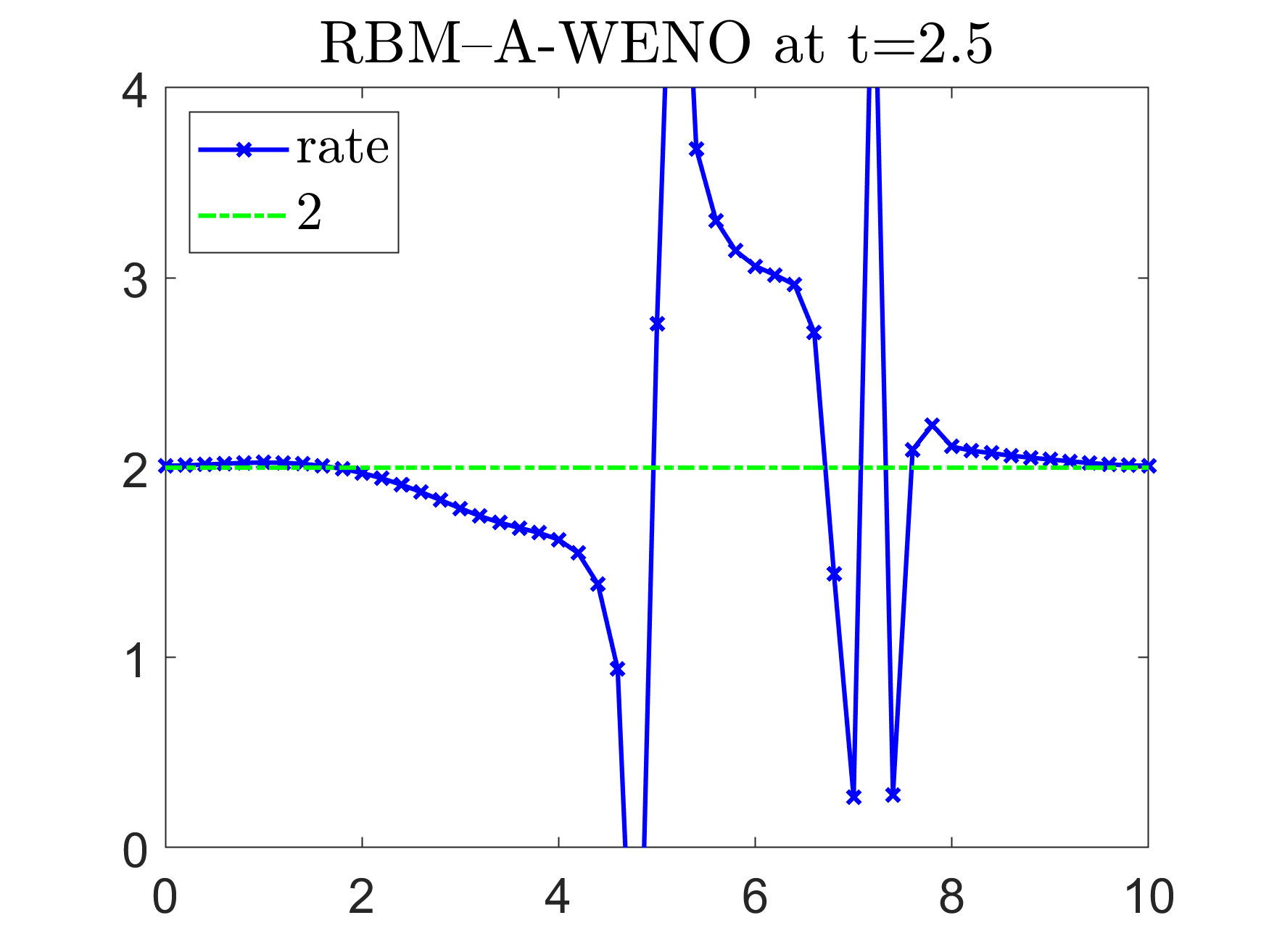}}
\caption{\sf Example 4: Average experimental rates of pointwise convergence for the RBM--CU (top row) and RBM--A-WENO (bottom row) schemes.
\label{Fig6.2}}
\end{figure}

We then check the pointwise convergence of the water depth $h$ by computing the average experimental pointwise convergence rates defined in
\eref{equ5.6} using three imbedded grids with $N=2000$ for the RBM--CU and RBM--A-WENO schemes. We plot the obtained values
$r^{\rm AVE}_{160k}$ in Figure \ref{Fig6.2}, where one can clearly see that the average convergence rates for the RBM--CU and RBM--A-WENO
schemes are close to those of the RBM scheme; compare Figure \ref{Fig6.2} and the middle row panels of Figure \ref{Fig4.5}.

\begin{figure}[ht!]
\centerline{\includegraphics[trim=1.3cm 0.3cm 1.0cm 0.1cm, clip, width=5.3cm]{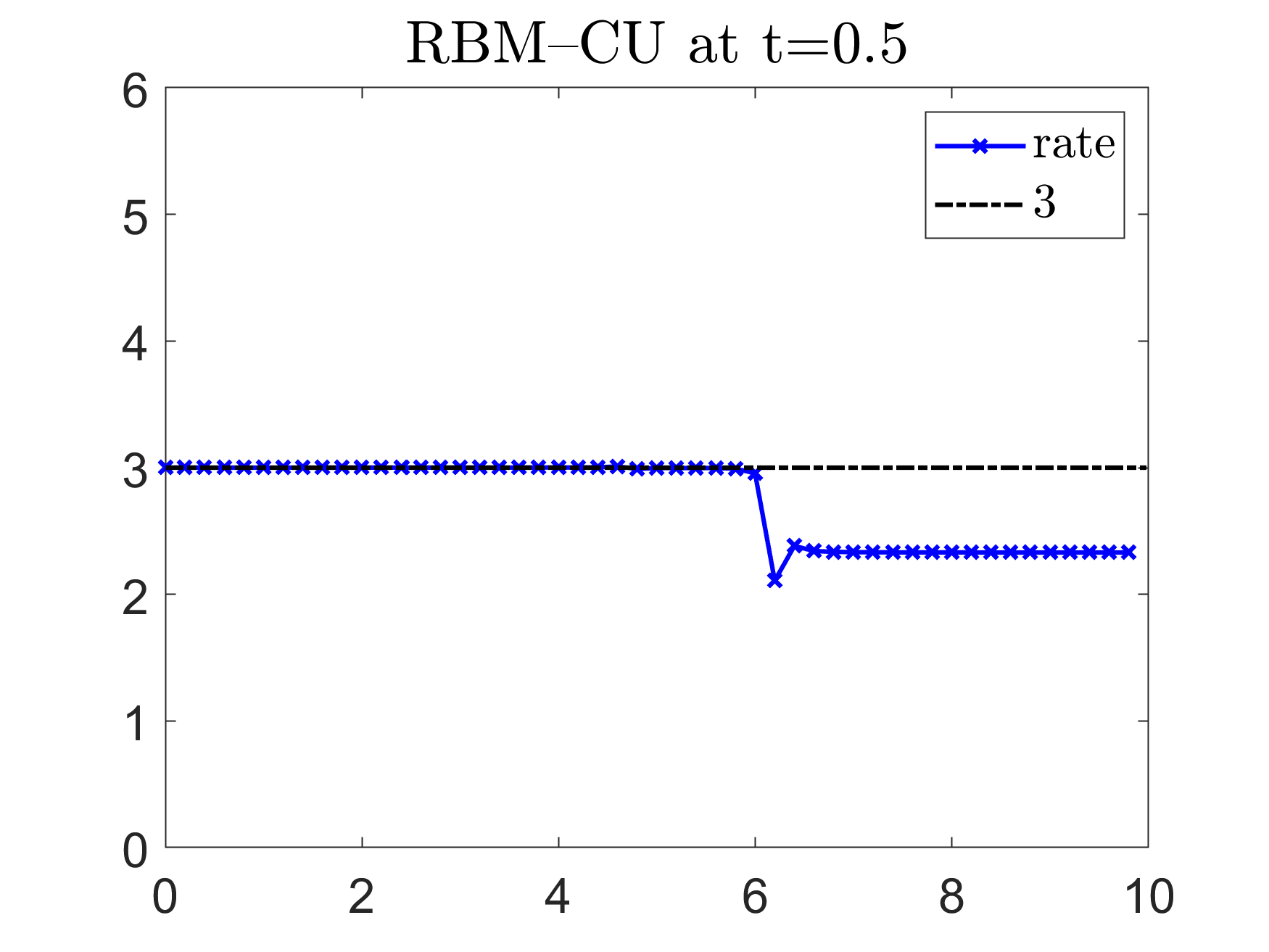}\hspace*{0.5cm}
            \includegraphics[trim=1.3cm 0.3cm 1.0cm 0.1cm, clip, width=5.3cm]{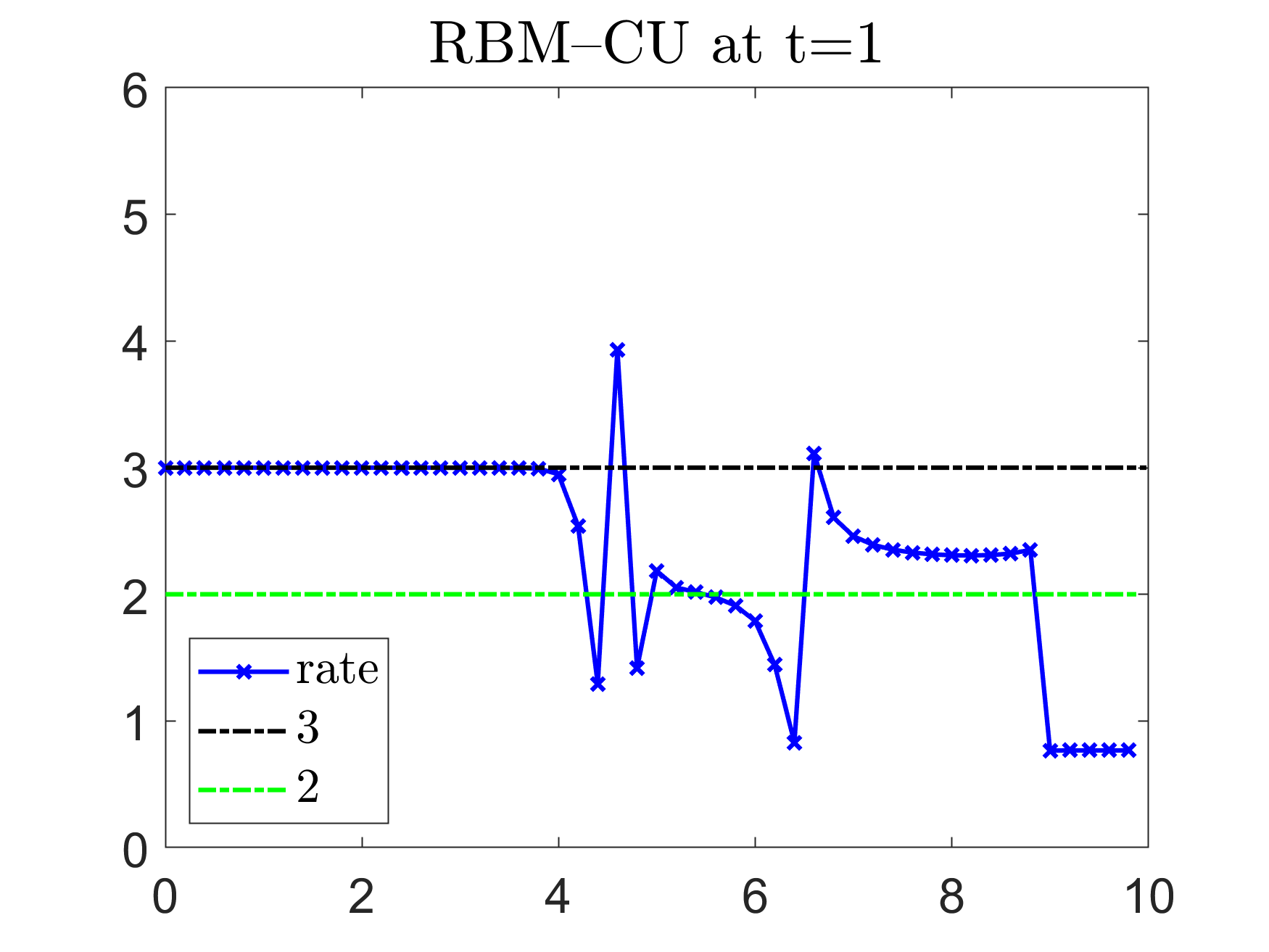}\hspace*{0.5cm}
            \includegraphics[trim=1.3cm 0.3cm 1.0cm 0.1cm, clip, width=5.3cm]{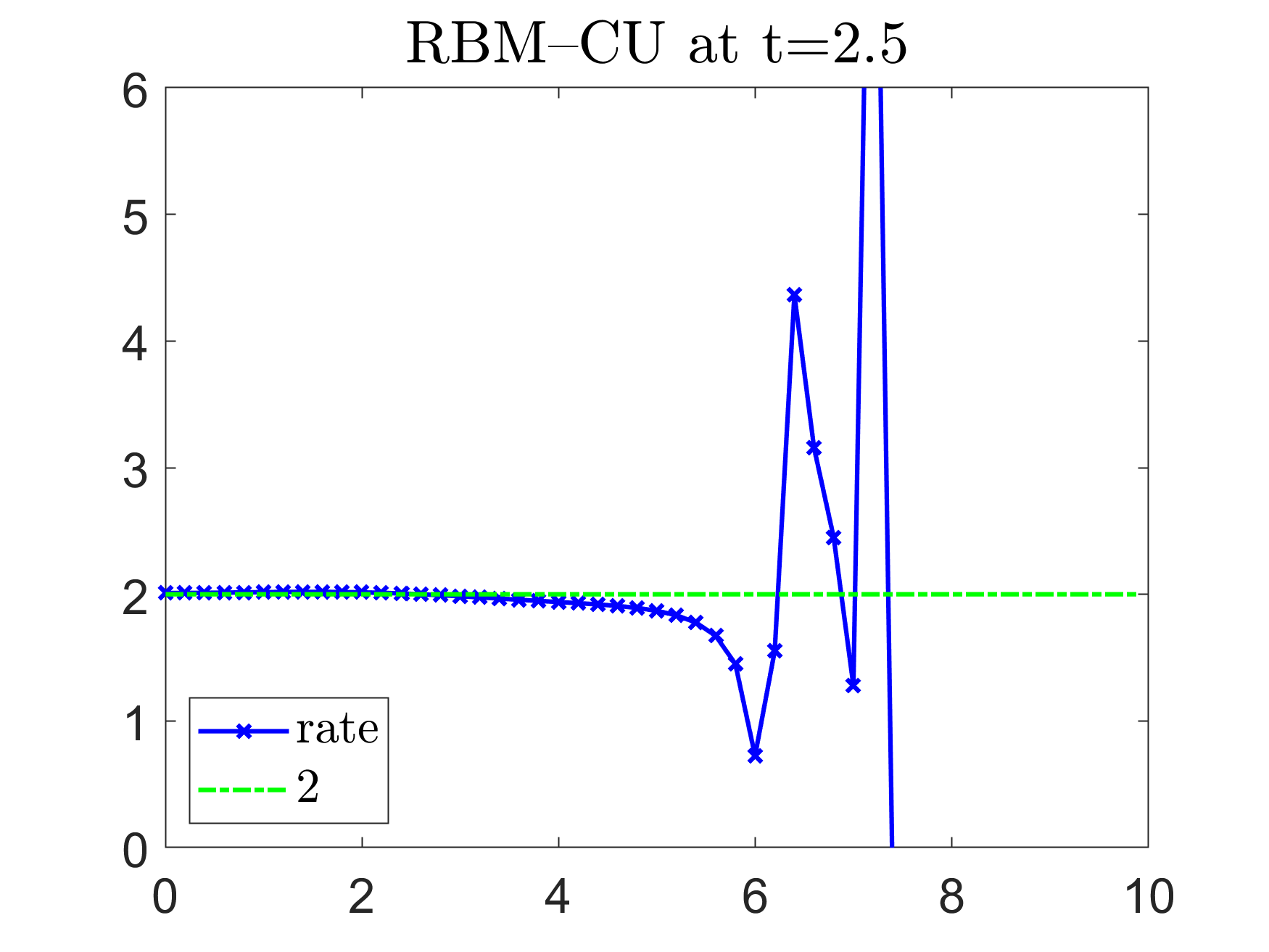}}
\vskip 8pt
\centerline{\includegraphics[trim=1.3cm 0.3cm 1.0cm 0.1cm, clip, width=5.3cm]{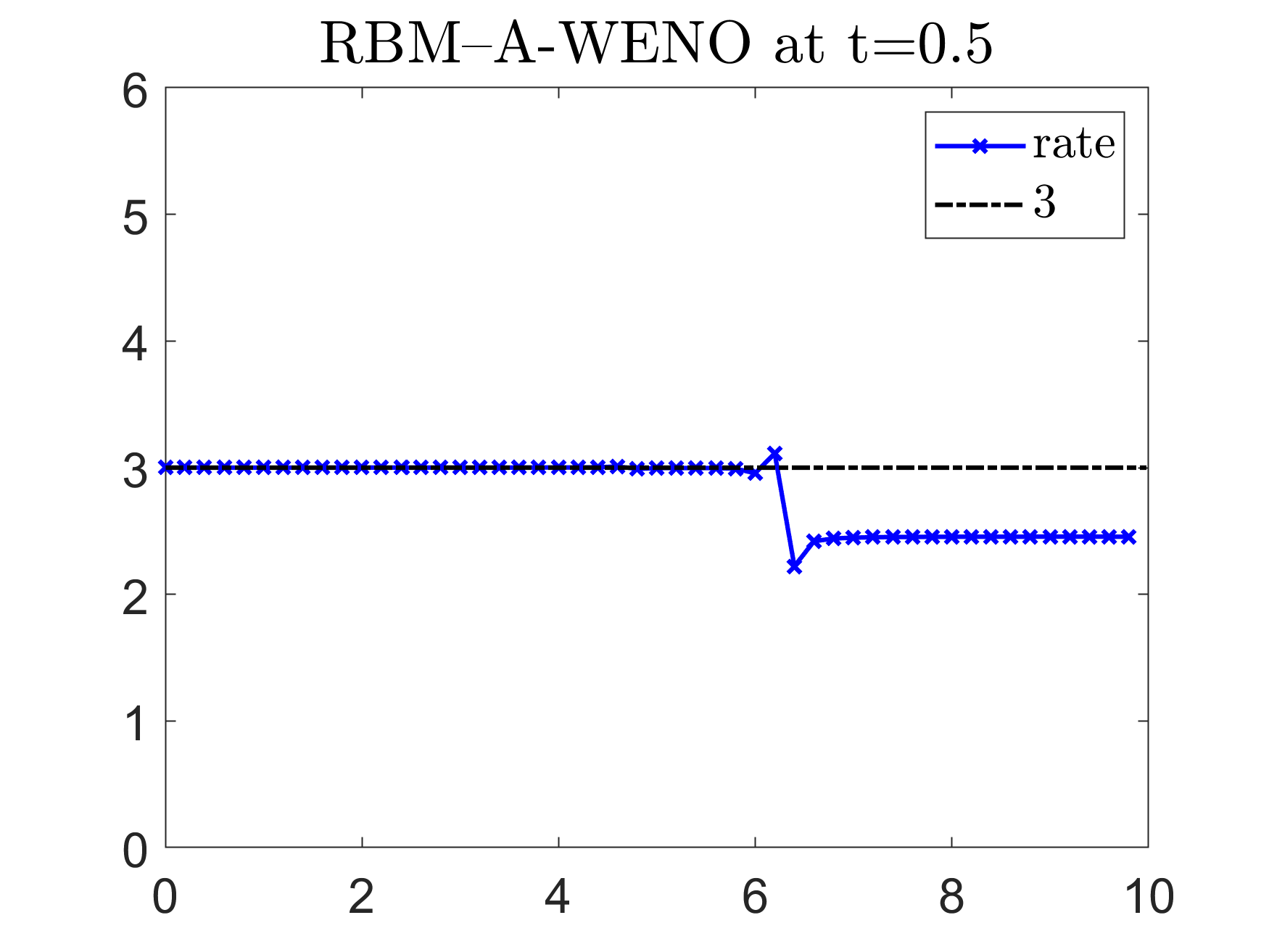}\hspace*{0.5cm}
            \includegraphics[trim=1.3cm 0.3cm 1.0cm 0.1cm, clip, width=5.3cm]{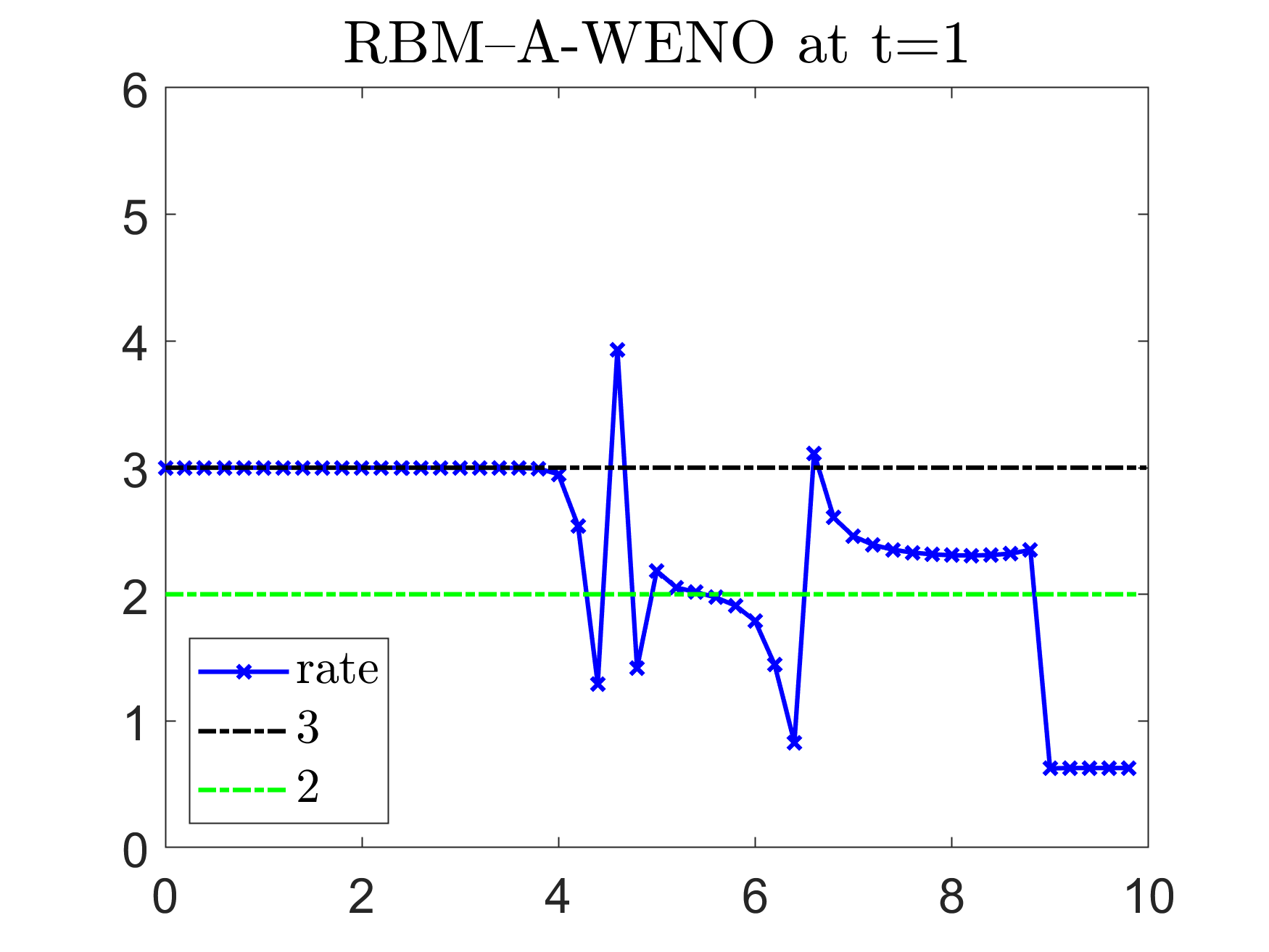}\hspace*{0.5cm}
            \includegraphics[trim=1.3cm 0.3cm 1.0cm 0.1cm, clip, width=5.3cm]{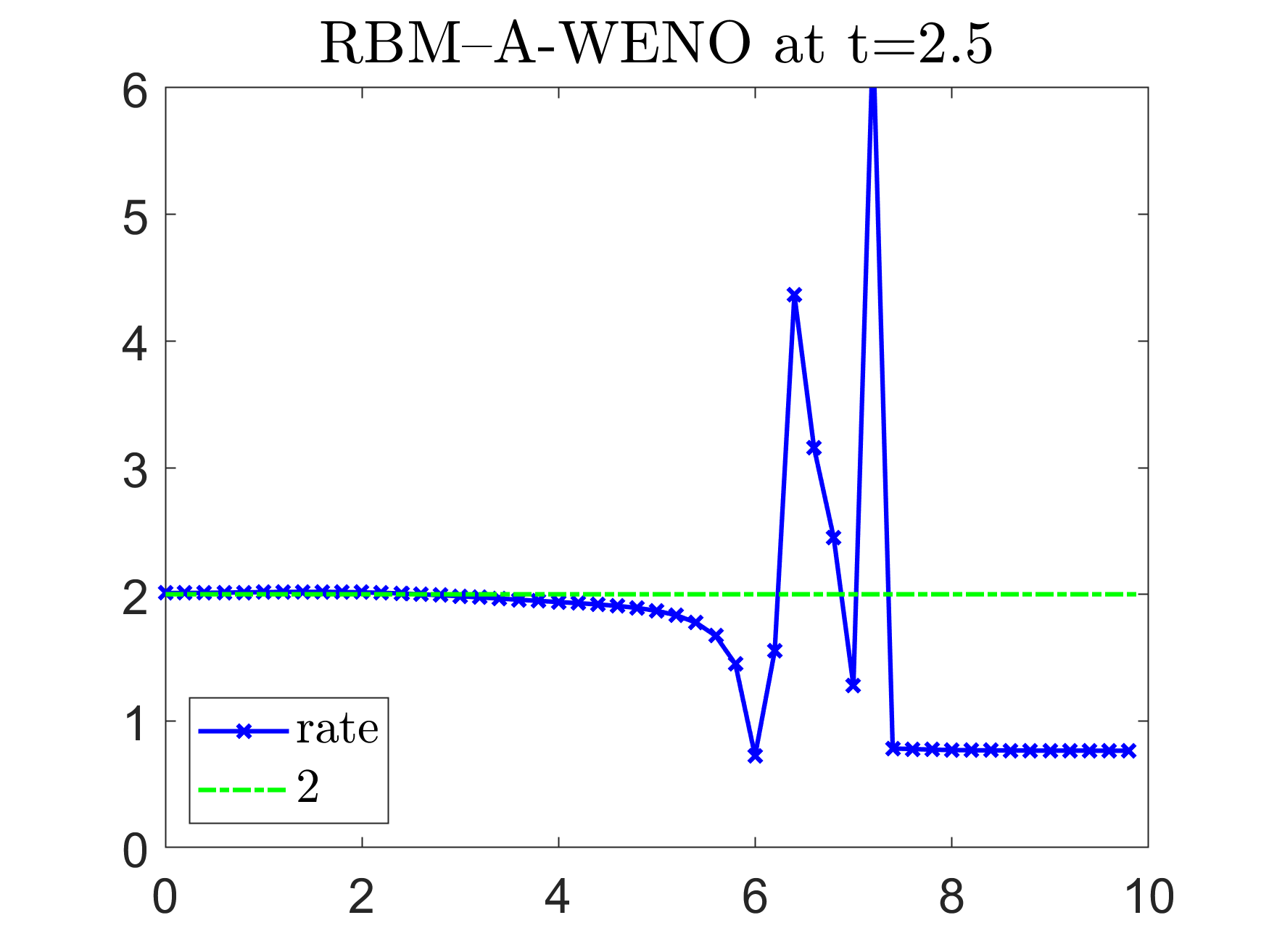}}
\caption{\sf Example 4: Experimental integral rates of convergence for the RBM--CU (top row) and RBM--A-WENO (bottom row) schemes at different times.\label{Fig6.3}}
\end{figure}

We also compute the integral rates of convergence using \eref{4.1} and present $r^{\rm INT}_{160k}$ for the RBM--CU and RBM--A-WENO schemes
with $k=0,\ldots,50$ in Figure \ref{Fig6.3}. As one can clearly see, the integral rates of convergence for both the RBM--CU and RBM--A-WENO
schemes substantially affected by combining two schemes of a different nature, namely, the RBM scheme and the CU/A-WENO scheme; compare
Figure \ref{Fig6.3} with the results reported in Figure \ref{Fig4.6}. This occurs even at time $t=0.5$, when the solution is still smooth,
but the integral convergence rates reduce to second order for the both of the two combined schemes. One can also notice that at time
$t=2.5$, the RBM--CU integral convergence rates deteriorate to the right of the shock front. This indicates that the integral convergence
rates \eref{4.1} may not offer a reliable tool for measuring the accuracy of the combined schemes. This is also confirmed by computing the
$W^{-1,1}$ convergence rates defined in \eref{4.4}. These rates are reported in Table \ref{Tab4}, where one can observe that the rates
reduce to second order at a small time $t=0.5$ and to the first order at a larger time $t=1$. However, at the final time $t=2.5$, by which
the shock has traveled throughout the entire computational domain, no $W^{-1,1}$ convergence is observed even though the pointwise
convergence rates at this time for the non-oscillatory combined schemes are as good as those for the RBM scheme.
\begin{table}[ht!]
\centering
\begin{tabular}{|c|c|c|c|c|c|c|c|}\hline
\multicolumn{2}{|c|}{$~$}&\multicolumn{2}{c|}{RBM--CU Scheme}&\multicolumn{2}{c|}{RBM--A-WENO Scheme}\\\hline
$t$&$N$&$||I^N-I^{2N}||_{L^1}$&$r^{\rm INT}$&$||I^N-I^{2N}||_{L^1}$&$r^{\rm INT}$\\\hline
\multirow{3}*{0.5} &1000&2.80e-4&1.69&8.95e-5&1.83\\
\cline{2-6}        &2000&8.78e-5&2.34&2.51e-5&2.44\\
\cline{2-6}        &4000&1.74e-5&--- &4.61e-6&--- \\\hline\hline
\multirow{3}*{1}   &1000&1.43e-3&0.79&7.55e-4&0.56\\
\cline{2-6}        &2000&8.27e-4&0.80&5.11e-4&0.68\\
\cline{2-6}        &4000&4.75e-4&--- &3.18e-4& ---\\\hline\hline
\multirow{3}*{2.5} &1000&2.31e-3&1.22&9.54e-4&-0.69\\
\cline{2-6}        &2000&9.94e-4&-0.321&1.53e-3&0.79\\
\cline{2-6}        &4000&1.24e-3&--- &8.89e-4&---  \\
\hline
\end{tabular}
\caption{\sf Example 4: $W^{-1,1}$ convergence rates at different times.\label{Tab4}}
\end{table}

\subsubsection*{Example 5---Test with Two Interacting Shocks}
\begin{figure}[ht!]
\centerline{\includegraphics[trim=1.3cm 0.3cm 1.0cm 0.8cm, clip, width=5.3cm]{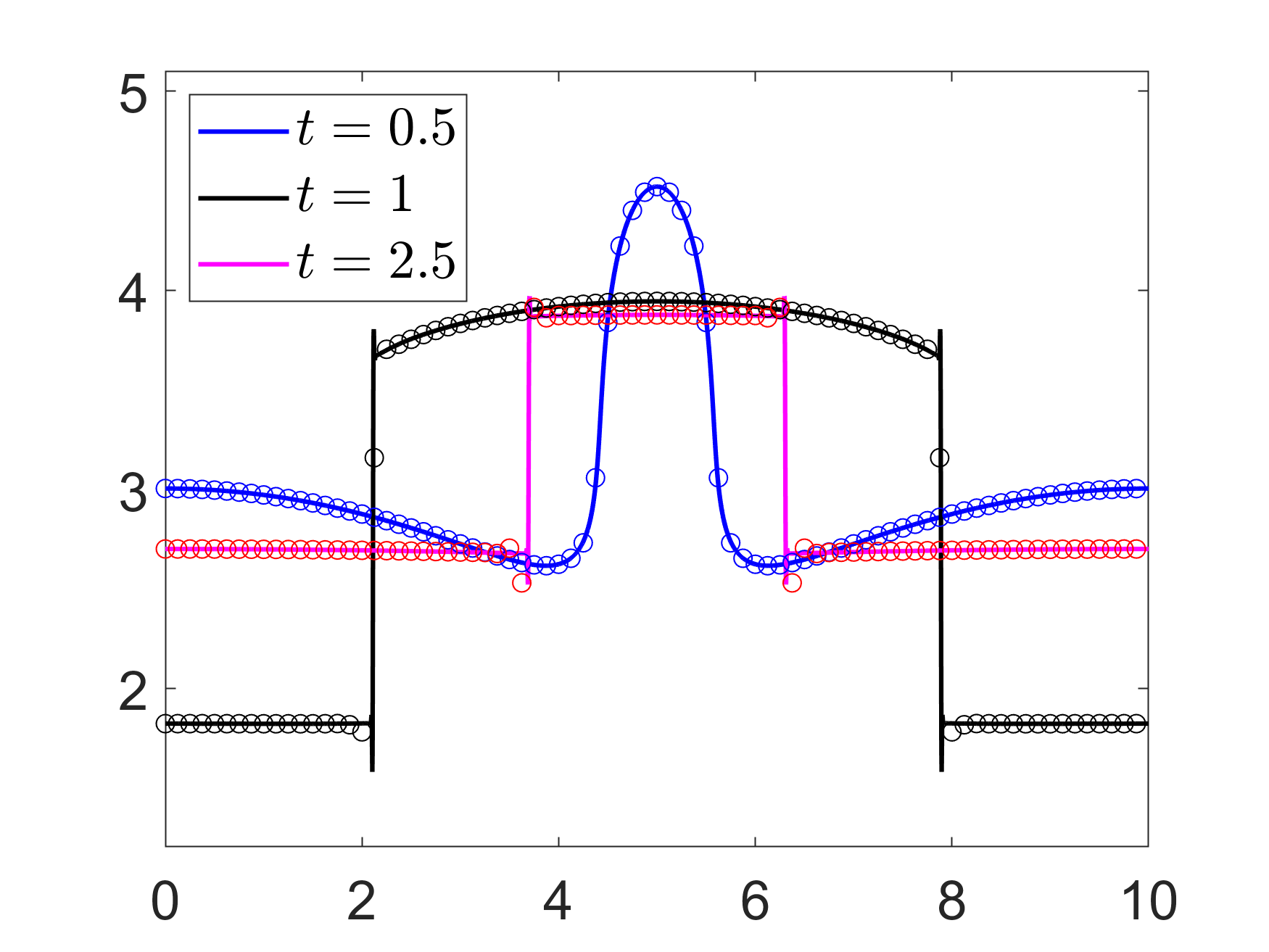}\hspace*{0.5cm}
            \includegraphics[trim=1.3cm 0.3cm 1.0cm 0.8cm, clip, width=5.3cm]{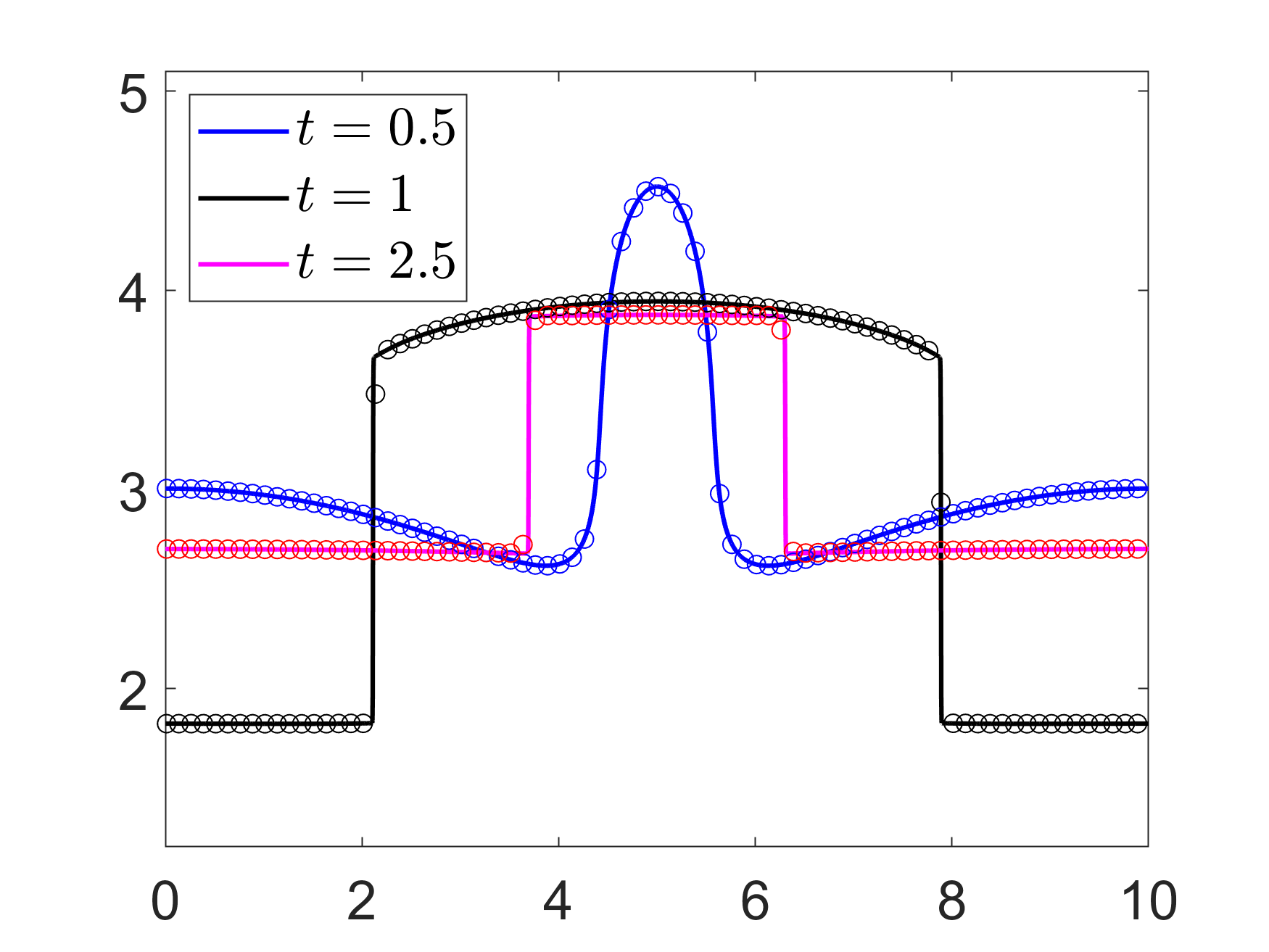}\hspace*{0.5cm}
            \includegraphics[trim=1.3cm 0.3cm 1.0cm 0.8cm, clip, width=5.3cm]{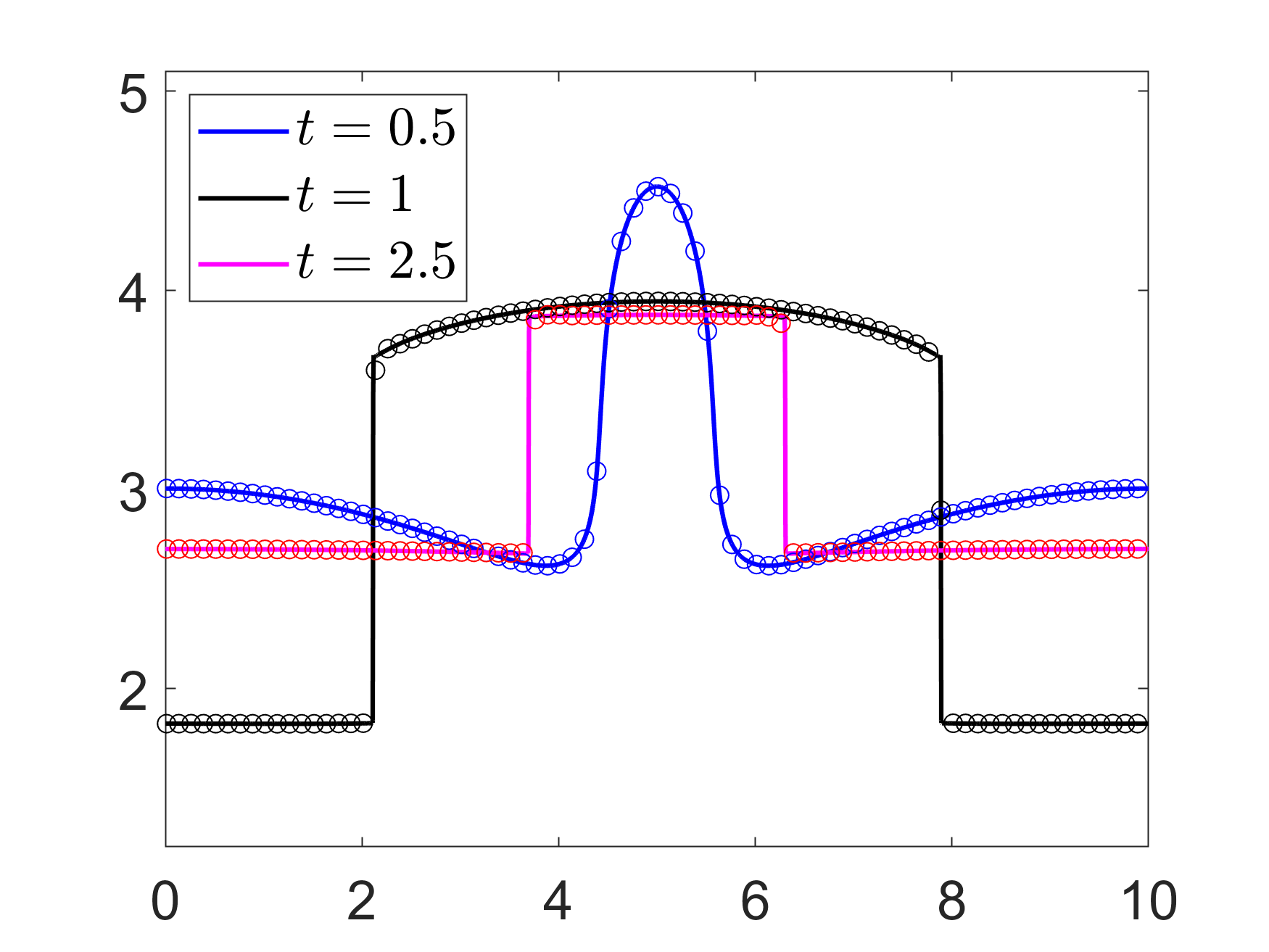}}
\caption{\sf Example 5: Water depth $h$ computed by the RBM (left), RBM--CU (middle), and RBM--A-WENO (right) schemes on the coarse
(circles) and fine (solid lines) grids.\label{fig72}}
\end{figure}

In this example, we use the same settings as in Example 2 and compute the solutions by the combined schemes at times $t=0.5$, 1, and 2.5 on
the computational domain $[0,10]$ using 1000, 2000, 4000, and 8000 uniform cells.  As in Example 4, we first tune the coefficient
$\texttt{C}$ on a coarse mesh with 400 uniform cells and then use it for the finer meshes. We take $\mu=0.1$ for both the RBM--CU and
RBM--A-WENO schemes. The results, computed on the grids with 400 and 4000 uniform cells are presented in Figures \ref{fig72} along with the
``pure'' RBM solutions. We plot the values $h_{5\jph}$, $j=0,\ldots,80$ for the numerical results computed on the coarser mesh. As one can
see from Figure \ref{fig72}, both combined schemes produce non-oscillatory solutions and the transitions between the RBM and CU/A-WENO parts
in the combined solutions are smooth.

\begin{figure}[ht!]
\centerline{\includegraphics[trim=1.3cm 0.3cm 1.0cm 0.1cm, clip, width=5.3cm]{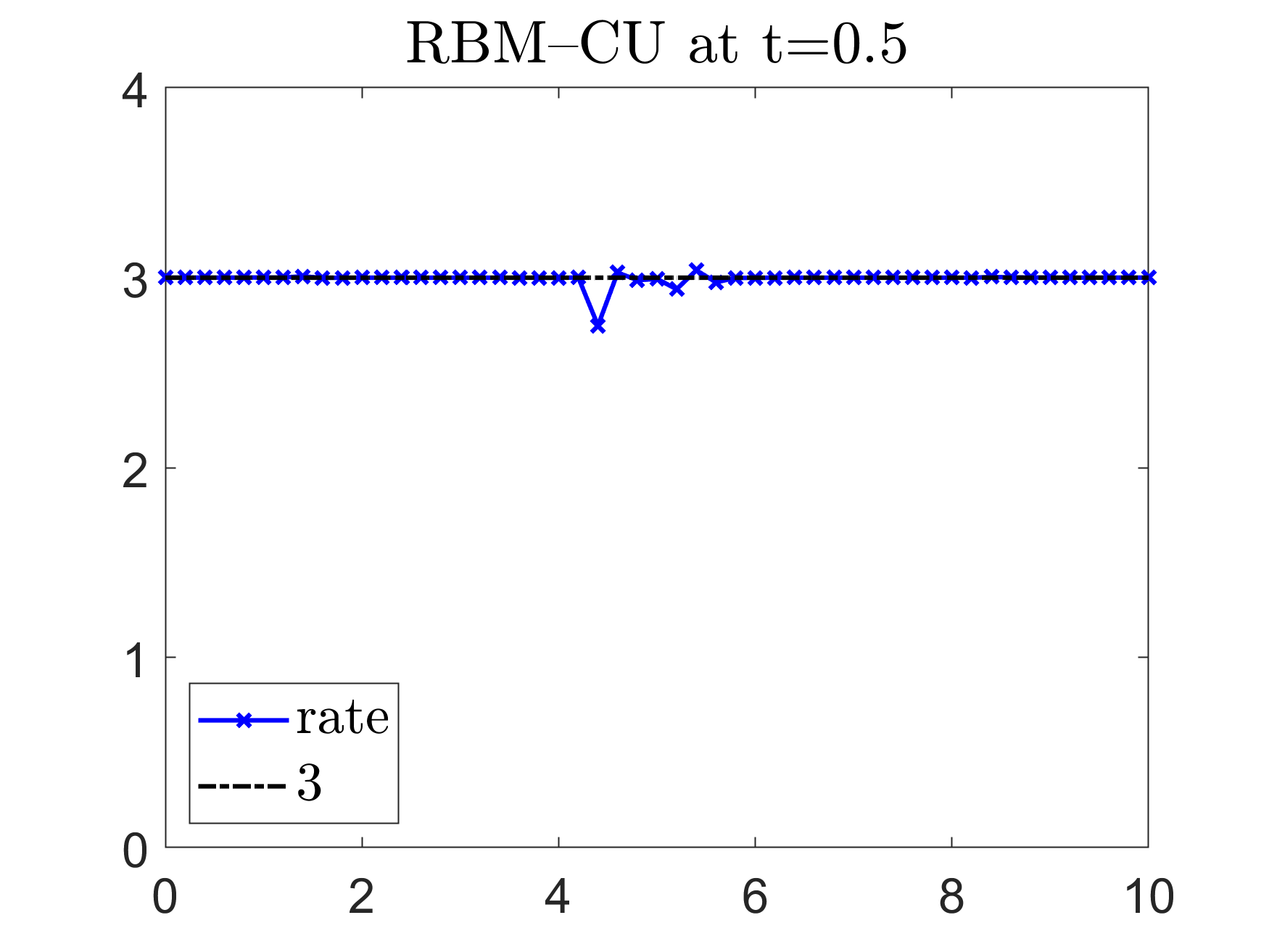}\hspace*{0.5cm}
            \includegraphics[trim=1.3cm 0.3cm 1.0cm 0.1cm, clip, width=5.3cm]{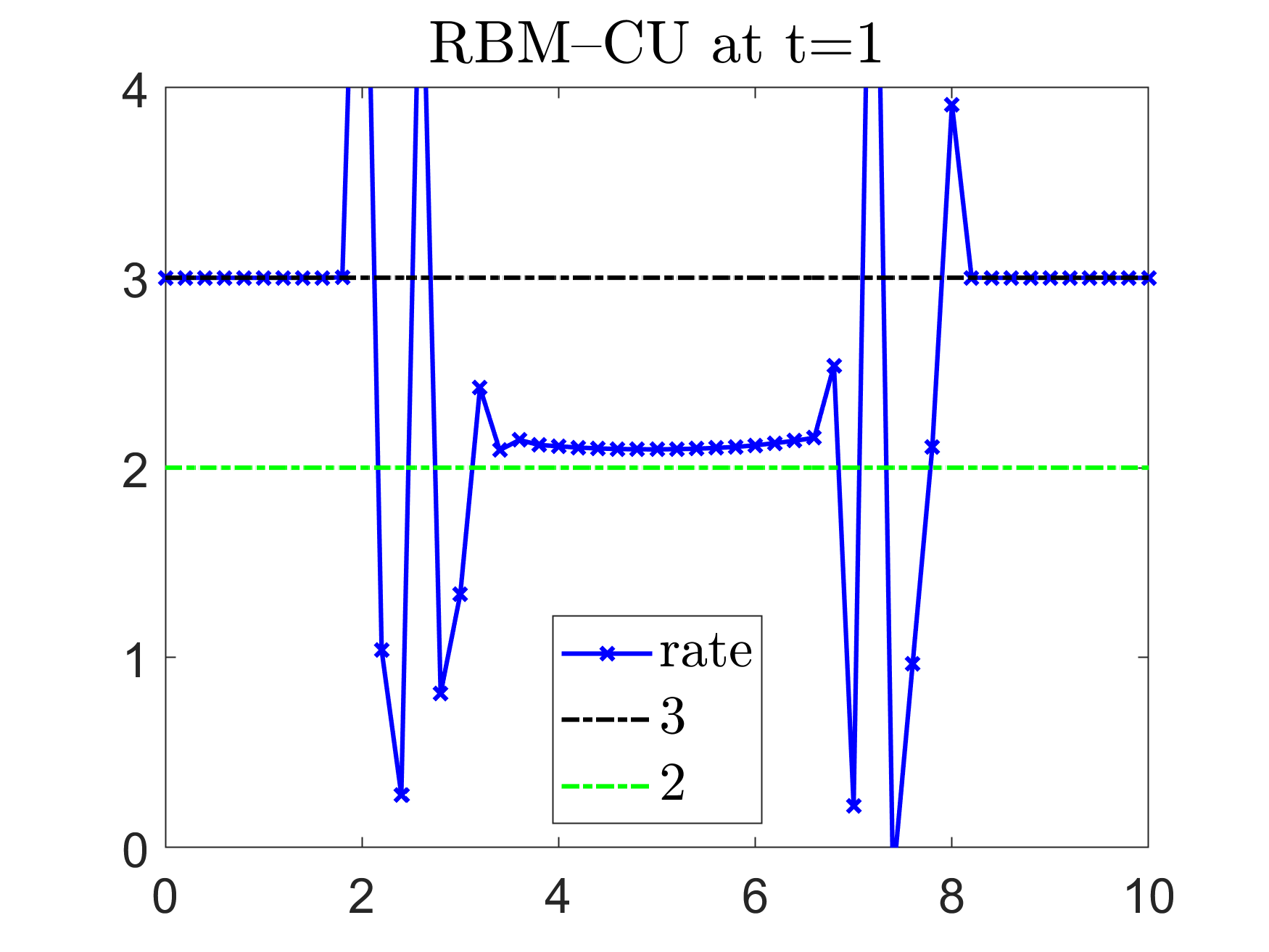}\hspace*{0.5cm}
            \includegraphics[trim=1.3cm 0.3cm 1.0cm 0.1cm, clip, width=5.3cm]{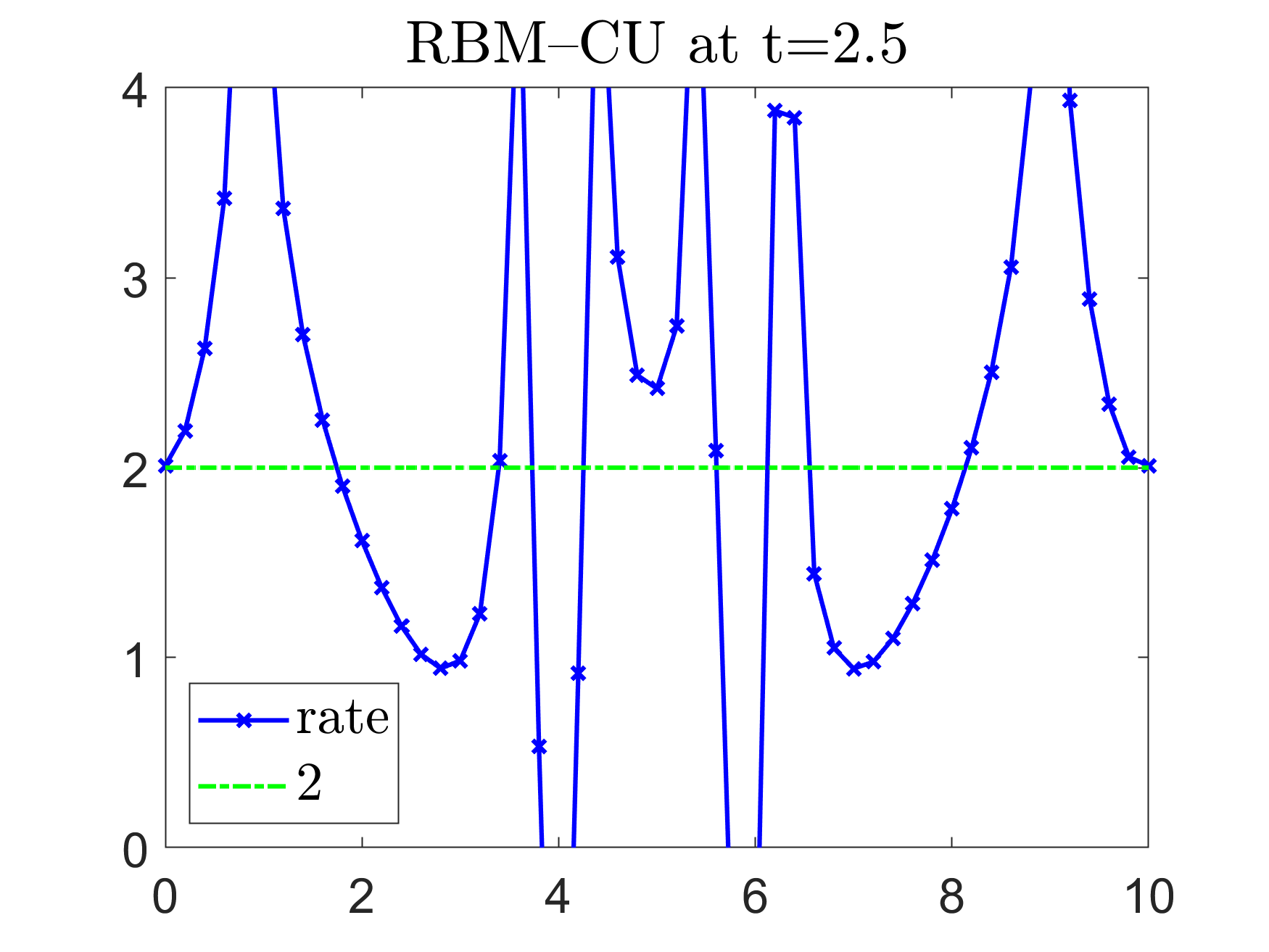}}
\vskip8pt
\centerline{\includegraphics[trim=1.3cm 0.3cm 1.0cm 0.1cm, clip, width=5.3cm]{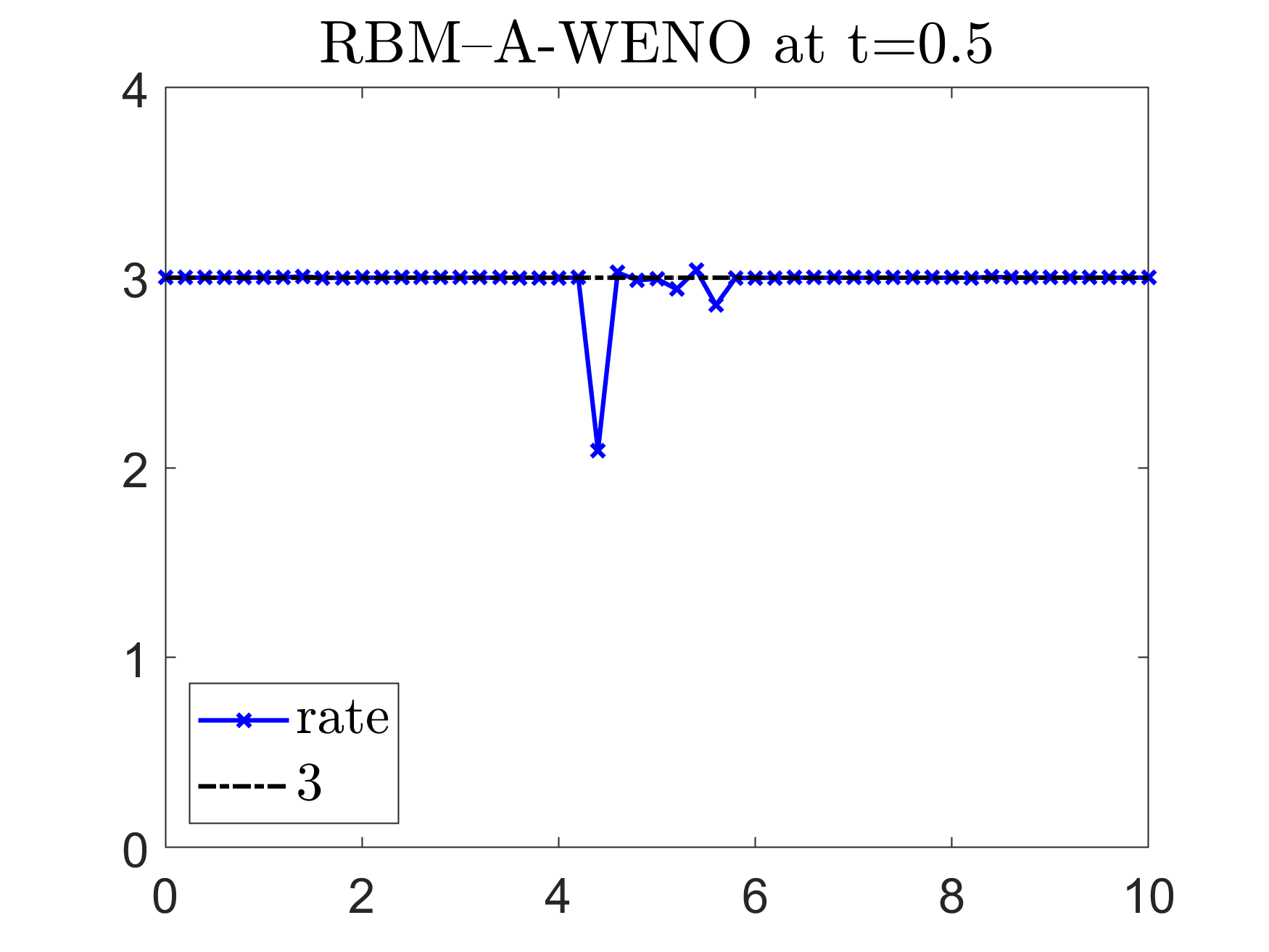}\hspace*{0.5cm}
            \includegraphics[trim=1.3cm 0.3cm 1.0cm 0.1cm, clip, width=5.3cm]{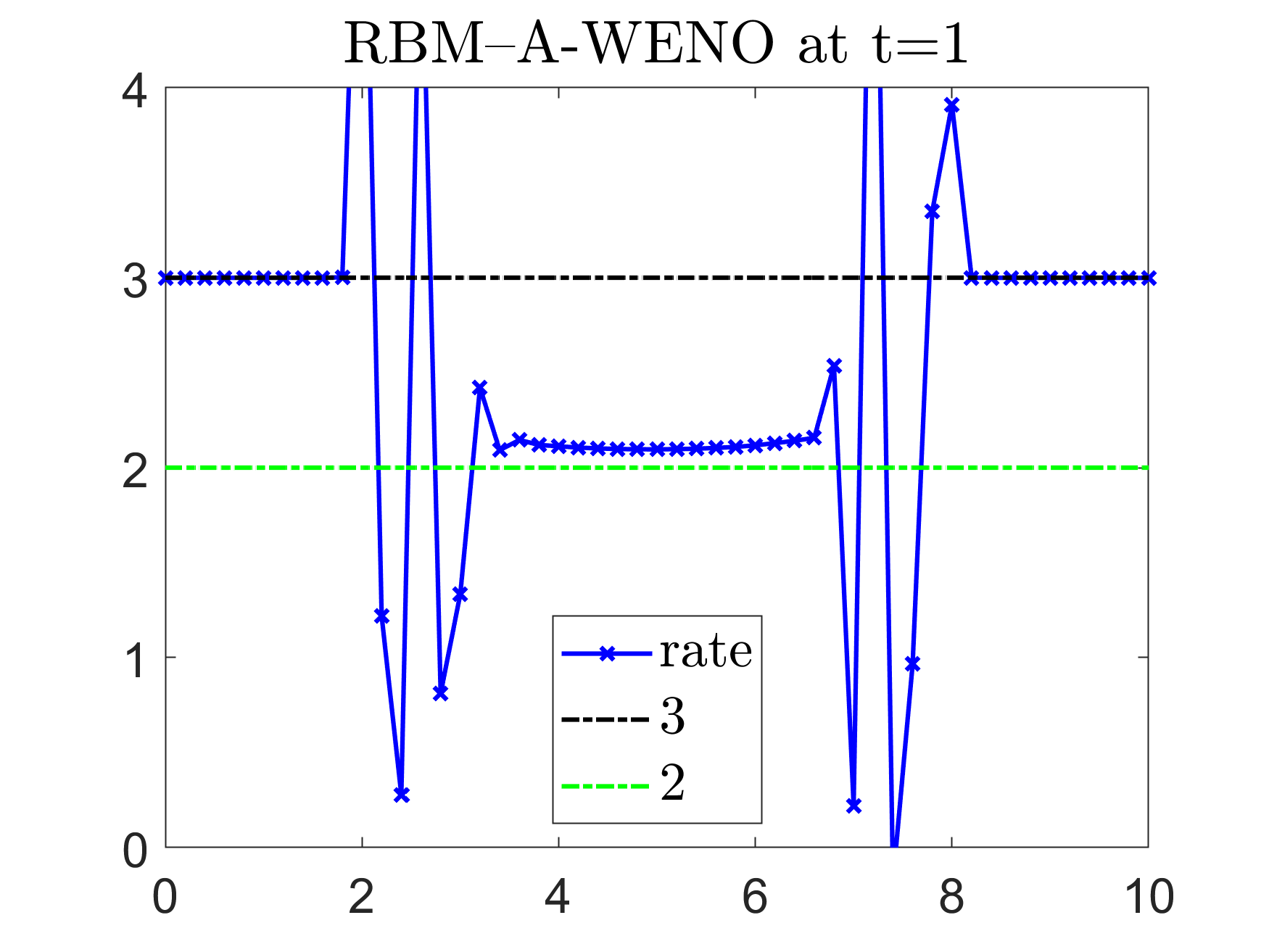}\hspace*{0.5cm}
            \includegraphics[trim=1.3cm 0.3cm 1.0cm 0.1cm, clip, width=5.3cm]{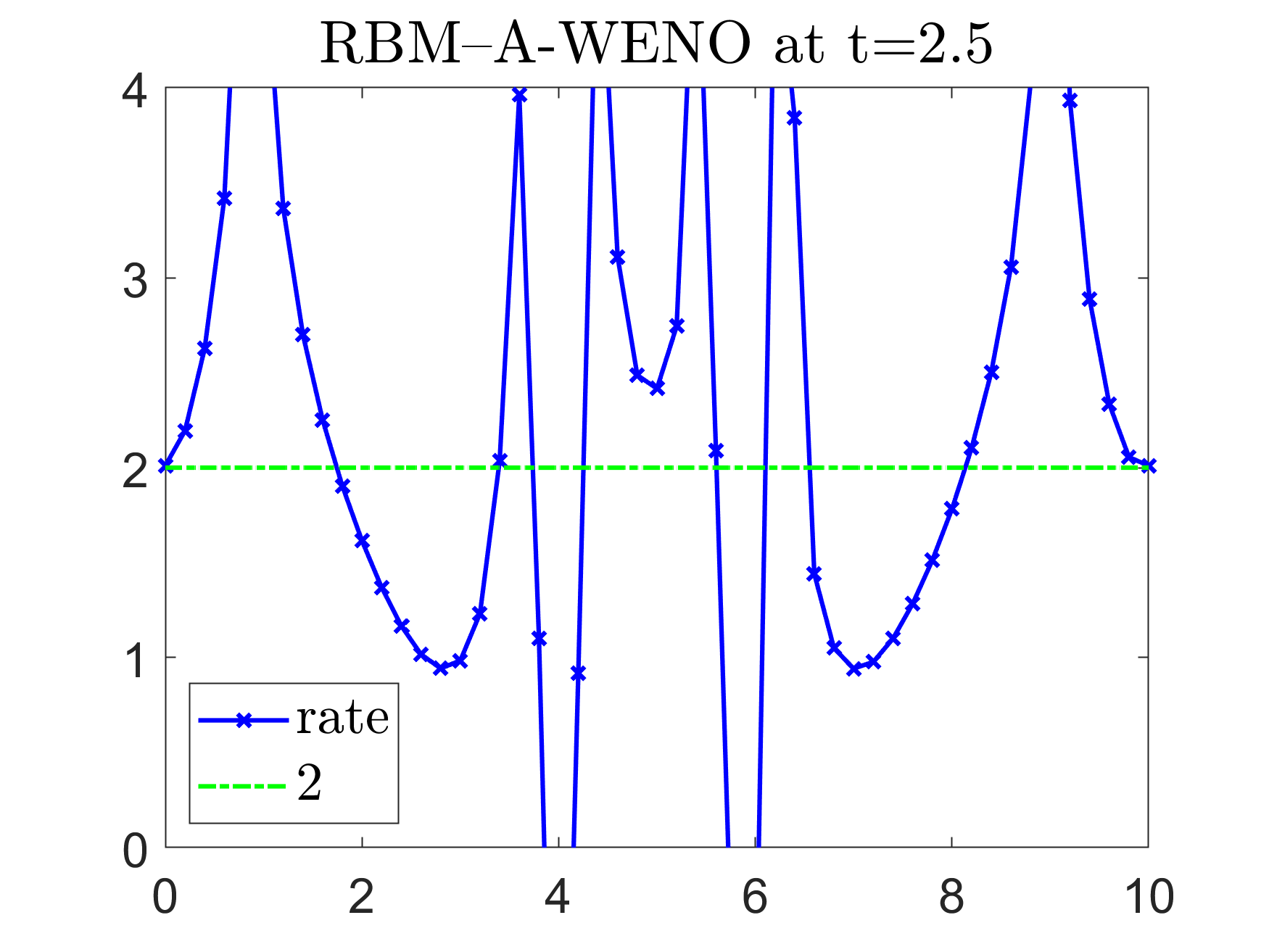}}
\caption{\sf Example 5: Average experimental rates of pointwise convergence for the RBM--CU (top row) and RBM--A-WENO (bottom row) schemes.
\label{Fig6.5}}
\end{figure}

We then check the pointwise convergence of the water depth $h$ by computing the average experimental pointwise convergence rates defined in
\eref{equ5.6} using three imbedded grids with $N=2000$ for the RBM--CU and RBM--A-WENO schemes. We plot the obtained values
$r^{\rm AVE}_{160k}$ in Figure \ref{Fig6.5}, where one can clearly see that the average convergence rates for both the RBM--CU and
RBM--A-WENO schemes are close to those of the RBM scheme; compare Figure \ref{Fig6.5} and the middle row panels of Figure \ref{Fig4.11}.

\begin{figure}[ht!]
\centerline{\includegraphics[trim=1.3cm 0.3cm 1.0cm 0.1cm, clip, width=5.3cm]{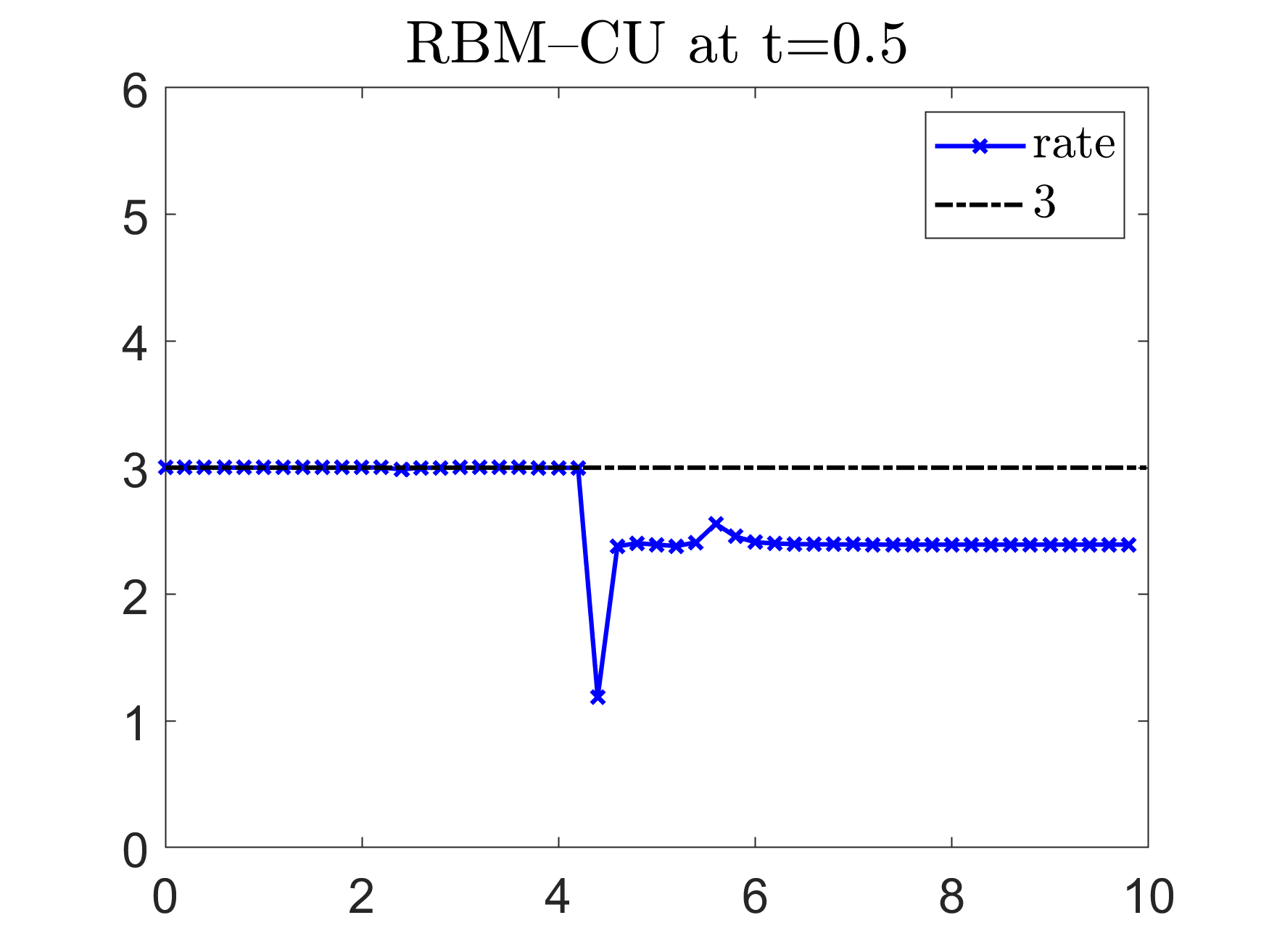}\hspace*{0.5cm}
            \includegraphics[trim=1.3cm 0.3cm 1.0cm 0.1cm, clip, width=5.3cm]{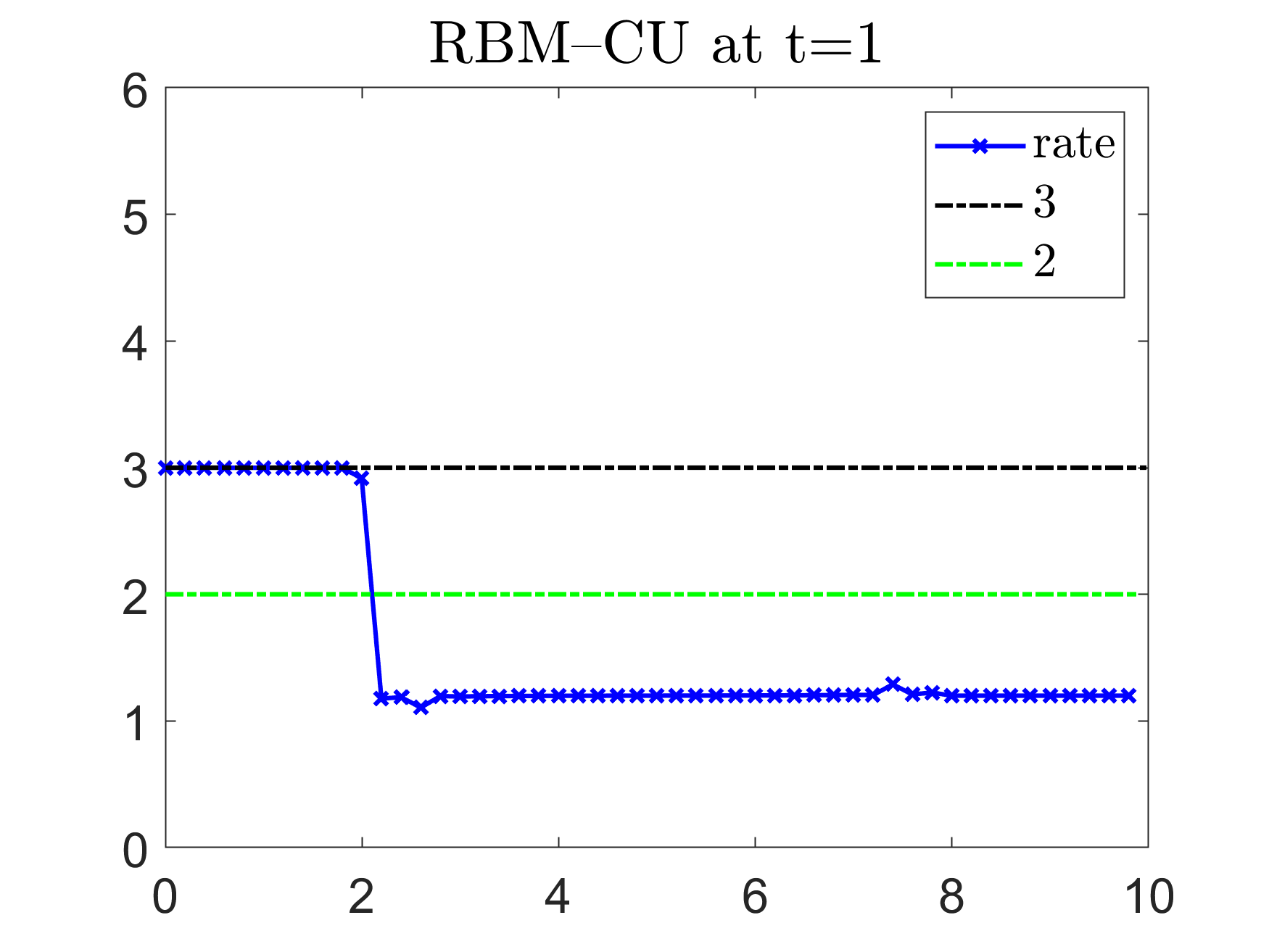}\hspace*{0.5cm}
            \includegraphics[trim=1.3cm 0.3cm 1.0cm 0.1cm, clip, width=5.3cm]{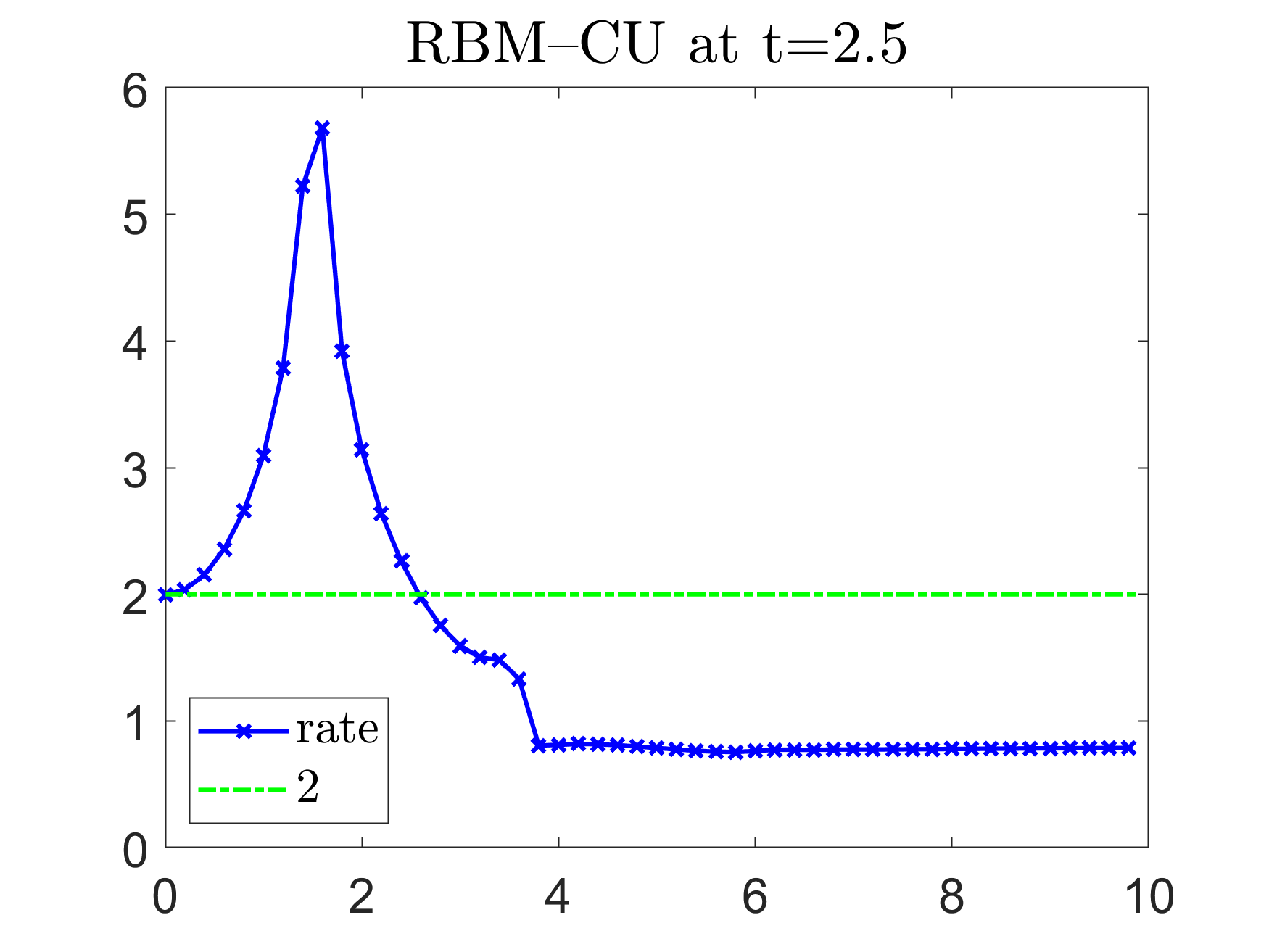}}
\vskip 8pt
\centerline{\includegraphics[trim=1.3cm 0.3cm 1.0cm 0.1cm, clip, width=5.3cm]{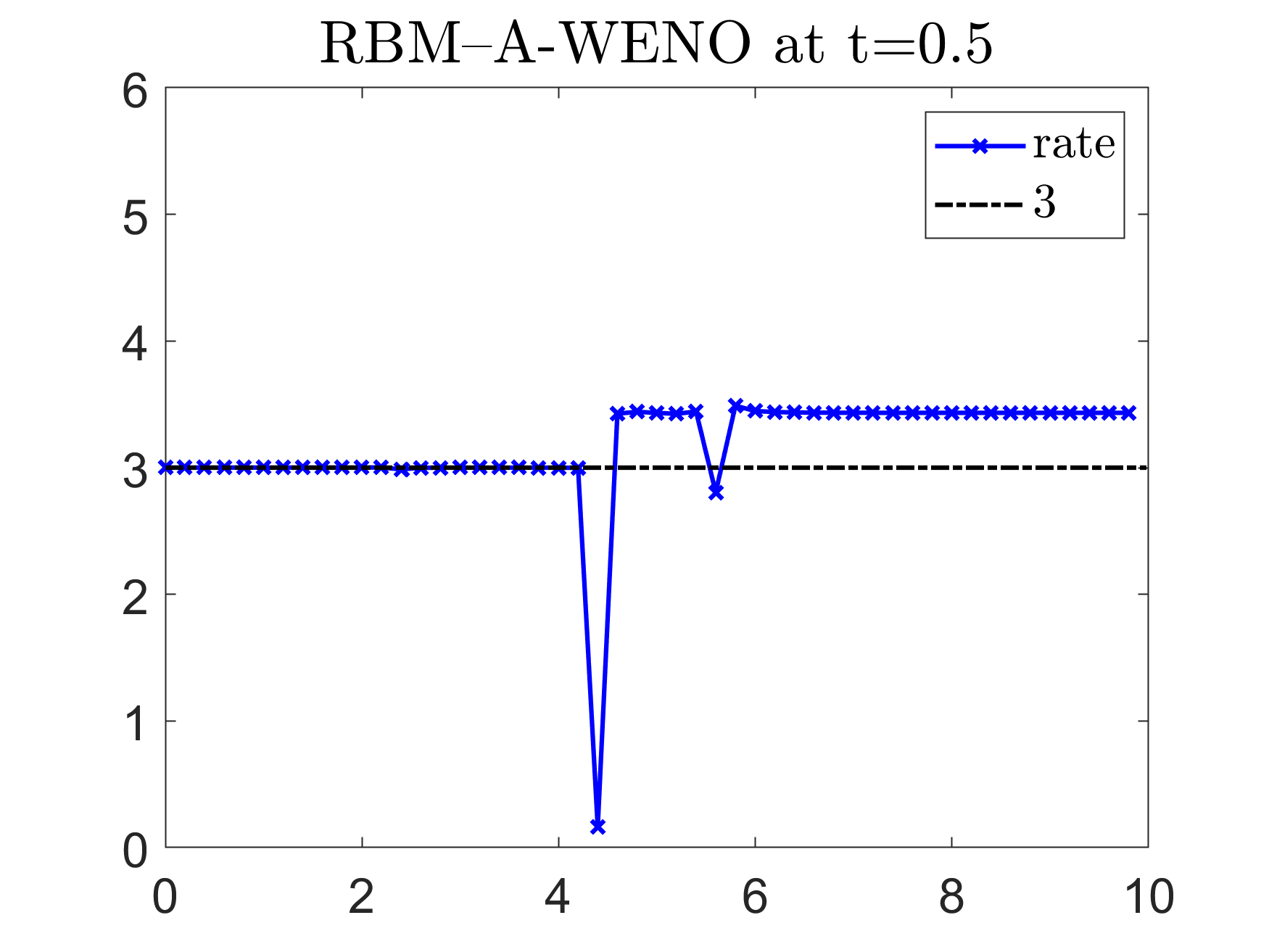}\hspace*{0.5cm}
            \includegraphics[trim=1.3cm 0.3cm 1.0cm 0.1cm, clip, width=5.3cm]{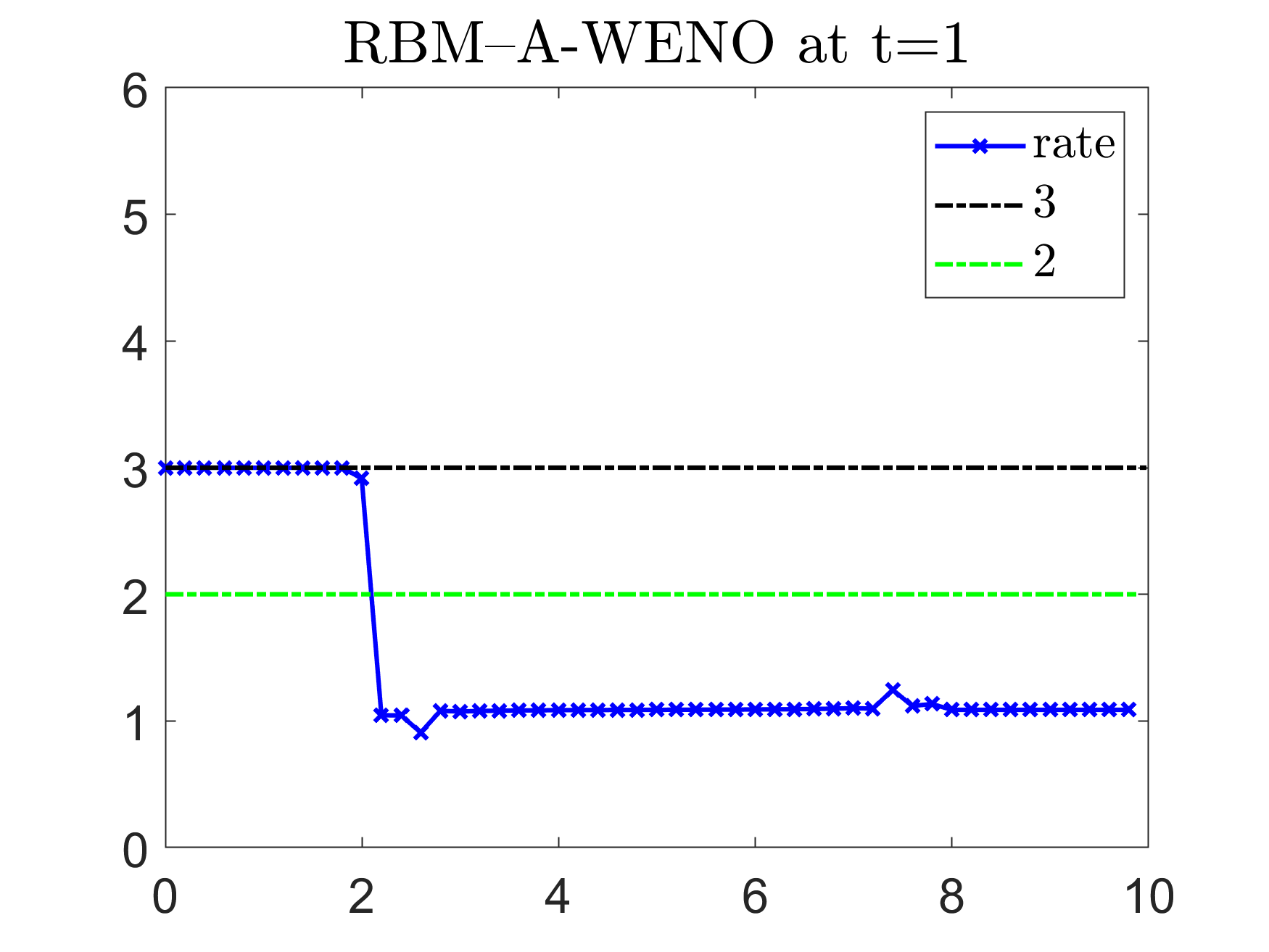}\hspace*{0.5cm}
            \includegraphics[trim=1.3cm 0.3cm 1.0cm 0.1cm, clip, width=5.3cm]{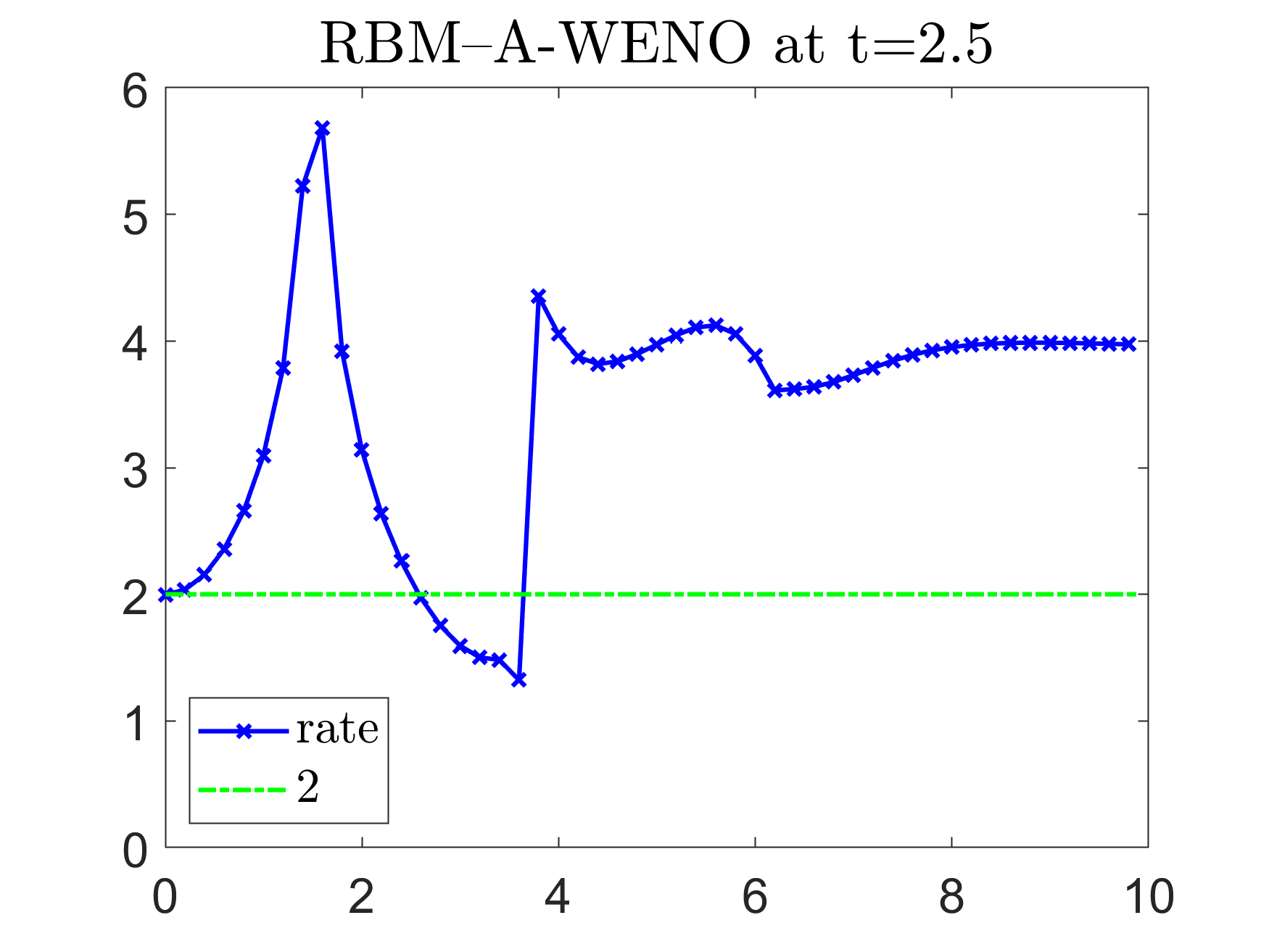}}
\caption{\sf Example 5: Experimental integral rates of convergence for the RBM--CU (top row) and RBM--A-WENO (bottom row) schemes at different times.\label{Fig6.3a}}
\end{figure}

We also compute the integral rates of convergence using \eref{4.1} and the $W^{-1,1}$ convergence rates defined in \eref{4.4}. The obtained
results are shown in Figure \ref{Fig6.3a} and Table \ref{Tab4a}. Unlike Example 4, one can now observe the integral convergence for both the
RBM--CU and RBM--A-WENO schemes at all of the times, but the convergence rates fluctuate and we believe that this is caused by combining
schemes of a different type.
\begin{table}[ht!]
\centering
\begin{tabular}{|c|c|c|c|c|c|c|c|}\hline
\multicolumn{2}{|c|}{$~$}&\multicolumn{2}{c|}{RBM--CU Scheme}&\multicolumn{2}{c|}{RBM--A-WENO Scheme}\\\hline
$t$&$N$&$||I^N-I^{2N}||_{L^1}$&$r^{\rm INT}$&$||I^N-I^{2N}||_{L^1}$&$r^{\rm INT}$\\\hline
\multirow{3}*{0.5} &1000&4.04e-4&3.12&2.55e-4&1.97\\
\cline{2-6}        &2000&4.64e-5&2.41&6.51e-5&3.40\\
\cline{2-6}        &4000&8.71e-6&--- &6.18e-6&--- \\\hline\hline
\multirow{3}*{1}   &1000&9.23e-3&0.48&4.39e-3&0.26\\
\cline{2-6}        &2000&6.60e-3&1.20&3.66e-3&1.09\\
\cline{2-6}        &4000&2.87e-3&--- &1.72e-3& ---\\\hline\hline
\multirow{3}*{2.5} &1000&7.92e-3&1.81&1.27e-2&3.61\\
\cline{2-6}        &2000&2.26e-3&0.80&1.04e-3&3.78\\
\cline{2-6}        &4000&1.30e-3&--- &7.58e-5&---  \\
\hline
\end{tabular}
\caption{\sf Example 5: $W^{-1,1}$ convergence rates at different times.\label{Tab4a}}
\end{table}

\subsubsection*{Example 6---Test with an Isolated Shock}
\begin{figure}[ht!]
\centerline{\includegraphics[trim=0.8cm 0.3cm 1.0cm 0.8cm, clip, width=5.3cm]{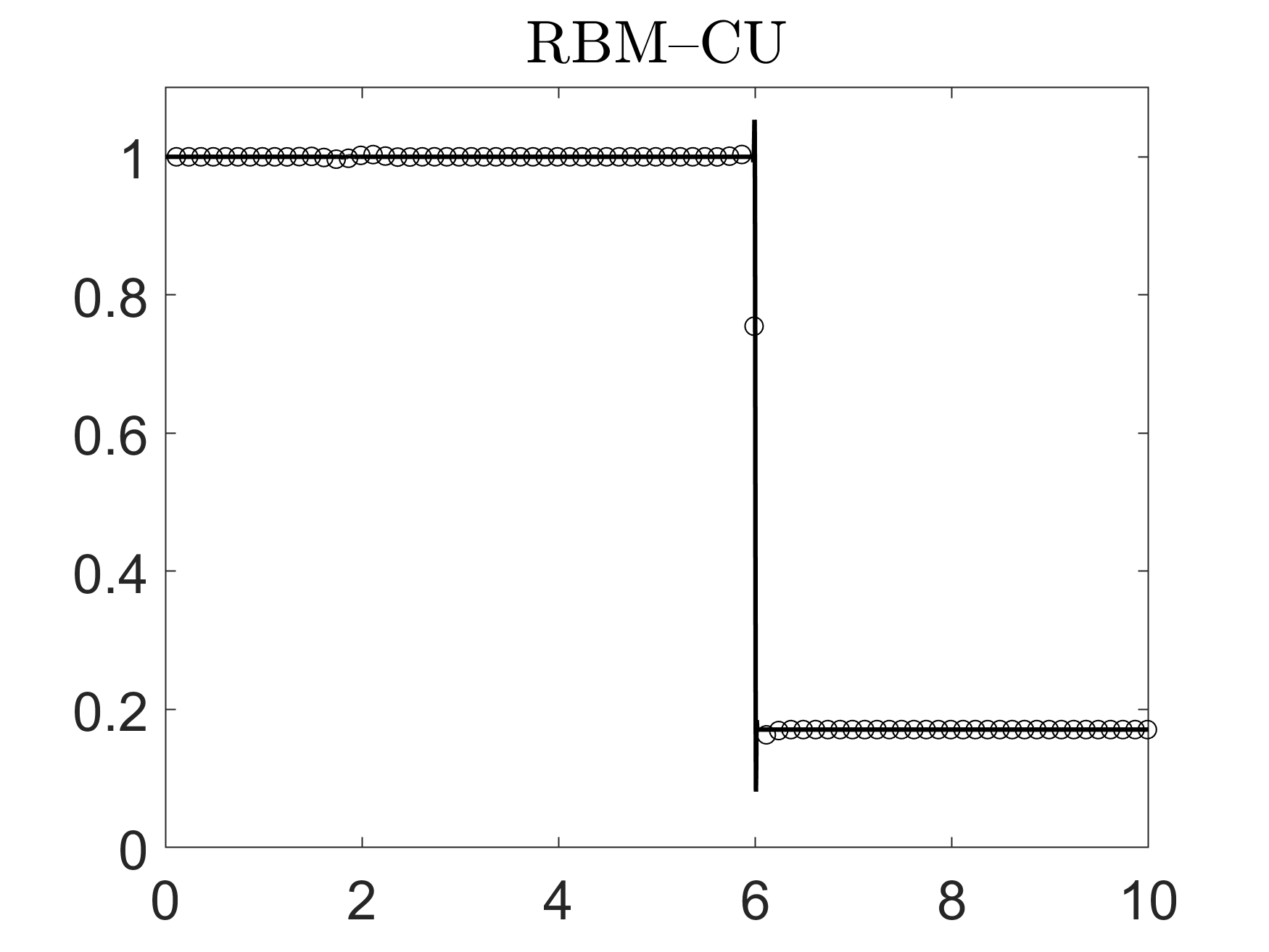}\hspace*{0.5cm}
            \includegraphics[trim=0.8cm 0.3cm 1.0cm 0.8cm, clip, width=5.3cm]{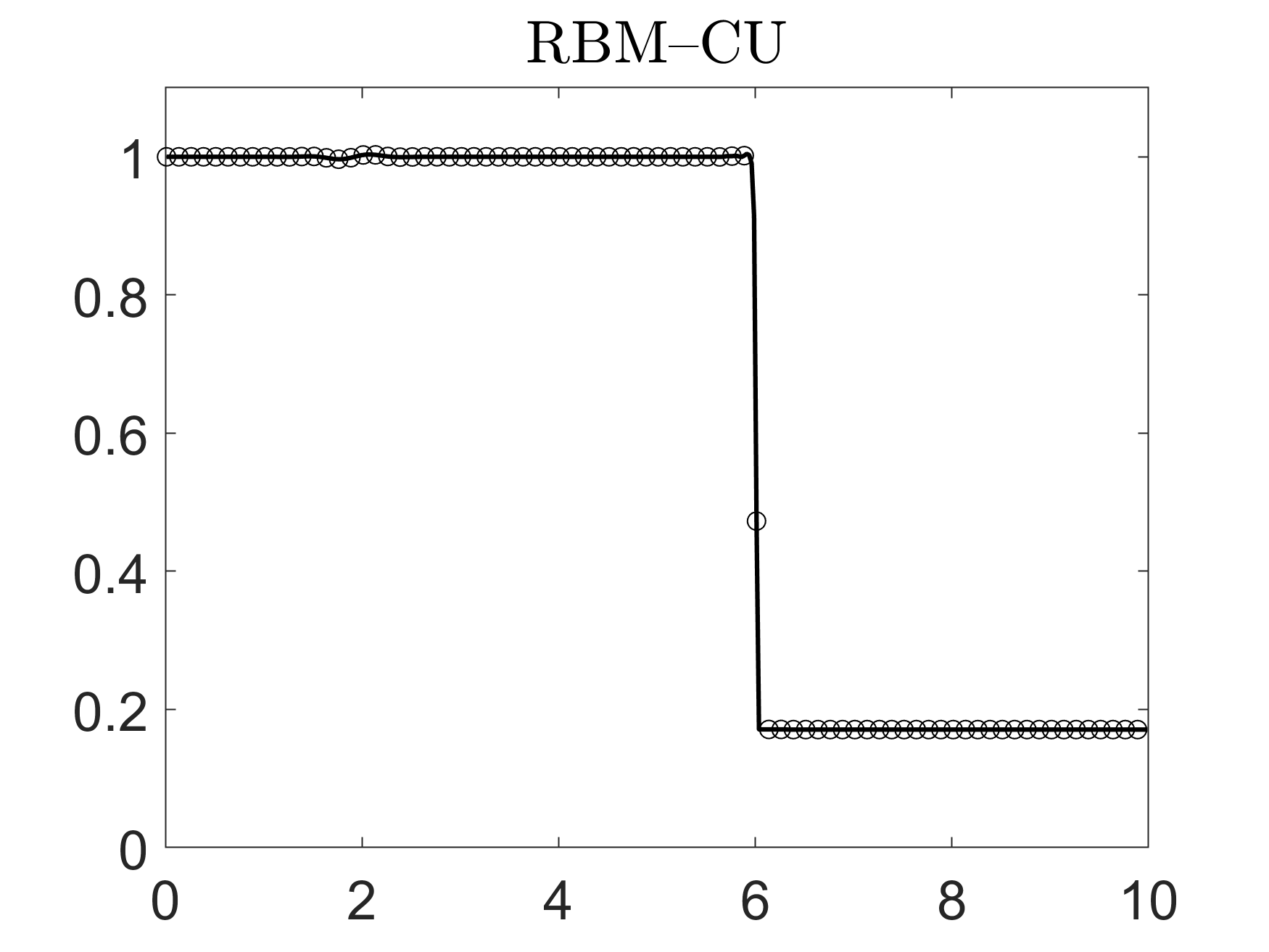}\hspace*{0.5cm}
            \includegraphics[trim=0.8cm 0.3cm 1.0cm 0.8cm, clip, width=5.3cm]{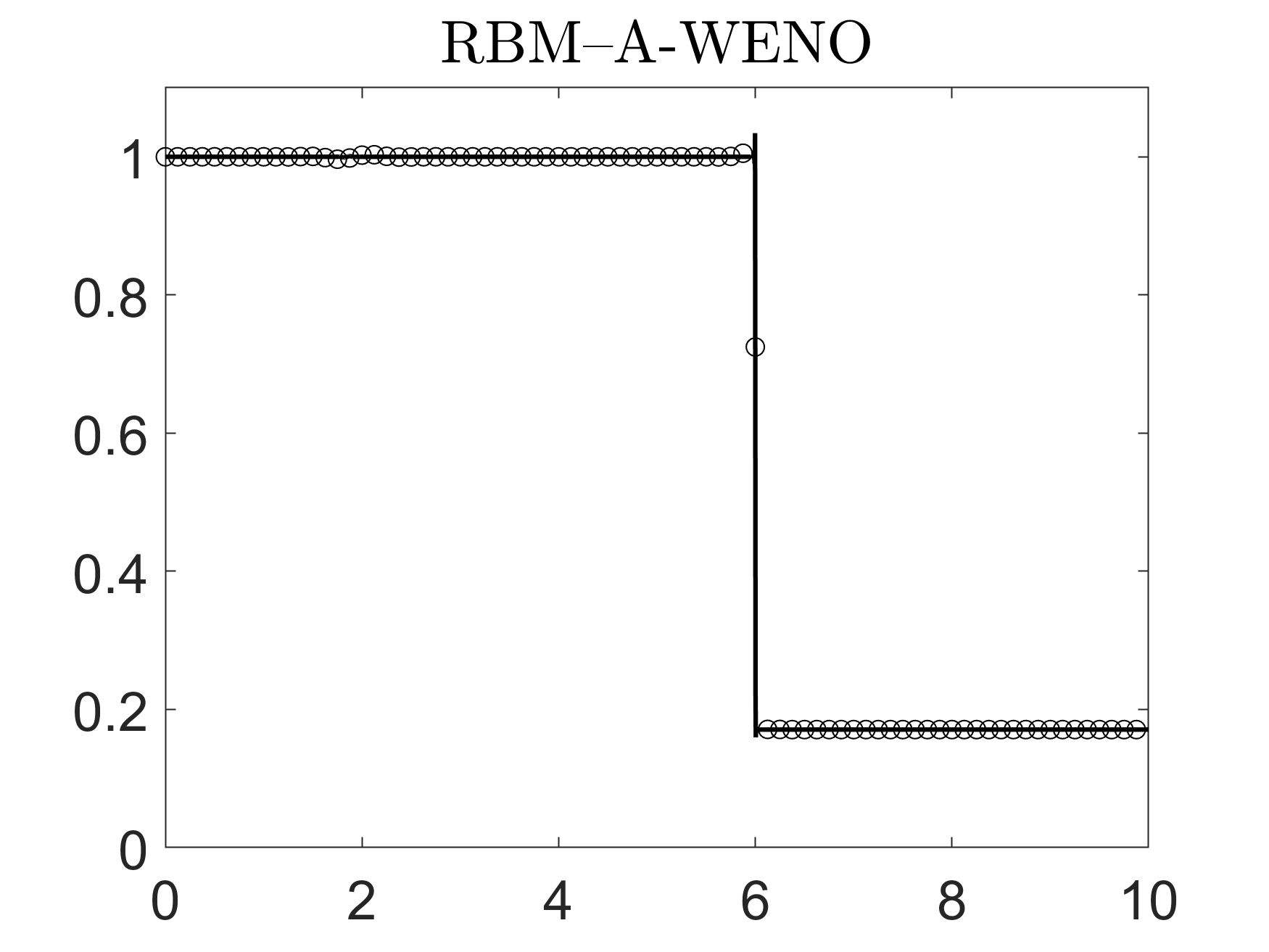}}
\caption{\sf Example 6: Water depth $h$ computed by the RBM (left), combined RBM--CU (middle), and combined RBM--A-WENO (right) schemes on
the coarse (circles) and fine (solid lines) grids.\label{fig67}}
\end{figure}

In the last example, we use the same settings as in Example 3 and compute the solutions by the studied RBM--CU and RBM--A-WENO schemes at
time $t=1$ on the computational domain $[0,10]$ using 1000, 2000, 4000, and 8000 uniform cells. We first tune the coefficient $\texttt{C}$
on a coarse mesh with 400 uniform cells and then use it for the finer mesh computations. In this example, we take $\mu=0.2$ for both the
RBM--CU and RBM--A-WENO schemes. The results, computed on the grids with 400 and 4000 uniform cells are presented in Figure \ref{fig67}
along with the ``pure'' RBM solutions. Here, we plot the values $h_{5\jph}$, $j=0,\ldots,80$ for the numerical results computed on the
coarser mesh. As one can see, both combined schemes produce non-oscillatory solutions and the transitions between the RBM and CU/A-WENO
parts in the combined solutions are smooth.

\begin{figure}[ht!]
\centerline{\includegraphics[trim=1.3cm 0.3cm 1.0cm 0.1cm, clip, width=5.3cm]{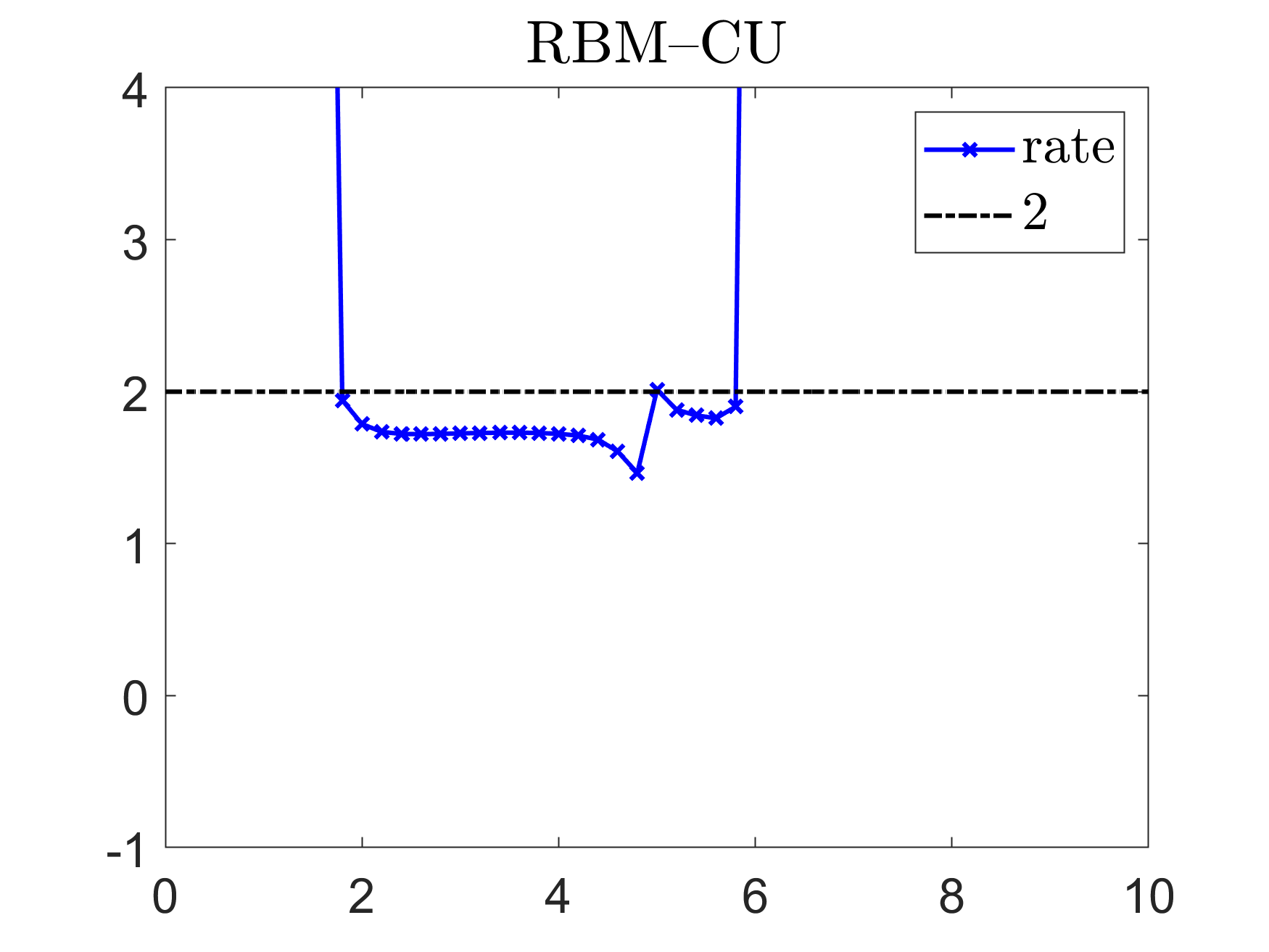}\hspace*{2.0cm}
            \includegraphics[trim=1.3cm 0.3cm 1.0cm 0.1cm, clip, width=5.3cm]{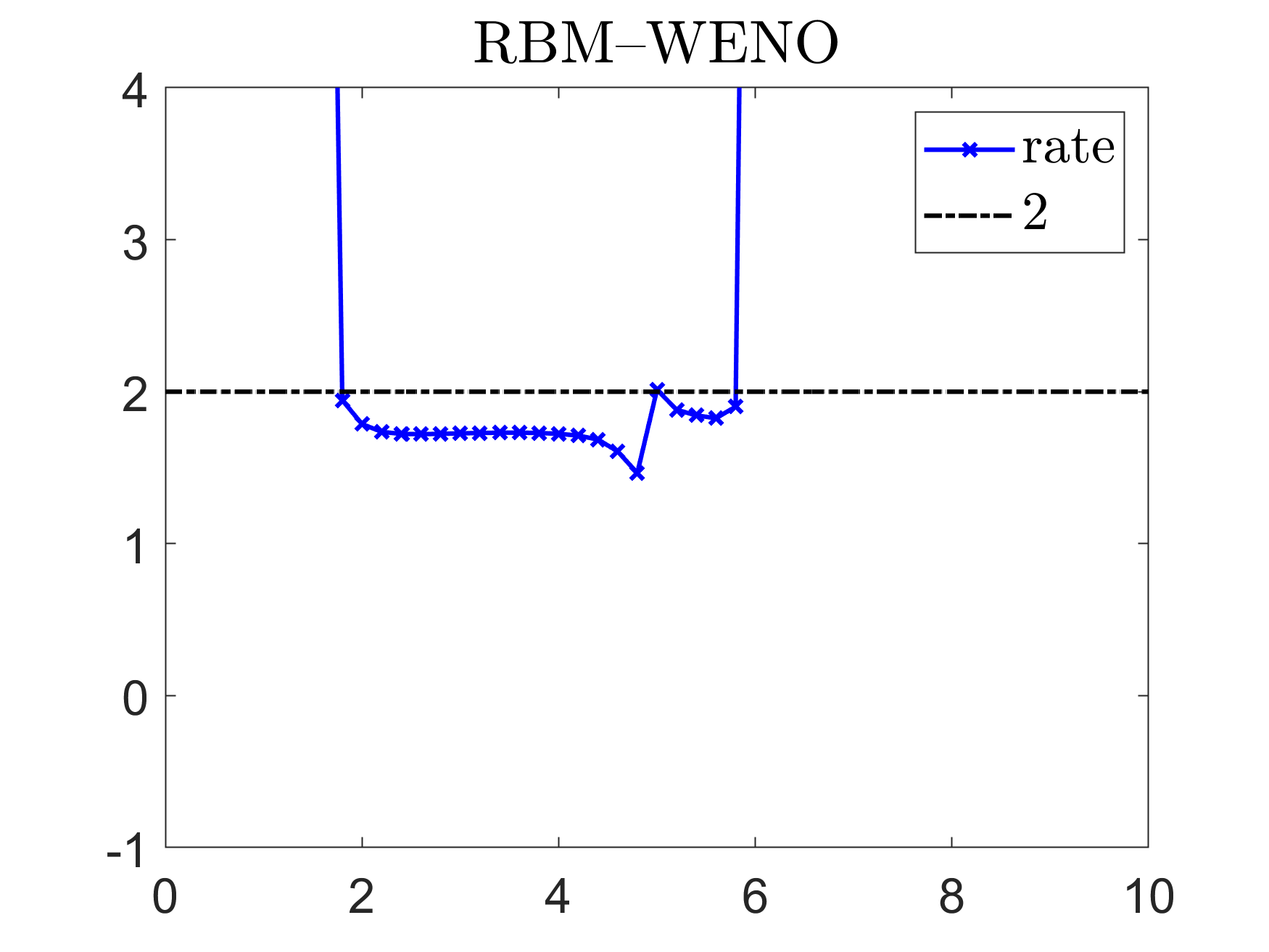}}
\caption{\sf Example 6: Average experimental rates of pointwise convergence for the RBM--CU (left) and RBM--A-WENO (right) schemes.
\label{Fig6.8}}
\end{figure}

We then check the pointwise convergence of the water depth $h$ by computing the average experimental pointwise convergence rates defined in
\eref{equ5.6} using three imbedded grids with $N=2000$ for the RBM--CU and RBM--A-WENO schemes. As in Example 3, we compute the values
$r^{\rm AVE}_{160k}$ using the exact solution and plot the obtained results in Figure \ref{Fig6.8}, where one can clearly see that the
average convergence rates for the RBM--CU and RBM--A-WENO schemes are close to the RBM scheme shown in Figure \ref{Fig5.13} (middle).

\begin{figure}[ht!]
\centerline{\includegraphics[trim=1.3cm 0.3cm 1.0cm 0.1cm, clip, width=5.3cm]{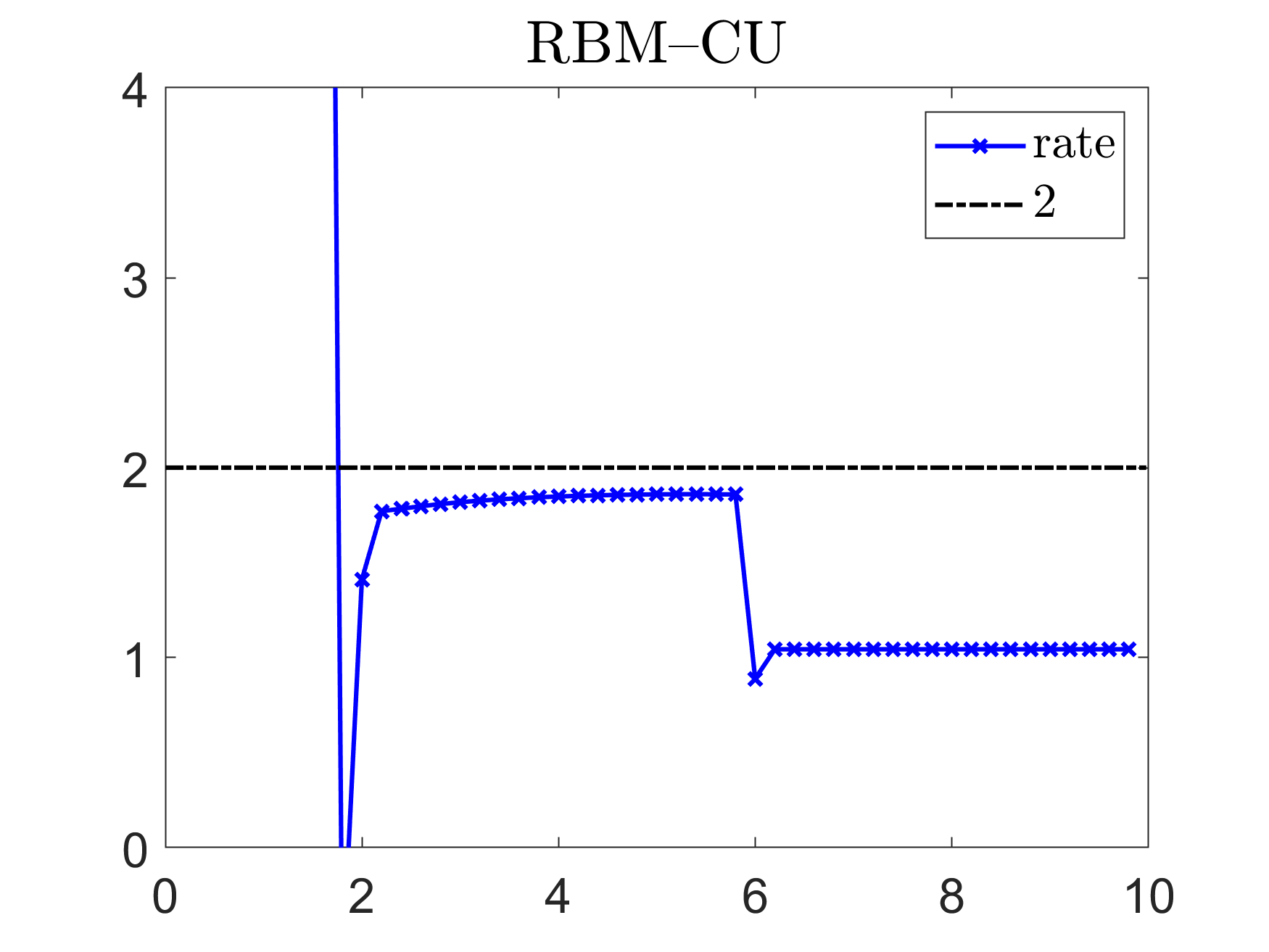}\hspace*{2.0cm}
            \includegraphics[trim=1.3cm 0.3cm 1.0cm 0.1cm, clip, width=5.3cm]{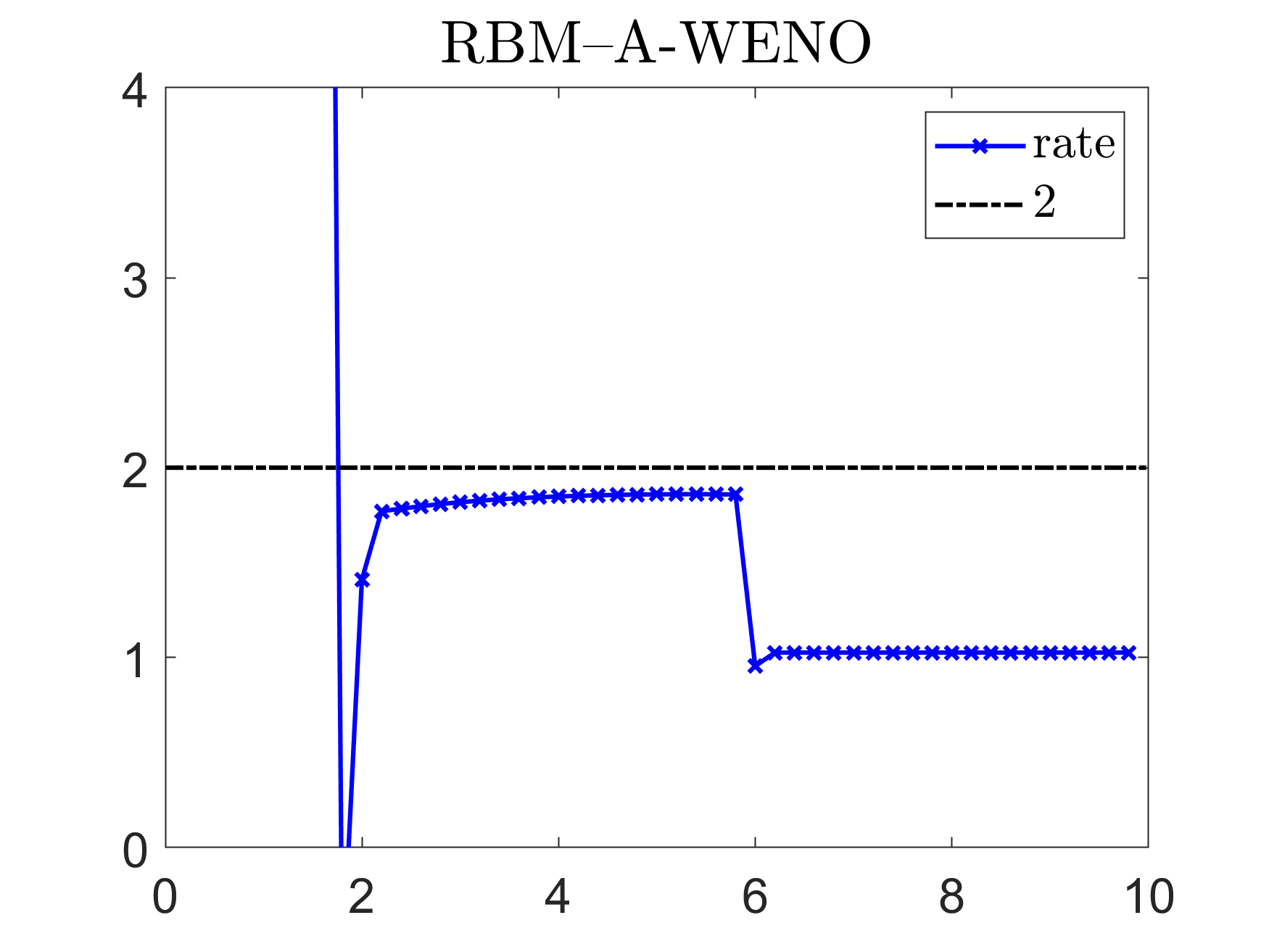}}
\caption{\sf Example 6: Experimental integral rates of convergence for the RBM--CU (left) and RBM--A-WENO (right) schemes.\label{Fig6.8a}}
\end{figure}

As in the previous examples, we also compute the integral and $W^{-1,1}$ convergence rates, which are computed precisely as in Example 3.
The obtained results are shown in Figure \ref{Fig6.8a} and Table \ref{Tab6a}. One can observe that the integral rates of convergence for
both the RBM--CU and RBM--A-WENO schemes reduce to first order even though the pointwise convergence rates for the non-oscillatory combined
schemes are as good as those for the RBM scheme. Once again, this occurs since we combine the schemes of a different nature.
\begin{table}[ht!]
\centering
\begin{tabular}{|c|c|c|c|c|c|c|c|}\hline
&\multicolumn{2}{c|}{RBM--CU Scheme}&\multicolumn{2}{c|}{RBM--A-WENO Scheme}\\\hline
$N$&$||I^N-I^{Exact}||_{L^1}$&$r^{\rm INT}$&$||I^N-I^{Exact}||_{L^1}$&$r^{\rm INT}$\\\hline
2000&1.44e-3&---&1.90e-3&---\\ \hline
4000&1.35e-3&0.09&1.58e-3&0.27\\ \hline
8000&6.56e-4&1.04&7.76e-4&1.03\\\hline
\end{tabular}
\caption{\sf Example 6: $W^{-1,1}$ convergence rates for the studied combined schemes.\label{Tab6a}}
\end{table}

\section{Conclusion}\label{sec8}
In this paper, we have conducted an experimental convergence rate study of the three high-order schemes for hyperbolic systems of
conservation laws: the second-order central-upwind (CU), the third-order Rusanov-Burstein-Mirin (RBM), and the fifth-order alternative
weighted essentially non-oscillatory (A-WENO) schemes. Our goal was to verify that both the pointwise and integral convergence rates of
formally high-order schemes based on nonlinear stabilization mechanism (the CU and A-WENO schemes) reduce to just the first order in the
smooth parts of the solution, which are affected by the propagating shock waves. At the same time, the accuracy of the RBM scheme, which
relies on the linear stabilization mechanism, reduces to the second order in the post shock areas. However, an obvious drawback of the RBM
scheme is that its linear stabilization mechanism is insufficient to prevent oscillations in the shock areas. This suggests that while the
RBM scheme cannot be considered as a good numerical method, it can be used as an excellent component for the development of combined
schemes, which capture shock discontinuities in a non-oscillatory manner, while preserving higher order of accuracy in the smooth parts of
the solution. This goal has been achieved in \S\ref{sec7}, where we have developed two new combined schemes based on the RBM as a basic
scheme and either the CU or A-WENO as an internal scheme. We have conducted numerical experiments to verify both the pointwise
and integral rates of convergence for the developed RBM--CU and RBM--A-WENO schemes. The obtained numerical results suggest that the
combined schemes can achieve similar rates of pointwise convergence to those achieved by the RBM scheme, while their integral rates of
convergence are substantially reduced. This, in turn, indicates that the integral convergence rates do not offer a reliable tool for
measuring the accuracy of combined scheme. We have also observed that the results obtained by the combined schemes with the second- and
fifth-order internal schemes are almost the same even when the computed solution is smooth. This suggests that combining the RBM scheme,
which is highly accurate in the smooth regions and the second-order CU scheme, which captures shocks in a sharp and non-oscillatory manner,
leads to a highly accurate, efficient, and robust combined scheme, and the use of a higher-order internal scheme may not be beneficial.

\begin{acknowledgment}
The reported study was funded in part by RFBR and NSFC, project number 21-51-53012 (RFBR) and 12111530004 (NSFC). The work of O. A.
Kovyrkina and V. V. Ostapenko on the development of a methodology for assessing the accuracy of shock-capturing schemes was supported by the
Russian Science Foundation (project number 22-11-00060). The work of A. Kurganov was supported in part by NSFC grant 12171226 and by the
fund of the Guangdong Provincial Key Laboratory of Computational Science and Material Design (No. 2019B030301001).
\end{acknowledgment}

\appendix
\section{Second-Order CU Scheme: A Brief Overview}\label{appa}
In this appendix, we review the second-order FV CU scheme from \cite{KNP} for the hyperbolic system \eref{1.1}. We describe the CU scheme
applied on the finest mesh of size $\dx$ introduced in \S\ref{sec2a}. Assume that the numerical solution realized in terms of its cell
averages,
\begin{equation*}
\xbar{\mV}_\jph\approx\frac{1}{\dx}\int\limits_{C^{4N}_\jph}\mU(x,t)\,{\rm d}x,
\end{equation*}
is available at a certain time level $t$. The numerical solution is then evolved in time by solving the following system of ODEs:
\begin{equation}
\frac{{\rm d}\xbar{\mV}_\jph}{{\rm d}t}=-\frac{\mH_{j+1}-\mH_j}{\dx},
\label{A2}
\end{equation}
where
\begin{equation}
\mH_j=\frac{a^+_j\bm F(\mV^-_j)-a^-_j\bm F(\mV^+_j)}{a^+_j-a^-_j}+\frac{a^+_ja^-_j}{a^+_j-a^-_j}\left(\mV^+_j-\mV^-_j\right)
\label{A3}
\end{equation}
are CU numerical fluxes derived in \cite{KNP}. Here, $\mV^\pm_j$ are the right- and left-sided point values of $\bm U$ at the cell interface
$x=x_j$ computed using a piecewise linear interpolation
\begin{equation}
\widetilde{\mV}(x)=\,\xbar{\mV}_\jph+(\mV_x)_\jph\big(x-x_\jph\big),\quad x\in C^{4N}_\jph,
\label{A4}
\end{equation}
which gives
\begin{equation}
{\mV}^-_j=\,\xbar{\mV}_\jmh+\frac{\dx}{2}(\mV_x)_\jmh,\quad\mV^+_j=\,\xbar{\mV}_\jph-\frac{\dx}{2}(\mV_x)_\jph.
\label{A5}
\end{equation}
In order to ensure a non-oscillatory nature of the reconstruction \eref{A4}, the slopes $(\mV_x)_\jph$ are to be computed using a nonlinear
limiter. In all of the numerical experiments reported in \S\ref{sec5} and \S\ref{sec7}, we have used the minmod limiter (see, e.g.,
\cite{Sweby84}):
\begin{equation*}
(\mV_x)_\jph={\rm minmod}\left(\frac{\xbar{\mV}_\jph-\,\xbar{\mV}_\jmh}{\dx},
\frac{\xbar{\mV}_{j+\frac{3}{2}}-\,\xbar{\mV}_\jph}{\dx}\right),
\end{equation*}
applied in the component-wise manner. Here, the minmod function is defined as
\begin{equation}
{\rm minmod}(z_1,z_2):=\frac{{\rm sign}(z_1)+{\rm sign}(z_2)}{2}\cdot\min\big(|z_1|,|z_2|\big).
\label{A7}
\end{equation}
Finally, the one-sided local speeds of propagation,
\begin{equation}
a^+_j=\max\Big\{\lambda_d(A\big(\mV^+_j)\big),\lambda_d\big(A(\mV^-_j)\big),0\Big\},\quad
a^-_j=\min\Big\{\lambda_1\big(A(\mV^+_j)\big),\lambda_1\big(A(\mV^-_j)\big),0\Big\},
\label{A8}
\end{equation}
are estimated using the eigenvalues $\lambda_1(A)\le\ldots\le\lambda_d(A)$ of the Jacobian $A=\nicefrac{\partial\mF}{\partial\mU}$.

\section{Third-Order RBM Scheme: A Brief Overview}\label{appb}
In this appendix, we briefly review the third-order RBM scheme introduced in \cite{Rusanov68,Burstein70} for the hyperbolic system
\eref{1.1}. We describe the RBM scheme applied on the finest mesh of size $\dx$ introduced in \S\ref{sec2a}.

Assume that the computed point values $\mV_\jph(t)$ are available at a certain time level $t$. According to the RBM scheme, the solution at
the next time level $(t+\dt)$ is obtained as follows:
\begin{equation*}
\begin{aligned}
&\mV^{(1)}_j=\frac{\mV_\jmh(t)+\mV_\jph(t)}{2}-\frac{\dt}{3\dx}\Big[\mF\big(\mV_\jph(t)\big)-\mF\big(\mV_\jmh(t)\big)\Big],\\
&\mV^{(2)}_\jph=\mV_\jph(t)-\frac{2\dt}{3\dx}\Big[\mF\big(\mV^{(1)}_{j+1}\big)-\mF\big(\mV^{(1)}_j\big)\Big],\\
&\mV_\jph(t+\dt)=\mV_\jph(t)-\frac{\dt}{24\dx}\Big[7\left\{(\mF\big(\mV_{j+\frac{3}{2}}(t)\big)-\mF\big(\mV_{j-\hf}(t)\big)\right\}-
2\left\{(\mF\big(\mV_{j+\frac{5}{2}}(t)\big)-\mF\big(\mV_{j-\frac{3}{2}}(t)\big)\right\}\Big]\\
&\hspace{2.35cm}-\frac{3\dt}{8\dx}\Big[\mF\big(\mV^{(2)}_{j+\frac{3}{2}}\big)-\mF\big(\mV^{(2)}_{j-\hf}\big)\Big]-
\frac{\texttt C}{24}\bm w_\jph,
\end{aligned}
\end{equation*}
where
\begin{equation*}
\bm w_\jph:=\mV_{j+\frac{5}{2}}(t)-4\mV_{j+\frac{3}{2}}(t)+6\mV_\jph(t)-4\mV_\jmh(t)+\mV_{j-\frac{3}{2}}(t)
\end{equation*}
is an artificial viscosity that approximates the fourth spatial derivative term $(\dx)^4\mU_{xxxx}$.

The constant $\texttt C$ in \eref{4.1} is the viscosity coefficient. At $\texttt C=0$ the RBM scheme has the order of accuracy
${\cal O}(h^4+\tau^3)$. However, in this case, it is linearly unstable and in order to stabilize the RBM scheme, it is necessary to take a
positive coefficient $\texttt C$, which satisfies the inequalities $z^2(4-z^2)\le\texttt C\le3$, where $z$ is the CFL number. This leads to
the third order of the RBM scheme both in time and space. Different values of $\texttt C$ were used in \cite{Burstein70} ($\texttt C=2.8$)
and \cite{Rusanov68} ($\texttt C=2.5$). In all of the examples reported in \S\ref{sec5} and \S\ref{sec7}, we have taken $\texttt C=2.8$.

\section{Fifth-Order Finite-Difference A-WENO Scheme: A Brief Overview}\label{appc}
In this appendix, we proceed along the lines of \cite{Kurganov20} and review the fifth-order A-WENO scheme for the system \eref{1.1}, which
reads as
\begin{equation}
\frac{{\rm d}{\bm{\mV}}_\jph}{{\rm d}t}=-\frac{\bm{\mH}_{j+1}-\bm{\mH}_j}{\dx}
+\frac{\dx}{24}\left[(\mF_{xx})_{j+1}-(\mF_{xx})_j\right]-\frac{7(\dx)^3}{5760}\left[(\mF_{xxxx})_{j+1}-(\mF_{xxxx})_j\right].
\label{C1}
\end{equation}
Here,
\begin{equation}
\mH_j=\frac{a^+_j\mF\big(\mV^-_j\big)-a^-_j\mF\big(\mV^+_j\big)}{a^+_j-a^-_j}+\frac{a^+_ja^-_j}{a^+_j-a^-_j}
\left(\mV^+_j-\mV^-_j-\bm Q_j\right),
\label{C2a}
\end{equation}
where
\begin{equation*}
\bm Q_j={\rm minmod}\left(\mV^+_j-\mV^*_j,\mV^*_j-\mV^-_j\right)
\end{equation*}
with
\begin{equation}
\mV^*_j=\frac{a^+_j\mV^+_j-a^-_j\mV^-_j-\left\{\mF\big(\mV^+_j\big)-\mF\big(\mV^-_j\big)\right\}}{a^+_j-a^-_j}
\label{C2c}
\end{equation}
is a built-in anti-diffusion term in the numerical flux \eref{C2a}. In \eref{C2a} and \eref{C2c}, $a^\pm_j$ are the local one-sided
propagation speeds defined in \eref{A8}. In \eref{C2a}--\eref{C2c}, $\bm{\mV}^\pm_j$ are the one-sided point values computed using the
fifth-order WENO-Z interpolant from \cite{Jiang13,liu17,wang18} applied to the local characteristic variables introduced in \S\ref{sec4}.
For the sake of brevity, we present the details on the computation of the left-sided value $\mV^-_j$, as $\mV^+_j$ can be obtained in a
mirror-symmetric way.

For simplicity of presentation, we proceed with a component-wise approach; an extension to the characteristic-wise case is quite
straightforward. For a certain component of $\mU$, the value $\cV_j^-$ is calculated using a weighted average of the three parabolic
interpolants ${\cal P}_0(x)$, ${\cal P}_1(x)$ and ${\cal P}_2(x)$ obtained using the stencils
$[x_{j-\frac{5}{2}},x_{j-\frac{3}{2}},x_\jmh]$, $[x_{j-\frac{3}{2}},x_\jmh,x_\jph]$ and $[x_\jmh,x_\jph,x_{j+\frac{3}{2}}]$, respectively:
\begin{equation}
\cV_j^-=\sum_{k=0}^2\omega^{(i)}_k{\cal P}^{(i)}_k(x_j),
\label{C2}
\end{equation}
where
\begin{equation*}
\begin{aligned}
&{\cal P}^{(i)}_0(x_j)=\frac{3}{8}\,V^{(i)}_{j-\frac{5}{2}}-\frac{5}{4}\,V^{(i)}_{j-\frac{3}{2}}+\frac{15}{8}\,V^{(i)}_\jmh,\\
&{\cal P}^{(i)}_1(x_j)=-\frac{1}{8}\,V^{(i)}_{j-\frac{3}{2}}+\frac{3}{4}\,V^{(i)}_\jmh+\frac{3}{8}\,V^{(i)}_\jph,\\
&{\cal P}^{(i)}_2(x_j)=\frac{3}{8}\,V^{(i)}_\jmh+\frac{3}{4}\,V^{(i)}_\jph-\frac{1}{8}\,V^{(i)}_{j+\frac{3}{2}},
\end{aligned}
\end{equation*}
and the weights $\omega^{(i)}_k$ are computed by
\begin{equation*}
\omega^{(i)}_k=\frac{\alpha^{(i)}_k}{\alpha^{(i)}_0+\alpha^{(i)}_1+\alpha^{(i)}_2},\quad
\alpha^{(i)}_k=d_k\left[1+\bigg(\frac{\tau^{(i)}_5}{\beta^{(i)}_k+\varepsilon}\bigg)^p\right],\quad k=0,1,2,
\end{equation*}
with $d_0=\frac{1}{16}$, $d_1=\frac{5}{8}$ and $d_2=\frac{5}{16}$. The smoothness indicators $\beta^{(i)}_k$ for the corresponding parabolic
interpolants ${\cal P}^{(i)}_k(x)$ are defined by
\begin{equation}
\beta^{(i)}_k=\sum_{\ell=1}^2(\dx)^{2\ell-1}\int\limits_{C^{4N}_j}\bigg(\frac{\partial^\ell{\cal P}^{(i)}_k}{\partial x^\ell}\bigg)^2
{\rm d}x,\quad k=0,1,2.
\label{C4}
\end{equation}
Evaluating the integrals in \eref{C4}, we obtain
\begin{equation}
\begin{aligned}
&\beta^{(i)}_0=\frac{13}{12}\big(V^{(i)}_{j-\frac{5}{2}}-2V^{(i)}_{j-\frac{3}{2}}+V^{(i)}_\jmh\big)^2+
\frac{1}{4}\big(V^{(i)}_{j-\frac{5}{2}}-4V^{(i)}_{j-\frac{3}{2}}+3V^{(i)}_\jmh\big)^2,\\
&\beta^{(i)}_1=\frac{13}{12}\big(V^{(i)}_{j-\frac{3}{2}}-2V^{(i)}_\jmh+V^{(i)}_\jph\big)^2+
\frac{1}{4}\big(V^{(i)}_{j-\frac{3}{2}}-V^{(i)}_\jph\big)^2,\\
&\beta^{(i)}_2=\frac{13}{12}\big(V^{(i)}_\jmh-2V^{(i)}_\jph+V^{(i)}_{j+\frac{3}{2}}\big)^2+
\frac{1}{4}\big(3V^{(i)}_\jmh-4V^{(i)}_\jph+V^{(i)}_{j+\frac{3}{2}}\big)^2.
\end{aligned}
\label{C5}
\end{equation}
Finally, in formula \eref{C4}, $\tau^{(i)}_5=\big|\beta^{(i)}_2-\beta^{(i)}_0\big|$, and in all of the numerical examples reported in
\S\ref{sec5} and \S\ref{sec7}, we have chosen $p=2$ and $\varepsilon=10^{-12}$.

Finally, $({\mF_{xx}})_j$ and $({\mF_{xxxx}})_j$ are the higher-order correction terms computed by the fourth- and second-order accurate
FDs, respectively:
\begin{equation*}
\begin{aligned}
&(\mF_{xx})_j=\frac{1}{48(\dx)^2}\left(-5\mF_{j-\frac{5}{2}}+39\mF_{j-\frac{3}{2}}-34\mF_\jmh-34\mF_\jph+39\mF_{j+\frac{3}{2}}-
5\mF_{j+\frac{5}{2}}\right),\\
&(\mF_{xxxx})_j=\frac{1}{2(\dx)^4}\left(\mF_{j-\frac{5}{2}}-3\mF_{j-\frac{3}{2}}+2\mF_\jmh+2\mF_\jph-3\mF_{j+\frac{3}{2}}+
\mF_{j+\frac{5}{2}}\right),
\end{aligned}
\end{equation*}
where $\mF_\jph:=\mF\big(\mV_\jph\big)$.

\bibliographystyle{siam}
\bibliography{ref}
\end{document}